\documentclass{article}

\usepackage[english]{babel}
\usepackage[T1]{fontenc}
\usepackage[utf8]{inputenc}
\usepackage{amsfonts}
\usepackage{amsmath}
\usepackage{dsfont}
\usepackage[thmmarks,amsmath]{ntheorem}
\usepackage{amssymb}
\usepackage{geometry}
\usepackage{mathtools}
\usepackage[backend=biber,sorting=nty,style=alphabetic,maxbibnames=6]{biblatex}
\usepackage{xcolor}

{
\newtheorem{Def}{Definition}
}

{
\newtheorem{Lem}[Def]{\textit{Lemma}}
}

{
\newtheorem{Cor}[Def]{\textit{Corollary}}
}

\newtheorem{Prop}[Def]{\textit{Proposition}}

\newtheorem{Theo}[Def]{Theorem}

{
\theoremstyle{break}
\theoremsymbol{$\square$}
\theorembodyfont{}
\newtheorem*{Dem}{Proof}
}

{
\theorembodyfont{}
\theoremindent 0.5cm
\newtheorem*{Intuit}{Strategy of the proof (informal)}
}

{
\theorembodyfont{}
\newtheorem{Rem}[Def]{\textit{Remark}}
}

\renewenvironment{subequations}[1][]{
  \refstepcounter{equation}%
  \setcounter{parentequation}{\value{equation}}
  \setcounter{equation}{0}
  \let\parentlabel\label
  \ifx\\#1\\\relax\else\label{#1}\fi
  \ignorespaces
}{%
  \setcounter{equation}{\value{parentequation}}
  \ignorespacesafterend
}

\newcommand*{\nextParentEquation}[1][]{
  \refstepcounter{parentequation}
  \setcounter{equation}{0}
  \ifx\\#1\\\relax\else\parentlabel{#1}\fi
}

\usepackage{hyperref}

\title{Stability of the constant states in the augmented Born-Infeld system}
\date{}
\author{Philippe Anjolras\footnote{Laboratoire de Mathématiques d'Orsay (LMO), philippe.anjolras@universite-paris-saclay.fr}}

\begin{document}

\allowdisplaybreaks

\maketitle

\begin{abstract} In this paper, we consider the Born-Infeld system, arising as a nonlinear model of electromagnetism, and its extension introduced by Brenier \cite{Brenier_ABI} the so-called "augmented Born-Infeld system". We show that this system enjoys a non-resonance structure and prove global existence and linear asymptotic behaviour of small (admissible) perturbations of arbitrary constant states. 
\end{abstract}


\tableofcontents

\section{Introduction}

We consider the Born-Infeld equation~: 
\begin{equation} \left\{ \begin{array}{l}
\partial_t B + \nabla \wedge E = 0 \\
\partial_t D - \nabla \wedge H = 0 \\
\nabla \cdot B = 0 \\
\nabla \cdot D = 0 \\
E = \frac{B \wedge V + D}{h}, ~~~ H = \frac{-D \wedge V + B}{h} \\
h = \sqrt{1 + B^2 + D^2 + |D \wedge B|^2}, ~~~ V = D \wedge B \end{array} \right.  \label{equ_BI} \end{equation}
which is a nonlinear model of electromagnetism. Our main goal is to prove stability of this system around constant states under small and localized perturbations. 

As shown by Brenier in \cite{Brenier_ABI}, the system \eqref{equ_BI} can be seen as a particular case of the so-called "augmented Born-Infeld system" (in short ABI), which consists in writing a conservation law satisfied by $h, V$ and then considering them as independant variables. Then, applying a suitable change of variables, one can recover a very elegant and much more practical version of the system \eqref{ABI_Brenier}, in which all nonlinearities are quadratic and the hyperbolicity is an immediate consequence of the symmetry~: 
\begin{equation}
\left\{ \begin{array}{l}
\partial_t \tau + v \cdot \nabla \tau - \tau \nabla \cdot v = 0 \\
\partial_t v + v \cdot \nabla v - b \cdot \nabla b - d \cdot \nabla d - \tau \nabla \tau = 0 \\
\partial_t b + v \cdot \nabla b - b \cdot \nabla v + \tau \nabla \wedge d = 0 \\
\partial_t d + v \cdot \nabla d - d \cdot \nabla v - \tau \nabla \wedge b = 0
\end{array} \right. \label{ABI_Brenier}
\end{equation}
In the system above, we call a constant solution \textit{admissible} if $\tau > 0$. We can then restate the divergence-free conditions in terms of this new set of variables, and it appears natural in our analysis to add a similar condition on the additional variable $V$~: 
\begin{equation}
\left\{ \begin{array}{l}
\tau \nabla \cdot b = b \cdot \nabla \tau \\
\tau \nabla \cdot d = d \cdot \nabla \tau \\
\tau \nabla \wedge v = b \cdot \nabla d - d \cdot \nabla b \label{equ_cont_explicite}
\end{array} \right. 
\end{equation}
that we call "constraint equations", since they do not depend on the time and are preserved by the time evolution, thus only adding a constraint on the initial data. Given a constant state $c \in \mathbb{R}^{10}$, we denote by $u : [0, T) \times \mathbb{R}^3 \to \mathbb{R}^{10}$ a smooth enough function such that $c + u = (\tau, v, b, d)$ is a solution to \eqref{ABI_Brenier} under the constraints \eqref{equ_cont_explicite} on $[0, T)$, $T \in \mathbb{R}^{+*} \cup \{ +\infty \}$. We can then state our main result~: 

\begin{Theo} The admissible constant solutions $c \in \mathbb{R}^{10}$ of the augmented Born-Infeld system consisting of \eqref{ABI_Brenier} and \eqref{equ_cont_explicite} are stable under small and localized perturbations $u_0 : \mathbb{R}^3 \to \mathbb{R}^{10}$, in the sense that if the following norms are small enough (with respect to the constant state only)~: 
\[ \Vert u_0 \Vert_{H^N}, ~~~ \Vert \langle x \rangle^2 u_0 \Vert_{H^6} \]
and $c+u_0$ satisfies the constraint equations \eqref{equ_cont_explicite}, then there exists a global solution to \eqref{ABI_Brenier} starting from $c+u_0$ such that the $L^{\infty}$ norm of the perturbation and its first derivative decays like $t^{-1}$. \end{Theo}

In particular, we get global existence and stability of constants for two physical subsystems~: the original Born-Infeld system \eqref{equ_BI} and the irrotational Chaplygin gas, obtained from \eqref{ABI_Brenier} + \eqref{equ_cont_explicite} when $b \equiv d \equiv 0$ (see section \ref{section_structure_BI_Chap}). 

\paragraph{Comparison with earlier work} Concerning the Born-Infeld system, a particular case of this result of stability, the stability of the zero constant state, was proved by Speck in \cite{Speck_BI}, using the vector fields method and the null condition~: stability was only proved for the initial system \eqref{equ_BI}, under different assumptions of smallness on the initial perturbations, and only around the trivial state $B \equiv D \equiv 0$. For the Chaplygin gas, the same kind of result using the method of vector fields with a null condition was also obtained in \cite{Chaplygin_null}, following Speck's approach. We are therefore able to group these two (seemingly unrelated) PDEs in a larger system. The method of vector fields with a null condition was introduced by Klainerman in \cite{Klainerman_null1, Klainerman_null2} for the study of the long-time behaviour of nonlinear wave equations with small initial data. It was then used to study Maxwell equations in the Minkowski space \cite{Christodoulou_Minkowski}, to prove global nonlinear stability results for the Einstein equation \cite{BieriZipser_Mink, CK_Mink, Dafermos_blackholes, Klainerman_evolutionpb, Lindblad_Einstein, Lindblad_Mink, Rodnianski_EulerEinstein, Speck_Mink}, global existence for nonlinear elastic waves with small initial data \cite{Sideris_nullcond}, the formation of shocks in solutions to the relativistic Euler equations \cite{Christodoulou_formationshock}, the study of decay of linear solutions on curved spaces \cite{Anderrson_hiddensym, Blue_decay, Dafermos_collapse, Dafermos_redshift, Holzegel_massivewave}, the local existence and classical limits of the relativistic Euler equations \cite{Speck_Eulerwp, Speck_Eulerlim, Speck_Boltz}. 

In our case, having a non-zero constant state as reference complicates the geometry of the vector fields method. Therefore, our analysis relies on the space-time resonance method initially developed by Germain, Masmoudi and Shatah in \cite{Germain_3DSchrod}, \cite{Germain_2DSchrod} and \cite{GMSWaterwaves} for the quadratic Schrödinger equation and gravity water waves, and simultaneously by Gustafson, Nakanishi and Tsai in \cite{NakanishiResonances} and \cite{NakanishiResonances2} for the Gross-Pitaevskii equation. More precisely, the Born-Infeld system around a constant is a quasilinear wave equation coupled with a non-linear non-dispersive part. In particular, for the purely wave part, we adapt the proof by Pusateri and Shatah in \cite{ShatahPusateri} of the global existence result for semilinear quadratic wave equations satisfying a null form condition, which extended the results of Klainerman \cite{Klainerman_null1, Klainerman_null2}. 

As we shall explain  below some new difficulties arise here due to the quasilinear structure of the Born-Infeld system (whereas \cite{ShatahPusateri} only deals with semilinear first order systems) and more importantly in order to handle the non-dispersive part by using the constraint equations. The space-time resonance methode has been widely used in the recent years to prove global regularity of nonlinear PDEs, e.g. water waves \cite{GMSWaterwaves,IonescuPusateri_waterwaves,GMS_waterwavescap}, the Klein-Gordon equation \cite{IonescuPausader_KG}, the Euler-Maxwell system \cite{Germain_EulerMaxwell, GuoIonescuPausader_plasma,GuoIonescuPausader_EulerMaxwell,DengIonescuPausader_EulerMaxwell}, the Euler-Poisson system \cite{GuoPausader_EulerPoisson,GuoIonescuPausader_plasma}. It has also been used to get scattering \cite{IonescuPausader_KG,IonescuPusateri_waterwaves,GMS_waterwavescap,GermainPusateriRousset_mKDV}. Similar results concerning global regularity and scattering have also been obtained in \cite{DelortAlazard_gravitywaterwaves2015,HunterIfrimTataru_ww,IfrimTataru_ww}. 

The result of \cite{ShatahPusateri} cannot apply directly to our case because of two differences. First, the nonlinearity handled in \cite{ShatahPusateri} is homogeneous of degree $0$ in derivatives, while the Born-Infeld system has nonlinearities of degree $1$. This leads to combine the quasilinear structure of the equations with the space-time resonance method, in particular to handle weighted estimates, in a similar fashion as what was done in \cite{Stewart_complexmkdv} for the complex mKdV equation. Furthermore, recall that for the traditional linear Maxwell equation in vacuum~: 
\[ \begin{aligned}
&\partial_t B + \nabla \wedge E = 0, ~~~~~ \nabla \cdot B = 0 \\
&\partial_t E - \nabla \wedge B = 0, ~~~~~ \nabla \cdot E = 0 \end{aligned} \]
it is only the set of constraint equations $\nabla \cdot B = \nabla \cdot E = 0$ that ensures that $B$ and $E$ follow a wave equation. The same occurs for the Born-Infeld system (note that the linear Maxwell equations are recovered, at least formally, as the weak-field limit of Born-Infeld), which satisfies two conditions of non-resonance~: one is encoded in the evolution equations, and another appears in the constraints equations. It is only the combination of these two structures that allows to prove stability. In comparison with earlier works treating equations coupling a dispersive PDE with an elliptic constraint with the space-time resonance method, the constraints were linear (or could be made linear by a suitable change of coordinates), usually written as a null-divergence of some field, thus making it possible to simply project all equations in Fourier space, while in the Born-Infeld system we deal with nonlinear constraints and analyze their resonances as well. We also point out the correct constraint equation to be added on $v$ (third line of \eqref{equ_cont_explicite}) in the augmented system, which was not known previously. 

Note that the non-resonance structure appearing in the Born-Infeld system is simpler than the general case treated in \cite{ShatahPusateri}~: indeed, our non-resonance structure only involves the so-called space-resonances. However, the time-resonances appear anyway in the a priori estimates in order to control the weighted estimates. We also get rid of a technical assumption in \cite{ShatahPusateri} (assumption $(2.4)$) which is not verified in the Born-Infeld system, by using Besov spaces instead of Lebesgue at one point of the proof. 

\section{Outline of the article}

\paragraph{} In section \ref{section_structure}, we analyze the structure of the augmented Born-Infeld equation \eqref{ABI_Brenier}. In subsection \ref{section_structure_BI}, we recall results proved by Brenier in \cite{Brenier_ABI} concerning the derivation of \eqref{ABI_Brenier} from \eqref{equ_BI}. In subsection \ref{section_structure_BI_cst}, we specialize this analysis by rewriting the equation around a constant state, separating the linear and non-linear terms. We also define and prove the non-resonance structure that we will use in the proof and introduce the following separation of the perturbation into three terms~: 
\begin{equation} u = u^{+} + u^{-} + u^0 \label{sep_termes} \end{equation}
They are defined by the fact that the linear part of the evolution equation satisfied by each is 
\[ \mathcal{L}^{+} = \partial_t + i \Lambda_0, ~~~ \mathcal{L}^{-} = \partial_t - i \Lambda_0, ~~~ \mathcal{L}^0 = \partial_t \]
for $u^{+}$, $u^{-}$ and $u^0$ respectively, where $\Lambda_0 = |\nabla|_0$ the linear operator of Fourier symbol $|\xi|_0$, and $|\cdot|_0$ is an (Euclidean) norm depending on the constant state. Here we see that there is no linear dispersion on $u^0$ from the evolution equation \eqref{ABI_Brenier} alone~: another key point is to prove that the constraint equations \eqref{equ_cont_explicite} exactly control $u^0$, and that they also exhibit a non-resonance structure. In subsection \ref{section_structure_BI_Chap}, we provide a short introduction to the Chaplygin gas system and show its relation to the augmented Born-Infeld system. In subsection \ref{section_structure_ABI}, we prove that the constraint equations \eqref{equ_cont_explicite} are preserved by the augmented Born-Infeld system \eqref{equ_ABI} in the general case. 

\paragraph{} Let us introduce briefly what we mean by "non-resonance structure". Using the separation \eqref{sep_termes}, we may decompose the non-linearity for each $u^{\epsilon_1}$, $\epsilon_1 \in \{ +, -, 0 \}$ into nine different interactions, characterized by $(\epsilon_2, \epsilon_3) \in \{ +, -, 0 \}^2$. Introducing the profile $f^{\epsilon}(t) = e^{\epsilon i t \Lambda} u^{\epsilon}(t)$, $f^{\epsilon_1}$ satisfies an equation without any linear term and with non-linearities that we may write as a sum over $(\epsilon_2, \epsilon_3)$ of 
\[ \mathcal{F}^{-1} \int_{\mathbb{R}^3} e^{i t \varphi^{\epsilon_1, \epsilon_2 \epsilon_3}(\xi, \eta)} m(\xi, \eta) \widehat{f}^{\epsilon_2}(t, \xi - \eta) \widehat{f}^{\epsilon_3}(t, \eta) ~ d\eta \]
for some symbol $m$, and where 
\[ \varphi^{\epsilon_1, \epsilon_2 \epsilon_3}(\xi, \eta) = \epsilon_1 |\xi|_0 - \epsilon_2 |\xi - \eta|_0 - \epsilon_3 |\eta|_0 \]
The non-resonance structure, which will be satisfied if $(\epsilon_1, \epsilon_2, \epsilon_3) \in \{ +, - \}^3$, states that 
\[ m(\xi, \eta) = |\xi - \eta| \mu_0(\xi, \eta) \nabla_{\eta} \varphi(\xi, \eta) \]
for some symbol $\mu_0$, homogeneous of order $0$, and satisfying appropriate smoothness conditions that allow to control the bilinear terms above. The presence of such a factorisation allows then to integrate by parts in $\eta$, thus gaining a factor $s^{-1}$, up to the apparition of terms $\nabla_{\eta} \widehat{f}$, which translates in the physical space to $x f$. 

Note that, when $\epsilon_2 = 0$ and $\epsilon_3 \neq 0$, or inversely, $\nabla_{\eta} \varphi$ never vanishes (outside the axis $\{ \xi = 0 \} \cup \{ \eta = 0 \} \cup \{ \xi = \eta \}$) and thus such a factorisation is always allowed. 

On the other hand, to control the nonlinearities when $\epsilon_2 = \epsilon_3 = 0$ or when $\epsilon_1 = 0$, we also show that the same type of factorisation appears in the constraint equations \eqref{equ_cont_explicite}. Since these equations do not involve any time derivative, we won't have to apply any Duhamel formula for them and thus will not have any time integration~: therefore, $u^0$ will actually satisfy stronger bounds than $u^{\pm}$, thus allowing to bound the interaction terms in which it is involved. However, since the constraint equations involve derivatives, we actually only get to control $\Lambda u^0$ as a quadratic non-linearity, which adds a singularity. 

\paragraph{} Section \ref{section_main_result} is devoted to the statement of our main result again, followed in section \ref{section_statement_apriori} by the statement of our main a priori estimate, which is then proved in sections \ref{section_energyest-0} to \ref{section_L2x2}. Namely, we will bound as follows~: 
\[ \left\{ 
\begin{array}{l}
\Vert u^{\pm} \Vert_{H^N} \lesssim t^{\varepsilon}, ~~~ \Vert u^{\pm} \Vert_{H^6} \lesssim 1, ~~~ \Vert u^{\pm} \Vert_{W^{1, \infty}} \lesssim t^{-1}, ~~~ \Vert x f^{\pm} \Vert_{H^5} \lesssim 1, ~~~ \Vert |x|^2 \Lambda f^{\pm} \Vert_{H^5} \lesssim t^{\gamma} \\
\Vert u^0 \Vert_{H^N} \lesssim t^{-1+\varepsilon+a}, ~~~ \Vert u^0 \Vert_{H^6} \lesssim t^{-2+a}, ~~~ \Vert u^0 \Vert_{W^{1, \infty-}} \lesssim t^{-2+a}, ~~~ \Vert x \Lambda u^0 \Vert_{H^4} \lesssim t^{-1+a}
\end{array} \right. \]
where $N$ is an integer large enough, $\varepsilon, \gamma, a$ are small enough parameters, $\infty-$ denotes a large enough real number. Here and in the whole article, we use the notation $A \lesssim B$ when there exists an universal constant $C > 0$ (which may depend on the parameters of the problem, like $N, \varepsilon, \gamma, ...$ assumed to be fixed once and for all, or on the reference constant state, but not on the added perturbation or the time $t$) such that $A \leq C B$. 

In section \ref{section_geninequ}, we prove or recall useful lemmas and inequalities that we will use repeatedly in the proof, as well as some identities specific to the wave equation and to our system. 

In subsection \ref{section_energyest}, we prove an energy estimate by the classical tools on hyperbolic systems. We do not need the resonance analysis in this section. In subsection \ref{section_est0}, we prove all the estimates on $u^0$, relying on the constraint equations. In subsection \ref{section_estH6}, we prove the $H^6$-estimate on $u^{\pm}$ in a very similar manner to the corresponding bound from \cite{ShatahPusateri}, up to control of terms involving $u^0$. 

In section \ref{section_L2x}, we prove the estimate on the Sobolev norm of $x f^{\pm}$. Once again, we follow a similar approach to the one from \cite{ShatahPusateri}. The proof is divided, as in the following, in subsections corresponding to the different types of interactions. 

In section \ref{section_Besov}, we prove a bound involving the Besov spaces $\dot{B}^0_{\infty, 1}$ and $\dot{B}^1_{\infty, 1}$. These norms control in particular the $W^{1, \infty}$ norm and are needed due to the lack of $L^2 \times L^{\infty}$ estimate for the bilinear interactions we consider. The analysis of the $\pm \pm$ interactions is similar to the one in \cite{ShatahPusateri} (up to dealing with Besov spaces instead of $L^{\infty}$) with the introduction of an angular repartition that separates between different domains depending on the space- and time-resonant sets. However, the analysis of the $\pm 0$ interactions involves a finer analysis by expanding $u^0$ again as a bilinear term, using the constraint equations, and the introduction of a new angular repartition. 

In subsection \ref{section_L2x2}, we prove the bound on $|x|^2 f^{\pm}$ in Sobolev spaces. Since our nonlinearities have one derivative, we need to combine both the quasi-linear hyperbolic structure of the equation to avoid losing derivatives and the space-time analysis. Again, the new $\pm 0$ interactions require new arguments in a similar fashion to what was done for the Besov norms. 

\section{Structure of the equations} \label{section_structure}

\paragraph{Notation} We denote by $\Lambda = |\nabla|$ the linear operator of Fourier symbol $|\xi|$. 

\subsection{Born-Infeld system} \label{section_structure_BI}

The Born-Infeld system, introduced by Max Born and Leopold Infeld in 1934 (\cite{Born_Infeld34}), comes from the Lagrangian~: 
\[ L = - \sqrt{1 - E^2 + B^2 - (E \cdot B)^2} \]
(expressed in renormalized units) under the differential constraints~: 
\[ \nabla \cdot B = 0, ~~~ \partial_t B + \nabla \wedge E = 0 \]
which express the fact that the Faraday tensor associated to $(B, E)$ is a closed form. This leads to the equations~: 
\[ \left\{ \begin{array}{l}
\partial_t B + \nabla \wedge E = 0 \\
\partial_t D - \nabla \wedge H = 0 \\
\nabla \cdot B = 0 \\
\nabla \cdot D = 0 \\
E = \frac{B \wedge V + D}{h}, ~~~ H = \frac{-D \wedge V + B}{h} \\
h = \sqrt{1 + B^2 + D^2 + |D \wedge B|^2}, ~~~ V = D \wedge B \end{array} \right. \]

In 2004, Brenier introduced the following augmented version of the Born-Infeld system (\cite{Brenier_ABI})~: 
\[ \left\{ \begin{array}{l}
\partial_t \tau + v \cdot \nabla \tau - \tau \nabla \cdot v = 0 \\
\partial_t v + v \cdot \nabla v - b \cdot \nabla b - d \cdot \nabla d - \tau \nabla \tau = 0 \\
\partial_t b + v \cdot \nabla b - b \cdot \nabla v + \tau \nabla \wedge d = 0 \\
\partial_t d + v \cdot \nabla d - d \cdot \nabla v - \tau \nabla \wedge b = 0
\end{array} \right. \]
by setting~: 
\[ \tau = h^{-1}, ~~~ v = \tau V, ~~~ b = \tau B, ~~~ d = \tau D \]
This system is symmetric and therefore well-posed, regardless of the values taken by $\tau, v, b, d$. We recover the initial Born-Infeld system by introducing the following constraints~: 
\[ \tau > 0, ~~~ \tau^2 + v^2 + b^2 + d^2 = 1, ~~~ \tau v = d \wedge b \]
and we call Born-Infeld manifold the set of initial data satisfying these. An important property of the Born-Infeld manifold is that it is preserved under the flow of the augmented Born-Infeld system.

Besides, if $\tau > 0$, the divergence constraints can be rewritten as~: 
\[ \nabla \cdot \frac{b}{\tau} = \nabla \cdot \frac{d}{\tau} = 0 \]
and are also preserved by the equation in the case of smooth solutions. 

\paragraph{Galilean invariance} If $(t, x) \mapsto (\tau, v, b, d)(t, x)$ is a solution of the augmented system, then 
\[ (t, x) \mapsto (\tau, v - v_0, b, d)(t, x + v_0 t) \]
is also a solution, where $v_0$ is a constant vector. However, the Born-Infeld manifold is not preserved by this galilean invariance. 

\subsubsection{The Born-Infeld system around a constant solution} \label{section_structure_BI_cst}

If $(\tau_0, v_0, b_0, d_0)$ is a constant (satisfying or not the constraints of the Born-Infeld manifold), then the unique smooth solution is global and defined by being constant in space and time. 

If we rewrite the system around such a solution, replacing $(\tau, v, b, d)$ by $(\tau_0 + \tau, v_0 + v, b_0 + b, d_0 + d)$, we obtain the following equations~: 
\[ \left\{ \begin{array}{l}
\partial_t \tau + v_0 \cdot \nabla \tau - \tau_0 \nabla \cdot v = - v \cdot \nabla \tau + \tau \nabla \cdot v \\
\partial_t v + v_0 \cdot \nabla v - b_0 \cdot \nabla b - d_0 \cdot \nabla d - \tau_0 \nabla \tau = - v \cdot \nabla v + b \cdot \nabla b + d \cdot \nabla d + \tau \nabla \tau \\
\partial_t b + v_0 \cdot \nabla b - b_0 \cdot \nabla v + \tau_0 \nabla \wedge d = - v \cdot \nabla b + b \cdot \nabla v - \tau \nabla \wedge d \\
\partial_t d + v_0 \cdot \nabla d - d_0 \cdot \nabla v - \tau_0 \nabla \wedge b = - v \cdot \nabla d + d \cdot \nabla v + \tau \nabla \wedge b
\end{array} \right. \]
that we can write under the form~: 
\begin{equation} \partial_t U + A_0(D) U = \mathcal{N}(\Lambda U, U) \label{equ_ABI} \end{equation}
where $U = (\tau, v, b, d)$, $A_0$ is a differential operator of order 1, linear and with constant coefficients depending on the constant $U_0 = (\tau_0, b_0, d_0, v_0)$, while $\mathcal{N}$ is a bilinear operator with constant coefficients (but potentially involving Riesz transforms). 

The differential constraints can be written~: 
\[ \nabla \cdot \frac{b_0 + b}{\tau_0 + \tau} = 0 ~~ \Longrightarrow ~~ \tau_0 \nabla \cdot b - b_0 \cdot \nabla \tau = - \tau \nabla \cdot b + b \cdot \nabla \tau \]
and likewise~: 
\[ \tau_0 \nabla \cdot d - d_0 \cdot \nabla \tau = - \tau \nabla \cdot b + b \cdot \nabla \tau \]

Finally, the constraints of the Born-Infeld manifold imply~: 
\[ \tau_0 \nabla \wedge v - b_0 \cdot \nabla d + d_0 \cdot \nabla b = - \tau \nabla \wedge v + b \cdot \nabla d - d \cdot \nabla b \]

We can therefore rewrite all these constraints under the form~: 
\begin{equation} L_0(D) U = \mathcal{N}'(\Lambda U, U) \label{equ_cont} \end{equation}
where $L_0$ is a differential operator of order 1, linear and with constant coefficients depending on the constant $U_0$, while $\mathcal{N}'$ is a bilinear operator with constant coefficients. 

\paragraph{Simplification using the galilean invariance} We saw that the augmented Born-Infeld equation enjoyed a galilean invariance, which does not preserve the Born-Infeld manifold a priori. However, this invariance preserves \eqref{equ_cont}, since no temporal derivatives appear in it and $v$ is present only through its derivatives. Therefore, we can study the simplified problem in which $v_0 = 0$, and we suppose that this condition is satisfied in the following. 

\begin{Prop} For any $\xi$, the matrix $A_0(\xi)$ is symmetric, with eigenvalues $0, |\xi|_0, - |\xi|_0$, where $| \cdot |_0$ is the euclidean norm associated to the scalar product~: 
\[ g_0(\xi, \eta) = \tau_0 \xi \cdot \eta + (b_0 \cdot \xi) (b_0 \cdot \eta) + (d_0 \cdot \xi) (d_0 \cdot \eta) \]
(It is a scalar product because $\tau_0 > 0$.) 

If $\xi \neq 0$, let $(e_1, e_2, e_3)$ be an orthonormal direct basis such that $e_1$ is positively collinear to $\xi$, and set~:
\[ \alpha = \frac{\tau_0 |\xi|}{|\xi|_0}, ~~~ \beta = \frac{b \cdot \xi}{|\xi|_0}, ~~~ \delta = \frac{d \cdot \xi}{|\xi|_0} \]
then the eigenspaces of $A_0(\xi)$ are, respectively~: 
\begin{align*}
E(0) = \mbox{Vect}&\left( \begin{pmatrix} -\beta \\ 0 \\ \alpha e_1 \\ 0 \end{pmatrix}, ~~~ \begin{pmatrix} -\delta \\ 0 \\ 0 \\ \alpha e_1 \end{pmatrix} , ~~~ \begin{pmatrix} 0 \\ \alpha e_2 \\ \delta e_3 \\ - \beta e_3 \end{pmatrix}, ~~~ \begin{pmatrix} 0 \\ \alpha e_3 \\ - \delta e_2 \\ \beta e_2 \end{pmatrix} \right) \\
E(|\xi|_0) = \mbox{Vect}&\left( \begin{pmatrix} \alpha \\ e_1 \\ \beta e_1 \\ \delta e_1 \end{pmatrix}, ~~~ \begin{pmatrix} 0 \\ \delta e_2 - \beta \alpha e_3 \\ \beta \delta e_2 - \alpha e_3 \\ (1 - \beta^2) e_2 \end{pmatrix}, ~~~ \begin{pmatrix} 0 \\ \beta \alpha e_2 + \delta e_3 \\ \alpha e_2 + \beta \delta e_3 \\ (1 - \beta^2) e_3 \end{pmatrix} \right) \\
E(-|\xi|_0) = \mbox{Vect}&\left( \begin{pmatrix} -\alpha \\ e_1 \\ -\beta e_1 \\ -\delta e_1 \end{pmatrix}, ~~~ \begin{pmatrix} 0 \\ -\delta e_2 - \beta \alpha e_3 \\ \beta \delta e_2 + \alpha e_3 \\ (1 - \beta^2) e_2 \end{pmatrix}, ~~~ \begin{pmatrix} 0 \\ \beta \alpha e_2 - \delta e_3 \\ - \alpha e_2 + \beta \delta e_3 \\ (1 - \beta^2) e_3 \end{pmatrix} \right)
\end{align*}
\label{def_alpha_beta_delta}
\end{Prop}

\begin{Rem} In the case $\xi = 0$, $\alpha, \beta, \delta$ are not well defined and $A_0 = 0$ only has one eigenvalue. However, we can consider any limit since the previous spaces are orthogonal to each other and generate $\mathbb{R}^{10}$ as long as $\alpha^2 + \beta^2 + \delta^2 = 1$, $\alpha > 0$. Note furthermore that, when considering for instance $L^2$ norms, removing the point $\xi = 0$ is harmless. \end{Rem} 

\begin{Dem}
We can check by a computation that, for each of the basis vector given in the proposition, $A_0(\xi) X = \lambda X$ with $\lambda \in \{ - |\xi|_0, 0, |\xi|_0 \}$ accordingly. Furthermore, these basis vectors are independant in each eigenspace, so they generate $\mathbb{R}^{10}$. Orthogonality follows from the symmetry of $A_0(\xi)$. 
\end{Dem}

We denote by $r(\xi)$ the operator associated with the cross product~: 
\[ r(D) f = \nabla \wedge f \]

\begin{Cor} The following operators are the projection operators on the eigenspaces of $A_0(\xi)$~: 
\begin{gather*}
P^0(\xi) = \begin{pmatrix} 1 - \alpha^2 & 0 & - \alpha \beta \frac{\xi^T}{|\xi|} & - \alpha \delta \frac{\xi^T}{|\xi|} \\
0 & \alpha^2 \left( I_3 - \frac{\xi \otimes \xi}{|\xi|^2} \right) & - \alpha \delta \frac{r(\xi)}{|\xi|} & \alpha \beta \frac{r(\xi)}{|\xi|} \\
-\alpha \beta \frac{\xi}{|\xi|} & \alpha \delta \frac{r(\xi)}{|\xi|} & \delta^2 I_3 + \alpha^2 \frac{\xi \otimes \xi}{|\xi|^2} & - \beta \delta I_3 \\
-\alpha \delta \frac{\xi}{|\xi|} & -\alpha \beta \frac{r(\xi)}{|\xi|} & -\beta \delta I_3 & \beta^2 + \alpha^2 \frac{\xi \otimes \xi}{|\xi|^2} \end{pmatrix} \\
 P^{+}(\xi) = \frac{1}{2} \begin{pmatrix} 
\alpha^2 & \alpha \frac{\xi^T}{|\xi|} & \alpha \beta \frac{\xi^T}{|\xi|} & \alpha \delta \frac{\xi^T}{|\xi|} \\
\alpha \frac{\xi}{|\xi|} & (1 - \alpha^2) I_3 + \alpha^2 \frac{\xi \otimes \xi}{|\xi|^2} & \beta I_3 + \alpha \delta \frac{r(\xi)}{|\xi|} & \delta I_3 - \alpha \beta \frac{r(\xi)}{|\xi|} \\
\alpha \beta \frac{\xi}{|\xi|} & \beta I_3 - \alpha \delta \frac{r(\xi)}{|\xi|} & (1 - \delta^2) I_3 - \alpha^2 \frac{\xi \otimes \xi}{|\xi|^2} & \beta \delta I_3 - \alpha \frac{r(\xi)}{|\xi|} \\
\alpha \delta \frac{\xi}{|\xi|} & \delta I_3 + \alpha \beta \frac{r(\xi)}{|\xi|} & \beta \delta I_3 + \alpha \frac{r(\xi)}{|\xi|} & (1 - \beta^2) I_3 - \alpha^2 \frac{\xi \otimes \xi}{|\xi|^2}
\end{pmatrix} \\
P^{-}(\xi) = \frac{1}{2} \begin{pmatrix} \alpha^2 & - \alpha \frac{\xi^T}{|\xi|} & \alpha \beta \frac{\xi^T}{|\xi|} & \alpha \delta \frac{\xi^T}{|\xi|} \\
-\alpha \frac{\xi}{|\xi|} & (1 - \alpha^2) I_3 + \alpha^2 \frac{\xi \otimes \xi}{|\xi|^2} & - \beta I_3 + \alpha \delta \frac{r(\xi)}{|\xi|} & - \delta I_3 - \alpha \beta \frac{r(\xi)}{|\xi|} \\
\alpha \beta \frac{\xi}{|\xi|} & - \beta I_3 - \alpha \delta \frac{r(\xi)}{|\xi|} & (1 - \delta^2) I_3 - \alpha^2 \frac{\xi \otimes \xi}{|\xi|^2} & \beta \delta I_3 + \alpha \frac{r(\xi)}{|\xi|} \\
\alpha \delta \frac{\xi}{|\xi|} & - \delta I_3 + \alpha \beta \frac{r(\xi)}{|\xi|} & \beta \delta I_3 - \alpha \frac{r(\xi)}{|\xi|} & (1 - \beta^2) I_3 - \alpha^2 \frac{\xi \otimes \xi}{|\xi|^2}
\end{pmatrix} 
\end{gather*}
expressed in the canonical basis of $\mathbb{R}^{10}$. 
\end{Cor} 

\begin{Dem}
Since we know a basis of eigenvectors for each eigenspace, if to a fixed eigenvalue $\lambda$ we associate $M$ the matrix with these eigenvectors as columns, we can apply the following formula to compute the projection matrix~: 
\[ P = M (M^T M)^{-1} M^T \]
This allows to compute the formulas of the corollary. To obtain the expression in the canonical basis, we can choose arbitrarily the basis $e_1, e_2, e_3$, for instance~: 
\begin{gather*}
e_1 = \frac{1}{|\xi|} (\xi_1, \xi_2, \xi_3) \\ 
e_2 = \frac{1}{\sqrt{\xi_1^2 + \xi_2^2}} (\xi_2, -\xi_1, 0) \\
e_3 = \frac{1}{|\xi| \sqrt{\xi_1^2 + \xi_2^2}} (\xi_1 \xi_3, \xi_2 \xi_3, -\xi_1^2 - \xi_2^2)
\end{gather*}
as long as none of the coordinates of $\xi$ vanish, and then extend by continuity.
\end{Dem}

\begin{Prop} The operator $L_0$ associated to the constraints \eqref{equ_cont} satisfies~: 
\[ Q(\xi) L_0(\xi) = |\xi| P^0(\xi) \]
for a certain invertible operator $Q(\xi)$ which is homogeneous of degree $0$ in $\xi$. 
\end{Prop}

\begin{Dem}
We know that~: 
\[ L_0(\xi) = \begin{pmatrix} b_0 \cdot \xi & 0 & - \tau_0 \xi^T & 0 \\
d_0 \cdot \xi & 0 & 0 & - \tau_0 \xi^T \\
0 & \tau_0 r(\xi) & d_0 \cdot \xi I_3 & - b_0 \cdot \xi I_3 \end{pmatrix} \]
in the canonical basis. 

The eigenvectors of $0$ are 
\[ \begin{aligned}
&~~~ \left(-\beta, 0, \alpha e_1, 0\right), ~~~ \left(-\delta, 0, 0, \alpha e_1\right), ~~~ \\
&\left(0 , \alpha e_2, \delta e_3, - \beta e_3\right), ~~~ \left(0, \alpha e_3, - \delta e_2, \beta e_2\right) 
\end{aligned} \]
that is they correspond to the following components of a given solution, in Fourier space~: 
\[ \begin{aligned}
&~~~ \tau_0 \xi \cdot \hat{b} - (b_0 \cdot \xi) \hat{\tau}, ~~~ \tau_0 \xi \cdot \hat{d} - (d_0 \cdot \xi) \hat{\tau}, ~~~ \\
&\tau_0 |\xi| \hat{v}_2 + (d_0 \cdot \xi) \hat{b}_3 - (b_0 \cdot \xi) \hat{d}_3, ~~~ \tau_0 |\xi| \hat{v}_3 - (d_0 \cdot \xi) \hat{b}_2 + (b_0 \cdot \xi) \hat{d}_2 
\end{aligned} \]
Going back to the physical space, we may combine these components into~: 
\[ \begin{aligned}
&~~~ \tau_0 \nabla \cdot b - b_0 \cdot \nabla \tau, ~~~ \tau_0 \nabla \cdot d - d_0 \cdot \nabla \tau, ~~~ \\
&~~~~~~\tau_0 \nabla \wedge v - b_0 \cdot \nabla d + d_0 \cdot \nabla b 
\end{aligned} \]

Therefore, $L_0(\xi) U$ corresponds exactly to the coordinates of the projection of $U$ on $E(0)$, expressed in the basis above. (More precisely, $L_0(\xi)$ has an image of dimension $5$, but there is a redundancy due to the fact that the Fourier transform of $\nabla \wedge v$ lives the space of dimension $2$ orthogonal to $\xi$.) In particular, there exists a matrix $Q(\xi)$, homogeneous of degree $0$ in $\xi$, such that $Q(\xi) L_0(\xi) = |\xi| P^0(\xi)$, and $Q$ is invertible. 
\end{Dem}

\begin{Def} We define~: 
\[ \mathcal{N}^{\epsilon_1, \epsilon_2 \epsilon_3}\left( \cdot, \cdot \right) = P^{\epsilon_1}(D) ~ \mathcal{N}\left(P^{\epsilon_2}(D) ~\cdot~, P^{\epsilon_3}(D) ~\cdot~ \right) \]
for any choice of $\epsilon_1, \epsilon_2, \epsilon_3 \in \{ +, 0, - \}$. 

Let us set
\[ \varphi^{\epsilon_1, \epsilon_2 \epsilon_3}(\xi, \eta) = \epsilon_1 |\xi|_0 - \epsilon_2 |\xi - \eta|_0 - \epsilon_3 |\eta|_0 \]

Finally, we define the space-resonant sets as 
\[ \mathcal{S}^{\epsilon_1, \epsilon_2 \epsilon_3} = \{ (\xi, \eta) \in \mathbb{R}^6, \nabla_{\eta} \varphi^{\epsilon_1, \epsilon_2 \epsilon_3}(\xi, \eta) = 0 \} \]
and the time-resonant sets as 
\[ \mathcal{T}^{\epsilon_1, \epsilon_2 \epsilon_3} = \{ (\xi, \eta) \in \mathbb{R}^6, \varphi^{\epsilon_1, \epsilon_2 \epsilon_3}(\xi, \eta) = 0 \} \]
and the space-time-resonant sets as 
\[ \mathcal{R}^{\epsilon_1, \epsilon_2 \epsilon_3} = \mathcal{S}^{\epsilon_1, \epsilon_2 \epsilon_3} \cap \mathcal{T}^{\epsilon_1, \epsilon_2 \epsilon_3} \]

Let us set
\[ u^{+} = P^{+} U, ~~~ u^{-} = P^{-} U, ~~~ u^{0} = P^{0} U \]
and
\[ f(t) = e^{t A_0(D)} U(t), ~~~ f^{\epsilon} = e^{t A_0(D)} u^{\epsilon} = P^{\epsilon} f = e^{t \epsilon \Lambda_0} u^{\epsilon} \]
Then the constraints can be written~: 
\[ u^0 = \Lambda^{-1} Q^{-1}(D) \mathcal{N}'(\Lambda U, U) \]
\label{def_nonlinearityN} 
\end{Def}

\begin{Rem}[Intuition] The choice of $f$ ensures that 
\[ \partial_t f = e^{t A_0(D)} \mathcal{N}(\Lambda e^{-t A_0(D)} f, e^{-t A_0(D)} f) \]
without linear part, so we can break the non-linearity down and consider uniquely terms of the form 
\[ e^{i t \epsilon_1 \Lambda_0 } \mathcal{N}^{\epsilon_1, \epsilon_2 \epsilon_3}(e^{-i t \epsilon_2 \Lambda_0} \Lambda f^{\epsilon_2}, e^{-i t \epsilon_3 \Lambda_0} f^{\epsilon_3}) \]
A Duhamel formula and a Fourier transform will make appear terms of the form~: 
\[ \int_1^t \int e^{i s \varphi(\xi, \eta)} b(\xi, \eta) \widehat{f}(s, \xi - \eta) |\xi - \eta| \widehat{f}(s, \eta) ~ d\eta ds \]
where we omit the superscripts $\epsilon$. Then, away from the space-resonant set, $\nabla_{\eta} \varphi \neq 0$ and we can perform an integration by parts in $\eta$ to win a factor $s^{-1}$. Away from the time-resonant set, $\varphi \neq 0$ so we can perfom an integration by parts in time and obtain a term of the form $\partial_s f$, so quadratic in $f$ and simpler to control. 

But we need, however, to control the nonlinearity close to the space-time-resonant set. 
\end{Rem}

\begin{Def} We say that a symbol $\mu_s : (\xi, \eta) \in \mathbb{R}^6 \mapsto \mu_s(\xi, \eta) \in \mathbb{R}$ is a symbol of order $s \in \mathbb{R}$ if $\mu_s$ is homogeneous of degree $s$, smooth outside of $\{ \xi = 0 \} \cup \{ \eta = 0 \} \cup \{ \xi = \eta \}$ and such that, if we write $(\xi_1, \xi_2, \xi_3)$ the three variables $(\xi, \eta, \xi - \eta)$ (in any order),  we have that $\mu_s(\xi, \eta) = \mathcal{A}\left( |\xi_1|, \frac{\xi_1}{|\xi_1|}, \xi_2 \right)$ for a smooth $\mathcal{A}$, as long as $|\xi_1| \ll |\xi_2|, |\xi_3| \sim 1$. (Note that $\xi_3$ is determined by $\xi_1$ and $\xi_2$.) 

We say that a nonlinearity satisfies the non-resonant condition of type $\epsilon_1, \epsilon_2 \epsilon_3$ (with $\epsilon_i \in \{ -, 0, + \}$) if it can be written under the form~: 
\[ \widehat{N(u, v)}(\xi) = \int b(\xi, \eta) \widehat{u}(\xi - \eta) \widehat{v}(\eta) ~ d\eta \]
where $b(\xi, \eta)$ is a symbol of order $0$ such that there exists a symbol of order $0$ denoted by $\mu_0$ such that  
\[ b(\xi, \eta) = \mu_0(\xi, \eta) \nabla_{\eta} \varphi^{\epsilon_1, \epsilon_2 \epsilon_3}(\xi, \eta) \]
(Note that this condition does not depend on $\epsilon_1$.) 
\end{Def}

\paragraph{Notation} We will write $\mu_0$ to denote any such symbol of order $0$. Since the precise expression of these symbols bears no importance whatsoever, we will denote them all with the same notation $\mu_0$ (but keeping in mind that they might actually differ) throughout the article. Let us also denote by $\mu_0^{n, m, l}$ for $n, m, l \in \mathbb{R}$ a symbol that can be written as
\[ \mu_0^{n, m, l}(\xi, \eta) = |\xi|^n |\eta|^m |\xi-\eta|^l \mu_0(\xi, \eta) \]
for some symbol $\mu_0$ of order $0$. Here again, we may use the same notation $\mu_0^{n, m, l}$ for two distinct symbols if we don't care about their precise expressions.

\begin{Lem} Let $F : S^2 \to F(S^2) \subset \mathbb{R}^3 \setminus \{ 0 \}$ be a $\mathcal{C}^{\infty}$-diffeomorphism. There exists $\mathcal{C}^{\infty}$ functions $m, m'$ (vector- resp. matrix-valued) such that for all $X, Y \in S^2$, 
\[ |F(X)| - |F(Y)| = m(X, Y) \cdot (X - Y), ~~~ X - Y = m'(X, Y) (F(X) - F(Y)) \]
\label{lem_factorisation_abstraite}
\end{Lem}

\begin{Dem}
We first prove the result locally. Set $G(X) = |F(X)|$, which is a smooth function. 

If $X \neq Y$, there exists a neighborhood on which $|X - Y|$ remains bounded from below by a strictly positive constant and we may set $m(X, Y) = (|F(X)| - |F(Y)|) \frac{X - Y}{|X - Y|}$ on this neighborhood. If $X = Y$, we use a local chart of $S^2$ to get back to open sets of $\mathbb{R}^2$, $F : U \to V$ a diffeomorphism, $G : U \to \mathbb{R}$ a smooth function, and thus~: 
\[ G(X) - G(Y) = (X - Y) \cdot \int_0^1 \nabla G (Y + t(X - Y)) ~ dt \]
and $\int_0^1 \nabla G(Y + t(X - Y)) ~ dt$ is well-defined and smooth, at least locally. 

We then use a partition of unity of $S^2 \times S^2$ to obtain $m$ on every piece. 

For the second statement, we procede the same~: if $X \neq Y$, $F(X) - F(Y)$ has at least one component that remains uniformly separated from $0$ locally and we may easily find $m'$~; if $X = Y$, we can use a local chart and use the same argument as before, noting that $F^{-1}$ satisfies the same properties. A partition of unity then allows to construct $m'$ on the whole $S^2 \times S^2$. 
\end{Dem}

From the previous lemma, we deduce~: 

\begin{Lem} There exists symbols $\mu_0, \mu_0'$ of order $0$ such that 
\[ \frac{|\eta|}{|\eta|_0} - \frac{|\xi - \eta|}{|\xi - \eta|_0} = \mu_0(\xi, \eta) \cdot \left( \frac{\eta}{|\eta|_0} - \frac{\xi - \eta}{|\xi - \eta|_0} \right) \]
and
\[ \frac{\eta}{|\eta|} - \frac{\xi - \eta}{|\xi - \eta|} = \mu_0'(\xi, \eta) \left( \frac{\eta}{|\eta|_0} - \frac{\xi - \eta}{|\xi - \eta|_0} \right) \]
\label{lem_factorisation}
\end{Lem}

\begin{Dem}
Let us first prove the second part. Set $F : X \in S^2 \mapsto \frac{X}{|X|_0} \in \mathbb{R}^3 \setminus \{ 0 \}$, which is a $\mathcal{C}^{\infty}$-diffeomorphism. By lemma \ref{lem_factorisation_abstraite}, we have $m' \in \mathcal{C}^{\infty}$ such that, for every $\eta, \xi$~: 
\[ \frac{\eta}{|\eta|} - \frac{\xi - \eta}{|\xi - \eta|} = m'\left( \frac{\eta}{|\eta|}, \frac{\xi - \eta}{|\xi - \eta|} \right) \left( \frac{\eta}{|\eta|_0} - \frac{\xi - \eta}{|\xi - \eta|_0} \right) \]
Set $\mu_0'(\xi, \eta) = m'\left( \frac{\eta}{|\eta|}, \frac{\xi - \eta}{|\xi - \eta|} \right)$ to get the desired result. 

The first part of lemma \ref{lem_factorisation_abstraite} then gives $m \in \mathcal{C}^{\infty}$ such that~: 
\[ \frac{|\eta|}{|\eta|_0} - \frac{|\xi - \eta|}{|\xi - \eta|_0} = m\left( \frac{\eta}{|\eta|}, \frac{\xi - \eta}{|\xi - \eta|} \right) \cdot \left( \frac{\eta}{|\eta|} - \frac{\xi - \eta}{|\xi - \eta|} \right) = m\left( \frac{\eta}{|\eta|}, \frac{\xi - \eta}{|\xi - \eta|} \right) \cdot m'\left( \frac{\eta}{|\eta|}, \frac{\xi - \eta}{|\xi - \eta|} \right) \left( \frac{\eta}{|\eta|_0} - \frac{\xi - \eta}{|\xi - \eta|_0} \right) \]
Set $\mu_0(\xi, \eta) = m\left( \frac{\eta}{|\eta|}, \frac{\xi - \eta}{|\xi - \eta|} \right) \cdot m'\left( \frac{\eta}{|\eta|}, \frac{\xi - \eta}{|\xi - \eta|} \right)$, which is a symbol of order $0$, to deduce the first part. 
\end{Dem}

\begin{Prop} The Born-Infeld system can be written under the form 
\[ \forall \epsilon_1 \in \{ -, 0, +\}, ~~~ \partial_t u^{\epsilon_1} + \epsilon_1 |D|_0 u^{\epsilon_1} = \sum_{\epsilon_2, \epsilon_3} \mathcal{N}^{\epsilon_1, \epsilon_2 \epsilon_3}(\Lambda u^{\epsilon_2}, u^{\epsilon_3}) \]
under the constraint
\[ u^0 = \Lambda^{-1} \sum_{\epsilon_2, \epsilon_3} \mathcal{N}^{', \epsilon_2 \epsilon_3}(\Lambda u^{\epsilon_2}, u^{\epsilon_3}) \]
Then $\mathcal{N}^{\epsilon_1, \epsilon_2 \epsilon_3}$ and $\mathcal{N}^{', \epsilon_2 \epsilon_3}$ satisfy the condition of non-resonance of type $\epsilon_1, \epsilon_2 \epsilon_3$ as long as $\epsilon_i \in \{ -, + \}$. 
\end{Prop}

\begin{Dem}
We already computed the projection operators and so we have formal expressions of $\mathcal{N}^{\epsilon_1, \epsilon_2 \epsilon_3}$ from definition \ref{def_nonlinearityN}. Each of these nonlinearities can be seen as tensors of size $10 \times 10 \times 10$ representing a bilinear operator $\mathbb{R}^{10} \times \mathbb{R}^{10} \to \mathbb{R}^{10}$. 

Furthermore, each entry of these tensors can be expressed as polynomials in the variables $\frac{\xi}{|\xi|}, \frac{\xi-\eta}{|\xi-\eta|}$, $\frac{\eta}{|\eta|}$, and $\alpha, \beta, \delta$ obtained from $\xi, \xi-\eta, \eta$ (see proposition \ref{def_alpha_beta_delta} for the definition of $\alpha, \beta, \delta$). Let us consider a fixed entry $n_0$ (from a nonlinearity satisfying the hypothesis of the proposition) and write it as
\[ n_0 \in \mathbb{R}[X_1, X_2, \dots, X_{18}] \]
where the only values of interest are when 
\[ X \in \mathcal{M} := \left\{ \left( \frac{\xi}{|\xi|}, \frac{\eta}{|\eta|}, \frac{\xi -\eta}{|\xi-\eta|}, \alpha(\xi), \beta(\xi), \delta(\xi), \alpha(\eta), \beta(\eta), \delta(\eta), \alpha(\xi-\eta), \beta(\xi-\eta), \delta(\xi-\eta) \right), ~ \xi, \eta \in \mathbb{R}^6 \setminus \{ 0 \}, \xi - \eta \neq 0 \right\} \]
Denote by 
\[ \iota(\xi, \eta) := \left( \frac{\xi}{|\xi|}, \frac{\eta}{|\eta|}, \frac{\xi -\eta}{|\xi-\eta|}, \alpha(\xi), \beta(\xi), \delta(\xi), \alpha(\eta), \beta(\eta), \delta(\eta), \alpha(\xi-\eta), \beta(\xi-\eta), \delta(\xi-\eta) \right) \]

Let us set
\[ \begin{aligned}
&P_1 = \left( X_1^2 + X_2^2 + X_3^2 - 1 \right), \quad P_2 = \left( X_4^2 + X_5^2 + X_6^2 - 1 \right), \quad P_3 = \left( X_7^2 + X_8^2 + X_9^2 - 1 \right), \\
&\quad P_4 = \left( X_{10}^2 + X_{11}^2 + X_{12}^2 - 1 \right), \quad P_5 = \left( X_{13}^2 + X_{14}^2 + X_{15}^2 - 1 \right), \quad P_6 = \left( X_{16}^2 + X_{17}^2 + X_{18}^2 - 1 \right), \\
&\quad P_7 = \left( X_4 \pm X_7 \right), \quad P_8 = \left( X_5 \pm X_8 \right), \quad P_9 = \left( X_6 \pm X_9 \right), \\
&\quad P_{10} = \left( X_{13} - X_{16} \right), \quad P_{11} = \left( X_{14} \pm X_{17} \right), \quad P_{12} = \left( X_{15} \pm X_{18} \right) 
\end{aligned} \]
where $\pm$ is a sign depending on $\epsilon_2, \epsilon_3$. Note that, when $X \in \mathcal{M}$, $X = \iota(\xi, \eta)$, then $P_i(X) = 0$ for all $i = 1, 2, ..., 6$ ; furthermore, up to choosing the sign correctly, $P_i(X)$ for $i = 7, 8, 9$ correspond to coordinates of $\nabla_{\eta} \varphi^{\epsilon_1, \epsilon_2 \epsilon_3}(\xi, \eta)$ ; and by lemma \ref{lem_factorisation}, $P_i(X)$ for $i = 10, 11, 12$ can also be factorized by $\nabla_{\eta} \varphi^{\epsilon_1, \epsilon_2 \epsilon_3}(\xi, \eta)$. 

This means that if $n_0$ can be written as 
\begin{equation} n_0 = \sum_{i = 1}^{12} Q_i P_i \label{prop_fact_ABI-fact-pol} \end{equation}
for some $Q_i \in \mathbb{R}[X]$, then one has 
\[ n_0(\iota(\xi, \eta)) = \sum_{i = 7}^{12} Q_i(\iota(\xi, \eta)) P_i(\iota(\xi, \eta)) = m(\xi, \eta) \nabla_{\eta} \varphi^{\epsilon_1, \epsilon_2 \epsilon_3}(\xi, \eta) \]
for some symbol $m$ of order $0$, thus proving the proposition. 

Although Euclidean division over $\mathbb{R}[X_1, ..., X_{18}]$ is generically ill-defined, one can here adapt it in a simple way given the expression of the $P_i$'s. Indeed, by basic Euclidean division, one can always write 
\[ n_0 = R_1 + \sum_{i = 7}^{12} Q_i P_i \]
for $R_1 \in \mathbb{R}[X]$ such that $R$ does not depend on $X_7, X_8, X_9, X_{16}, X_{17}, X_{18}$. Then, one can decompose again 
\[ R_1 = R_2 + \sum_{i = 1}^6 Q_i P_i \]
for $R_2$ having only monomials where $X_3, X_6, X_9, X_{12}, X_{15}$ are of order $0$ or $1$. 

We claim that, for our specific choice of $n_0$, $R_2 = 0$. 

Since the computations are cumbersome (one needs to compute composition of tensors of size $10 \times 10$ and $10 \times 10 \times 10$, and this for each of the cases covered by the proposition), we prove this cancelation by running the following algorithm on a computer : we first compute each entry as a formal polynomial (note that each entry, in this setting, will be an integer), then replace formally $X_1^2 + X_2^2 + X_3^2$ by $1$, $X_7$ by $\mp X_4$, ... and so on whenever it is possible, and then check that the polynomial obtained that way is zero. Note that there is no hardware error since we only manipulate integer coefficients. 
\end{Dem}

\begin{Prop}[Quasi-linear structure] The equation \eqref{equ_ABI} is a symmetric quasilinear hyperbolic system. 

More precisely, we may write equivalently~: 
\[ A_0(D) U + \mathcal{N}(\Lambda U, U) = \sum_{i = 1}^3 (M_{0i} + M_{1i}(U)) \partial_i U \]
where $M_{0i}, M_{1i}(U)$ are symmetric matrices for any $i = 1, 2, 3$ and any $U$, and $M_{1i}(U)$ is linear in $U$. 
\label{prop_quasi_lin}
\end{Prop}

\begin{Dem}
This is immediate from the system written explicitly \eqref{ABI_Brenier}. 
\end{Dem}

\subsubsection{Chaplygin gas} \label{section_structure_BI_Chap}

\paragraph{} The Chaplygin gas is a gas model for which the pression is determined by 
\[ p = -\frac{A}{\rho} \]
where $\rho$ is the density and $A$ a constant. Up to a change of units, we fix $A = 1$. The equations can therefore be written 
\[ \left\{ \begin{array}{l}
\partial_t \rho + \nabla \cdot (\rho v) = 0 \\
\partial_t v + v \cdot \nabla v = \frac{1}{\rho} \nabla \frac{1}{\rho} \end{array} \right. \]
If we set $\tau = \frac{1}{\rho}$, we get 
\[ \left\{ \begin{array}{l}
\partial_t \tau + v \cdot \nabla \tau - \tau \nabla \cdot v = 0 \\
\partial_t v + v \cdot v = \tau \nabla \tau
\end{array} \right. \]
and we recognize exactly the augmented Born-Infeld system where we set $b \equiv d \equiv 0$, which is an initial condition preserved by \eqref{equ_ABI}. Then, studying the Chaplygin gas around any constant solution comes back to studying the augmented Born-Infeld system around $(\tau_0, v_0, 0, 0)$, enforcing $b \equiv d \equiv 0$. 

In the constraints \eqref{equ_cont}, the divergence equations on $b, d$ are automatically satisfied. The final constraint is 
\[ \tau_0 \nabla \wedge v + \tau \nabla \wedge v = 0 \]
Therefore, if we assume we are in the case of irrotational velocity fields, $\nabla \wedge v \equiv 0$ and this constraint is satisfied as well. Furthermore, $\nabla \wedge v \equiv 0$ is also preserved by the equation~: 
\[ \partial_t \nabla \wedge v = - \nabla \wedge (v \cdot \nabla v + \tau \nabla \tau) = - v \cdot \nabla (\nabla \wedge v) - (\nabla \cdot v) \nabla \wedge v + (\nabla \wedge v) \cdot \nabla v \]

So the same structure properties can be applied directly to the Chaplygin gas with irrotational velocity fields. 

\subsubsection{General case} \label{section_structure_ABI}

We now consider the general augmented Born-Infeld system \eqref{equ_ABI} under the constraint equations \eqref{equ_cont}. 

\begin{Prop} \eqref{equ_cont} is preserved under the flow of \eqref{equ_ABI}. 
\end{Prop}

\begin{Dem}
\textbf{Case of the divergence constraints} We differentiate~: 
\[ \begin{aligned}
\partial_t \left( \tau \nabla \cdot b - b \cdot \nabla \tau \right)
&= - (v \cdot \nabla \tau) (\nabla \cdot b) + \tau (\nabla \cdot v) (\nabla \cdot b) + \tau \nabla \cdot \left( - v \cdot \nabla b + b \cdot \nabla v - \tau \nabla \wedge d \right) \\
&+ v \cdot \nabla b \cdot \nabla \tau - b \cdot \nabla v \cdot \nabla \tau + \tau \nabla \wedge d \cdot \nabla \tau + b \cdot \nabla \left( v \cdot \nabla \tau - \tau \nabla \cdot v \right) \\
&= - v \cdot \nabla \left( \tau \nabla \cdot b - b \cdot \nabla \tau \right) + (\nabla \cdot v) (\tau \nabla \cdot b - b \cdot \nabla \tau)
\end{aligned} \]
In particular, if this quantity vanishes at the initial time, it remains zero for any time as long as the solution is smooth. Likewise, we can obtain a similar equation for the divergence of $d$. 

\textbf{Case of the rotational constraint} We compute 
\[ \begin{aligned}
\partial_t \left( \tau \nabla \wedge v - b \cdot \nabla d + d \cdot \nabla b \right)
&= - (v \cdot \nabla \tau) \nabla \wedge v + \tau (\nabla \cdot v) \nabla \wedge v + \tau \nabla \wedge \left( - v \cdot \nabla v + b \cdot \nabla b + d \cdot \nabla d \right) \\
&+ v \cdot \nabla b \cdot \nabla d - b \cdot \nabla v \cdot \nabla d + \tau \nabla \wedge d \cdot \nabla d - b \cdot \nabla \left( - v \cdot \nabla d + d \cdot \nabla v + \tau \nabla \wedge b \right) \\
&- v \cdot \nabla d \cdot \nabla b + d \cdot \nabla v \cdot \nabla b + \tau \nabla \wedge b \cdot \nabla b + d \cdot \nabla \left( - v \cdot \nabla b + b \cdot \nabla v - \tau \nabla \wedge d \right) \\
&= - v \cdot \nabla \left( \tau \nabla \wedge v - b \cdot \nabla d + d \cdot \nabla b \right) + (\tau \nabla \wedge v - b \cdot \nabla d + d \cdot \nabla b) \cdot \nabla v \\
&+ (\nabla \wedge b) \left( \tau \nabla \cdot b - b \cdot \nabla \tau \right) + (\nabla \wedge d) \left( \tau \nabla \cdot d - d \cdot \nabla \tau \right)
\end{aligned} \]
In particular, if the divergence constraints are satisfied, if this quantity vanishes at the initial time, then it remains zero for any time when the solutions are smooth. 
\end{Dem}

\begin{Def} Given a constant state $U_0$ and an initial perturbation $u : \mathbb{R}^3 \to \mathbb{R}^{10}$, we say that $u = (\tau, v, b, d)$ is an \textit{admissible initial data} if $\tau + \tau_0 > 0$ on all $\mathbb{R}^3$, and $u$ satisfies the constraint equations \eqref{equ_cont}. 
\end{Def}

\section{Main result} \label{section_main_result}

From now on, we consider the equation \eqref{equ_ABI} + \eqref{equ_cont}, so the augmented Born-Infeld system around a constant solution, assuming the initial data satisfies the constraint equations. 

\begin{Theo}
Let $U_0$ be a constant state with $\tau_0 > 0$. There exists an integer $N$ and constants $\delta, C > 0$ depending continuously on $U_0$, such that for any admissible initial data satisfying
\[ \Vert u_0 \Vert_{H^N} \leq \delta, ~~~ \Vert \langle x \rangle^2 u_0 \Vert_{H^6} \leq \delta \]
that we write $\Vert u_0 \Vert_{X_0} \leq \delta$, then the solution is global and 
\[ \begin{aligned}
& \Vert u(t) \Vert_{W^{1, \infty}} \lesssim C \Vert u_0 \Vert_{X_0} \langle t \rangle^{-1}
\end{aligned} \]
\end{Theo}

\begin{Rem}[Scattering of the solution] Even though the decay of the $L^{\infty}$ norm is only $t^{-1}$, our method of proof will immediately imply scattering~: indeed, in section \ref{section_estH6}, we will show $\Vert f \Vert_{H^6} \lesssim \Vert u_0 \Vert_{X_0}$ by proving that $\Vert \partial_t f \Vert_{H^6} \lesssim \langle t \rangle^{-1-c}$ for some $c > 0$, which also ensures that $f(t)$ has a limit in $H^6$ when $t \to \infty$. 
\end{Rem}

\subsection{Statement of the a priori estimate} \label{section_statement_apriori}

We now consider the system introduced in section \ref{section_structure}, with the constraints preserved by the evolution, and the decomposition $u^{+} + u^{-} + u^0$. The system is symmetric, so hyperbolic, and we have the existence and uniqueness of local in time smooth solutions. We denote by $u$ the solution and $f = e^{t A_0(D)} U$, with the same decomposition on the eigenspaces in Fourier. 

Furthermore, without any loss of generality, we assume that the initial time is $t = 1$ in place of $t = 0$, and we try to extend the solution for any time $t \geq 1$. 

If $T > 1$ is an existence time, we can introduce the norm  
\begin{equation} \begin{aligned}
\Vert u \Vert_{X, T} = \sup_{1 \leq t \leq T} &\left\{ \Vert u^{\pm} \Vert_{H^6}, t^{-\varepsilon} \Vert u^{\pm} \Vert_{H^N}, t \Vert u^{\pm} \Vert_{\dot{B}^1_{\infty, 1}}, t \Vert u^{\pm} \Vert_{\dot{B}^0_{\infty, 1}}, \Vert x f^{\pm} \Vert_{H^5}, t^{-b} \Vert x f^{\pm} \Vert_{H^6}, \sup_{0 \leq k \leq 5} t^{-\gamma_k} \Vert \Lambda^k |x|^2 \Lambda f \Vert_{L^2}, \right. \\
&\left. ~~~~~ t \Vert u^0 \Vert_{H^7}, t^{1-\varepsilon} \Vert u^0 \Vert_{H^N}, t^{2-a} \Vert u^0 \Vert_{W^{1, \infty-}}, t^{1-a} \Vert \Lambda x u^0 \Vert_{H^4}, t^{1-a-\gamma_5/2} \Vert \Lambda x u^0 \Vert_{H^6} \right\} 
\end{aligned} \end{equation}
for $\varepsilon, a, b, \gamma_k > 0$ small enough, $N$ a large enough integer, $\infty-$ a large enough number. 

Recall that $\dot{B}^s_{p, q}$ for $p, q \in [1, \infty], s \in \mathbb{R}$ is the homogeneous Besov space, with norm~: 
\[ \Vert v \Vert_{\dot{B}^s_{p, q}} := \left( \sum_{j \in \mathbb{Z}} 2^{sj} \Vert \varphi_j(D) v \Vert_{L^p}^q \right)^{1/q} \]
where $\varphi_j(\xi) = \varphi(2^{-j} \xi)$ is an appropriate Littlewood-Paley localisation. In particular, if $v \in L^2$, then $\Vert v \Vert_{W^{1, \infty}} \lesssim \Vert v \Vert_{\dot{B}^1_{\infty, 1}} + \Vert v \Vert_{\dot{B}^0_{\infty, 1}}$. Recall also that $\dot{B}^0_{2, 2}$ and $L^2$ are the same spaces. See \cite{BahouriFourierbook} for more details on the construction of Besov spaces. 

\begin{Rem}[Choice of the parameters] Here above, the $X$-norm depends on various parameters, that need to be chosen in a particular way for the long-time control of the $X$-norm to hold. Let us explain how they have to depend on each other. 

We first fix $\gamma_k = \frac{k \gamma_5 + (5-k) \gamma_0}{5}$ so that only $\gamma_0$ and $\gamma_5$ are free parameters. We then choose $(\gamma_0, \gamma_5)$ small (with respect to $1$) and satisfying $5 \gamma_0 < \gamma_5$. We then choose $b > 0$ such that $\frac{\gamma_4}{2} < b \leq \frac{\gamma_5}{2}$. Eventually, we choose $a > 0$ small enough with respect to the previous parameters, $N$ large enough with respect to the previous parameters, $\varepsilon > 0$ small enough with respect to the previous parameters, and $\infty-$ large enough with respect to the previous parameters (in this order), with in particular $\frac{\gamma_4 + \varepsilon}{2} \leq b < \frac{\gamma_5 + \varepsilon}{2}$. 

None of these parameters will depend on the $(\tau_0, v_0, b_0, d_0)$ constant state. 
\end{Rem}

\begin{Rem}[Position of the derivatives] Note that the position of the derivatives $\Lambda$ plays no role in the weighted $L^2$ estimates. Indeed, one has the commutation identity for any function $g$ : 
\[ \Lambda x g = x \Lambda g - R g \]
for $R$ a Riesz transform, so that only lower-order terms appear when commuting $\Lambda^j$ and $x$. This means that $\Lambda^j x f^{\epsilon} = x \Lambda^j f^{\epsilon} + \mbox{remainder}$ where the remainder have better decay ; the same holds for the $x^2$-weighted estimate. 
\end{Rem}

The a priori estimate we will prove will be~: 

\begin{Prop}[A priori estimate] If the parameters of the $X$-norm are well-chosen, 
\[ \Vert u \Vert_{X, T} \leq C \Vert u_0 \Vert_{X} + C \Vert u \Vert_{X, T}^{3/2} \left( 1 + \Vert u \Vert_{X, T}^3 \right) \]
where $C$ is independant of the size of the data, but may depend on $(\tau_0, v_0, b_0, d_0)$ and on the parameters of the $X$-norm. 
\end{Prop}

Therefore, if we apply the smallness hypothesis on $u_0$ and we choose $T > 0$ maximal such that $\Vert u \Vert_{X, T} \leq 2 C \delta$, with $0 < \delta < 1$, we have that
\[ \Vert u \Vert_{X, T} \leq C \delta (1 + 2^{3/2} C^{3/2} \delta^{1/2} + 2^{9/2} C^{9/2} \delta^{7/2} ) \]
In particular, if $\delta$ is chosen small enough, 
\[ \Vert u \Vert_{X, T} < 2 C \delta \]
But if $T < \infty$, we could extend the solution and this would contradict the maximality of $T$. Therefore, there exists a global solution that remains close to the constant solution in $X$-norm. 

In the following, we only write $\Vert \cdot \Vert_X$ and omit the dependance in $T$. We will write $\gamma = \gamma_5$, which controls all the other $\gamma_k$. 

\paragraph{Notation} Recall that we use the notation $\Lambda = |\nabla|$, ie the linear operator of Fourier symbol $|\xi|$. The Hörmander-Mikhlin theorem implies that it is equivalent to $\Lambda_0$ of symbol $|\xi|_0$ (and their quotient is a symbol of order $0$ with our definitions), in the sense that, for any $1 < p < \infty$, $\Vert \Lambda v \Vert_{L^p} \lesssim \Vert \Lambda_0 v \Vert_{L^p} \lesssim \Vert \Lambda v \Vert_{L^p}$ for all $v$. Furthermore, we will write $U = u^{+} + u^{-} + u^0$, and sometimes simply $u$ in place of $u^{+}$ or $u^{-}$ (but not $u^0$). $\mu_0$ will name a generic symbol of order $0$ and will be authorized to vary at each line. 

\section{General inequalities} \label{section_geninequ}

In this section, we prove useful lemmas for the a priori estimates and general inequalities. 

\subsection{Functional inequalities} \label{section_funcinequ}

\begin{Lem}[Hardy's inequality] We have the following $L^2$ Hardy inequalities~:  
\[ \Vert g \Vert_{L^2} \lesssim \Vert \nabla (|x| g) \Vert_{L^2} ~~~~~~~~~~ \mbox{ and } ~~~~~~~~~~ \Vert g \Vert_{L^2} \lesssim \Vert x \Lambda g \Vert_{L^2} \] \label{lem_Hardy}
\end{Lem}

\begin{Dem}
The second inequality comes from the first by applying Parseval's inequality~: 
\[ \Vert g \Vert_{L^2} = \Vert \widehat{g} \Vert_{L^2} \lesssim \Vert \nabla ( |\xi| \widehat{g} ) \Vert_{L^2} = \Vert x \Lambda g \Vert_{L^2} \]
For the first inequality, we use a polar decomposition and set $h(x) = |x| g(x)$. Assuming $g$ is a $C^{\infty}$ function with compact support, we compute~: 
\[ \begin{aligned}
\int |g(x)|^2 ~ dx &= \int_0^{\infty} \int_{S^2} \frac{|h(ry)|^2}{r^2} r^2 ~ d\sigma(y) dr = - \int_0^{\infty} \int_{S^2} \int_r^{\infty} 2 \frac{\partial h}{\partial r}(sy) h(sy) ~ ds d\sigma(y) dr \\
&= - \int_0^{\infty} \int_{S^2} 2 r \frac{\partial h}{\partial r}(ry) h(ry) ~ d\sigma(y) dr = - \int_{\mathbb{R}^3} 2 \frac{\partial h}{\partial r}(x) \frac{h(x)}{|x|} ~ dx \\
&\leq 2 \Vert \partial_r h \Vert_{L^2} \Vert g \Vert_{L^2} \lesssim \Vert g \Vert_{L^2} \Vert \nabla (|x| g) \Vert_{L^2}
\end{aligned} \]
by Hölder's inequality. So we get~:  
\[ \Vert g \Vert_{L^2} \lesssim \Vert \nabla (|x| g) \Vert_{L^2} \]

We conclude by a density argument. 
\end{Dem}

%

\begin{Lem}[Moser estimate] Let $g, h$ be two functions and $\alpha \in \mathbb{N}^3$ with $|\alpha| = k \in \mathbb{N}$. Then~: 
\[ \Vert D^{\alpha} (g h) - g D^{\alpha} h \Vert_{L^2} \lesssim \Vert g \Vert_{H^k} \Vert h \Vert_{L^{\infty}} + \Vert \nabla g \Vert_{L^{\infty}} \Vert h \Vert_{H^{k-1}} \]
\label{Moser_est}
\end{Lem}

See for instance \cite{Taylor_PDE3}, Proposition 3.7. The next lemma shows which norms are controlled by $\Vert u \Vert_X$. 

\begin{Lem} For any integer $k < N$ and any $2 < q < \infty-$, there exists $\delta = \delta(N, \varepsilon, k, \infty-, a, q)$ such that
\[ \Vert u^{\pm} \Vert_{W^{k, q}} \lesssim t^{-1+2/q+\delta} \Vert u \Vert_X, ~~~ \Vert u^0 \Vert_{W^{k, q}} \lesssim t^{-2+2/q+\delta} \Vert u \Vert_X \]
with $\delta$ going to $0$, when $k, q$ are fixed, if $N \to \infty$, $\varepsilon \to 0$, $a \to 0$. 
\label{lem_infini}
\end{Lem}

The demonstration can be found in \cite{GMSWaterwaves}, lemma 5.1, with a slight variation for $u^0$ (since we do not control $\Vert u^0 \Vert_{L^{\infty}}$ but only the $L^{\infty-}$ norm). 

\begin{Lem} We have that 
\[ \Vert \Lambda x g \Vert_{L^2} \lesssim \Vert g \Vert_{H^2}^{1/2} \Vert \langle x \rangle^2 g \Vert_{L^2}^{1/2} \]
In particular, 
\[ \Vert x f^{\pm} \Vert_{\dot{H}^7} \lesssim t^{\gamma/2+\varepsilon/2} \Vert u \Vert_X \] \label{lem_xH7}
\end{Lem}

\begin{Dem}
We write that 
\[ \Vert \Lambda x g \Vert_{L^2} \lesssim \Vert |x| \nabla g \Vert_{L^2} + \Vert g \Vert_{L^2} \]
Then
\[ \Vert |x| \nabla g \Vert_{L^2}^2 = \int |x|^2 \nabla g \cdot \overline{\nabla g} ~ dx = - \int 2 g x \cdot \overline{\nabla g} ~ dx - \int |x|^2 g \nabla \cdot \overline{\nabla g} ~ dx \lesssim \Vert x g \Vert_{L^2} \Vert \nabla g \Vert_{L^2} + \Vert |x|^2 g \Vert_{L^2} \Vert \Lambda^2 g \Vert_{L^2} \]
So summing each contribution we have that
\[ \Vert \Lambda x g \Vert_{L^2} \lesssim \Vert \langle x \rangle^2 g \Vert_{L^2}^{1/2} \Vert g \Vert_{H^2}^{1/2} \]

By setting $g = \Lambda^6 f^{\pm}$, we get
\[ \Vert x f^{\pm} \Vert_{\dot{H}^7} \lesssim \Vert \Lambda x g \Vert_{L^2} + \Vert f \Vert_{H^6} \lesssim \Vert \langle x \rangle^2 g \Vert_{L^2}^{1/2} \Vert g \Vert_{H^2}^{1/2} + \Vert u \Vert_X \lesssim t^{\gamma/2+\varepsilon/2} \Vert u \Vert_X \]
as desired. 
\end{Dem}

\begin{Lem}[Dispersion inequality for the wave equation] Let $2 \leq p \leq \infty$ and $p'$ be its conjugate exponent. We have~: 
\[ 
\Vert e^{i t \Lambda} g \Vert_{L^p} \lesssim t^{-1+2/p} \Vert \Lambda^{2 - 4/p} g \Vert_{L^{p'}}
\] \label{lem_disp}
\end{Lem}

The previous theorem can be found in \cite{Shatah_geometricwave}. The following two lemmas come from \cite{ShatahPusateri} (appendix). 

\begin{Lem}[Inequality $L^1 - L^2$] We have that
\[ \Vert g \Vert_{L^1} \lesssim \Vert |x| g \Vert_{L^2}^{1/2} \Vert |x|^2 g \Vert_{L^2}^{1/2} \] \label{lem_L1L2}
\end{Lem}

\begin{Lem}[Fractional integration] For any $\alpha > 0$ and any $1 < p, q < \infty$ such that $\alpha = \frac{3}{p} - \frac{3}{q}$, we have
\[ \Vert g \Vert_{L^q} \lesssim \Vert \Lambda^{\alpha} g \Vert_{L^p} \]
If furthermore $p \leq 2 \leq q$, 
\[ \Vert e^{i t \Lambda} g \Vert_{L^q} \lesssim \Vert \Lambda^{\alpha} g \Vert_{L^p} \] \label{lem_intfrac}
\end{Lem}

Finally, let us recall the following lemma (see \cite{BahouriFourierbook}, Proposition 2.22)~: 

\begin{Lem}[Interpolation of Besov spaces] Let $\kappa > 0$. There exists $\theta = \theta(\kappa) \in (0, 1)$ such that, for any $g$~: 
\[ \Vert g \Vert_{\dot{B}^0_{2, 1}} \lesssim \Vert \Lambda^{-\kappa} g \Vert_{L^2}^{\theta} \Vert \Lambda^{1/2} g \Vert_{L^2}^{1-\theta} \] \label{interp_Besov}
\end{Lem}

\subsection{Identities for the wave symbol} \label{section_waveid}

\begin{Lem} Let us consider  
\[ \varphi(\xi, \eta) = \varphi^{\epsilon_1, \epsilon_2 \epsilon_3}(\xi, \eta) = \epsilon_1 |\xi|_0 - \epsilon_2 |\xi - \eta|_0 - \epsilon_3 |\eta|_0 \]
with $\epsilon_i \in \{ +, - \}$, $i = 1, 2, 3$. Denote by $g_0$ the matrix associated to the norm $| \cdot |_0$, that is~: 
\[ |\xi|_0^2 = \xi \cdot g_0 \xi, ~~~ g_0 = \tau_0^2 I_3 + b_0 \otimes b_0 + d_0 \otimes d_0 \]
Then~: 
\[ \epsilon_1 |\xi|_0 \nabla_{\xi} \varphi(\xi, \eta) = - \epsilon_3 |\eta|_0 \nabla_{\eta} \varphi(\xi, \eta) - \epsilon_2 g_0 \frac{\xi - \eta}{|\xi - \eta|_0} \varphi(\xi, \eta) \] \label{identite_fond}
\end{Lem}

\begin{Dem}
We compute~:  
\[ \begin{aligned}
\epsilon_1 |\xi|_0 \nabla_{\xi} \varphi(\xi, \eta) &= g_0 \xi - \epsilon_1 \epsilon_2 \frac{|\xi|_0}{|\xi - \eta|_0} g_0 (\xi - \eta) \\
&= - \epsilon_2 \frac{g_0 (\xi - \eta)}{|\xi - \eta|_0} \left( \epsilon_1 |\xi|_0 - \epsilon_2 |\xi - \eta|_0 - \epsilon_3 |\eta|_0 \right) - g_0 (\xi - \eta) + g_0 \xi - \epsilon_2 \epsilon_3 \frac{|\eta|_0}{|\xi - \eta|_0} g_0 (\xi - \eta) \\
&= - \epsilon_2 \frac{g_0 (\xi - \eta)}{|\xi - \eta|_0} \varphi(\xi, \eta) - \epsilon_3 |\eta|_0 \nabla_{\eta} \varphi(\xi, \eta)
\end{aligned} \]

Hence the result. 
\end{Dem}

\begin{Lem} Let $g_0$ be the matrix associated to the norm $|\cdot |_0$ as in the previous lemma, and $|\cdot |_{0'}$ the norm associated to $g_0^{-1}$. The following identities hold~: 
\[ \begin{aligned}
&\varphi^{+, ++}(\xi, \eta) = - \frac{|\xi - \eta|_0 |\eta|_0}{|\xi|_0 + |\eta|_0 + |\xi - \eta|_0} |\nabla_{\eta} \varphi^{+, ++}(\xi, \eta)|_{0'}^2 \\
&2 \varphi^{+, +-}(\xi, \eta) = \frac{(|\xi|_0 + |\eta|_0 + |\xi - \eta|_0)|\xi - \eta|_0 |\eta|_0 }{|\xi|_0 |\xi - \eta|_0 + \xi \cdot g_0 (\xi - \eta)} |\nabla_{\eta} \varphi^{+, +-}(\xi, \eta)|_{0'}^2 \\
&2 \varphi^{+, -+}(\xi, \eta) = \frac{(|\xi|_0 + |\eta|_0 + |\xi - \eta|_0) |\xi - \eta|_0 |\eta|_0}{|\xi|_0 |\eta|_0 + \xi \cdot g_0 \eta} |\nabla_{\eta} \varphi^{+, -+}(\xi, \eta)|_{0'}^2 \\
&2 \varphi^{-, +-}(\xi, \eta) = -\frac{(|\xi|_0 + |\eta|_0 + |\xi - \eta|_0) |\xi - \eta|_0 |\eta|_0}{|\xi|_0 |\eta|_0 + \xi \cdot g_0 \eta} |\nabla_{\eta} \varphi^{-, +-}(\xi, \eta)|_{0'}^2 \\
&2 \varphi^{-, -+}(\xi, \eta) = - \frac{(|\xi|_0 + |\eta|_0 + |\xi - \eta|_0)|\eta|_0 |\xi - \eta|_0}{|\xi|_0 |\xi - \eta|_0 + \xi \cdot g_0 (\xi - \eta)} |\nabla_{\eta} \varphi^{-, -+}(\xi, \eta)|_{0'}^2 \\
&\varphi^{-, --}(\xi, \eta) = \frac{|\xi - \eta|_0 |\eta|_0}{|\xi|_0 + |\eta|_0 + |\xi - \eta|_0} |\nabla_{\eta} \varphi^{-, --}(\xi, \eta)|_{0'}^2
\end{aligned} \]
at any point where the denominator does not vanish. 
\label{identite_Linf}
\end{Lem}

\begin{Lem} The following identities hold~: 
\[ \begin{aligned}
&|\xi|_0 \nabla_{\xi} \varphi^{+, +0}(\xi, \eta) = \left( |\xi|_0 - |\xi - \eta|_0 \right) \nabla_{\eta} \varphi^{+, +0}(\xi, \eta) + g_0 \eta \\
&|\xi|_0 \nabla_{\xi} \varphi^{-, -0}(\xi, \eta) = \left( |\xi|_0 - |\xi - \eta|_0 \right) \nabla_{\eta} \varphi^{-, -0}(\xi, \eta) - g_0 \eta
\end{aligned} \] \label{identite_fond_0}
\end{Lem}

%

The proof of Lemmas \ref{identite_Linf} and \ref{identite_fond_0} is a straightforward computation. 

\subsection{Inequalities with a symbol} \label{section_symbinequ}

Let $\mu$ be a symbol. We denote by $T_{\mu}$ the operator defined by~: 
\[ \widehat{T_{\mu}(g, h)}(\xi) = \int \mu(\xi, \eta) \widehat{g}(\xi - \eta) \widehat{h}(\eta) ~ d\eta \]

\begin{Lem}[Symbols] Let $1 < p, q, r < \infty$ be such that $\frac{1}{p} + \frac{1}{q} = \frac{1}{r}$, and $\mu_0$ a symbol of order $0$. Then
\[ \Vert T_{\mu_0}(g, h) \Vert_{L^r} \lesssim \Vert g \Vert_{L^p} \Vert h \Vert_{L^q} \]
If $M > 3$, we also have
\[ \Vert T_{\mu_0}(g, h) \Vert_{L^2} \lesssim \Vert g \Vert_{L^2} \Vert h \Vert_{W^{1, M}} \] \label{lem_symb}
\end{Lem}

Again, see the appendix of \cite{ShatahPusateri}. Note that the second inequality here above does not quite match the one given by \cite{ShatahPusateri} ; however, it will follow from the embedding $W^{1, M} \hookrightarrow \dot{B}_{\infty, 1}^0$ when $M > 3$ and the next lemma, which is a refinement involving Besov spaces needed in our analysis~: 

\begin{Lem} Let $\mu_0$ be a symbol of order $0$. For any $f \in L^2, g \in \dot{B}_{\infty, 1}^0$, we have~: 
\[ \Vert T_{\mu_0}(f, g) \Vert_{L^2} \lesssim \Vert f \Vert_{L^2} \Vert g \Vert_{\dot{B}_{\infty, 1}^0} \]
\label{lem_symb_Besov}
\end{Lem}

\begin{Dem}
We follow the proof of lemma \ref{lem_symb} in \cite{ShatahPusateri}, adapting only at few steps where estimates fail on $L^{\infty}$ but not on $\dot{B}_{\infty, 1}^0$. 

Away from the coordinate axes $\{ \xi = 0 \} \cup \{ \eta = 0 \} \cup \{ \xi - \eta = 0 \}$, the Coifman-Meyer theorem \cite{CoifmanMeyer} applies to give a control by $\Vert f \Vert_{L^2} \Vert g \Vert_{L^{\infty}}$, and the $\dot{B}_{\infty, 1}^0$-norm is stronger than the $L^{\infty}$-norm on the intersection of these spaces. 

Consider now the case $|\eta| \ll |\xi - \eta|, |\xi| \sim 1$. Notice that the case $|\xi - \eta| \ll |\eta|, |\xi| \sim 1$ is symmetric and the case $|\xi| \ll |\eta|, |\xi - \eta|$ can be reduced to the previous ones by duality. This allows to consider only~: 
\[ \sum_j T_{\mu_0}(P_{< j - 100} f, P_{j} g) \]
where $P_{< k}$ denotes the Littlewood-Paley projection on low frequencies of order $2^{k}$ and $P_k = P_{<k+1} - P_{<k}$. The definition of the class of symbols considered here ensures that, close to the axis $\{ \eta = 0 \}$, we can write $\mu_0(\xi, \eta) = \mathcal{A}\left( |\eta|, \frac{\eta}{|\eta|}, \xi \right)$ for some smooth $\mathcal{A}$. By homogeneity, we write $\mu_0(\xi, \eta) = \mathcal{A}\left( \frac{|\eta|}{|\xi|}, \frac{\eta}{|\eta|}, \frac{\xi}{|\xi|} \right)$ and expand in $\frac{|\eta|}{|\xi|}$~: 
\[ \mu_0(\xi, \eta) = \sum_{k = 0}^L \frac{|\eta|^k}{|\xi|^k} m_k\left( \frac{\eta}{|\eta|}, \frac{\xi}{|\xi} \right) + ~~~ \mbox{remainder} \]
where the $m_k$ are smooth. If $L$ is large enough, the singularity of the remainder at $\eta = 0$ becomes weak enough to allow the use of the Coifman-Meyer theorem and we may focus on the finite sum in $k$. 

We then expand each $m_k$ in spherical harmonics~: 
\[ \mu_0(\xi, \eta) = \sum_{k = 0}^L \sum_{l, l'} a_{k, l, l'} \frac{|\eta|^k}{|\xi|^k} Z_l\left( \frac{\eta}{|\eta|} \right) Z_{l'} \left( \frac{\xi}{|\xi|} \right) \]
The spherical harmonics $Z_l$ are bounded on $L^2$, and on $\dot{B}_{\infty, 1}^0$ by Bernstein's lemma. Furthermore, the bounds grow polynomially in $l, l'$ while the $a_{k, l, l'}$ decay faster than any polynomial by smoothness of $m_k$. Therefore, up to adding a constant eventually, we may disregard the summation in $l, l'$ in what follows. Since the summation in $k$ is finite, we also only consider the case of only one fixed $k$. But then~: 
\[ \left\Vert \sum_j \Lambda^{-k} \left( \left( P_{<j-100} \Lambda^k f \right) \left( P_{j} g \right) \right) \right\Vert_{L^2} \lesssim \left\Vert \left( \sum_j 2^{-2jk} \left( P_{<j-100} \Lambda^k f \right)^2 \left( P_{j} g \right)^2 \right)^{1/2} \right\Vert_{L^2} \]
by the Littlewood-Paley square function estimate. Let us now distinguish two cases, depending on whether we prefer having $f$ in $L^2$ and $g$ in $\dot{B}_{\infty, 1}^0$ or the other way around. In the first case, we write by Hölder's inequality~: 
\[ \left\Vert \left( \sum_j 2^{-2jk} \left( P_{<j-100} \Lambda^k f \right)^2 \left( P_{j} g \right)^2 \right)^{1/2} \right\Vert_{L^2} \lesssim \left\Vert \sup_j \left| 2^{-jk} P_{<j-100} \Lambda^k f \right| \right\Vert_{L^2} \left\Vert \left( \sum_j (P_{j} g )^2 \right)^{1/2} \right\Vert_{L^{\infty}} \]
By the Littlewood-Paley maximal function estimate, 
\[ \left\Vert \sup_j \left| 2^{-jk} P_{<j-100} \Lambda^k f \right| \right\Vert_{L^2} \lesssim \Vert f \Vert_{L^2} \]
On the other hand, 
\[ \left\Vert \left( \sum_j (P_{j} g )^2 \right)^{1/2} \right\Vert_{L^{\infty}} \leq \left( \sum_j \Vert P_{j} g \Vert_{L^{\infty}}^2 \right)^{1/2} = \Vert g \Vert_{\dot{B}_{\infty, 2}^0} \lesssim \Vert g \Vert_{\dot{B}_{\infty, 1}^0} \]
In the second case, we apply again Hölder's inequality~: 
\[  \left\Vert \left( \sum_j 2^{-2jk} \left( P_{<j-100} \Lambda^k f \right)^2 \left( P_{j} g \right)^2 \right)^{1/2} \right\Vert_{L^2} \lesssim \left\Vert \sup_j \left| 2^{-jk} P_{<j-100} \Lambda^k f \right| \right\Vert_{L^{\infty}} \left\Vert \left( \sum_j (P_{j} g )^2 \right)^{1/2} \right\Vert_{L^2} \]
and by the Littlewood-Paley square function estimate~: 
\[ \left\Vert \left( \sum_j (P_{j} g )^2 \right)^{1/2} \right\Vert_{L^2} \lesssim \Vert g \Vert_{L^2} \]
While for $f$ we estimate by~: 
\[ \left\Vert \sup_j \left| 2^{-jk} P_{<j-100} \Lambda^k f \right| \right\Vert_{L^{\infty}} \leq \sup_j 2^{-jk} \sum_{l < j} \Vert P_l \Lambda^k f \Vert_{L^{\infty}} \leq \sup_j \sum_{l < j} 2^{-(j-l) k} \Vert P_l f \Vert_{L^{\infty}} \leq \Vert f \Vert_{\dot{B}_{\infty, 1}^0} \]

Putting everything together, we get the desired bound. 
\end{Dem}

The next lemma is a computation that will be useful several times. 

\begin{Lem} Let $f, g$ be two smooth enough functions and recall that $\mu_0$ denote any symbol of order $0$ (that can vary in each term), and $\mu_0^{n, m, l}$ any symbol that can be written as $|\xi|^n |\eta|^m |\xi-\eta|^l \mu_0$ for some $\mu_0$ a symbol of order $0$. Then one has~: 
\begin{subequations}
\begin{align}
&\nabla_{\xi} \int_{\mathbb{R}^d} e^{i t \varphi(\xi, \eta)} \mu_0^{0, 0, 1}(\xi, \eta) \nabla_{\eta} \varphi(\xi, \eta) \widehat{f}(t, \xi - \eta) \widehat{g}(t, \eta) ~ d\eta \label{lemm_computation_termeinit} \\
&= \int_{\mathbb{R}^d} e^{i t \varphi(\xi, \eta)} \mu_0(\xi, \eta) \widehat{f}(t, \xi - \eta) \widehat{g}(t, \eta) ~ d\eta 
+ \int_{\mathbb{R}^d} e^{i t \varphi(\xi, \eta)} \mu_0^{0, -1, 1}(\xi, \eta) \widehat{f}(t, \xi - \eta) \widehat{g}(t, \eta) ~ d\eta \\
&+ \int_{\mathbb{R}^d} e^{i t \varphi(\xi, \eta)} \mu_0^{0, 0, 1}(\xi, \eta) \nabla_{\eta} \widehat{f}(t, \xi - \eta) \widehat{g}(t, \eta) ~ d\eta + \int_{\mathbb{R}^d} e^{i t \varphi(\xi, \eta)} \mu_0^{0, 0, 1}(\xi, \eta) \widehat{f}(t, \xi - \eta) \nabla_{\eta} \widehat{g}(t, \eta) ~ d\eta \\
&+ \int_{\mathbb{R}^d} e^{i t \varphi(\xi, \eta)} \nabla_{\eta} \varphi(\xi, \eta) \mu_0^{-1, 0, 1}(\xi, \eta) \widehat{f}(t, \xi - \eta) \widehat{g}(t, \eta) ~ d\eta
\end{align}
\end{subequations}
where $\varphi$ denotes any of the $\varphi^{\epsilon_1, \epsilon_2 \epsilon_3}$, $\epsilon_1, \epsilon_2, \epsilon_3 \in \{ +, -, 0 \}$. \label{lemm_computation}
\end{Lem}

\begin{Dem}
We first develop the expression by applying $\nabla_{\xi}$~: 
\begin{subequations}
\begin{align}
\eqref{lemm_computation_termeinit} &= \nabla_{\xi} \int_{\mathbb{R}^d} e^{i t \varphi(\xi, \eta)} \mu_0^{0, 0, 1}(\xi, \eta) \nabla_{\eta} \varphi(\xi, \eta) \widehat{f}(t, \xi - \eta) \widehat{g}(t, \eta) ~ d\eta \nonumber \\
&= \int_{\mathbb{R}^d} e^{i t \varphi(\xi, \eta)} i t \nabla_{\xi} \varphi(\xi, \eta) \mu_0^{0, 0, 1}(\xi, \eta) \nabla_{\eta} \varphi(\xi, \eta) \widehat{f}(t, \xi - \eta) \widehat{g}(t, \eta) ~ d\eta \label{lemm_computation_termeipp} \\
&+ \int_{\mathbb{R}^d} e^{i t \varphi(\xi, \eta)} \nabla_{\xi} \left( \mu_0^{0, 0, 1}(\xi, \eta) \nabla_{\eta} \varphi(\xi, \eta) \right) \widehat{f}(t, \xi - \eta) \widehat{g}(t, \eta) ~ d\eta \\
&+ \int_{\mathbb{R}^d} e^{i t \varphi(\xi, \eta)} \mu_0^{0, 0, 1}(\xi, \eta) \nabla_{\eta} \varphi(\xi, \eta) \nabla_{\xi} \widehat{f}(t, \xi - \eta) \widehat{g}(t, \eta) ~ d\eta 
\end{align}
\end{subequations}
Note now that
\begin{equation} \nabla_{\xi} \mu_0^{n, m, l} = \mu_0^{n-1, m, l} + \mu_0^{n, m, l-1}, \quad \nabla_{\eta} \mu_0^{n, m, l} = \mu_0^{n, m-1, l} + \mu_0^{n, m, l-1} \label{derivation_symbols} \end{equation}
and
\begin{equation} \mu_0^{n_1, m_1, l_1} \mu_0^{n_2, m_2, l_2} = \mu_0^{n_1+n_2, m_1+m_2, l_1+l_2} \label{summation_symbols} \end{equation}
with our notation, while $\nabla_{\eta} \varphi = \mu_0$. More precisely, one has $\nabla_{\xi} \nabla_{\eta} \varphi = \mu_0^{0, 0, -1}$ and therefore, 
\[ \nabla_{\xi} \left( \mu_0^{0, 0, 1} \nabla_{\eta} \varphi \right) = \mu_0^{-1, 0, 1} \nabla_{\eta} \varphi + \mu_0 \] 
Up to a change of sign, which reduces to changing the symbol, we may also replace $\nabla_{\xi} \widehat{f}(t, \xi - \eta)$ by $\nabla_{\eta} \widehat{f}(t, \xi - \eta)$. Finally, on the first term above \eqref{lemm_computation_termeipp}, we may apply an integration by parts in frequency to get~: 
\[ \begin{aligned}
\eqref{lemm_computation_termeinit}
&= \int_{\mathbb{R}^d} e^{i t \varphi(\xi, \eta)} \mu_0^{0, -1, 1}(\xi, \eta) \widehat{f}(t, \xi - \eta) \widehat{g}(t, \eta) ~ d\eta + \int_{\mathbb{R}^d} e^{i t \varphi(\xi, \eta)} \mu_0(\xi, \eta) \widehat{f}(t, \xi - \eta) \widehat{g}(t, \eta) ~ d\eta \\
&+ \int_{\mathbb{R}^d} e^{i t \varphi(\xi, \eta)} \mu_0^{0, 0, 1}(\xi, \eta) \nabla_{\eta} \widehat{f}(t, \xi - \eta) \widehat{g}(t, \eta) ~ d\eta + \int_{\mathbb{R}^d} e^{i t \varphi(\xi, \eta)} \mu_0^{0, 0, 1}(\xi, \eta) \widehat{f}(t, \xi - \eta) \nabla_{\eta} \widehat{g}(t, \eta) ~ d\eta \\
&+ \int_{\mathbb{R}^d} e^{i t \varphi(\xi, \eta)} \nabla_{\eta} \varphi(\xi, \eta) \mu_0^{-1, 0, 1}(\xi, \eta) \widehat{f}(t, \xi - \eta) \widehat{g}(t, \eta) ~ d\eta 
\end{aligned} \]

This is the desired form. 
\end{Dem}

\begin{Lem} Let $\mathcal{N}$ be a nonlinearity satisfying the nonresonance conditions, and set  
\[ h = e^{\pm i t \Lambda} \mathcal{N}(\Lambda u, v) \]
Then
\[ \begin{aligned}
&\Vert h \Vert_{H^6} \lesssim t^{-2+\tau+\gamma/2} \Vert u \Vert_X \Vert v \Vert_X , ~~~ & \Vert h \Vert_{H^5} \lesssim t^{-2+\tau+b} \Vert u \Vert_X \Vert v \Vert_X, ~~~ &\Vert h \Vert_{H^4} \lesssim t^{-2+\tau} \Vert u \Vert_X \Vert v \Vert_X \\
&\Vert x h \Vert_{H^6} \lesssim t^{-1+\tau+\gamma/2} \Vert u \Vert_X \Vert v \Vert_X, ~~~ &\Vert x h \Vert_{H^5} \lesssim t^{-1+\tau+b} \Vert u \Vert_X \Vert v \Vert_X, ~~~ &\Vert x h \Vert_{H^4} \lesssim t^{-1+\tau} \Vert u \Vert_X \Vert v \Vert_X
\end{aligned} \]
for $\tau > 0$ that can be made arbitrarily small by choosing adequately $N$ big enough, $\varepsilon$ small enough. 

If we set instead
\[ h = \Lambda^{-1} \mathcal{N}(\Lambda u, v) \]
Then
\[ \Vert h \Vert_{H^7} \lesssim t^{-4/3+\tau+b} \Vert u \Vert_X \Vert v \Vert_X \] \label{lem_decnonres}
\end{Lem}

\begin{Dem}
Any nonlinearity $\mathcal{N}$ can be written as a sum of interactions~: 
\[ \mathcal{N} = \sum_{\epsilon_2, \epsilon_3 \in \{ -, 0, + \}} \mathcal{N}^{\epsilon_2 \epsilon_3} \]
where 
\[ \mathcal{N}^{\epsilon_2 \epsilon_3}(\cdot, \cdot) = \mathcal{N}(P^{\epsilon_2}(D) \cdot, P^{\epsilon_3}(D) \cdot) \]
that is, $\mathcal{N}^{\epsilon_2 \epsilon_3}$ only deals with the interaction between the $\epsilon_2$ part and the $\epsilon_3$ part. In the following, we call "$\pm \pm$ interaction" any such term with $\epsilon_2, \epsilon_3 \in \{ -, + \}$, "$0 \pm$ interaction" any such term with $\epsilon_2 = 0, \epsilon_3 \in \{ -, + \}$, and likewise "$\pm 0$ interaction", "$0 0$ interaction". 

\textbf{Norm }$H^6$ We decompose the nonlinearity and start by controlling the $\pm \pm$ interactions, that is terms of the form 
\begin{equation} \mathcal{F}^{-1} \int e^{i t \varphi^{\pm, \pm \pm}(\xi, \eta)} \nabla_{\eta} \varphi^{\pm, \pm \pm}(\xi, \eta) \mu_0^{0, 0, 1}(\xi, \eta) \widehat{f^{\pm}}(t, \xi - \eta) \widehat{g^{\pm}}(t, \eta) d\eta \label{lem_decnonres_termegen} \end{equation}
where we used the notation $f = e^{t A(D)} u, g = e^{t A(D)} v$. In the following, we omit the superscripts on $\varphi$ and on $u, v, f, g$. $\mu_0$ is, as before, a symbol of order $0$. 

We can apply an integration by parts~: 
\begin{equation*}
\eqref{lem_decnonres_termegen} = t^{-1} \int e^{i t \varphi(\xi, \eta)} \nabla_{\eta} \left( \frac{1}{i} \mu_0^{0, 0, 1}(\xi, \eta) \widehat{f}(t, \xi - \eta) \widehat{g}(t, \eta) \right) ~ d\eta 
\end{equation*}
Now, we distribute the derivative, using \eqref{derivation_symbols}. Therefore, we get~: 
\begin{subequations}
\begin{align}
\eqref{lem_decnonres_termegen} &= t^{-1} \int e^{i t \varphi(\xi, \eta)} \mu_0^{0, -1, 1}(\xi, \eta) \widehat{f}(t, \xi - \eta) \widehat{g}(t, \eta) d\eta \label{lem_decnonres_H61} \\
&+ t^{-1} \int e^{i t \varphi(\xi, \eta)} \mu_0^{0, 0, 1}(\xi, \eta) \nabla_{\eta} \widehat{f}(t, \xi - \eta) \widehat{g}(t, \eta) d\eta \label{lem_decnonres_H62} \\
&+ t^{-1} \int e^{i t \varphi(\xi, \eta)} \mu_0(\xi, \eta) \widehat{f}(t, \xi - \eta) \widehat{g}(t, \eta) d\eta \label{lem_decnonres_H63} \\
&+ t^{-1} \int e^{i t \varphi(\xi, \eta)} \mu_0^{0, 0, 1}(\xi, \eta) \widehat{f}(t, \xi - \eta) \nabla_{\eta} \widehat{g}(t, \eta) d\eta \label{lem_decnonres_H64}
\end{align}
\end{subequations}
We estimate them by using in particular the lemmas \ref{lem_symb}, \ref{lem_infini} and \ref{lem_Hardy}~: 
\begin{align*}
\Vert \eqref{lem_decnonres_H61} \Vert_{H^6} &= t^{-1} \Vert T_{\mu_0}(\Lambda u, \Lambda^{-1} v) \Vert_{H^6} \lesssim t^{-1} \Vert u \Vert_{W^{2, \infty-}} \Vert v \Vert_{H^5} + t^{-1} \Vert u \Vert_{W^{8, \infty-}} \Vert \Lambda^{-1} v \Vert_{L^2} \\
&\lesssim t^{-2+\delta} \Vert u \Vert_X \Vert v \Vert_X + t^{-2+\delta} \Vert u \Vert_X \Vert v \Vert_X \\
\Vert \eqref{lem_decnonres_H62} \Vert_{H^6} &= t^{-1} \Vert T_{\mu_0}(e^{\pm i t \Lambda} \Lambda x f, v) \Vert_{H^6} \lesssim t^{-1} \Vert x f \Vert_{H^7} \Vert v \Vert_{W^{7, \infty-}} \lesssim t^{-2+\delta+\varepsilon+\frac{\gamma}{2}} \Vert u \Vert_X \Vert v \Vert_X \\
\Vert \eqref{lem_decnonres_H63} \Vert_{H^6} &= t^{-1} \Vert T_{\mu_0}(u, v) \Vert_{H^6} \lesssim t^{-1} \Vert u \Vert_{H^6} \Vert v \Vert_{W^{8, \infty-}} \lesssim t^{-2+\delta} \Vert u \Vert_X \Vert v \Vert_X \\
\Vert \eqref{lem_decnonres_H64} \Vert_{H^6} &= t^{-1} \Vert T_{\mu_0}(\Lambda u, e^{\pm i t \Lambda} x g) \Vert_{H^6} \lesssim t^{-1} \Vert u \Vert_{W^{8, \infty-}} \Vert x g \Vert_{H^6} \lesssim t^{-2+\delta+b} \Vert u \Vert_X \Vert v \Vert_X
\end{align*}
For the second inequality, we used the fact that $\Vert x f \Vert_{H^6} \lesssim t^b \Vert u \Vert_X$ and the lemma \ref{lem_xH7}. Finally, since $b \leq \frac{\gamma+\varepsilon}{2}$, we get the desired result. Note that, strictly speaking, we do not estimate, for instance, \eqref{lem_decnonres_H61} in $H^6$, but $\mathcal{F}^{-1} \eqref{lem_decnonres_H61}$ in $H^6$~: however, since we will generically make computations in Fourier space and then estimate in norms defined on the physical space, we omit the $\mathcal{F}^{-1}$ in our notation. 

If we add a $\Lambda^{-1}$ in front, the only additional terms to control are $\Vert \Lambda^{-1} (\ref{lem_decnonres_termegen}.i) \Vert_{L^2}$, with $i = 1, 2, 3, 4$. We then apply the lemma \ref{lem_intfrac} to recover a $L^{6/5}$ norm, then we apply again the lemmas \ref{lem_symb}, \ref{lem_infini} and \ref{lem_Hardy}~: 
\begin{align*}
\Vert \eqref{lem_decnonres_H61} \Vert_{L^{6/5}} &\lesssim t^{-1} \Vert u \Vert_{W^{1, 3}} \Vert \Lambda^{-1} v \Vert_{L^2} \lesssim t^{-4/3+\delta} \Vert u \Vert_X \Vert v \Vert_X \\
\Vert \eqref{lem_decnonres_H62} \Vert_{L^{6/5}} &\lesssim t^{-1} \Vert u \Vert_{W^{1, 3}} \Vert xg \Vert_{L^2} \lesssim t^{-4/3+\delta} \Vert u \Vert_X \Vert v \Vert_X \\
\Vert \eqref{lem_decnonres_H63} \Vert_{L^{6/5}} &\lesssim t^{-1} \Vert u \Vert_{L^2} \Vert v \Vert_{W^{1, 3}} \lesssim t^{-4/3+\delta} \Vert u \Vert_X \Vert v \Vert_X \\
\Vert \eqref{lem_decnonres_H64} \Vert_{L^{6/5}} &\lesssim t^{-1} \Vert x f \Vert_{H^1} \Vert v \Vert_{L^3} \lesssim t^{-4/3+\delta} \Vert u \Vert_X \Vert v \Vert_X
\end{align*} 

Let us now consider the case of a $0 \pm$, $\pm 0$ or $0 0$ interaction. This time, we can directly estimate~: 
\[ \Vert T_{\mu_0}(\Lambda u, v) \Vert_{H^6} \lesssim \Vert u \Vert_{H^7} \Vert v \Vert_{W^{7, \infty-}} \]
but since at least one of $u$ or $v$ is of type $0$, the $\Vert \cdot \Vert_X$ norm gives a stronger decrease~:  
\[ \Vert T_{\mu_0}(\Lambda u, v) \Vert_{H^6} \lesssim t^{-2+\delta} \Vert u \Vert_X \Vert v \Vert_X \]

In presence of a $\Lambda^{-1}$, we recover as above a $L^{6/5}$ norm and obtain the same gain~: 
\[ \Vert T_{\mu_0}(\Lambda u, v) \Vert_{L^{6/5}} \lesssim t^{-4/3+\delta} \Vert u \Vert_X \Vert v \Vert_X \]

\textbf{Norm }$H^5$ The $H^5$ norm can be estimated the same way as the $H^6$ norm, by noticing that the term involving a $t^{\gamma/2}$ isn't present here~; however, in \eqref{lem_decnonres_H63}, $\Vert x f \Vert_{H^7}$ becomes $\Vert x f \Vert_{H^6}$ and one exponent $t^b$ remains. 

\textbf{Norm }$H^4$ The $H^4$ norm can be estimated just like the $H^6$ and $H^5$ norms by noticing that we dropped the derivatives enough to avoid the presence of $\gamma$ or $b$. 

\textbf{Norm with weight }$x$ We start by considering the $\pm \pm$ interactions. In Fourier space, we write, using lemma \ref{lemm_computation}~: 
\begin{subequations}
\begin{align}
&\nabla_{\xi} \int e^{i t \varphi(\xi, \eta)} \nabla_{\eta} \varphi(\xi, \eta) \mu_0^{0, 0, 1}(\xi, \eta) \widehat{f}(t, \xi - \eta) \widehat{g}(t, \eta) d\eta \nonumber \\
&= \int e^{i t \varphi(\xi, \eta)} \mu_0(\xi, \eta) \widehat{f}(t, \xi - \eta) \widehat{g}(t, \eta) d\eta \label{lem_decnonres_xH1} \\
&+ \int e^{i t \varphi(\xi, \eta)} \mu_0^{0, -1, 0}(\xi, \eta) \widehat{f}(t, \xi - \eta) \widehat{g}(t, \eta) d\eta \label{lem_decnonres_xH2} \\
&+ \int e^{i t \varphi(\xi, \eta)} \mu_0^{0, 0, 1}(\xi, \eta) \nabla_{\eta} \widehat{f}(t, \xi - \eta) \widehat{g}(t, \eta) d\eta \label{lem_decnonres_xH3} \\
&+ \int e^{i t \varphi(\xi, \eta)} \mu_0^{0, 0, 1}(\xi, \eta) \widehat{f}(t, \xi - \eta) \nabla_{\eta} \widehat{g}(t, \eta) d\eta \label{lem_decnonres_xH4} \\
&+ \int e^{i t \varphi(\xi, \eta)} \nabla_{\eta} \varphi \mu_0^{-1, 0, 1}(\xi, \eta) \widehat{f}(t, \xi - \eta) \widehat{g}(t, \eta) d\eta \label{lem_decnonres_xH5}
\end{align}
\end{subequations}
Recall that we allowed again $\mu_0^{n, m, l}$ to denote general symbols of order $n+m+l$ (which may be different at each line and even in each term). Notice that $\eqref{lem_decnonres_xH1} = t \eqref{lem_decnonres_H63}$, $\eqref{lem_decnonres_xH2} = t \eqref{lem_decnonres_H61}$, $\eqref{lem_decnonres_xH3} = t \eqref{lem_decnonres_H62}$, $\eqref{lem_decnonres_xH4} = t \eqref{lem_decnonres_H64}$, so that these terms have already been estimated above. Finally, $\eqref{lem_decnonres_xH5} = \Lambda^{-1} \eqref{lem_decnonres_termegen}$, and therefore we can procede the same way. 

For the $0 \pm$ interactions, we can procede the same way by noticing that, in this case, $\nabla_{\eta} \varphi = \pm g_0 \frac{\eta}{|\eta|_0}$ (with the notation of lemma \ref{identite_fond}) and therefore never vanishes, which allows to write~: 
\[ \mu_0(\xi, \eta) = \nabla_{\eta} \varphi \mu_0'(\xi, \eta) \]
for some other symbol of order $0$ $\mu_0'$, and then apply again lemma \ref{lemm_computation}~: 
\begin{align*}
&\nabla_{\xi} \int e^{i t \varphi(\xi, \eta)} \mu_0^{0, 0, 1}(\xi, \eta) \widehat{u^0}(t, \xi - \eta) \widehat{g}(t, \eta) ~ d\eta \\
&= \int e^{i t \varphi(\xi, \eta)} \mu_0(\xi, \eta) \widehat{u^0}(t, \xi - \eta) \widehat{g}(t, \eta) d\eta + \int e^{i t \varphi(\xi, \eta)} \mu_0^{0, 0, 1}(\xi, \eta) \widehat{u^0}(t, \xi - \eta) \nabla_{\eta} \widehat{g}(t, \eta) d\eta \\
&+ \int e^{i t \varphi(\xi, \eta)} \mu_0^{-1, 0, 1}(\xi, \eta) \widehat{u^0}(t, \xi - \eta) \widehat{g}(t, \eta) d\eta + \int e^{i t \varphi(\xi, \eta)} \mu_0^{0, -1, 1}(\xi, \eta) \widehat{u^0}(t, \xi - \eta) \widehat{g}(t, \eta) d\eta \\
&+ \int e^{i t \varphi(\xi, \eta)} \mu_0^{0, 0, 1}(\xi, \eta) \nabla_{\eta} \widehat{u^0}(t, \xi - \eta) \widehat{g}(t, \eta) d\eta
\end{align*}
All these terms can be estimated the same way as before, by noticing that $\Lambda u^0$ behaves at least as well as $f^{\pm}$. 

For the $\pm 0$ integrations, we can also apply lemma \ref{lemm_computation} for the same reason~: 
\begin{align*}
&\nabla_{\xi} \int e^{i t \varphi(\xi, \eta)} \mu_0^{0, 0, 1}(\xi, \eta) \widehat{f}(t, \xi - \eta) \widehat{v^0}(t, \eta) d\eta \\
&= \int e^{i t \varphi(\xi, \eta)} \mu_0(\xi, \eta) \widehat{f}(t, \xi - \eta) \widehat{v^0}(t, \eta) d\eta + \int e^{i t \varphi(\xi, \eta)} \mu_0^{0, 0, 1}(\xi, \eta) \widehat{f}(t, \xi - \eta) \nabla_{\eta} \widehat{v^0}(t, \eta) d\eta \\
&+ \int e^{i t \varphi(\xi, \eta)} \mu_0^{-1, 0, 1}(\xi, \eta) \widehat{f}(t, \xi - \eta) \widehat{v^0}(t, \eta) d\eta + \int e^{i t \varphi(\xi, \eta)} \mu_0^{0, -1, 1}(\xi, \eta) \widehat{f}(t, \xi - \eta) \widehat{v^0}(t, \eta) d\eta \\
&+ \int e^{i t \varphi(\xi, \eta)} \mu_0^{0, 0, 1}(\xi, \eta) \nabla_{\eta} \widehat{f}(t, \xi - \eta) \widehat{v^0}(t, \eta) d\eta
\end{align*}
This time, we may however have more singular terms to control, like $\Lambda^{-1} v^0$ or $x v^0$. But we can estimate these terms in $L^6$, apply the lemma \ref{lem_intfrac} and recover terms with a strong decay.  

For the $0 0$ interactions, we cannot apply lemma \ref{lemm_computation} but it is not needed~: 
\begin{subequations}
\begin{align}
&\nabla_{\xi} \int e^{i t \varphi(\xi, \eta)} \mu_0^{0, 0, 1}(\xi, \eta) \widehat{u^0}(t, \xi - \eta) \widehat{v^0}(t, \eta) d\eta \nonumber \\
&= \int e^{i t \varphi(\xi, \eta)} t \mu_0^{0, 0, 1}(\xi, \eta) \widehat{u^0}(t, \xi - \eta) \widehat{v^0}(t, \eta) d\eta \label{lem_decnonres_xH00} \\
&+ \int e^{i t \varphi(\xi, \eta)} \mu_0^{-1, 0, 1}(\xi, \eta) \widehat{u^0}(t, \xi - \eta) \widehat{v^0}(t, \eta) d\eta + \int e^{i t \varphi(\xi, \eta)} \mu_0(\xi, \eta) \widehat{u^0}(t, \xi - \eta) \widehat{v^0}(t, \eta) d\eta \nonumber \\
&+ \int e^{i t \varphi(\xi, \eta)} \mu_0^{0, 0, 1}(\xi, \eta) \nabla_{\eta} \widehat{u^0}(t, \xi - \eta) \widehat{v^0}(t, \eta) d\eta \nonumber
\end{align}
\end{subequations}
All these terms are simpler to estimate than before, except \eqref{lem_decnonres_xH00} for which we write~: 
\[ \Vert \eqref{lem_decnonres_xH00} \Vert_{H^k} = t \Vert T_{\mu_0}(\Lambda u^0, v^0) \Vert_{H^k} \lesssim t \Vert u^0 \Vert_{W^{k+2, \infty-}} \Vert v^0 \Vert_{H^k} \lesssim t^{-2+2a+\delta} \Vert u \Vert_X^2 \]
as desired. 

We conclude that, with a weight $x$, we have the same estimate with a factor $t^{-1+\tau}$ instead of $t^{-2+\tau}$. 
\end{Dem}

\begin{Cor} Consider now the system \eqref{equ_ABI} + \eqref{equ_cont}. Then~: 
\[ \begin{aligned}
&\Vert \partial_t f^{\pm} \Vert_{H^6} \lesssim t^{-2+\tau+\gamma/2} \Vert u \Vert_X^2, ~~~ &\Vert \partial_t f^{\pm} \Vert_{H^5} \lesssim t^{-2+\tau+b} \Vert u \Vert_X^2, ~~~ &\Vert \partial_t f^{\pm} \Vert_{H^4} \lesssim t^{-2+\tau} \Vert u \Vert_X^2 \\
&\Vert \partial_t x f^{\pm} \Vert_{H^6} \lesssim t^{-1+\tau+\gamma/2} \Vert u \Vert_X^2, ~~~ &\Vert \partial_t x f^{\pm} \Vert_{H^5} \lesssim t^{-1+\tau+b} \Vert u \Vert_X^2, ~~~ &\Vert \partial_t x f^{\pm} \Vert_{H^4} \lesssim t^{-1+\tau} \Vert u \Vert_X^2
\end{aligned} \]
where $\tau > 0$ depends only on $\varepsilon, N$. 
\end{Cor}

\section{Energy and \texorpdfstring{$u^0$}{u0} estimates} \label{section_energyest-0}

\subsection{Energy estimate} \label{section_energyest}

The goal of this section is to prove the following estimate~: 

\begin{Prop}[Energy estimate] For any $\epsilon \in \{ -, 0, + \}$, we have~: 
\[ \Vert u^{\epsilon}(t) \Vert_{H^N}^2 \lesssim \Vert U(t = 1) \Vert_{H^N}^2 + t^{2\varepsilon} \Vert u \Vert_X^3 \]
\end{Prop}

For this, let us write the equation satisfied by $U$~:
\[ \partial_t U + A_0(D) U = \mathcal{N}(\Lambda U, U) \]
By proposition \ref{prop_quasi_lin}, we may write this as~: 
\[ \partial_t U = \sum_{i = 1}^3 (M_{0i} + M_{1i}(U)) \partial_i U \]
for some symmetric matrices $M_{0i}$ and some linear, symmetric-matrices-valued applications $M_{1i}$. Therefore, for $|\alpha| \leq N$~: 
\[ \partial_t \partial^{\alpha} U = \sum_{i = 1}^3 (M_{0i} + M_{1i}(U)) \partial_i \partial^{\alpha} U 
+ \sum_{i = 1}^3 \left( \partial^{\alpha} (M_{1i}(U) \partial_i U) - M_{1i}(U) \partial_i \partial^{\alpha} U \right) \]
Hence~: 
\[ \partial_t \Vert \partial^{\alpha} U \Vert_{L^2}^2 = \int \partial^{\alpha} U \cdot \sum_{i = 1}^3 (M_{0i} + M_{1i}(U)) \partial_i \partial^{\alpha} U ~ dx + \int \partial^{\alpha} U \cdot \sum_{i = 1}^3 \left( \partial^{\alpha} (M_{1i}(U) \partial_i U) - M_{1i}(U) \partial_i \partial^{\alpha} U \right) ~ dx \]
By symmetry of the matrices, the first term can be rewritten as~: 
\[ \int \partial^{\alpha} U \cdot \sum_{i = 1}^3 (M_{0i} + M_{1i}(U)) \partial_i \partial^{\alpha} U ~ dx
= -\frac{1}{2} \int \partial^{\alpha} U \cdot \sum_{i = 1}^3 \partial_i M_{1i}(U) \partial^{\alpha} U ~ dx \]
Finally, applying the Moser estimate of lemma \ref{Moser_est} and Hölder inequalities~: 
\[ \partial_t \Vert \partial^{\alpha} U \Vert_{L^2}^2 \lesssim \Vert U \Vert_{H^N}^2 \Vert \nabla U \Vert_{L^{\infty}} \]

We then notice that, since $U = u^{+} + u^{-} + u^0$, we have that
\[ \Vert U \Vert_{W^{1, \infty}} \leq \Vert u^{+} \Vert_{W^{1, \infty}} + \Vert u^{-} \Vert_{W^{1, \infty}} + \Vert u^0 \Vert_{W^{1, \infty}} \lesssim t^{-1} \Vert u \Vert_X + \Vert u^0 \Vert_{H^6} \lesssim t^{-1} \Vert u \Vert_X \]

Therefore, we showed~: 
\[ \partial_t \Vert U \Vert_{H^N}^2 \lesssim \Vert U \Vert_{H^N}^2 \Vert U \Vert_{W^{1, \infty}} \lesssim t^{2 \varepsilon - 1} \Vert u \Vert_X^3 \]

This implies~: 
\[ \Vert U \Vert_{H^N}^2 \lesssim \Vert U(t = 1) \Vert_{H^N}^2 + t^{2\varepsilon} \Vert u \Vert_X^3 \]

We deduce the same inequality for $u^{\pm}$ instead of $U$ by applying the projection operators and Parseval's identity~: since these projection operators are pointwise orthogonal, we have that 
\[ |\widehat{u^{\pm}}(\xi)| \leq |\widehat{U}(\xi)| \]
and thus $\Vert u^{\pm} \Vert_{H^N} \leq \Vert U \Vert_{H^N}$. 

\paragraph{Conclusion} $\varepsilon, N$ can be chosen independently of any other parameter. 

\subsection{Estimates for \texorpdfstring{$u^0$}{u0}} \label{section_est0}

In this section, we establish all the estimates concerning $u^0$ by making use of the constraint equation \eqref{equ_cont} and its structure, summarized in the following proposition. In the remaining sections, we will only consider $u^{+}$ and $u^{-}$. 

\begin{Prop}[Estimates on $u^0$] $u^0$ satisfies the following bounds~: 
\[ \begin{aligned}
&\Vert u^0 \Vert_{H^7} \lesssim t^{-1} \Vert u \Vert_X^2, ~~~~~ \Vert u^0 \Vert_{H^N} \lesssim t^{\varepsilon - 1} \Vert u \Vert_X^2, ~~~~~ \Vert u^0 \Vert_{W^{1, \infty-}} \lesssim t^{-2+a} \Vert u \Vert_X^2, \\
&~~~~~ \Vert \Lambda x u^0 \Vert_{H^4} \lesssim t^{-1+a} \Vert u \Vert_X^2, ~~~~~ \Vert \Lambda x u^0 \Vert_{H^6} \lesssim t^{-1+a+\gamma/2} \Vert u \Vert_X^2
\end{aligned} \]
\end{Prop}

\paragraph{$H^7$ estimate} It is a consequence of lemma \ref{lem_decnonres}~: 
\[ \Vert u^0 \Vert_{H^7} \lesssim t^{-4/3+\gamma/2+\tau} \Vert u \Vert_X^2 \lesssim t^{-1} \Vert u \Vert_X^2 \]
provided $\gamma, \tau$ are small enough in front of $1$. 

\paragraph{$H^N$ estimate} By using the fact that the nonlinearity is a product~: 
\[ \begin{aligned} 
\Vert \Lambda^N u^0 \Vert_{L^2} &\lesssim \Vert \Lambda^{N-1} \mathcal{N}'(\nabla U, U) \Vert_{L^2} \lesssim \Vert U \Vert_{W^{1, \infty}} \Vert U \Vert_{H^N} \lesssim t^{\varepsilon-1} \Vert u \Vert_X^2 
\end{aligned} \]
This controls the norm $\dot{H}^N$, but since we already controlled the $L^2$ norm above, we deduce by interpolation that $\Vert u^0 \Vert_{H^N} \lesssim t^{\varepsilon - 1} \Vert u \Vert_X^2$. 

\paragraph{$L^{\infty-}$ estimate} Likewise, 
\[ \Vert \Lambda u^0 \Vert_{L^{\infty-}} \lesssim \Vert \mathcal{N}'(\nabla U, U) \Vert_{L^{\infty-}} \lesssim t^{-2+} \Vert u \Vert_X^2 \]
and
\[ \Vert u^0 \Vert_{L^{\infty-}} = \Vert \Lambda^{-1} \mathcal{N}'(\Lambda U, U) \Vert_{L^{\infty-}} \lesssim \Vert \mathcal{N}(\Lambda U, U) \Vert_{L^{3-}} \]
Let us now use the structure of the nonlinearity. For the $\pm \pm$ interactions, we apply the same integrations by parts as in the proof of lemma \ref{lem_decnonres} (but this time we will need to control this in $L^3$)~: 
\begin{subequations}
\begin{align}
\widehat{\Lambda u^0}(t, \xi) &= \int e^{i t \varphi(\xi, \eta)} \nabla_{\eta} \varphi(\xi, \eta) b(\xi, \eta) |\xi - \eta| \widehat{f}(t, \xi - \eta) \widehat{f}(t, \eta) ~ d\eta \nonumber \\
&= t^{-1} \int e^{i t \varphi(\xi, \eta)} \mu_0(\xi, \eta) \widehat{f}(t, \xi - \eta) \widehat{f}(t, \eta) ~ d\eta \label{est_u0_Linf1} \\
&+ t^{-1} \int e^{i t \varphi(\xi, \eta)} \mu_0^{0, -1, 1}(\xi, \eta) \widehat{f}(t, \xi - \eta) \widehat{f}(t, \eta) ~ d\eta \label{est_u0_Linf2} \\
&+ t^{-1} \int e^{i t \varphi(\xi, \eta)} \mu_0^{0, 0, 1}(\xi, \eta) \nabla_{\eta} \widehat{f}(t, \xi - \eta) \widehat{f}(t, \eta) ~ d\eta \label{est_u0_Linf3} \\
&+ t^{-1} \int e^{i t \varphi(\xi, \eta)} \mu_0^{0, 0, 1}(\xi, \eta) \widehat{f}(t, \xi - \eta) \nabla_{\eta} \widehat{f}(t, \eta) ~ d\eta \label{est_u0_Linf4}
\end{align} 
\end{subequations}
Then~: 
\begin{align*}
\Vert \eqref{est_u0_Linf1} \Vert_{L^{3-}} &\lesssim t^{-1} \Vert u \Vert_{L^{6-}}^2 \lesssim t^{-7/3+} \Vert u \Vert_X^2 \\
\Vert \eqref{est_u0_Linf2} \Vert_{L^{3-}} &\lesssim t^{-1} \Vert \Lambda u \Vert_{L^{12-}} \Vert \Lambda^{-1} e^{\pm i t \Lambda} f \Vert_{L^4} \lesssim t^{-7/3+} \Vert u \Vert_X \Vert f \Vert_{L^{4/3}} \lesssim t^{-7/3+} \Vert u \Vert_X \Vert \langle x \rangle f \Vert_{L^2} \lesssim t^{-7/3+b+} \Vert u \Vert_X^2 \\
\Vert \eqref{est_u0_Linf3} \Vert_{L^{3-}} &\lesssim t^{-1} \Vert \Lambda e^{\pm i t \Lambda} x f \Vert_{L^3} \Vert u \Vert_{W^{1, \infty-}} \lesssim t^{-7/3+} \Vert \Lambda \langle x \rangle^2 f \Vert_{H^1} \Vert u \Vert_X^2 \lesssim t^{-7/3+\gamma+} \Vert u \Vert_X^2 \\
\Vert \eqref{est_u0_Linf4} \Vert_{L^{3-}} &\lesssim t^{-1} \Vert u \Vert_{W^{1, 12-}} \Vert e^{\pm i t \Lambda} x f \Vert_{L^4} \lesssim t^{-7/3+} \Vert u \Vert_X \Vert \Lambda x f \Vert_{L^{4/3}} \lesssim t^{-7/3+\gamma+} \Vert u \Vert_X^2
\end{align*}
In particular, if $\gamma > 0$ is small enough, we indeed have a $t^{-2+}$ decay above, with $-2+$ arbitrarily close to $-2$ (in a way depending only on $\infty-$). For the $0 \pm$ or $\pm 0$ interactions, we can also apply an integration by parts and the estimate on \eqref{est_u0_Linf1} is identical~; for the $0 \pm$ interaction, the estimates of \eqref{est_u0_Linf2}, \eqref{est_u0_Linf4} are identical while 
\[ \begin{aligned}
\Vert \eqref{est_u0_Linf3} \Vert_{L^{3-}} &\lesssim t^{-1} \Vert \Lambda x u^0 \Vert_{L^3} \Vert u \Vert_{L^{\infty-}} \lesssim t^{-2+} \Vert \Lambda x u^0 \Vert_{L^2} \Vert u \Vert_X \lesssim t^{-3+} \Vert u \Vert_X^2
\end{aligned} \]
Finally, for a $\pm 0$ interaction, the estimate on \eqref{est_u0_Linf3} is identical and 
\begin{align*}
\Vert \eqref{est_u0_Linf2} \Vert_{L^{3-}} &\lesssim t^{-1} \Vert u \Vert_{W^{1, 6-}} \Vert \Lambda^{-1} u^0 \Vert_{L^6} \lesssim t^{-5/3+\delta+} \Vert u \Vert_X \Vert u^0 \Vert_{L^2} \lesssim t^{-8/3+\delta+} \Vert u \Vert_X^2 \\
\Vert \eqref{est_u0_Linf4} \Vert_{L^{3-}} &\lesssim t^{-1} \Vert u \Vert_{W^{1, 6-}} \Vert x u^0 \Vert_{L^6} \lesssim t^{-5/3+\delta+} \Vert u \Vert_X \Vert \Lambda x u^0 \Vert_{L^2} \lesssim s^{-8/3+\delta+} \Vert u \Vert_X^2 
\end{align*}
Only the $00$ interaction remains, for which we do not need to integrate by parts~: 
\[ \Vert T_{\mu_0}(\Lambda u^0, u^0) \Vert_{L^{3-}} \lesssim \Vert u^0 \Vert_{W^{1, \infty-}} \Vert u^0 \Vert_{W^{1, 3-}} \lesssim t^{-7/3+} \Vert u \Vert_X^2 \]
Therefore, summing every contribution we get~: 
\[ \Vert u^0 \Vert_{W^{1, \infty-}} \lesssim t^{-2+a} \Vert u \Vert_X^2 \]
with $a > 0$ that we can choose arbitrarily small, depending only on $\infty-$. 

\paragraph{Estimate with weight $x$} In contrast with others, we only estimate $\Vert \Lambda x u^0 \Vert_{H^6}$ and not $\Vert x u^0 \Vert_{L^2}$, because the presence of a $\Lambda^{-1}$ is too singular otherwise. Notice that 
\[ \Vert \Lambda x u^0 \Vert_{H^6} \lesssim \Vert x \Lambda u^0 \Vert_{H^6} + \Vert u^0 \Vert_{H^6} \]
The estimate on $\Vert u^0 \Vert_{H^6}$ is a consequence of lemma \ref{lem_decnonres}. Furthermore, since $\Lambda u^0 = \mathcal{N}(U, \nabla U)$, we can also apply lemma \ref{lem_decnonres} and deduce that
\[ \Vert \Lambda x u^0 \Vert_{H^6} \lesssim t^{-1+\tau+\gamma/2} \Vert u \Vert_X^2 \]
and
\[ \Vert \Lambda x u^0 \Vert_{H^4} \lesssim t^{-1+\tau} \Vert u \Vert_X^2 \]
Therefore, here, we obtain the condition $a \geq \tau = \tau(\varepsilon, N)$. 

\paragraph{Conclusion} We showed the desired inequalities for all the norms involving $u^0$. The exponent $a$ has only to satisfy $a \geq \tau(\varepsilon, N)$, and we need to choose $\infty-$ close enough of $+\infty$ (once $\varepsilon, N$ are fixed). 

\subsection{Estimate of the \texorpdfstring{$H^6$}{H6} norm} \label{section_estH6}

In this section, we want to control $\Vert u \Vert_{H^6}$, where $u = u^{+}$ or $u = u^{-}$. One has that $\Vert u \Vert_{H^6} = \Vert f \Vert_{H^6}$ (with the same convention of omitting the $\pm$ superscript). 

Let us write 
\[ f(t) = f(1) + \int_1^t \partial_s f(s) ~ ds \]
But by lemma \ref{lem_decnonres}, 
\[ \Vert \partial_s f(s) \Vert_{H^6} \lesssim s^{-2+\gamma/2+\tau} \Vert u \Vert_X^2 \]
In particular, if $\tau = \tau(\varepsilon, N), \gamma$ are small enough,  
\[ \Vert f(t) \Vert_{H^6} \lesssim \Vert f(1) \Vert_{H^6} + \int_1^t s^{-2+\gamma/2+\tau} \Vert u \Vert_X^2 ~ ds \lesssim \Vert f(1) \Vert_{H^6} + \Vert u \Vert_X^2 \]
as desired. 

\section{Estimate of the \texorpdfstring{$L^2$}{L2} norm with weight \texorpdfstring{$x$}{x}} \label{section_L2x}

The goal of this section is to prove the following estimate~: 

\begin{Prop}[First weighted estimate] For $\epsilon = +$ or $\epsilon = -$, we have 
\[ \Vert x f^{\epsilon}(t) \Vert_{H^5} \leq \Vert x f^{\epsilon}(t = 1) \Vert_{H^5} + C \Vert u \Vert_X^2 \left( 1 + \Vert u \Vert_X \right) , ~~~~~ \Vert x f^{\epsilon}(t) \Vert_{H^6} \leq \Vert x f^{\epsilon}(t = 1) \Vert_{H^6} + C t^b \Vert u \Vert_X^2 \left( 1 + \Vert u \Vert_X \right) \]
\end{Prop}

\begin{Intuit}
We saw in lemma \ref{lem_decnonres} that $\partial_t x f$ was estimated in $H^6$ with a $t^{-1+\tau+\gamma/2}$ decay. If we integrate this relation, we obtain a $t^{\tau+\gamma/2}$ growth. 

In order to get a better estimate, we will use the fact that we have here an integral in time. More precisely, the weight $x$ acts as a derivative in Fourier, and if this derivative hits the exponential $e^{i t \varphi}$, we can use lemma \ref{identite_fond} to transform $\nabla_{\xi} \varphi$ into $\nabla_{\eta} \varphi$ and $\varphi$. On the first term, we may apply an integration by parts in frequency to win an additional decay $1/t$~; on the second one, we can integrate by parts in time and win an additional decay $1/t$ through the use of lemma \ref{lem_decnonres}. 
\end{Intuit}

In the first subsection, we will only consider $\pm \pm$ interactions. 

In order to use lemma \ref{identite_fond}, we need the (artificial) presence of a factor $|\xi|_0$. Therefore, in order to control the $L^2$ norm without derivative, we will need to apply a $\Lambda^{-1}$ in front~; however, in order to estimate the $\dot{H}^6$ norm, it will be enough to estimate the $H^5$ norm of the computed terms.

\subsection{\texorpdfstring{$\pm \pm$}{+/- +/-} interactions}

\paragraph{Simplification of the terms} Let us write in Fourier~: 
\begin{subequations}\label{est_Hx++_init}
\begin{align}
&|\xi|_0 \nabla_{\xi} \int_1^t \int e^{i s \varphi(\xi, \eta)} \nabla_{\eta} \varphi(\xi, \eta) \mu_0^{0, 0, 1}(\xi, \eta) \widehat{f}(s, \xi - \eta) \widehat{f}(s, \eta) ~ d\eta ds \tag{\ref{est_Hx++_init}} \\
&= \int_1^t \int e^{i s \varphi(\xi, \eta)} s |\xi|_0 \nabla_{\xi} \varphi(\xi, \eta) \nabla_{\eta} \varphi(\xi, \eta) \mu_0^{0, 0, 1}(\xi, \eta) \widehat{f}(s, \xi - \eta) \widehat{f}(s, \eta) ~ d\eta ds \label{est_Hx++_init_id} \\
&+ \int_1^t \int e^{i s \varphi(\xi, \eta)} |\xi|_0 \nabla_{\xi} \nabla_{\eta} \varphi(\xi, \eta) \mu_0^{0, 0, 1}(\xi, \eta) |\xi - \eta| \widehat{f}(s, \xi - \eta) \widehat{f}(s, \eta) ~ d\eta ds \label{est_Hx++_initb} \\
&+ \int_1^t \int e^{i s \varphi(\xi, \eta)} \nabla_{\eta} \varphi(\xi, \eta) \mu_0^{0, 0, 1}(\xi, \eta) \widehat{f}(s, \xi - \eta) \widehat{f}(s, \eta) ~ d\eta ds \label{est_Hx++_initc} \\
&+ \int_1^t \int e^{i s \varphi(\xi, \eta)} \nabla_{\eta} \varphi(\xi, \eta) \mu_0^{1, 0, 0}(\xi, \eta) \widehat{f}(s, \xi - \eta) \widehat{f}(s, \eta) ~ d\eta ds \label{est_Hx++_initd} \\
&+ \int_1^t \int e^{i s \varphi(\xi, \eta)} \nabla_{\eta} \varphi(\xi, \eta) \mu_0^{1, 0, 1}(\xi, \eta) \nabla_{\xi} \widehat{f}(s, \xi - \eta) \widehat{f}(s, \eta) ~ d\eta ds \label{est_Hx++_inite} \\
&= \int_1^t \int e^{i s \varphi(\xi, \eta)} s \varphi(\xi, \eta) \nabla_{\eta} \varphi(\xi, \eta) \mu_0^{0, 0, 1}(\xi, \eta) \widehat{f}(s, \xi - \eta) \widehat{f}(s, \eta) ~ d\eta ds \tag{\ref{est_Hx++_init_id}1} \label{est_Hx++_init_ippt} \\
&- \int_1^t \int e^{i s \varphi(\xi, \eta)} s \epsilon_1 \epsilon_2 |\eta|_0 \nabla_{\eta} \varphi(\xi, \eta) \nabla_{\eta} \varphi(\xi, \eta) \mu_0^{0, 0, 1}(\xi, \eta) \widehat{f}(s, \xi - \eta) \widehat{f}(s, \eta) ~ d\eta ds \tag{\ref{est_Hx++_init_id}2} \\
&+ \eqref{est_Hx++_initb} + \eqref{est_Hx++_initc} + \eqref{est_Hx++_initd} + \eqref{est_Hx++_inite} \nonumber
\end{align} 
\end{subequations}
by lemma \ref{identite_fond}. On these terms, we can apply an integration by parts in frequency or in time, to obtain~: 
\begin{subequations} 
\label{est_Hx++_termes_ipp}
\begin{align}
\eqref{est_Hx++_init}
&= \int_1^t \int e^{i s \varphi(\xi, \eta)} s \nabla_{\eta} \varphi(\xi, \eta) \mu_0^{0, 0, 1} (\xi, \eta) \partial_s \widehat{f}(s, \xi - \eta) \widehat{f}(s, \eta) ~ d\eta ds \label{est_Hx++_ipp1} \\
&+ \int_1^t \int e^{i s \varphi(\xi, \eta)} s \nabla_{\eta} \varphi(\xi, \eta) \mu_0^{0, 0, 1}(\xi, \eta) \widehat{f}(s, \xi - \eta) \partial_s \widehat{f}(s, \eta) ~ d\eta ds \label{est_Hx++_ipp2} \\
&+ \int e^{i t \varphi(\xi, \eta)} t \nabla_{\eta} \varphi(\xi, \eta) \mu_0^{0, 0, 1}(\xi, \eta) \widehat{f}(t, \xi - \eta) \widehat{f}(t, \eta) ~ d\eta \label{est_Hx++_ipp3} \\
&+ \int e^{i \varphi(\xi, \eta)} \nabla_{\eta} \varphi(\xi, \eta) \mu_0^{0, 0, 1}(\xi, \eta) \widehat{f}(1, \xi - \eta) \widehat{f}(1, \eta) \eta ~ d\eta \label{est_Hx++_ipp4} \\
&+ \int_1^t \int e^{i s \varphi(\xi, \eta)} \nabla_{\eta} \varphi(\xi, \eta) \mu_0^{0, 1, 0}(\xi, \eta) \widehat{f}(s, \xi - \eta) \widehat{f}(s, \eta) ~ d\eta ds \label{est_Hx++_sym1a} \\
&+ \int_1^t \int e^{i s \varphi(\xi, \eta)} \nabla_{\eta} \varphi(\xi, \eta) \mu_0^{0, 1, 1}(\xi, \eta) \nabla_{\eta} \widehat{f}(s, \xi - \eta) \widehat{f}(s, \eta) ~ d\eta ds \label{est_Hx++_sym2a} \\
&+ \int_1^t \int e^{i s \varphi(\xi, \eta)} \nabla_{\eta} \varphi(\xi, \eta) \mu_0^{0, 1, 1}(\xi, \eta) \widehat{f}(s, \xi - \eta) \nabla_{\eta} \widehat{f}(s, \eta) ~ d\eta ds \label{est_Hx++_sym2b} \\
&+ \int_1^t \int e^{i s \varphi(\xi, \eta)} \nabla_{\eta} \left( |\xi|_0 \nabla_{\xi} + \epsilon_1 \epsilon_2 |\eta|_0 \nabla_{\eta} \right) \varphi(\xi, \eta) \mu_0^{0, 0, 1}(\xi, \eta) \widehat{f}(s, \xi - \eta) \widehat{f}(s, \eta) ~ d\eta ds \\
&+ \int_1^t \int e^{i s \varphi(\xi, \eta)} \nabla_{\eta} \varphi(\xi, \eta) \mu_0^{0, 0, 1}(\xi, \eta) \widehat{f}(s, \xi - \eta) \widehat{f}(s, \eta) ~ d\eta ds \label{est_Hx++_sym1b} \\
&+ \int_1^t \int e^{i s \varphi(\xi, \eta)} \nabla_{\eta} \varphi(\xi, \eta) \mu_0^{1, 0, 0}(\xi, \eta) \widehat{f}(s, \xi - \eta) \widehat{f}(s, \eta) ~ d\eta ds \label{est_Hx++_sym_der} \\
&+ \int_1^t \int e^{i s \varphi(\xi, \eta)} \nabla_{\eta} \varphi(\xi, \eta) \mu_0^{1, 0, 1}(\xi, \eta) \nabla_{\xi} \widehat{f}(s, \xi - \eta) \widehat{f}(s, \eta) ~ d\eta ds \label{est_Hx++_ipp5}
\end{align}
\end{subequations}
where we used the fact that $|\eta|_0 = |\eta| \mu_0$, and likewise for $\xi$. We can start simplifying the expression by noticing that some of the terms are symmetric if we exchange the role of $\xi - \eta$ and $\eta$~: in particular, \eqref{est_Hx++_sym1a} and \eqref{est_Hx++_sym1b}, or \eqref{est_Hx++_sym2a} and \eqref{est_Hx++_sym2b}. Furthermore, we also have
\begin{equation} \mu_0^{n+N, m, l} = \mu_0^{n, m+N, l} + \mu_0^{n, m, l+N} \label{distrib_deriv_symbols} \end{equation}
(which corresponds to distributing the derivatives), we can express \eqref{est_Hx++_sym_der} as a combination of \eqref{est_Hx++_sym1a} and \eqref{est_Hx++_sym1b}. Finally, by lemma \ref{identite_fond}, we have that
\[ \nabla_{\eta} \left( |\xi|_0 \nabla_{\xi} + \epsilon_1 \epsilon_2 |\eta|_0 \nabla_{\eta} \right) \varphi(\xi, \eta) = \mu_0 \nabla_{\eta} \varphi + \mu_0^{0, 0, -1} \varphi \]
It thus only remains to study 
\begin{subequations}
\begin{align}
\eqref{est_Hx++_init}
&= \eqref{est_Hx++_ipp1} + \eqref{est_Hx++_ipp2} + \eqref{est_Hx++_ipp3} + \eqref{est_Hx++_ipp4} + \eqref{est_Hx++_sym1a} + \eqref{est_Hx++_sym2a} + \eqref{est_Hx++_ipp5} \nonumber \\
&+ \int_1^t \int e^{i s \varphi(\xi, \eta)} \varphi(\xi, \eta) \mu_0(\xi, \eta) \widehat{f}(s, \xi - \eta) \widehat{f}(s, \eta) ~ d\eta ds \tag{\ref{est_Hx++_termes_ipp}l}
\end{align} 
\end{subequations}
On each of these terms, we apply an integration by parts and get~: 
\begin{subequations}\label{est_Hx++_term1}
\begin{align}
\eqref{est_Hx++_init}
&= \int_1^t \int e^{i s \varphi(\xi, \eta)} \mu_0(\xi, \eta) \partial_s \widehat{f}(s, \xi - \eta) \widehat{f}(s, \eta) ~ d\eta ds \label{est_Hx++_term1-1} \\
&+ \int_1^t \int e^{i s \varphi(\xi, \eta)} \mu_0^{0, -1, 1}(\xi, \eta) \partial_s \widehat{f}(s, \xi - \eta) \widehat{f}(s, \eta) ~ d\eta ds \label{est_Hx++_term1-2} \\
&+ \int_1^t \int e^{i s \varphi(\xi, \eta)} \mu_0^{0, 0, 1}(\xi, \eta) \partial_s \nabla_{\eta} \widehat{f}(s, \xi - \eta) \widehat{f}(s, \eta) ~ d\eta ds \label{est_Hx++_term1-3} \\
&+ \int_1^t \int e^{i s \varphi(\xi, \eta)} \mu_0^{0, 0, 1}(\xi, \eta) \partial_s \widehat{f}(s, \xi - \eta) \nabla_{\eta} \widehat{f}(s, \eta) ~ d\eta ds \label{est_Hx++_term1-4} \\
&+ \int_1^t \int e^{i s \varphi(\xi, \eta)} \mu_0^{0, -1, 1}(\xi, \eta) \widehat{f}(s, \xi - \eta) \partial_s \widehat{f}(s, \eta) ~ d\eta ds \label{est_Hx++_term1-5} \\
&+ \int_1^t \int e^{i s \varphi(\xi, \eta)} \mu_0^{0, 0, 1}(\xi, \eta) \nabla_{\eta} \widehat{f}(s, \xi - \eta) \partial_s \widehat{f}(s, \eta) ~ d\eta ds \label{est_Hx++_term1-6} \\
&+ \int_1^t \int e^{i s \varphi(\xi, \eta)} \mu_0^{0, 0, 1}(\xi, \eta) \widehat{f}(s, \xi - \eta) \partial_s \nabla_{\eta} \widehat{f}(s, \eta) ~ d\eta ds \label{est_Hx++_term1-7} \\
\nextParentEquation[est_Hx++_term2]
&+ \int e^{i t \varphi(\xi, \eta)} \mu_0(\xi, \eta) \widehat{f}(t, \xi - \eta) \widehat{f}(t, \eta) ~ d\eta \label{est_Hx++_term2-1} \\
&+ \int e^{i t \varphi(\xi, \eta)} \mu_0^{0, -1, 1}(\xi, \eta) \widehat{f}(t, \xi - \eta) \widehat{f}(t, \eta) ~ d\eta \label{est_Hx++_term2-2} \\
&+ \int e^{i t \varphi(\xi, \eta)} \mu_0^{0, 0, 1}(\xi, \eta) \nabla_{\eta} \widehat{f}(t, \xi - \eta) \widehat{f}(t, \eta) ~ d\eta \label{est_Hx++_term2-3} \\
&+ \int e^{i t \varphi(\xi, \eta)} \mu_0^{0, 0, 1}(\xi, \eta) \widehat{f}(t, \xi - \eta) \nabla_{\eta} \widehat{f}(t, \eta) ~ d\eta \label{est_Hx++_term2-4} \\
&+ \int e^{i \varphi(\xi, \eta)} \left( \mu_0(\xi, \eta) + \mu_0^{0, 0, 1}(\xi, \eta) \right) \widehat{f}(1, \xi - \eta) \widehat{f}(1, \eta) \eta ~ d\eta \label{est_Hx++_term2-init} \\
\nextParentEquation[est_Hx++_term3]
&+ \int_1^t \int e^{i s \varphi(\xi, \eta)} s^{-1} \mu_0(\xi, \eta) \widehat{f}(s, \xi - \eta) \widehat{f}(s, \eta) ~ d\eta ds \\
\nextParentEquation[est_Hx++_term4]
&+ \int_1^t \int e^{i s \varphi(\xi, \eta)} s^{-1} \mu_0^{0, 1, -1}(\xi, \eta) \widehat{f}(s, \xi - \eta) \widehat{f}(s, \eta) ~ d\eta ds \label{est_Hx++_term4-1} \\
&+ \int_1^t \int e^{i s \varphi(\xi, \eta)} s^{-1} \mu_0^{0, 1, 0}(\xi, \eta) \nabla_{\eta} \widehat{f}(s, \xi - \eta) \widehat{f}(s, \eta) ~ d\eta ds \label{est_Hx++_term4-2} \\
&+ \int_1^t \int e^{i s \varphi(\xi, \eta)} s^{-1} \mu_0^{0, 1, 0}(\xi, \eta) \widehat{f}(s, \xi - \eta) \nabla_{\eta} \widehat{f}(s, \eta) ~ d\eta ds \label{est_Hx++_term4-3} \\
\nextParentEquation[est_Hx++_term5]
&+ \int_1^t \int e^{i s \varphi(\xi, \eta)} s^{-1} \mu_0^{0, 1, 1}(\xi, \eta) \nabla^2_{\eta} \widehat{f}(s, \xi - \eta) \widehat{f}(s, \eta) ~ d\eta ds \label{est_Hx++_term5-1} \\
&+ \int_1^t \int e^{i s \varphi(\xi, \eta)} s^{-1} \mu_0^{0, 1, 1}(\xi, \eta) \nabla_{\eta} \widehat{f}(s, \xi - \eta) \nabla_{\eta} \widehat{f}(s, \eta) ~ d\eta ds \label{est_Hx++_term5-2} \\
&+ \int_1^t \int e^{i s \varphi(\xi, \eta)} s^{-1} \mu_0^{0, -1, 2}(\xi, \eta) \nabla_{\xi} \widehat{f}(s, \xi - \eta) \widehat{f}(s, \eta) ~ d\eta ds \label{est_Hx++_term5-3} \\
&+ \int_1^t \int e^{i s \varphi(\xi, \eta)} s^{-1} \mu_0^{0, 0, 2}(\xi, \eta) \nabla^2_{\xi} \widehat{f}(s, \xi - \eta) \widehat{f}(s, \eta) ~ d\eta ds \label{est_Hx++_term5-4} \\
&+ \int_1^t \int e^{i s \varphi(\xi, \eta)} s^{-1} \mu_0^{0, 0, 2}(\xi, \eta) \nabla_{\xi} \widehat{f}(s, \xi - \eta) \nabla_{\eta} \widehat{f}(s, \eta) ~ d\eta ds \label{est_Hx++_term5-5} 
\end{align} 
\end{subequations}
where we simplified again by taking symmetries and redundancies into account.

This repartition corresponds to~: 
\begin{itemize}
\item \eqref{est_Hx++_term1} contains the terms with a time derivative. 
\item \eqref{est_Hx++_term2} contains the terms without time integration. 
\item \eqref{est_Hx++_term4} contains the terms with only one frequency derivative or only one $\Lambda^{-1}$. 
\item \eqref{est_Hx++_term5} contains the terms with two frequency derivatives or one frequency derivative and one $\Lambda^{-1}$. 
\end{itemize}
Note that one term remains above, \eqref{est_Hx++_term3}, but it is easier to estimate. The term at initial time \eqref{est_Hx++_term2-init} is estimated by the hypothesis of our main theorem. 

\subsubsection{Estimate of \eqref{est_Hx++_term1}}

We apply lemma \ref{lem_decnonres}~: 
\begin{align*}
\Vert \Lambda^{-1} \eqref{est_Hx++_term1-2} \Vert_{L^2} &\lesssim \int_1^t \Vert T_{\mu_0}(\Lambda e^{\pm i s \Lambda} \partial_s f, \Lambda^{-1} u) \Vert_{L^{6/5}} ~ ds \lesssim \int_1^t \Vert e^{\pm i s \Lambda}  \partial_s f \Vert_{W^{1, 3}} \Vert \Lambda^{-1} u \Vert_{L^2} ~ ds \lesssim \int_1^t \Vert \partial_s f \Vert_{H^2} \Vert x f \Vert_{L^2} ~ ds \\
&\lesssim s^{-2+\tau} \Vert u \Vert_X^3 ~ ds \lesssim \Vert u \Vert_X^3 \\
\Vert \Lambda^{-1} \eqref{est_Hx++_term1-3} \Vert_{L^2} &\lesssim \int_1^t \Vert T_{\mu_0}(\Lambda e^{\pm i s \Lambda} \partial_s x f, u) \Vert_{L^{6/5}} ~ ds \lesssim \int_1^t \Vert \partial_s x f \Vert_{L^2} \Vert u \Vert_{L^3} ~ ds \\
&\lesssim \int_1^t s^{-4/3+\tau} \Vert u \Vert_X^3 ~ ds \lesssim \Vert u \Vert_X^3 \\
\Vert \Lambda^{-1} \eqref{est_Hx++_term1-4} \Vert_{L^2} &\lesssim \int_1^t \Vert T_{\mu_0}(\Lambda e^{\pm i s \Lambda} \partial_s f, e^{\pm i s \Lambda} x f) \Vert_{L^{6/5}} ~ ds \lesssim \int_1^t \Vert e^{\pm i s \Lambda} \partial_s f \Vert_{W^{1, 3}} \Vert x f \Vert_{L^2} ~ ds \\
&\lesssim \int_1^t \Vert \partial_s f \Vert_{H^2} \Vert x f \Vert_{L^2} ~ ds \lesssim \int_1^t s^{-2+\tau} \Vert u \Vert_X^3 \lesssim \Vert u \Vert_X^3 \\
\Vert \eqref{est_Hx++_term1-2} \Vert_{H^5} &\lesssim \int_1^t \Vert T_{\mu_0}(\Lambda e^{\pm i s \Lambda}  \partial_s f, \Lambda^{-1} u) \Vert_{H^5} ~ ds \lesssim \int_1^t \Vert \partial_s f \Vert_{H^6} \Vert \Lambda^{-1} u \Vert_{W^{6, \infty-}} ~ ds \\
&\lesssim \int_1^t s^{-2+\tau+\gamma/2} \Vert u \Vert_X^2 \Vert u \Vert_{W^{6, 3-}} ~ ds \lesssim \int_1^t s^{-2-1/3+\tau+\gamma/2+\delta} \Vert u \Vert_X^3 ~ ds \lesssim \Vert u \Vert_X^3 \\
\Vert \eqref{est_Hx++_term1-3} \Vert_{H^5} &\lesssim \int_1^t \Vert T_{\mu_0}(\Lambda e^{\pm i s \Lambda} \partial_s x f, u) \Vert_{H^5} ~ ds \lesssim \int_1^t \Vert \partial_s x f \Vert_{H^6} \Vert u \Vert_{W^{6, \infty-}} ~ ds \\
&\lesssim \int_1^t s^{-1+\tau+\gamma/2-2/3+\delta} \Vert u \Vert_X^3 ~ ds \lesssim \Vert u \Vert_X^3 \\
\Vert \eqref{est_Hx++_term1-4} \Vert_{H^5} &\lesssim \int_1^t \Vert T_{\mu_0}(\Lambda e^{\pm i s \Lambda} \partial_s f, e^{\pm i s \Lambda} x f) \Vert_{H^5} ~ ds \lesssim \int_1^t \left( \Vert \partial_s f \Vert_{H^6} \Vert e^{\pm i s \Lambda} x f \Vert_{W^{1, 6}} + \Vert e^{\pm i s \Lambda} \partial_s f \Vert_{W^{1, 6}} \Vert x f \Vert_{H^5} \right) ~ ds \\
&\lesssim \int_1^t \Vert \partial_s f \Vert_{H^6} \Vert x f \Vert_{H^5} ~ ds \lesssim \int_1^t s^{-2+\tau+\gamma/2} \Vert u \Vert_X^3 ~ ds \lesssim \Vert u \Vert_X^3
\end{align*}
\eqref{est_Hx++_term1-1} is easier to estimate~; \eqref{est_Hx++_term1-5}, \eqref{est_Hx++_term1-6}, \eqref{est_Hx++_term1-7} are very similar to previous ones. 

\subsubsection{Estimate of \eqref{est_Hx++_term2}} 
\begin{align*}
\Vert \Lambda^{-1} \eqref{est_Hx++_term2-3} \Vert_{L^2} &\lesssim \Vert T_{\mu_0}(\Lambda e^{\pm i s \Lambda} x f, u) \Vert_{L^{6/5}} \lesssim \Vert x f \Vert_{H^1} \Vert u \Vert_{L^3} \lesssim \Vert u \Vert_X^2 \\
\Vert \eqref{est_Hx++_term2-3} \Vert_{H^5} &\lesssim \Vert x f \Vert_{H^6} \Vert u \Vert_{W^{6, \infty-}} \lesssim \Vert u \Vert_X^2
\end{align*}
\eqref{est_Hx++_term2-1} and \eqref{est_Hx++_term2-init} are easier to estimate~; \eqref{est_Hx++_term2-2}, \eqref{est_Hx++_term2-4} can be estimated in a similar way as \eqref{est_Hx++_term2-3} (by using Hardy's inequality for \eqref{est_Hx++_term2-2}). 

\subsubsection{Estimate of \eqref{est_Hx++_term4}} 
\begin{align*}
\Vert \Lambda^{-1} \eqref{est_Hx++_term4-2} \Vert_{L^2} &\lesssim \int_1^t s^{-1} \Vert T_{\mu_0}(e^{\pm i s \Lambda} xf, \Lambda u) \Vert_{L^{6/5}} ~ ds \lesssim \int_1^t s^{-1} \Vert x f \Vert_{L^2} \Vert u \Vert_{W^{1, 3}} ~ ds \lesssim \int_1^t s^{-4/3} \Vert u \Vert_X^2 ~ ds \lesssim \Vert u \Vert_X^2 \\
\Vert \eqref{est_Hx++_term4-2} \Vert_{H^5} &\lesssim \int_1^t s^{-1} \Vert T_{\mu_0}(e^{\pm i s \Lambda} xf, \Lambda u) \Vert_{H^5} ~ ds \lesssim \int_1^t s^{-1} \Vert x f \Vert_{H^5} \Vert u \Vert_{W^{7, \infty-}} ~ ds \lesssim \int_1^t s^{-2+\delta} \Vert u \Vert_X^2 ~ ds \lesssim \Vert u \Vert_X^2
\end{align*}
Again, \eqref{est_Hx++_term4-1} and \eqref{est_Hx++_term4-3} are similar to \eqref{est_Hx++_term4-2}. 

\subsubsection{Estimate of \eqref{est_Hx++_term5}} 
\begin{align*}
\Vert \Lambda^{-1} \eqref{est_Hx++_term5-1} \Vert_{L^2} &\lesssim \int_1^t s^{-1} \Vert T_{\mu_0}(\Lambda e^{\pm i s \Lambda} x^2 f, \Lambda u) \Vert_{L^{6/5}} ~ ds \lesssim \int_1^t s^{-1} \Vert \Lambda |x|^2 f \Vert_{L^2} \Vert u \Vert_{W^{1, 3}} ~ ds \lesssim \int_1^t s^{-4/3+\gamma} \Vert u \Vert_X^2 ~ ds \lesssim \Vert u \Vert_X^2 \\
\Vert \Lambda^{-1} \eqref{est_Hx++_term5-3} \Vert_{L^2} &\lesssim \int_1^t s^{-1} \Vert T_{\mu_0}(\Lambda^2 e^{\pm i s \Lambda} x f, \Lambda^{-1} u) \Vert_{L^{6/5}} ~ ds \lesssim \int_1^t s^{-1} \Vert \Lambda^2 e^{\pm i s \Lambda} x f \Vert_{L^3} \Vert \Lambda^{-1} u \Vert_{L^2} ~ ds \\
&\lesssim \int_1^t s^{-4/3} \Vert \Lambda^{8/3} x f \Vert_{L^{3/2}} \Vert x f \Vert_{L^2} ~ ds \lesssim \int_1^t s^{-4/3} \Vert \langle x \rangle \Lambda x f \Vert_{H^2} \Vert x f \Vert_{L^2} ~ ds \lesssim \int_1^t s^{-4/3+\gamma} \Vert u \Vert_X^2 ~ ds \lesssim \Vert u \Vert_X^2 \\
\Vert \Lambda^{-1} \eqref{est_Hx++_term5-4} \Vert_{L^2} &\lesssim \int_1^t s^{-1} \Vert T_{\mu_0}(\Lambda^2 e^{\pm i s \Lambda} x^2 f, u) \Vert_{L^{6/5}} ~ ds \lesssim \int_1^t s^{-1} \Vert \Lambda |x|^2 f \Vert_{H^1} \Vert u \Vert_{L^3} ~ ds \lesssim \int_1^t s^{-4/3+\gamma} \Vert u \Vert_X^2 ~ ds \lesssim \Vert u \Vert_X^2 \\
\Vert \Lambda^{-1} \eqref{est_Hx++_term5-5} \Vert_{L^2} &\lesssim \int_1^t s^{-1} \Vert T_{\mu_0}(\Lambda^2 e^{\pm i s \Lambda} x f, e^{\pm i s \Lambda} x f) \Vert_{L^{6/5}} ~ ds \lesssim \int_1^t s^{-1} \Vert e^{\pm i s \Lambda} x f \Vert_{W^{2, 12/5}}^2 ~ ds \\
&\lesssim \int_1^t s^{-4/3} \Vert \Lambda x f \Vert_{W^{2, 12/7}}^2 ~ ds \lesssim \int_1^t s^{-4/3} \Vert \langle x \rangle \Lambda x f \Vert_{H^2}^2 ~ ds \lesssim \int_1^t s^{-4/3+2\gamma} \Vert u \Vert_X^2 ~ ds \lesssim \Vert u \Vert_X^2 \\
~ \\
\Vert \eqref{est_Hx++_term5-1} \Vert_{H^5} &\lesssim \int_1^t s^{-1} \Vert T_{\mu_0}(\Lambda e^{\pm i s \Lambda} x^2 f, \Lambda u) \Vert_{H^5} ~ ds \lesssim \int_1^t s^{-1} \Vert \Lambda |x|^2 f \Vert_{H^5} \Vert u \Vert_{W^{7, \infty-}} ~ ds \lesssim \int_1^t s^{-2+\gamma+\delta} \Vert u \Vert_X^2 ~ ds \lesssim \Vert u \Vert_X^2 \\
\Vert \eqref{est_Hx++_term5-3} \Vert_{H^5} &\lesssim \int_1^t s^{-1} \Vert T_{\mu_0}(\Lambda^2 e^{\pm i s \Lambda} x f, \Lambda^{-1} u) \Vert_{H^5} ~ ds \lesssim \int_1^t s^{-1} \Vert x f \Vert_{H^7} \Vert \Lambda^{-1} u \Vert_{W^{6, \infty-}} ~ ds \\
&\lesssim \int_1^t s^{-1+b+\gamma/2+\varepsilon/2} \Vert u \Vert_X \Vert u \Vert_{W^{6, 3-}} ~ ds \lesssim \int_1^t s^{-4/3+b+\gamma/2+\varepsilon/2} \Vert u \Vert_X^2 ~ ds \lesssim \Vert u \Vert_X^2 \\
\Vert \eqref{est_Hx++_term5-5} \Vert_{H^5} &\lesssim \int_1^t s^{-1} \Vert T_{\mu_0}(\Lambda^2 e^{\pm i s \Lambda} x f, e^{\pm i s \Lambda} x f) \Vert_{H^5} ~ ds \lesssim \int_1^t s^{-1} \Vert x f \Vert_{H^7} \Vert e^{\pm i s \Lambda} x f \Vert_{W^{2, 4}} ~ ds \\
&\lesssim \int_1^t s^{-3/2+\gamma/2+\varepsilon/2} \Vert u \Vert_X \Vert \Lambda x f \Vert_{W^{2, 4/3}} ~ ds \lesssim \int_1^t s^{-3/2+\gamma/2+\varepsilon/2} \Vert \langle x \rangle \Lambda x f \Vert_{H^2} \Vert u \Vert_X ~ ds \\
&\lesssim \int_1^t s^{-3/2+3\gamma/2+\varepsilon/2} \Vert u \Vert_X^2 ~ ds \lesssim \Vert u \Vert_X^2
\end{align*}
where we used the dispersion inequality of lemma \ref{lem_disp} and Hardy's inequality of lemma \ref{lem_Hardy}, as well as lemma \ref{lem_xH7}. Then, \eqref{est_Hx++_term5-2} is similar to \eqref{est_Hx++_term5-5}. 

The term \eqref{est_Hx++_term5-4} is a little particular, because it is responsible for the slight growth when we apply all the derivatives. Indeed, when we distribute $|\xi|^5$ using \eqref{distrib_deriv_symbols}, we have to estimate~: 
\begin{subequations}
\begin{align}
\eqref{est_Hx++_term5-4} &= \int_1^t \int e^{i s \varphi(\xi, \eta)} s^{-1} \mu_0^{0, 0, 7}(\xi, \eta) \nabla^2_{\xi} \widehat{f}(s, \xi - \eta) \widehat{f}(s, \eta) ~ d\eta ds \label{est_Hx++_term5-4-1} \\
&+ \int_1^t \int e^{i s \varphi(\xi, \eta)} s^{-1} \mu_0^{0, 5, 2}(\xi, \eta) \nabla^2_{\xi} \widehat{f}(s, \xi - \eta) \widehat{f}(s, \eta) ~ d\eta ds \label{est_Hx++_term5-4-2}
\end{align}
\end{subequations}
\eqref{est_Hx++_term5-4-2} is easy to estimate~: 
\[ \Vert \eqref{est_Hx++_term5-4-2} \Vert_{L^2} \lesssim \int_1^t s^{-1} \Vert T_{\mu_0}(\Lambda^2 e^{\pm i s \Lambda} x^2 f, \Lambda^5 u) \Vert_{L^2} ~ ds \lesssim \int_1^t s^{-1} \Vert \Lambda |x|^2 f \Vert_{L^2} \Vert u \Vert_{W^{6, \infty-}} ~ ds \lesssim \int_1^t s^{-2+\gamma+\delta} \Vert u \Vert_X^2 ~ ds \lesssim \Vert u \Vert_X^2 \]
However, \eqref{est_Hx++_term5-4-1} contains a factor $\Lambda^7 |x|^2 f$, that we cannot control a priori. Thus, we apply an integration by parts in frequency, noting that $\nabla_{\xi}^2 \widehat{f}(s, \xi - \eta) = - \nabla_{\eta} \nabla_{\xi} \widehat{f}(s, \xi - \eta)$~: 
\begin{subequations}
\begin{align}
\eqref{est_Hx++_term5-4-2}
&= \int_1^t \int e^{i s \varphi(\xi, \eta)} \mu_0^{0, 0, 7}(\xi, \eta) \nabla_{\xi} \widehat{f}(s, \xi - \eta) \widehat{f}(s, \eta) ~ d\eta ds \label{est_Hx++_term5-4-2-1} \\
&+ \int_1^t \int e^{i s \varphi(\xi, \eta)} s^{-1} \mu_0^{0, 0, 6}(\xi, \eta) \nabla_{\xi} \widehat{f}(s, \xi - \eta) \widehat{f}(s, \eta) ~ d\eta ds \label{est_Hx++_term5-4-2-2} \\
&+ \int_1^t \int e^{i s \varphi(\xi, \eta)} s^{-1} \mu_0^{0, -1, 7}(\xi, \eta) \nabla_{\xi} \widehat{f}(s, \xi - \eta) \widehat{f}(s, \eta) ~ d\eta ds \label{est_Hx++_term5-4-2-3} \\
&+ \int_1^t \int e^{i s \varphi(\xi, \eta)} s^{-1} \mu_0^{0, 0, 7}(\xi, \eta) \nabla_{\xi} \widehat{f}(s, \xi - \eta) \nabla_{\xi} \widehat{f}(s, \eta) ~ d\eta ds\label{est_Hx++_term5-4-2-4} 
\end{align}
\end{subequations}
\eqref{est_Hx++_term5-4-2-2} is similar to \eqref{est_Hx++_term4-3}, \eqref{est_Hx++_term5-4-2-3} to \eqref{est_Hx++_term5-3}, \eqref{est_Hx++_term5-4-2-4} to \eqref{est_Hx++_term5-5} and can be treated the same way. The most sensitive term is the one where the frequency derivative hit the exponential and we lost the $s^{-1}$ factor, ie \eqref{est_Hx++_term5-4-2-1}. We then write~: 
\[ \Vert \eqref{est_Hx++_term5-4-2-1} \Vert_{L^2} \lesssim \int_1^t \Vert T_{\mu_0}(\Lambda^7 e^{\pm i s \Lambda} x f, u) \Vert_{L^2} ~ ds \lesssim \int_1^t \Vert \Lambda^7 x f \Vert_{L^2} \Vert u \Vert_{\dot{B}^0_{\infty, 1}} ~ ds \lesssim \int_1^t s^{-1+\gamma/2+\varepsilon/2} \Vert u \Vert_X^2 ~ ds \lesssim t^{\gamma/2+\varepsilon/2} \Vert u \Vert_X^2 \]
Therefore, it is enough to have $b \geq \frac{\gamma+\varepsilon}{2}$ to conclude. 

\begin{Rem}[Amelioration of the growth exponent] By applying lemma \ref{lem_xH7}, we have~: 
\[ \Vert x f \Vert_{\dot{H}^6} \lesssim \Vert u \Vert_{H^7}^{1/2} \Vert \langle x \rangle^2 \Lambda f \Vert_{H^4}^{1/2} \lesssim t^{\varepsilon/2+\gamma_4/2} \Vert u \Vert_X \]
Therefore, we may actually choose $b = \frac{\varepsilon+\gamma_4}{2} < \frac{\varepsilon+\gamma_5}{2}$ ($\gamma_5 = \gamma$ being the maximal growth exponent). \end{Rem}

\subsection{\texorpdfstring{$\pm 0$}{+/- 0} and \texorpdfstring{$0 \pm$}{0 +/-} interactions} 

In both of these cases, we know that $\nabla_{\eta} \varphi$ never vanishes and we can always apply integrations by parts. Therefore, the terms we need to estimate are~: 
\begin{subequations} \label{est_Hx+0_init}
\begin{align}
&\nabla_{\xi} \int_1^t \int e^{i s \varphi(\xi, \eta)} \mu_0^{0, 0, 1}(\xi, \eta) \widehat{f}(s, \xi - \eta) \widehat{f}(s, \eta) ~ d\eta ds \tag{\ref{est_Hx+0_init}} \\
&= \int_1^t \int e^{i s \varphi(\xi, \eta)} s \mu_0^{0, 0, 1}(\xi, \eta) \widehat{f}(s, \xi - \eta) \widehat{f}(s, \eta) ~ d\eta ds \label{est_Hx+0_init_ipp} \\
&+ \int_1^t \int e^{i s \varphi(\xi, \eta)} \mu_0^{-1, 0, 1}(\xi, \eta) \widehat{f}(s, \xi - \eta) \widehat{f}(s, \eta) ~ d\eta ds \label{est_Hx+0_term3} \\
&+ \int_1^t \int e^{i s \varphi(\xi, \eta)} \mu_0(\xi, \eta) \widehat{f}(s, \xi - \eta) \widehat{f}(s, \eta) ~ d\eta ds \label{est_Hx+0_term4} \\
&+ \int_1^t \int e^{i s \varphi(\xi, \eta)} \mu_0^{0, 0, 1}(\xi, \eta) \nabla_{\xi} \widehat{f}(s, \xi - \eta) \widehat{f}(s, \eta) ~ d\eta ds \label{est_Hx+0_term5}
\end{align} 
where one of the $f$ is of type $\pm$, and the other of type $0$. We apply an integration by parts on the term containing the $s$ factor \eqref{est_Hx+0_init_ipp}, so that we get~: 
\begin{align}
\eqref{est_Hx+0_init}
&= \int_1^t \int e^{i s \varphi(\xi, \eta)} \mu_0^{0, -1, 1}(\xi, \eta) \widehat{f}(s, \xi - \eta) \widehat{f}(s, \eta) ~ d\eta ds \label{est_Hx+0_term1} \\
&+ \int_1^t \int e^{i s \varphi(\xi, \eta)} \mu_0^{0, 0, 1}(\xi, \eta) \widehat{f}(s, \xi - \eta) \nabla_{\eta} \widehat{f}(s, \eta) ~ d\eta ds \label{est_Hx+0_term2} \\
&+ \eqref{est_Hx+0_term3} + \eqref{est_Hx+0_term4} + \eqref{est_Hx+0_term5} \nonumber
\end{align}
\end{subequations}

\subsubsection{\texorpdfstring{$\pm 0$}{+/- 0} interactions} 
In this case, we will estimate \eqref{est_Hx+0_term1}, \eqref{est_Hx+0_term2}, \eqref{est_Hx+0_term3}~; the other are simpler. 
\begin{align*}
\Vert \eqref{est_Hx+0_term1} \Vert_{H^6} &\lesssim \int_1^t \Vert T_{\mu_0}(\Lambda u, \Lambda^{-1} u^0) \Vert_{H^6} ~ ds \lesssim \int_1^t \Vert u \Vert_{W^{7, 3}} \Vert \Lambda^{-1} u^0 \Vert_{W^{6,6}} ~ ds \\
&\lesssim \int_1^t s^{-1/3+\delta} \Vert u \Vert_X \Vert u^0 \Vert_{H^6} ~ ds \lesssim \int_1^t s^{-4/3+\delta} \Vert u \Vert_X^2 ~ ds \lesssim \Vert u \Vert_X^2 \\
\Vert \eqref{est_Hx+0_term2} \Vert_{H^6} &\lesssim \int_1^t \Vert T_{\mu_0}(\Lambda u, x u^0) \Vert_{H^6} ~ ds \lesssim \int_1^t \Vert u \Vert_{W^{7, 3}} \Vert x u^0 \Vert_{W^{6, 6}} ~ ds \\
&\lesssim \int_1^t s^{-1/3+\delta} \Vert u \Vert_X \Vert \Lambda x u^0 \Vert_{H^5} ~ ds \lesssim \int_1^t s^{-4/3+\delta+a+b} \Vert u \Vert_X^2 ~ ds \lesssim \Vert u \Vert_X^2 \\
\Vert \eqref{est_Hx+0_term3} \Vert_{L^2} &\lesssim \int_1^t \Vert T_{\mu_0}(\Lambda u, u^0) \Vert_{L^{6/5}} ~ ds \lesssim \int_1^t \Vert u \Vert_{H^1} \Vert u^0 \Vert_{L^3} ~ ds \lesssim \int_1^t s^{-4/3+a} \Vert u \Vert_X^2 ~ ds \lesssim \Vert u \Vert_X^2
\end{align*}
($\Vert \Lambda \eqref{est_Hx+0_term3} \Vert_{H^5}$ is also simpler.) 

\subsubsection{\texorpdfstring{$0 \pm$}{0 +/-} interactions} 
In this case, the estimate of \eqref{est_Hx+0_term3} is similar to the one we just did above. The presence of the $\Lambda$ in front of the term of type $0$ allows a nice control of every term except \eqref{est_Hx+0_term5} which has more derivatives. For \eqref{est_Hx+0_term5}, we recover~: 
\begin{subequations}
\begin{align}
&\xi^6 \int_1^t \int e^{i s \varphi(\xi, \eta)} \mu_0^{0, 0, 1}(\xi, \eta) \nabla_{\xi} \widehat{u^0}(s, \xi - \eta) \widehat{f}(s, \eta) ~ d\eta ds \nonumber \\
&= \int_1^t \int e^{i s \varphi(\xi, \eta)} \mu_0^{0, 0, 7}(\xi, \eta) \nabla_{\xi} \widehat{u^0}(s, \xi - \eta) \widehat{f}(s, \eta) ~ d\eta ds \label{est_Hx+0_term5-1} \\
&+ \int_1^t \int e^{i s \varphi(\xi, \eta)} \mu_0^{0, 6, 1}(\xi, \eta) \nabla_{\xi} \widehat{u^0}(s, \xi - \eta) \widehat{f}(s, \eta) ~ d\eta ds \label{est_Hx+0_term5-2}
\end{align}
\end{subequations}
\eqref{est_Hx+0_term5-2} is simpler. For \eqref{est_Hx+0_term5-1}, we have one derivative too much, so we apply an integration by parts in frequency~:
\begin{subequations}
\begin{align}
\eqref{est_Hx+0_term5-1}
&= \int_1^t \int e^{i s \varphi(\xi, \eta)} s \mu_0^{0, 0, 7}(\xi, \eta) \widehat{u^0}(s, \xi - \eta) \widehat{f}(s, \eta) ~ d\eta ds \label{est_Hx+0_term5-1-1} \\
&+ \int_1^t \int e^{i s \varphi(\xi, \eta)} \mu_0^{0, 0, 6}(\xi, \eta) \widehat{u^0}(s, \xi - \eta) \widehat{f}(s, \eta) ~ d\eta ds \\
&+ \int_1^t \int e^{i s \varphi(\xi, \eta)} \mu_0^{0, -1, 7}(\xi, \eta) \widehat{u^0}(s, \xi - \eta) \widehat{f}(s, \eta) ~ d\eta ds \\
&+ \int_1^t \int e^{i s \varphi(\xi, \eta)} \mu_0^{0, 0, 7}(\xi, \eta) \widehat{u^0}(s, \xi - \eta) \nabla_{\eta} \widehat{f}(s, \eta) ~ d\eta ds
\end{align}
\end{subequations}
All these terms have already been estimated in \eqref{est_Hx+0_term1}, \eqref{est_Hx+0_term2} or \eqref{est_Hx+0_term3}, except the first one \eqref{est_Hx+0_term5-1-1} that has a $s$ factor. For this one, we thus estimate directly~:  
\[ \Vert \eqref{est_Hx+0_term5-1-1} \Vert_{L^2} \lesssim \int_1^t s \Vert T_{\mu_0}(\Lambda^7 u^0, u) \Vert_{L^2} ~ ds \lesssim \int_1^t s \Vert u^0 \Vert_{W^{8, \infty-}} \Vert u \Vert_{L^2} ~ ds \lesssim \int_1^t s^{-1+a+\delta} \Vert u \Vert_X^2 ~ ds \lesssim t^{a+\delta} \Vert u \Vert_X^2 \]
Therefore, we also need $b$ to satisfy $b \geq a + \delta$. However, since $a$ only depends on $\varepsilon, N$, and $\delta$ as well, we may choose $\varepsilon$ small enough and $N$ big enough so that this second condition be weaker than $b \geq \frac{\gamma_4 + \varepsilon}{2}$. 

\subsection{\texorpdfstring{$0 0$}{0 0} interactions}

This time, the nonlinearity has no particular structure, but both $u^0$ factors have a strong decay. We write~: 
\begin{subequations} \label{est_Hx00_init}
\begin{align}
&\nabla_{\xi} \int_1^t \int e^{i s \varphi(\xi, \eta)} \mu_0^{0, 0, 1}(\xi, \eta) \widehat{u^0}(s, \xi - \eta) \widehat{u^0}(s, \eta) ~ d\eta ds \tag{\ref{est_Hx00_init}} \\
&= \int_1^t \int e^{i s \varphi(\xi, \eta)} s \mu_0^{0, 0, 1}(\xi, \eta) \widehat{u^0}(s, \xi - \eta) \widehat{u^0}(s, \eta) ~ d\eta ds \label{est_Hx00_term1} \\
&+ \int_1^t \int e^{i s \varphi(\xi, \eta)} \mu_0(\xi, \eta) \widehat{u^0}(s, \xi - \eta) \widehat{u^0}(s, \eta) ~ d\eta ds \label{est_Hx00_term2}  \\
&+ \int_1^t \int e^{i s \varphi(\xi, \eta)} \mu_0^{-1, 0, 1}(\xi, \eta) \widehat{u^0}(s, \xi - \eta) \widehat{u^0}(s, \eta) ~ d\eta ds \label{est_Hx00_term3} \\
&+ \int_1^t \int e^{i s \varphi(\xi, \eta)} \mu_0^{0, 0, 1}(\xi, \eta) \nabla_{\xi} \widehat{u^0}(s, \xi - \eta) \widehat{u^0}(s, \eta) ~ d\eta ds \label{est_Hx00_term4} 
\end{align}
\end{subequations} 
We then estimate directly~: 
\begin{align*}
\Vert \eqref{est_Hx00_term1} \Vert_{H^6} &\lesssim \int_1^t s \Vert T_{\mu_0}(\Lambda u^0, u^0) \Vert_{H^6} ~ ds \lesssim \int_1^t s \Vert u^0 \Vert_{H^7} \Vert u^0 \Vert_{W^{7, \infty-}} ~ ds \lesssim \int_1^t s^{-2+2a+\delta+\varepsilon} \Vert u \Vert_X^2 ~ ds \lesssim \Vert u \Vert_X^2 \\
\Vert \eqref{est_Hx00_term3} \Vert_{L^2} &\lesssim \int_1^t \Vert T_{\mu_0}(\Lambda u^0, u^0) \Vert_{L^{6/5}} ~ ds \lesssim \int_1^t \Vert u^0 \Vert_{H^1} \Vert u^0 \Vert_{L^3} ~ ds \lesssim \int_1^t s^{-7/3+2a} \Vert u \Vert_X^2 ~ ds \lesssim \Vert u \Vert_X^2
\end{align*}
\eqref{est_Hx00_term2} is simpler, as well as \eqref{est_Hx00_term4} in $H^5$ or \eqref{est_Hx00_term3} in $\dot{H}^6$~; however, if we apply $6$ derivatives on \eqref{est_Hx00_term4}, one bothersome term appears that we need to control~: 
\begin{subequations}
\begin{align}
\xi^6 \eqref{est_Hx00_term4}
&= \int_1^t \int e^{i s \varphi(\xi, \eta)} \mu_0^{0, 0, 7}(\xi, \eta) \nabla_{\xi} \widehat{u^0}(s, \xi - \eta) \widehat{u^0}(s, \eta) ~ d\eta ds \label{est_Hx00_term4-1} \\
&+ \int_1^t \int e^{i s \varphi(\xi, \eta)} \mu_0^{0, 6, 1}(\xi, \eta) \nabla_{\xi} \widehat{u^0}(s, \xi - \eta) \widehat{u^0}(s, \eta) ~ d\eta ds \label{est_Hx00_term4-2}
\end{align}
\end{subequations}
\eqref{est_Hx00_term4-2} is simpler to estimate. As for \eqref{est_Hx00_term4-1}, we apply an integration by parts in frequency~: 
\begin{subequations}
\begin{align}
\eqref{est_Hx00_term4-1}
&= \int_1^t \int e^{i s \varphi(\xi, \eta)} s \mu_0^{0, 0, 7}(\xi, \eta) \widehat{u^0}(s, \xi - \eta) \widehat{u^0}(s, \eta) ~ d\eta ds \label{est_Hx00_term4-1-1} \\
&+ \int_1^t \int e^{i s \varphi(\xi, \eta)} \mu_0^{0, 0, 6}(\xi, \eta) \widehat{u^0}(s, \xi - \eta) \widehat{u^0}(s, \eta) ~ d\eta ds \label{est_Hx00_term4-1-2} \\
&+ \int_1^t \int e^{i s \varphi(\xi, \eta)} \mu_0^{0, -1, 7}(\xi, \eta) \widehat{u^0}(s, \xi - \eta) \widehat{u^0}(s, \eta) ~ d\eta ds \label{est_Hx00_term4-1-3} \\
&+ \int_1^t \int e^{i s \varphi(\xi, \eta)} \mu_0^{0, 0, 7}(\xi, \eta) \widehat{u^0}(s, \xi - \eta) \nabla_{\eta} \widehat{u^0}(s, \eta) ~ d\eta ds \label{est_Hx00_term4-1-4}
\end{align}
\end{subequations}
We now can estimate~: 
\begin{align*}
\Vert \eqref{est_Hx00_term4-1-1} \Vert_{L^2} &\lesssim \int_1^t s \Vert T_{\mu_0}(\Lambda^7 u^0, u^0) \Vert_{L^2} ~ ds \lesssim \int_1^t s \Vert u^0 \Vert_{H^7} \Vert u^0 \Vert_{W^{1, \infty-}} ~ ds \lesssim \int_1^t s^{-2+\varepsilon+2a} \Vert u \Vert_X^2 ~ ds \lesssim \Vert u \Vert_X^2 \\
\Vert \eqref{est_Hx00_term4-1-3} \Vert_{L^2} &\lesssim \int_1^t \Vert T_{\mu_0}(\Lambda^7 u^0, \Lambda^{-1} u^0) \Vert_{L^2} ~ ds \lesssim \int_1^t \Vert u^0 \Vert_{W^{7, 3}} \Vert \Lambda^{-1} u^0 \Vert_{L^6} ~ ds \lesssim \int_1^t \Vert u^0 \Vert_{W^{7, 3}} \Vert u^0 \Vert_{L^2} ~ ds \\
&\lesssim \int_1^t s^{-7/3+\delta+2a} \Vert u \Vert_X^2 ~ ds \lesssim \Vert u \Vert_X^2 \\
\Vert \eqref{est_Hx00_term4-1-4} \Vert_{L^2} &\lesssim \int_1^t \Vert T_{\mu_0}(\Lambda^7 u^0, x u^0) \Vert_{L^2} ~ ds \lesssim \int_1^t \Vert u^0 \Vert_{W^{7, 3}} \Vert x u^0 \Vert_{L^6} ~ ds \lesssim \int_1^t \Vert u^0 \Vert_{W^{7, 3}} \Vert \Lambda x u^0 \Vert_{L^2} ~ ds \\
&\lesssim \int_1^t s^{-7/3+\delta+2a} \Vert u \Vert_X^2 ~ ds \lesssim \Vert u \Vert_X^2 
\end{align*}
\eqref{est_Hx00_term4-1-2} is simpler. 

\paragraph{Conclusion} We showed that
\[ \Vert x f^{\pm}(t) \Vert_{H^5} \leq \Vert x f^{\pm}(1) \Vert_{H^5} + C \Vert u \Vert_X^2 \]
and
\[ \Vert x f^{\pm}(t) \Vert_{H^6} \leq \Vert x f^{\pm}(1) \Vert_{H^6} + C t^b \Vert u \Vert_X^2 \]

We need to satisfy~: 
\[ a + \delta \leq \frac{\gamma_4+\varepsilon}{2} \leq b \]
In particular, we can choose $b$ such that $b < \frac{\gamma}{2}$ if $\varepsilon > 0$ is small enough with respect to $\gamma_5 - \gamma_4$. 

\section{Estimate of the \texorpdfstring{$\dot{B}^1_{\infty, 1}$ and $\dot{B}^0_{\infty, 1}$}{Besov} norm} \label{section_Besov}

\begin{Intuit} The dispersive inequality gives the intuition of a $t^{-1}$ decay for the $L^{\infty}$ norm of $f^{\pm}$, and thus of the Besov $\dot{B}^0_{\infty, 1}$ norm, or likewise for the $\dot{B}^1_{\infty, 1}$ norm. In order to control the $L^1$ we obtain by the dispersive inequality, we apply lemma \ref{lem_L1L2} that controls the $L^1$ norm by $L^2$ norms with weight $x$ or $x^2$, so exactly quantities controlled by the $\Vert u \Vert_X$ norm. 

However, once we applied lemma \ref{lem_L1L2}, we cannot apply any integration by parts in time (since the time integral is outside the $L^2$ norm). We therefore separate the $(\xi, \eta)$ space into two parts, one away from $\mathcal{S}$ the time-resonant set (so that we can apply integrations by parts in time for free before applying lemma \ref{lem_L1L2}), and one in the neighborhood of $\mathcal{S}$ where we can use the identities of lemma \ref{identite_fond} to link $\varphi$ to $\nabla_{\eta} \varphi$ in order to replace integrations by parts in time by integrations by parts in frequency. 
\end{Intuit}

The result we prove is the following~: 

\begin{Prop} For $\epsilon = +$ or $\epsilon = -$, $m = 0$ or $m = 1$, we have 
\[ t \Vert u^{\epsilon} \Vert_{\dot{B}^m_{\infty, 1}} \lesssim \Vert \langle x \rangle^2 f(t = 1) \Vert_{H^4} + \Vert u \Vert_X^2 \left( 1 + \Vert u \Vert_X^3 \right) \]
\label{prop-est-Linf}
\end{Prop} 

Let us write~: 
\[ u^{\pm}(t) = e^{\pm i t \Lambda} f^{\pm}(1) + e^{\pm i t \Lambda} \int_1^t \partial_s f^{\pm}(s) ~ ds \]
We fix $m = 0$ or $m = 1$ in the following. 

\paragraph{Decay of the initial term} The initial time term decays at a rate $t^{-1}$ by lemma \ref{lem_disp}~: 
\[ \begin{aligned}
\Vert e^{\pm i t \Lambda} f^{\pm}(1) \Vert_{\dot{B}^m_{\infty, 1}} &\lesssim t^{-1} \Vert \Lambda^2 f^{\pm}(1) \Vert_{\dot{B}^m_{1, 1}} \lesssim t^{-1} \sum_{j \in \mathbb{Z}} \Vert |x|^2 \varphi_j(D) \Lambda^2 f^{\pm}(1) \Vert_{H^1}^{1/2} \Vert |x| \varphi_j(D) f^{\pm}(1) \Vert_{H^1}^{1/2} \\
&\lesssim t^{-1} \left\Vert \left( \Vert |x|^2 \varphi_j(D) \Lambda^2 f^{\pm}(1) \Vert_{H^1} \right) \right\Vert_{l^1}^{1/2} \left\Vert \left( \Vert |x| \varphi_j(D) \Lambda^2 f^{\pm}(1) \Vert_{H^1} \right) \right\Vert_{l^1}^{1/2} 
\end{aligned} \]
Now, we have that 
\[ \Vert |x|^2 \varphi_j(D) \Lambda^2 f^{\pm}(1) \Vert_{H^1} \lesssim \Vert \varphi_j(D) |x|^2 \Lambda^2 f^{\pm}(1) \Vert_{H^1} + \Vert \widetilde{\varphi}_j(D) f^{\pm}(1) \Vert_{H^1} \]
where $\widetilde{\varphi} = \nabla^2 \varphi$, and $\widetilde{\varphi}_j(\xi) = \widetilde{\varphi}(2^j \xi)$ as usual. Since $\widetilde{\varphi}$ remains localized, we recover (up to an universal constant) Besov spaces for each term. The same kind of computation applies to the term with weight $|x|$. Therefore, by also applying lemma \ref{interp_Besov}~: 
\begin{equation} \Vert e^{\pm i t \Lambda} f^{\pm}(1) \Vert_{\dot{B}^m_{\infty, 1}} \lesssim t^{-1} \Vert \langle x \rangle^2 f^{\pm}(1) \Vert_{H^4} \label{decay_initial_infty} \end{equation}
which gives one term of Proposition \ref{prop-est-Linf}. 

\subsection{\texorpdfstring{$\pm \pm$}{+/- +/-} interactions} 

We start by introducing an angular repartition, depending on the type of interaction considered. 

\paragraph{Angular repartition} The time-resonant set is
\[ \begin{aligned}
&\mathcal{T}^{+, --} = \mathcal{T}^{-, ++} = \{ \xi = \eta = 0 \} \\
&\mathcal{T}^{+, -+} = \mathcal{T}^{-, +-} = \{ \eta = \lambda \xi, \lambda \geq 1 \} \cup \{ \xi = 0 \} \\
&\mathcal{T}^{+, ++} = \mathcal{T}^{-, --} = \{ \eta = \lambda \xi, 0 \leq \lambda \leq 1 \} \cup \{ \xi = 0 \} \\
&\mathcal{T}^{+, +-} = \mathcal{T}^{-, -+} = \{ \eta = \lambda \xi, \lambda \leq 0 \} \cup \{ \xi = 0 \}
\end{aligned} \]
Let us choose the following cutoff functions~: 
\[ \begin{aligned}
&\chi^{+, --}(\xi, \eta) = \chi^{-, ++}(\xi, \eta) = 1 \\
&\chi^{+, -+}(\xi, \eta) = \chi^{-, +-}(\xi, \eta) = \widetilde{\chi}\left( \frac{\xi}{|\xi|_0} \cdot g_0 \frac{\xi - \eta}{|\xi - \eta|_0} \right) \\
&\chi^{+, ++}(\xi, \eta) = \chi^{-, --}(\xi, \eta) = 0 \\
&\chi^{+, +-}(\xi, \eta) = \chi^{-, -+}(\xi, \eta) = \widetilde{\chi}\left( - \frac{\xi}{|\xi|_0} \cdot g_0 \frac{\xi - \eta}{|\xi - \eta|_0} \right) 
\end{aligned} \]
where $\widetilde{\chi}$ is a smooth function, taking values in $[0, 1]$, and such that
\[ \widetilde{\chi}(x) = \left\{ \begin{array}{ll} 0 & \mbox{ if } x \leq -\frac{1}{4} \\
1 & \mbox{ if } x \geq \frac{1}{4} \end{array} \right. \]
and $g_0$ is the matrix associated to the $| \cdot |_0$ norm (cf lemmas \ref{identite_fond} and \ref{identite_Linf}). 

From now on, we denote simply by $\chi$ the angular repartition, and set
\[ \chi_{+} = \chi, ~~~~~ \chi_{-} = 1 - \chi \]
In particular, $\chi_{+}$ and $\chi_{-}$ are symbols of order $0$ (but are not polynomials). 

This choice ensures that, on the support of $\chi_{+}$, $\varphi$ doesn't vanish (except possibly on the axes). Indeed~: 
\begin{itemize}
\item in the $+--$ or $-++$ cases, the time-resonant set is reduced to a single point~; 
\item in the $+-+$ or $-+-$ cases, on the support of $\chi_{+}$, $\xi \cdot g_0 (\xi - \eta) \geq - \frac{1}{4} |\xi|_0 |\xi - \eta|_0$ and since $|\cdot|_0$ is an euclidean norm, this implies that we stay away (angularly) of $\{ \xi - \eta = \lambda \xi, \lambda \leq 0 \}$, that is of $\{ \eta = \lambda \xi, \lambda \geq 1 \}$~; 
\item in the $+++$ or $---$ cases, the support of $\chi_{+}$ is empty~; 
\item in the $++-$ or $--+$ cases, on the support of $\chi_{+}$, $\xi \cdot g_0 (\xi - \eta) \leq \frac{1}{4} |\xi|_0 |\xi -\eta|_0$ so it implies that we stay away (angularly) of $\{ \xi - \eta = \lambda \xi, \lambda \geq 1 \}$, that is of $\{ \eta = \lambda \xi, \lambda \leq 0 \}$. 
\end{itemize}

On the other hand, on the support of $\chi_{-}$, we can apply the identities of lemma \ref{identite_Linf} because their denominators do not vanish (outside possibly the axes). Indeed~: 
\begin{itemize}
\item in the $+--$ or $-++$ cases, the support of $\chi_{-}$ is empty~; 
\item in the $+-+$ or $-+-$ cases, on the support of $\chi_{-}$, we have that $\xi \cdot g_0 (\xi - \eta) \leq \frac{1}{4} |\xi|_0 |\xi - \eta|_0$ so $\frac{\xi}{|\xi|_0} \cdot g_0 \frac{\eta}{|\eta|_0} \geq - \frac{1}{4} |\xi|_0 |\eta|_0$ while the denominator of the identity of lemma \ref{identite_Linf} is
\[ |\xi|_0 |\eta|_0 + \xi \cdot g_0 \eta \geq \frac{3}{4} |\xi|_0 |\eta|_0 \]
\item in the $+++$ or $---$ cases, the denominator of the identity of lemma \ref{identite_Linf} is
\[ |\xi|_0 + |\xi - \eta|_0 + |\eta|_0 \]
which does not vanish outside the axes~; 
\item in the $++-$ or $--+$ cases, on the support of $\chi_{-}$, we have that $\xi \cdot g_0 (\xi - \eta) \geq -\frac{1}{4} |\xi|_0 |\xi - \eta|_0$ and the denominator of the identity of lemma \ref{identite_Linf} is
\[ |\xi|_0 |\xi - \eta|_0 + \xi \cdot g_0 (\xi - \eta) \geq \frac{3}{4} |\xi|_0 |\xi - \eta|_0 \]
\end{itemize}
(NB~: This comes from the fact that we chose the angular repartition by taking into account that we always have $\mathcal{T} \subset \mathcal{S}$, so that we can separate a into an area containing $\mathcal{T}^c \cap \mathcal{S}$ and away from $\mathcal{T}$ (the support of $\chi_{+}$) and an area containing $\mathcal{T}$ and away from $\mathcal{T}^c \cap \mathcal{S}$ (the support of $\chi_{-}$). In this second area, $\nabla_{\eta} \varphi$ doesn't vanish outside of the vanishing set of $\varphi$, which allows the use of the identities of lemma \ref{identite_Linf} without making any singularity appear.) 

We now consider an interaction $\pm \pm$ which we separate into two terms~: 
\[ \begin{aligned}
&\int_1^t \int e^{i s \varphi(\xi, \eta)} \nabla_{\eta} \varphi(\xi, \eta) \mu_0^{0, 0, 1}(\xi, \eta) \chi_{+}(\xi, \eta) \widehat{f}(s, \xi - \eta) \widehat{f}(s, \eta) ~ d\eta ds \\
&+ \int_1^t \int e^{i s \varphi(\xi, \eta)}\nabla_{\eta} \varphi(\xi, \eta) \mu_0^{0, 0, 1}(\xi, \eta) \chi_{-}(\xi, \eta) \widehat{f}(s, \xi - \eta) \widehat{f}(s, \eta) ~ d\eta ds 
\end{aligned} \]

\subsubsection{Estimate away from the time-resonant space} 
We need to estimate in $\dot{B}^1_{\infty, 1}$ the following contributions~: 
\begin{equation} t e^{\pm i t \Lambda} \int_1^t \mathcal{F}^{-1} \int e^{i s \varphi(\xi, \eta)} \nabla_{\eta} \varphi(\xi, \eta) \mu_0^{0, 0, 1}(\xi, \eta) \chi_{+}(\xi, \eta) \widehat{f}(s, \xi - \eta) \widehat{f}(s, \eta) ~ d\eta \label{est_B++_away_init} \end{equation}
Since the support of $\chi_{+}$ stays away from the time-resonant set, we can integrate by parts in time. 
\begin{subequations}
\begin{align}
\eqref{est_B++_away_init} &= t e^{\pm i t \Lambda} \mathcal{F}^{-1} \int e^{i t \varphi(\xi, \eta)} \frac{\nabla_{\eta} \varphi(\xi, \eta)}{\varphi(\xi, \eta)} \mu_0^{0, 0, 1}(\xi, \eta) \chi_{+}(\xi, \eta) \widehat{f}(t, \xi - \eta) \widehat{f}(t, \eta) ~ d\eta \\
&+ t e^{\pm i t \Lambda} \mathcal{F}^{-1} \int e^{i \varphi(\xi, \eta)} \frac{\nabla_{\eta} \varphi(\xi, \eta)}{\varphi(\xi, \eta)} \mu_0^{0, 0, 1}(\xi, \eta) \chi_{+}(\xi, \eta) \widehat{f}(1, \xi - \eta) \widehat{f}(1, \eta) ~ d\eta \\
&+ t e^{\pm i t \Lambda} \mathcal{F}^{-1} \int_1^t \int e^{i \varphi(\xi, \eta)} \frac{\nabla_{\eta} \varphi(\xi, \eta)}{\varphi(\xi, \eta)} \mu_0^{0, 0, 1}(\xi, \eta) \chi_{+}(\xi, \eta) \partial_s \left( \widehat{f}(s, \xi - \eta) \widehat{f}(s, \eta) \right) ~ d\eta \\
&= t T_{\psi}(u, u) - t e^{\pm i t \Lambda} T_{\psi}(u(1), u(1)) \label{est_B++_away_term1} \\
&+ t e^{\pm i t \Lambda} \int_1^t e^{\mp i s \Lambda} T_{\psi}(e^{\pm i s \Lambda} \partial_s f, u) ~ ds \label{est_B++_away_term2} \\
&+ t e^{\pm i t \Lambda} \int_1^t e^{\mp i s \Lambda} T_{\psi}(u, e^{\pm i s \Lambda} \partial_s f) ~ ds \label{est_B++_away_term3}
\end{align}
\end{subequations}
where $\psi = \frac{\nabla_{\eta} \varphi}{\varphi} \mu_0^{0, 0, 1}(\xi, \eta) \chi_{+}$. 

\begin{Lem} We have that
\[ \psi(\xi, \eta) = \mu_0^{-1, 0, 1}(\xi, \eta) + \mu_0(\xi, \eta) \]

We also have
\[ \psi(\xi, \eta) = \mu_0^{0, -1, 1}(\xi, \eta) + \mu_0(\xi, \eta) \]
\label{lem_psi}
\end{Lem}

\begin{Dem}
\textbf{Case }$+--$ (or $-++$)~: $|\varphi(\xi, \eta)| = |\xi|_0 + |\eta|_0 + |\xi - \eta|_0$. Thus, 
\[ \psi(\xi, \eta) = \frac{|\xi - \eta|}{\varphi(\xi, \eta)} \frac{1}{|\xi - \eta|} \mu_0^{0, 0, 1}(\xi, \eta) \]
and $\frac{|\xi - \eta|}{\varphi(\xi, \eta)}$ is a symbol of order $0$. This proves both identities. 

\textbf{Case }$+-+$ (or $-+-$)~: On the support of $\chi_{+}$, $\xi \cdot g_0 (\xi - \eta) \geq - \frac{1}{4} |\xi|_0 |\xi - \eta|_0$. Yet~: 
\[ \epsilon_1 \varphi(\xi, \eta) \left( |\xi|_0 + |\xi - \eta|_0 + |\eta|_0 \right) = |\xi|_0^2 + |\xi - \eta|_0^2 + 2 |\xi|_0 |\xi - \eta|_0 - |\eta|_0^2 = 2 \left( \xi \cdot g_0 (\xi - \eta) + |\xi|_0 |\xi - \eta|_0 \right) \geq \frac{3}{2} |\xi|_0 |\xi - \eta|_0 \]
Thus, 
\[ \frac{\chi_{+}(\xi, \eta)}{\varphi(\xi, \eta)} = \frac{\mu_0 (|\xi|_0 + |\xi - \eta|_0 + |\eta|_0)}{|\xi|_0 |\xi - \eta|_0} = \mu_0^{0, 0, -1} + \mu_0^{-1, 0, 0} + \mu_0^{-1, 1, -1} \]
But the last term is similar to the previous ones by noting that the identity \eqref{distrib_deriv_symbols} holds if we switch the order of the terms, and here in particular : 
\begin{equation} \mu_0^{n, m+M, l} = \mu_0^{n+M, m, l} + \mu_0^{n, m, l+M} \label{distrib_deriv_symbols_back} \end{equation} 
We then obtain the first identity. For the second one, we write that $\epsilon_2 \frac{\xi - \eta}{|\xi - \eta|_0} - \epsilon_3 \frac{\eta}{|\eta|_0}$ vanishes at $\xi = 0$ (because $\epsilon_2 = - \epsilon_3$), so that $\nabla_{\eta} \varphi$ as well and 
\[ \frac{\nabla_{\eta} \varphi}{|\xi|} = \mu_0^{0, -1, 0} + \mu_0^{0, 0, -1} \]
which gives the second identity. 

\textbf{Case }$+++$ (or $---$)~: $\psi \equiv 0$. 

\textbf{Case }$++-$ (or $--+$)~: Again, we write~: 
\[ \epsilon_1 \varphi(\xi, \eta) \left( |\xi|_0 - |\xi - \eta|_0 - |\eta|_0 \right) = |\xi|_0^2 + |\xi - \eta|_0^2 - |\eta|_0^2 - 2 |\xi|_0 |\xi - \eta|_0 = - 2 \left( |\xi|_0 |\xi - \eta|_0 - \xi \cdot g_0 (\xi - \eta) \right) \geq \frac{3}{2} |\xi|_0 |\eta|_0 \]
because on the support of $\chi_{+}$ we have $\xi \cdot g_0 (\xi - \eta) \leq \frac{1}{4} |\xi|_0 |\xi - \eta|_0$. We conclude as in the $+-+$ case. For the second identity, $\epsilon_2 = - \epsilon_3$ so we can procede as in the $+-+$ case. 
\end{Dem}

In \eqref{est_B++_away_term1}, the term depending only on $u(1)$ at initial time can be estimated in a similar way as what was done in \eqref{decay_initial_infty}. In the following, we only consider the other term of \eqref{est_B++_away_term1}. 

To estimate it, we apply again lemma \ref{lem_psi} to write~: 
\[ \eqref{est_B++_away_term1} = t T_{\mu_0}(u, u) + t T_{\mu_0^{0, -1, 1}}(u, u) \]
Therefore~: 
\[ \begin{aligned}
\Vert \eqref{est_B++_away_term1} \Vert_{\dot{B}^m_{\infty, 1}} &\lesssim t \Vert T_{\mu_0}(u, u) \Vert_{\dot{B}^m_{\infty, 1}} + t \Vert T_{\mu_0}(\Lambda^{-1} u, \Lambda u) \Vert_{\dot{B}^m_{\infty, 1}} \lesssim t \Vert T_{\mu_0}(u, u) \Vert_{W^{2, \infty-}} + t \Vert T_{\mu_0}(\Lambda^{-1} u, \Lambda u) \Vert_{W^{2, \infty-}} \\
&\lesssim t \Vert u \Vert_{W^{2, \infty-}} \Vert u \Vert_{L^{\infty-}} + t \Vert \Lambda^{-1} u \Vert_{W^{2, \infty-}} \Vert u \Vert_{W^{3, \infty-}} \lesssim \Vert u \Vert_X^2 + t \Vert u \Vert_{W^{2, 3-}} \Vert u \Vert_{W^{3, \infty-}} 
\lesssim \Vert u \Vert_X^2 \end{aligned} \]

For \eqref{est_B++_away_term2} and \eqref{est_B++_away_term3}, we have terms of the form~: 
\[ t e^{\pm i t \Lambda} \int_1^t e^{\mp i s \Lambda} T_{\psi}(g, h) ~ ds \]
where $g = e^{\pm i s \Lambda} \partial_s f, h = u$ or $g = u, h = e^{\pm i s \Lambda} \partial_s f$. 

We apply the dispersive inequality and lemma \ref{lem_L1L2} to get~: 
\begin{equation} \begin{aligned}
&\left\Vert t e^{\pm i t \Lambda} \int_1^t e^{\mp i s \Lambda} T_{\psi}(g, h) ~ ds \right\Vert_{\dot{B}^m_{\infty, 1}} 
\lesssim \left\Vert \Lambda^2 \int_1^t e^{\mp i s \Lambda} T_{\psi}(g, h) ~ ds \right\Vert_{\dot{B}^m_{1, 1}} \\
&\lesssim \int_1^t \sum_j \Vert |x| \varphi_j(D) e^{\pm i s \Lambda} \Lambda^2 T_{\psi}(g, h) \Vert_{H^1}^{1/2} \Vert |x|^2 \varphi_j(D) e^{\pm i s \Lambda} \Lambda^2 T_{\psi}(g, h) \Vert_{H^1}^{1/2} ~ ds 
\end{aligned} \label{est_B++_away_term23} \end{equation}
As in the computation \eqref{decay_initial_infty}, we may exchange $|x|$ (or $|x|^2$) with the symbol $\varphi_j$, up to constant terms and the introduction of terms without weight~: 
\[ \begin{aligned}
\eqref{est_B++_away_term23} &\lesssim \int_1^t \left( \Vert |x| e^{\pm i s \Lambda} \Lambda^2 T_{\psi}(g, h) \Vert_{\dot{B}^m_{2, 1}} + \Vert e^{\pm i s \Lambda} \Lambda T_{\psi}(g, h) \Vert_{\dot{B}^m_{2, 1}} \right)^{1/2} \\
&\quad \quad \left( \Vert |x|^2 e^{\pm i s \Lambda} \Lambda^2 T_{\psi}(g, h) \Vert_{\dot{B}^m_{2, 1}} + \Vert e^{\pm i s \Lambda} T_{\psi}(g, h) \Vert_{\dot{B}^m_{2, 1}} \right)^{1/2} ~ ds 
\end{aligned} \]
We now apply lemma \ref{interp_Besov} to recover Sobolev spaces~: 
\[ \begin{aligned}
&\eqref{est_B++_away_term23} \lesssim \int_1^t \left( \Vert \Lambda^{-\kappa} |x| e^{\pm i s \Lambda} \Lambda^2 T_{\psi}(g, h) \Vert_{H^1}^{\theta} \Vert |x| e^{\pm i s \Lambda} \Lambda^2 T_{\psi}(g, h) \Vert_{H^{3/2}}^{1 - \theta} + \Vert T_{\psi}(g, h) \Vert_{H^{5/2}} \right)^{1/2} \\
&~~~~~\left( \Vert \Lambda^{-\kappa} |x|^2 e^{\pm i s \Lambda} \Lambda^2 T_{\psi}(g, h) \Vert_{H^1}^{\theta} \Vert |x|^2 e^{\pm i s \Lambda} \Lambda^2 T_{\psi}(g, h) \Vert_{H^{3/2}}^{1 - \theta} + \Vert \Lambda^{-\kappa} T_{\psi}(g, h) \Vert_{H^2} \right)^{1/2} ~ ds \\
&\lesssim \int_1^t \left( \Vert \Lambda^{-\kappa} |x| e^{\pm i s \Lambda} \Lambda^2 T_{\psi}(g, h) \Vert_{H^2} + \Vert \Lambda^{-\kappa} T_{\psi}(g, h) \Vert_{H^3} \right)^{1/2} \left( \Vert \Lambda^{-\kappa} |x|^2 e^{\pm i s \Lambda} \Lambda^2 T_{\psi}(g, h) \Vert_{H^2} + \Vert \Lambda^{-\kappa} T_{\psi}(g, h) \Vert_{H^3} \right)^{1/2} ~ ds
\end{aligned} \]
for a small enough parameter $\kappa$ (in particular $\kappa < 1/2$ here above), and $\theta = \theta(\kappa) \in (0, 1)$. 

We will show that the term with weight $|x|$ decays at a rate $s^{-2+}$ where $2+$ is close to $-2$, and the term with weight $|x|^2$ at a rate $s^{-1+}$ where $-1+$ is close to $-1$, while the terms without any weight decay at a rate $s^{-2+}$, so that we have a total decay of $s^{-3/2+}$ which is integrable if the parameters are chosen small enough. 

By symmetry, we only treat \eqref{est_B++_away_term2} (\eqref{est_B++_away_term3} is very analogous). Denote 
\begin{subequations}
\begin{align}
&\Lambda^{-\kappa} |x| e^{\pm i s \Lambda} \Lambda^2 T_{\psi}(g, h) \label{est_B++_away_term2-a} \\
&\Lambda^{-\kappa} |x|^2 e^{\pm i s \Lambda} \Lambda^2 T_{\psi}(g, h) \label{est_B++_away_term2-b} \\
&\Lambda^{-\kappa} T_{\psi}(g, h) \label{est_B++_away_term2-c} 
\end{align}
\end{subequations}


For \eqref{est_B++_away_term2-a}, we start by separating the terms by using lemma \ref{lem_psi} and distributing the derivatives~: 
\begin{align*} \eqref{est_B++_away_term2-a} &= |\xi|^{-\kappa} \nabla_{\xi} |\xi|^2 \int e^{i s \varphi} \mu_0^{0, 0, 1} \chi_{+} \frac{\nabla_{\eta} \varphi}{\varphi} \partial_s \widehat{f}(s, \xi - \eta) \widehat{f}(s, \eta) ~ d\eta \\
&= \int e^{i s \varphi} \mu_0^{0, -1, 2-\kappa} \partial_s \widehat{f}(s, \xi - \eta) \widehat{f}(s, \eta) ~ d\eta + \int e^{i s \varphi} \mu_0^{0, 2-\kappa, -1} \partial_s \widehat{f}(s, \xi - \eta) \widehat{f}(s, \eta) ~ d\eta \\
&+ \int e^{i s \varphi} \mu_0^{1-\kappa, 0, 1} \nabla_{\xi} \partial_s \widehat{f}(s, \xi - \eta) \widehat{f}(s, \eta) ~ d\eta + \int e^{i s \varphi} \mu_0^{2-\kappa, 0, 0} \nabla_{\xi} \partial_s \widehat{f}(s, \xi - \eta) \widehat{f}(s, \eta) ~ d\eta \\
&+ \int e^{i s \varphi} s \nabla_{\xi} \varphi \mu_0^{2-\kappa, 0, 1} \chi_{+} \frac{\nabla_{\eta} \varphi}{\varphi} \partial_s \widehat{f}(s, \xi - \eta) \widehat{f}(s, \eta) ~ d\eta 
\end{align*}
Then, by the identity of lemma \ref{identite_fond}, we can replace $|\xi| \nabla_{\xi} \varphi$ by $\varphi$ and $|\eta| \nabla_{\eta} \varphi$ (up to symbols of order $0$), thus
\begin{align*} \eqref{est_B++_away_term2-a}
&= \int e^{i s \varphi} \mu_0^{0, -1, 2-\kappa} \partial_s \widehat{f}(s, \xi - \eta) \widehat{f}(s, \eta) ~ d\eta + \int e^{i s \varphi} \mu_0^{0, 2-\kappa, -1} \partial_s \widehat{f}(s, \xi - \eta) \widehat{f}(s, \eta) ~ d\eta \\
&+ \int e^{i s \varphi} \mu_0^{1-\kappa, 0, 1} \nabla_{\xi} \partial_s \widehat{f}(s, \xi - \eta) \widehat{f}(s, \eta) ~ d\eta + \int e^{i s \varphi} \mu_0^{2-\kappa, 0, 0} \nabla_{\xi} \partial_s \widehat{f}(s, \xi - \eta) \widehat{f}(s, \eta) ~ d\eta \\
&+ \int e^{i s \varphi} s \mu_0^{1-\kappa, 0, 1} \chi_{+} \nabla_{\eta} \varphi \partial_s \widehat{f}(s, \xi - \eta) \widehat{f}(s, \eta) ~ d\eta + \int e^{i s \varphi} s \nabla_{\eta} \varphi \mu_0^{1-\kappa, 1, 1} \psi \partial_s \widehat{f}(s, \xi - \eta) \widehat{f}(s, \eta) ~ d\eta 
\end{align*}
On the last two terms carrying a factor $s$, we apply an integration by parts in frequency to get
\begin{subequations}
\begin{align} \eqref{est_B++_away_term2-a}
&= \int e^{i s \varphi} \mu_0^{0, 2-\kappa, -1} \partial_s \widehat{f}(s, \xi - \eta) \widehat{f}(s, \eta) ~ d\eta + \int e^{i s \varphi} \mu_0^{0, -1, 2-\kappa} \partial_s \widehat{f}(s, \xi - \eta) \widehat{f}(s, \eta) ~ d\eta \label{est_B++_away_term2-a-1} \\
&+ \int e^{i s \varphi} \mu_0^{0, 0, 2-\kappa} \nabla_{\xi} \partial_s \widehat{f}(s, \xi - \eta) \widehat{f}(s, \eta) ~ d\eta + \int e^{i s \varphi} \mu_0^{0, 2-\kappa, 0} \nabla_{\xi} \partial_s \widehat{f}(s, \xi - \eta) \widehat{f}(s, \eta) ~ d\eta \label{est_B++_away_term2-a-2} \\
&+ \int e^{i s \varphi} \mu_0^{0, 0, 2-\kappa} \partial_s \widehat{f}(s, \xi - \eta) \nabla_{\eta} \widehat{f}(s, \eta) ~ d\eta + \int e^{i s \varphi} \mu_0^{0, 2-\kappa, 0} \partial_s \widehat{f}(s, \xi - \eta) \nabla_{\eta} \widehat{f}(s, \eta) ~ d\eta \label{est_B++_away_term2-a-3}
\end{align} 
\end{subequations}
where we distributed the derivatives. The right-hand side terms are essentially symmetric to the left-hand side ones. 
Then, by lemma \ref{lem_decnonres}, 
\begin{align*}
\Vert \eqref{est_B++_away_term1} \Vert_{H^2} &= \Vert T_{\mu_0}(\Lambda^{-1} e^{\pm i s \Lambda} \partial_s f, \Lambda^{2-\kappa} u) \Vert_{H^2} \lesssim \Vert \Lambda^{-1} e^{\pm i s \Lambda} \partial_s f \Vert_{W^{2, 6}} \Vert \Lambda^{2-\kappa} u \Vert_{W^{2, 3}} \lesssim \Vert \partial_s f \Vert_{H^2} \Vert u \Vert_{W^{4, 3}} \\
&\lesssim s^{-7/3+\tau+\delta} \Vert u \Vert_X^3 \\
\Vert \eqref{est_B++_away_term2-a-2} \Vert_{H^2} &= \Vert T_{\mu_0}(\Lambda^{2 - \kappa} e^{\pm i s \Lambda} x \partial_s f, u) \Vert_{H^2} \lesssim \Vert x \partial_s f \Vert_{H^4} \Vert u \Vert_{W^{3, \infty-}} \lesssim s^{-2+\tau+\delta} \Vert u \Vert_X^3 \\
\Vert \eqref{est_B++_away_term2-a-3} \Vert_{H^2} &= \Vert T_{\mu_0}(\Lambda^{2-\kappa} e^{\pm i s \Lambda} \partial_s f, e^{\pm i s \Lambda} x f) \Vert_{H^2} \lesssim \Vert \partial_s f \Vert_{H^4} \Vert e^{\pm i s \Lambda} x f \Vert_{W^{3, 6}} \lesssim \Vert \partial_s f \Vert_{H^4} \Vert x f \Vert_{H^4} \lesssim s^{-2+\tau} \Vert u \Vert_X^2
\end{align*}

In order to treat \eqref{est_B++_away_term2-b}, let us notice that $\eqref{est_B++_away_term2-b} = |x| \eqref{est_B++_away_term2-a} + \Lambda^{-1} \eqref{est_B++_away_term2-a}$ (up to symbols of order $0$). Since we want to estimate these in norms based on $L^2$, by lemma \ref{lem_Hardy}, we only need to consider $|x| \eqref{est_B++_away_term2-a}$ and we can reuse the previous decomposition
\begin{align*} \eqref{est_B++_away_term2-b}
&= \nabla_{\xi} \int e^{i s \varphi} \mu_0^{0, 2-\kappa, -1} \partial_s \widehat{f}(s, \xi - \eta) \widehat{f}(s, \eta) ~ d\eta + \nabla_{\xi} \int e^{i s \varphi} \mu_0^{0, -1, 2-\kappa} \partial_s \widehat{f}(s, \xi - \eta) \widehat{f}(s, \eta) ~ d\eta \\
&+ \nabla_{\xi} \int e^{i s \varphi} \mu_0^{0, 0, 2-\kappa} \nabla_{\xi} \partial_s \widehat{f}(s, \xi - \eta) \widehat{f}(s, \eta) ~ d\eta + \nabla_{\xi} \int e^{i s \varphi} \mu_0^{0, 2-\kappa, 0} \nabla_{\eta} \partial_s \widehat{f}(s, \xi - \eta) \widehat{f}(s, \eta) ~ d\eta \\
&+ \nabla_{\xi} \int e^{i s \varphi} \mu_0^{0, 0, 2-\kappa} \partial_s \widehat{f}(s, \xi - \eta) \nabla_{\eta} \widehat{f}(s, \eta) ~ d\eta + \nabla_{\xi} \int e^{i s \varphi} \mu_0^{0, 0, 2-\kappa} \partial_s \widehat{f}(s, \xi - \eta) \nabla_{\eta} \widehat{f}(s, \eta) ~ d\eta
\end{align*}
If the derivative hits the exponential, only a factor $s$ appears (up to a symbol of order $0$) and we can apply the same estimates as in \eqref{est_B++_away_term2-a} with a decay $s^{-1+\tau}$ instead of $s^{-2+\tau}$, as desired, where $\tau = \tau(\varepsilon, N)$. It then only remains~: 
\begin{subequations} \label{est_B++_away_term2-b-1}
\begin{align}
&\int e^{i s \varphi} \mu_0^{-1, 2-\kappa, -1} \partial_s \widehat{f}(s, \xi - \eta) \widehat{f}(s, \eta) ~ d\eta \\
&+ \int e^{i s \varphi} \mu_0^{-1, -1, 2-\kappa} \partial_s \widehat{f}(s, \xi - \eta) \widehat{f}(s, \eta) ~ d\eta \\
&+ \int e^{i s \varphi} \mu_0^{-1, 0, 2-\kappa} \nabla_{\xi} \partial_s \widehat{f}(s, \xi - \eta) \widehat{f}(s, \eta) ~ d\eta \\
&+ \int e^{i s \varphi} \mu_0^{-1, 2-\kappa, 0} \nabla_{\eta} \partial_s \widehat{f}(s, \xi - \eta) \widehat{f}(s, \eta) ~ d\eta \\
&+ \int e^{i s \varphi} \mu_0^{-1, 0, 2-\kappa} |\xi - \eta| \partial_s \widehat{f}(s, \xi - \eta) \nabla_{\eta} \widehat{f}(s, \eta) ~ d\eta \\
&+ \int e^{i s \varphi} \mu_0^{-1, 2-\kappa, 0} \partial_s \widehat{f}(s, \xi - \eta) \nabla_{\eta} \widehat{f}(s, \eta) ~ d\eta \\
\nextParentEquation[est_B++_away_term2-b-2]
&+ \int e^{i s \varphi} \mu_0^{0, 2-\kappa, -2} \partial_s \widehat{f}(s, \xi - \eta) \widehat{f}(s, \eta) ~ d\eta \\
&+ \int e^{i s \varphi} \mu_0^{0, -1, 1-\kappa} \partial_s \widehat{f}(s, \xi - \eta) \widehat{f}(s, \eta) ~ d\eta \\
&+ \int e^{i s \varphi} \mu_0^{0, 0, 1-\kappa} \nabla_{\xi} \partial_s \widehat{f}(s, \xi - \eta) \widehat{f}(s, \eta) ~ d\eta \\
&+ \int e^{i s \varphi} \mu_0^{0, 2-\kappa, -1} \nabla_{\eta} \partial_s \widehat{f}(s, \xi - \eta) \widehat{f}(s, \eta) ~ d\eta \\
&+ \int e^{i s \varphi} \mu_0^{0, 0, 1-\kappa} \partial_s \widehat{f}(s, \xi - \eta) \nabla_{\eta} \widehat{f}(s, \eta) ~ d\eta \\
&+ \int e^{i s \varphi} \mu_0^{0, 2-\kappa, -1} \partial_s \widehat{f}(s, \xi - \eta) \nabla_{\eta} \widehat{f}(s, \eta) ~ d\eta \\
\nextParentEquation[est_B++_away_term2-b-3]
&+ \int e^{i s \varphi} \mu_0^{0, 2-\kappa, -1} \partial_s \nabla_{\xi} \widehat{f}(s, \xi - \eta) \widehat{f}(s, \eta) ~ d\eta \label{est_B++_away_term2-b-3-1} \\
&+ \int e^{i s \varphi} \mu_0^{0, -1, 2-\kappa} \partial_s \nabla_{\xi} \widehat{f}(s, \xi - \eta) \widehat{f}(s, \eta) ~ d\eta \label{est_B++_away_term2-b-3-2} \\
&+ \int e^{i s \varphi} \mu_0^{0, 0, 2-\kappa} \nabla^2_{\xi} \partial_s \widehat{f}(s, \xi - \eta) \widehat{f}(s, \eta) ~ d\eta \label{est_B++_away_term2-b-3-3} \\
&+ \int e^{i s \varphi} \mu_0^{0, 2-\kappa, 0} \nabla^2_{\eta} \partial_s \widehat{f}(s, \xi - \eta) \widehat{f}(s, \eta) ~ d\eta \label{est_B++_away_term2-b-3-4} \\
&+ \int e^{i s \varphi} \mu_0^{0, 0, 2-\kappa} \partial_s \nabla_{\xi} \widehat{f}(s, \xi - \eta) \nabla_{\eta} \widehat{f}(s, \eta) ~ d\eta \label{est_B++_away_term2-b-3-5} \\
&+ \int e^{i s \varphi} \mu_0^{0, 2-\kappa, 0} \partial_s \nabla_{\xi} \widehat{f}(s, \xi - \eta) \nabla_{\eta} \widehat{f}(s, \eta) ~ d\eta \label{est_B++_away_term2-b-3-6}
\end{align}
\end{subequations} 
\eqref{est_B++_away_term2-b-1} corresponds to the case when the derivative hits the symbol and a $|\xi|^{-1}$ appears. We can then apply the same estimates as for \eqref{est_B++_away_term2-a}, but in $L^{6/5}$ (using lemma \ref{lem_intfrac}), which gives a decay of order $s^{-4/3}$. \eqref{est_B++_away_term2-b-2} corresponds to the case when the derivative hits the symbol and a $|\xi - \eta|^{-1}$ appears. All the terms are simple to control (in particular when there already is a $|\xi - \eta|$ to absorb the singularity), except the following two~: 
\[ \int e^{i s \varphi} \mu_0^{0, 2-\kappa, -2} \partial_s \widehat{f}(s, \xi - \eta) \widehat{f}(s, \eta) ~ d\eta + \int e^{i s \varphi} \mu_0^{0, 2-\kappa, -1} \nabla_{\eta} \partial_s \widehat{f}(s, \xi - \eta) \widehat{f}(s, \eta) ~ d\eta  \]
But~: 
\begin{align*}
\Vert T_{\mu_0}(\Lambda^{-2} e^{\pm i s \Lambda} \partial_s f, \Lambda^{2 - \kappa} u) \Vert_{H^2} &\lesssim \Vert \Lambda^{-2} e^{\pm i s \Lambda} \partial_s f \Vert_{W^{2, 6}} \Vert u \Vert_{W^{4, 3}} \lesssim \Vert \Lambda^{-1} \partial_s f \Vert_{H^2} \Vert u \Vert_{W^{4, 3}} \lesssim \Vert \partial_s x f \Vert_{H^2} \Vert u \Vert_{W^{4, 3}} \\
&\lesssim s^{-4/3+\tau+\delta} \Vert u \Vert_X^3 \\
\Vert T_{\mu_0}(\Lambda^{-1} e^{\pm i s \Lambda} x \partial_s f, \Lambda^{2 - \kappa} u) \Vert_{H^2} &\lesssim \Vert \Lambda^{-1} e^{\pm i s \Lambda} x \partial_s f \Vert_{W^{2, 6}} \Vert u \Vert_{W^{4, 3}} \lesssim \Vert \partial_s x f \Vert_{H^2} \Vert u \Vert_{W^{4, 3}} \lesssim s^{-4/3+\tau+\delta} \Vert u \Vert_X^3
\end{align*}
Finally, in \eqref{est_B++_away_term2-b-3}, \eqref{est_B++_away_term2-b-3-1} and \eqref{est_B++_away_term2-b-3-2} are similar to terms from \eqref{est_B++_away_term2-a} or \eqref{est_B++_away_term2-b-2}. For \eqref{est_B++_away_term2-b-3-5}, we write that 
\[ \Vert T_{\mu_0}(\Lambda^{2 - \kappa} e^{\pm i s \Lambda} x \partial_s f, e^{\pm i s \Lambda} x f) \Vert_{H^2} \lesssim \Vert x \partial_s f \Vert_{H^4} \Vert e^{\pm i s \Lambda} x f \Vert_{W^{3, 4}} \lesssim s^{-3/2+\tau} \Vert u \Vert_X^2 \Vert \Lambda^2 x f \Vert_{W^{3, 4/3}} \lesssim s^{-3/2+\tau+\gamma} \Vert u \Vert_X^3 \]
and symmetrically for \eqref{est_B++_away_term2-b-3-6}. Finally, for \eqref{est_B++_away_term2-b-3-3} and \eqref{est_B++_away_term2-b-3-4}, we notice that $\nabla_{\xi} \widehat{f}(s, \xi - \eta) = - \nabla_{\eta} \widehat{f}(s, \xi - \eta)$ which allows one integration by parts in frequency on each of the terms. We then obtain terms that were already treated. 

Concerning \eqref{est_B++_away_term2-c}, we have~: 
\[ \begin{aligned}
\eqref{est_B++_away_term2-c} &= \int e^{i s \varphi} \mu_0^{-\kappa, 0, 1} \frac{\nabla_{\eta} \varphi}{\varphi} \partial_s \widehat{f}(s, \xi - \eta) \widehat{f}(s, \eta) ~ d\eta \\
&= \int e^{i s \varphi} \mu_0^{-1-\kappa, 0, 1} \partial_s \widehat{f}(s, \xi - \eta) \widehat{f}(s, \eta) ~ d\eta + \int e^{i s \varphi} \mu_0^{-\kappa, 0, 0} \partial_s \widehat{f}(s, \xi - \eta) \widehat{f}(s, \eta) ~ d\eta 
\end{aligned} \]
Therefore, by lemma \ref{lem_intfrac}, denoting $q_0 = q_0(\kappa)$ such that $\frac{1}{q_0} = \frac{1}{2} + \frac{\kappa}{3}$ and $q_1 = q_1(\kappa)$ such that $\frac{1}{q_1} = \frac{1}{2} + \frac{1+\kappa}{3}$~: 
\[ \Vert \eqref{est_B++_away_term2-c} \Vert_{H^3} \lesssim \Vert T_{\mu_0}(\Lambda e^{\pm i s \Lambda} \partial_s f, u) \Vert_{W^{q_1, 3}} + \Vert T_{\mu_0}(e^{\pm i s \Lambda} \partial_s f, u) \Vert_{W^{q_0, 3}} \lesssim \Vert \partial_s f \Vert_{H^4} \left( \Vert u \Vert_{W^{3/\kappa, 3}} + \Vert u \Vert_{W^{3/(1+\kappa), 3}} \right) \lesssim s^{-2+\tau} \Vert u \Vert_X^3 \]
for small $\kappa$. 

\subsubsection{Estimate close to the time-resonant set} 
On the time-resonant set, we need to estimate the following term in $\dot{B}^m_{\infty, 1}$~: 
\[ t e^{\pm i t \Lambda} \int_1^t e^{\mp i s \Lambda} T_{\nabla_{\eta} \varphi \mu_0 \chi_{-}}(\Lambda u, u) ~ ds \]
Again, we can apply the dispersive inequality to get 
\[ \begin{aligned}
&\int_1^t \left( \Vert \Lambda^{-\kappa} x e^{\mp i s \Lambda} \Lambda^2 T_{\nabla_{\eta} \varphi \mu_0 \chi_{-}}(\Lambda u, u) \Vert_{H^2} + \Vert \Lambda^{-\kappa} T_{\nabla_{\eta} \varphi \mu_0 \chi_{-}}(\Lambda u, u) \Vert_{H^2} \right)^{1/2} \\
&~~~~~\left( \Vert \Lambda^{-\kappa} |x|^2 e^{\mp i s \Lambda} \Lambda^2 T_{\nabla_{\eta} \varphi \mu_0 \chi_{-}}(\Lambda u, u) \Vert_{H^2} + \Vert \Lambda^{-\kappa} T_{\nabla_{\eta} \varphi \mu_0 \chi_{-}}(\Lambda u, u) \Vert_{H^2} \right)^{1/2} ~ ds 
\end{aligned} \]
Again, we denote : 
\begin{subequations}
\begin{align}
&x e^{\mp i s \Lambda} \Lambda^2 T_{\nabla_{\eta} \varphi \mu_0 \chi_{-}}(\Lambda u, u) \label{est_B++_close_term1} \\
&|x|^2 e^{\mp i s \Lambda} \Lambda^2 T_{\nabla_{\eta} \varphi \mu_0 \chi_{-}}(\Lambda u, u) \label{est_B++_close_term2} \\
&T_{\nabla_{\eta} \varphi \mu_0 \chi_{-}}(\Lambda u, u) \label{est_B++_close_term3} 
\end{align}
\end{subequations} 
and we need to control these terms with a $\Lambda^{-\kappa}$ in front. 

\begin{Intuit} The idea is the following~: the most troublesome term will be the one where the $x$ weight, which acts as a derivative in frequency, will hit the exponential and a factor $s$ appears, but with $\nabla_{\xi} \varphi$. Consequently, using lemma \ref{identite_fond}, we can replace it by $\nabla_{\eta} \varphi$ and $\varphi$~; furthermore, by lemma \ref{identite_Linf}, we can replace $\varphi$ by $\nabla_{\eta} \varphi$. Finally, we have enough $\nabla_{\eta} \varphi$ factors to apply integrations by parts and obtain a decay close to $s^{-2}$ on the term with weight $x$, and close to $s^{-1}$ on the term with weight $|x|^2$, as before. 
\end{Intuit}

\paragraph{Weight $x$} Let us begin with the term with weight $x$ \eqref{est_B++_close_term1}. We write it in Fourier~:
\begin{align*}
\eqref{est_B++_close_term1}
&= \int e^{i s \varphi(\xi, \eta)} \nabla_{\eta} \varphi(\xi, \eta) \mu_0^{1, 0, 1}(\xi, \eta) \widehat{f}(s, \xi - \eta) \widehat{f}(s, \eta) ~ d\eta \\
&+ \int e^{i s \varphi(\xi, \eta)} \nabla_{\eta} \varphi(\xi, \eta) \mu_0^{2, 0, 0}(\xi, \eta) \widehat{f}(s, \xi - \eta) \widehat{f}(s, \eta) ~ d\eta \\
&+ \int e^{i s \varphi(\xi, \eta)} |\xi|_0 \nabla_{\xi} \nabla_{\eta} \varphi(\xi, \eta) \mu_0^{1, 0, 1}(\xi, \eta) \chi_{-}(\xi, \eta) \widehat{f}(s, \xi - \eta) \widehat{f}(s, \eta) ~ d\eta \\
&+ \int e^{i s \varphi(\xi, \eta)} \nabla_{\eta} \varphi(\xi, \eta) \mu_0^{2, 0, 1}(\xi, \eta) \chi_{-}(\xi, \eta) \nabla_{\xi} \widehat{f}(s, \xi - \eta) \widehat{f}(s, \eta) ~ d\eta \\
&+ \int e^{i s \varphi(\xi, \eta)} s |\xi|_0 \nabla_{\xi} \varphi(\xi, \eta) \nabla_{\eta} \varphi(\xi, \eta) \mu_0^{1, 0, 1}(\xi, \eta) \chi_{-}(\xi, \eta) \widehat{f}(s, \xi - \eta) \widehat{f}(s, \eta) ~ d\eta 
\end{align*}
Let us apply lemma \ref{identite_fond} on the last line and apply integrations by parts~: 
\begin{align*}
&\int e^{i s \varphi(\xi, \eta)} s |\xi|_0 \nabla_{\xi} \varphi(\xi, \eta) \nabla_{\eta} \varphi(\xi, \eta) \mu_0^{1, 0, 1}(\xi, \eta) \chi_{-}(\xi, \eta) \widehat{f}(s, \xi - \eta) \widehat{f}(s, \eta) ~ d\eta \\
&= - \int e^{i s \varphi(\xi, \eta)} s \epsilon_1 \epsilon_3 |\eta|_0 \nabla_{\eta} \varphi(\xi, \eta) \nabla_{\eta} \varphi \mu_0^{1, 0, 1}(\xi, \eta) \chi_{-}(\xi, \eta) \widehat{f}(s, \xi - \eta) \widehat{f}(s, \eta) ~ d\eta \\
&+ \int e^{i s \varphi(\xi, \eta)} s \varphi(\xi, \eta) \nabla_{\eta} \varphi(\xi, \eta) \mu_0^{1, 0, 1}(\xi, \eta) \chi_{-}(\xi, \eta) \widehat{f}(s, \xi - \eta) \widehat{f}(s, \eta) ~ d\eta \\
&= \int e^{i s \varphi(\xi, \eta)} \epsilon_1 \epsilon_3 \nabla_{\eta} \left( |\eta|_0 \nabla_{\eta} \varphi(\xi, \eta) \right) \mu_0^{1, 0, 1}(\xi, \eta) \chi_{-}(\xi, \eta) \widehat{f}(s, \xi - \eta) \widehat{f}(s, \eta) ~ d\eta \\
&+ \int e^{i s \varphi(\xi, \eta)} \nabla_{\eta} \varphi(\xi, \eta) \mu_0^{1, 0, 1}(\xi, \eta) \widehat{f}(s, \xi - \eta) \widehat{f}(s, \eta) ~ d\eta \\
&+ \int e^{i s \varphi(\xi, \eta)} \nabla_{\eta} \varphi(\xi, \eta) \mu_0^{1, 1, 1}(\xi, \eta) \nabla_{\eta} \widehat{f}(s, \xi - \eta) \widehat{f}(s, \eta) ~ d\eta \\
&+ \int e^{i s \varphi(\xi, \eta)} \varphi(\xi, \eta) \mu_0^{1, 0, 0}(\xi, \eta) \chi_{-}(\xi, \eta) \widehat{f}(s, \xi - \eta) \widehat{f}(s, \eta) ~ d\eta \\
&+ \int e^{i s \varphi(\xi, \eta)} \varphi(\xi, \eta) \mu_0^{1, -1, 1}(\xi, \eta) \chi_{-}(\xi, \eta) \widehat{f}(s, \xi - \eta) \widehat{f}(s, \eta) ~ d\eta \\
&+ \int e^{i s \varphi(\xi, \eta)} \varphi(\xi, \eta) \mu_0^{1, 0, 1}(\xi, \eta) \chi_{-}(\xi, \eta) \nabla_{\eta} \widehat{f}(s, \xi - \eta) \widehat{f}(s, \eta) ~ d\eta \\
&+ \int e^{i s \varphi(\xi, \eta)} \varphi(\xi, \eta) \mu_0^{1, 0, 1}(\xi, \eta) \chi_{-}(\xi, \eta) \widehat{f}(s, \xi - \eta) \nabla_{\eta} \widehat{f}(s, \eta) ~ d\eta
\end{align*}
We can regroup the first line here above with the corresponding term in the total expression to obtain~:
\[ \int e^{i s \varphi(\xi, \eta)} \nabla_{\eta} \left( |\xi|_0 \nabla_{\xi} + \epsilon_1 \epsilon_3 |\eta|_0 \nabla_{\eta} \right) \varphi(\xi, \eta) \mu_0^{1, 0, 1}(\xi, \eta) \chi_{-}(\xi, \eta) \widehat{f}(s, \xi - \eta) \widehat{f}(s, \eta) ~ d\eta \]
and we apply again lemma \ref{identite_fond} to simplify~: 
\[ \nabla_{\eta} \left( |\xi|_0 \nabla_{\xi} + \epsilon_1 \epsilon_3 |\eta|_0 \nabla_{\eta} \right) \varphi(\xi, \eta) = \mu_0 \nabla_{\eta} \varphi + \mu_0^{0, 0, -1} \varphi \]
On the other hand, for all terms containing $\varphi \chi_{-}$, we can apply lemma \ref{identite_Linf} to obtain instead a factor $\nabla_{\eta} \varphi$. More precisely, let us exhaust the different cases. In the $+++$ or $---$ cases, we have that 
\[ |\xi| \varphi(\xi, \eta) = \mu_0^{0, 1, 1}(\xi, \eta) \nabla_{\eta} \varphi(\xi, \eta) \]
because $\frac{|\xi|}{|\eta| + |\xi - \eta|}$ is a symbol of order $0$. In the $++-$ or $-++$ cases, lemma \ref{identite_Linf} and the fact that we are on the support of $\chi_{-}$, where $\xi \cdot g_0 (\xi - \eta) \geq - \frac{1}{4} |\xi|_0 |\xi - \eta|_0$ and thus $|\xi|_0 |\xi - \eta|_0 + \xi \cdot g_0 (\xi - \eta) \geq \frac{3}{4} |\xi|_0 |\xi - \eta|_0$, means that
\[ \begin{aligned}
|\xi| \varphi(\xi, \eta) \chi_{-}(\xi, \eta) &= \mu_0(\xi, \eta) \chi_{-}(\xi, \eta) \frac{(|\xi|_0 + |\eta|_0 + |\xi - \eta|_0)|\xi| |\xi - \eta|_0 |\eta|_0 }{|\xi|_0 |\xi - \eta|_0 + \xi \cdot g_0 (\xi - \eta)} \nabla_{\eta} \varphi(\xi, \eta) \\
&= \left( \mu_0^{1, 1, 0} + \mu_0^{0, 1, 1} + \mu_0^{0, 2, 0} \right) \nabla_{\eta} \varphi(\xi, \eta) 
\end{aligned}\]
Finally, in the $+-+$ or $-+-$ cases, we have in a symmetric way~:
\[ |\xi| \varphi(\xi, \eta) \chi_{-}(\xi, \eta) = \chi_{-}(\xi, \eta) \left( \mu_0^{1, 0, 1} + \mu_0^{0, 0, 2} + \mu_0^{0, 1, 1} \right) \nabla_{\eta} \varphi(\xi, \eta) \]
but we can notice that, on the support of $\chi_{-}$, $\xi \cdot g_0 (\xi - \eta) \leq \frac{1}{4} |\xi|_0 |\xi - \eta|_0$, and therefore
\[ |\eta|_0^2 = |\xi|_0^2 + |\xi - \eta|_0^2 - 2 \xi \cdot g_0 (\xi - \eta) \geq |\xi|_0^2 + |\xi - \eta|_0^2 - \frac{1}{2} |\xi|_0 |\xi - \eta|_0 \geq \frac{3}{4} \left( |\xi|_0^2 + |\xi - \eta|_0^2 \right) \]
so that $|\eta|_0 \geq \frac{\sqrt{3}}{2} \max\left( |\xi|_0, |\xi - \eta|_0 \right)$. In particular, we can write that 
\[ |\xi| \varphi(\xi, \eta) \chi_{-}(\xi, \eta) = \left( \mu_0^{1, 1, 0} + \mu_0^{0, 1, 1} + \mu_0^{0, 2, 0} \right) \nabla_{\eta} \varphi(\xi, \eta) \]
The $+--$ and $-++$ cases are empty ($\chi_{-} \equiv 0$). In any case, we can therefore write that
\[ |\xi| \varphi(\xi, \eta) \chi_{-}(\xi, \eta) = \left( \mu_0^{0, 1, 1} + \mu_0^{0, 2, 0} \right) \nabla_{\eta} \varphi(\xi, \eta) \]
by distributing via \eqref{distrib_deriv_symbols}. Summing up, we get
\[ \begin{aligned}
\eqref{est_B++_close_term1}
&= \int e^{i s \varphi(\xi, \eta)} \nabla_{\eta} \varphi(\xi, \eta) \mu_0^{0, 0, 2}(\xi, \eta) \widehat{f}(s, \xi - \eta) \widehat{f}(s, \eta) ~ d\eta \\
&+ \int e^{i s \varphi(\xi, \eta)} \nabla_{\eta} \varphi(\xi, \eta) \mu_0^{0, 0, 3}(\xi, \eta) \nabla_{\xi} \widehat{f}(s, \xi - \eta) \widehat{f}(s, \eta) ~ d\eta \\
&+ \int e^{i s \varphi(\xi, \eta)} \nabla_{\eta} \varphi(\xi, \eta) \mu_0^{0, 2, 1}(\xi, \eta) \nabla_{\xi} \widehat{f}(s, \xi - \eta) \widehat{f}(s, \eta) ~ d\eta
\end{aligned} \]
by taking symmetries, redundancy and distributions of derivatives into account. We can then apply a last integration by parts in frequency~: 
\begin{subequations}
\begin{align}
\eqref{est_B++_close_term1}
&= s^{-1} \int e^{i s \varphi(\xi, \eta)} \mu_0^{0, 0, 1}(\xi, \eta) \widehat{f}(s, \xi - \eta) \widehat{f}(s, \eta) ~ d\eta \label{est_B++_close_term1-1} \\
&+ s^{-1} \int e^{i s \varphi(\xi, \eta)} \mu_0^{0, -1, 2}(\xi, \eta) \widehat{f}(s, \xi - \eta) \widehat{f}(s, \eta) ~ d\eta \label{est_B++_close_term1-2} \\
&+ s^{-1} \int e^{i s \varphi(\xi, \eta)} \mu_0^{0, 0, 2}(\xi, \eta) \nabla_{\eta} \widehat{f}(s, \xi - \eta) \widehat{f}(s, \eta) ~ d\eta \label{est_B++_close_term1-3} \\
&+ s^{-1} \int e^{i s \varphi(\xi, \eta)} \mu_0^{0, 2, 0}(\xi, \eta) \nabla_{\eta} \widehat{f}(s, \xi - \eta) \widehat{f}(s, \eta) ~ d\eta \label{est_B++_close_term1-4} \\
&+ s^{-1} \int e^{i s \varphi(\xi, \eta)} \mu_0^{0, -1, 3}(\xi, \eta) \nabla_{\xi} \widehat{f}(s, \xi - \eta) \widehat{f}(s, \eta) ~ d\eta \label{est_B++_close_term1-5} \\
&+ s^{-1} \int e^{i s \varphi(\xi, \eta)} \mu_0^{0, 0, 3}(\xi, \eta) \nabla_{\xi} \widehat{f}(s, \xi - \eta) \nabla_{\eta} \widehat{f}(s, \eta) ~ d\eta \label{est_B++_close_term1-6} \\
&+ s^{-1} \int e^{i s \varphi(\xi, \eta)} \mu_0^{0, 0, 3}(\xi, \eta) \nabla^2_{\xi} \widehat{f}(s, \xi - \eta) \widehat{f}(s, \eta) ~ d\eta \label{est_B++_close_term1-7} \\
&+ s^{-1} \int e^{i s \varphi(\xi, \eta)} \mu_0^{0, 2, 1}(\xi, \eta) \nabla^2_{\xi} \widehat{f}(s, \xi - \eta) \widehat{f}(s, \eta) ~ d\eta \label{est_B++_close_term1-8} 
\end{align}
\end{subequations}
We need to control these terms with $\Lambda^{-\kappa}$ in front, so we will apply lemma \ref{lem_intfrac} each time~: denote by $q_0 = q_0(\kappa)$ such that $\frac{1}{q_0} = \frac{1}{2} + \frac{\kappa}{3}$. Denote also by $q_1 = q_1(\kappa)$ such that $\frac{1}{q_1} = \frac{1}{4} + \frac{\kappa}{3}$. 
\begin{align*}
\Vert \Lambda^{-\kappa} \eqref{est_B++_close_term1-3} \Vert_{H^2} &\lesssim s^{-1} \Vert T_{\mu_0}(\Lambda^2 e^{\pm i s \Lambda} x f, u) \Vert_{W^{2, q_0}} \lesssim s^{-1} \Vert x f \Vert_{H^3} \Vert u \Vert_{W^{2, 3/\kappa}} \lesssim s^{-2+\delta+2\kappa/3} \Vert u \Vert_X^2 \\
\Vert \Lambda^{-\kappa} \eqref{est_B++_close_term1-5} \Vert_{H^2} &\lesssim s^{-1} \Vert T_{\mu_0}(\Lambda^3 e^{\pm i s \Lambda} x f, \Lambda^{-1} u) \Vert_{W^{2, q_0}} \lesssim s^{-1} \Vert \Lambda^3 e^{\pm i s \Lambda} x f \Vert_{W^{2, q_1}} \Vert \Lambda^{-1} e^{\pm i s \Lambda} f \Vert_{W^{2, 4}} \\
&\lesssim s^{-2+2\kappa/3} \Vert \Lambda^3 x f \Vert_{W^{3, q_1'}} \Vert f \Vert_{W^{2, 4/3}} \lesssim s^{-2+2\kappa/3} \Vert \langle x \rangle \Lambda x f \Vert_{H^5} \Vert x f \Vert_{H^2} \lesssim s^{-2+\gamma+2\kappa/3} \Vert u \Vert_X^2 \\
\Vert \Lambda^{-\kappa} \eqref{est_B++_close_term1-6} \Vert_{H^2} &\lesssim s^{-1} \Vert T_{\mu_0}(\Lambda^3 e^{\pm i s \Lambda} x f, e^{\pm i s \Lambda} x f) \Vert_{W^{2, q_0}} \lesssim s^{-1} \Vert e^{\pm i s \Lambda} \Lambda x f \Vert_{W^{4, q_1}} \Vert e^{\pm i s \Lambda} x f \Vert_{W^{2, 4}} \\
&\lesssim s^{-2+2\kappa/3} \Vert \Lambda x f \Vert_{W^{5, q_1'}} \Vert \Lambda x f \Vert_{W^{2, 4/3}} \lesssim s^{-2+2\kappa/3} \Vert \langle x \rangle \Lambda x f \Vert_{H^5} \lesssim s^{-2+2\gamma+2\kappa/3} \Vert u \Vert_X^2 \\
\Vert \Lambda^{-\kappa} \eqref{est_B++_close_term1-7} \Vert_{H^2} &\lesssim s^{-1} \Vert T_{\mu_0}(\Lambda^3 e^{\pm i s \Lambda} x^2 f, u) \Vert_{W^{2, q_0}} \lesssim s^{-1} \Vert \Lambda |x|^2 f \Vert_{H^4} \Vert u \Vert_{W^{2, 3/\kappa}} \lesssim s^{-2+\gamma+\delta+2\kappa/3} \Vert u \Vert_X^2
\end{align*}
Then, \eqref{est_B++_close_term1-1} is simpler, \eqref{est_B++_close_term1-2} is similar to \eqref{est_B++_close_term1-3} by Hardy's inequality, \eqref{est_B++_close_term1-4} is similar to \eqref{est_B++_close_term1-3}, \eqref{est_B++_close_term1-8} is similar to \eqref{est_B++_close_term1-7}. 

\paragraph{Weight $|x|^2$} For the term with weight $|x|^2$ \eqref{est_B++_close_term2}, we apply another derivative $\nabla_{\xi}$ to the above terms from \eqref{est_B++_close_term1}. More precisely, we can apply it before the last integration by parts~: 
\begin{align*}
\eqref{est_B++_close_term2}
&= \nabla_{\xi} \int e^{i s \varphi(\xi, \eta)} \nabla_{\eta} \varphi(\xi, \eta) \mu_0^{0, 0, 2}(\xi, \eta) \widehat{f}(s, \xi - \eta) \widehat{f}(s, \eta) ~ d\eta \\
&+ \nabla_{\xi} \int e^{i s \varphi(\xi, \eta)} \nabla_{\eta} \varphi(\xi, \eta) \mu_0^{0, 0, 3}(\xi, \eta) \nabla_{\xi} \widehat{f}(s, \xi - \eta) \widehat{f}(s, \eta) ~ d\eta \\
&+ \nabla_{\xi} \int e^{i s \varphi(\xi, \eta)} \nabla_{\eta} \varphi(\xi, \eta) \mu_0^{0, 2, 1}(\xi, \eta) \nabla_{\xi} \widehat{f}(s, \xi - \eta) \widehat{f}(s, \eta) ~ d\eta
\end{align*}
Again, if the additional derivative hits the exponential, we only have a factor $s$ appearing and we can apply the exact same estimates to obtain a $s^{-1+2\gamma+\delta}$ decay. If the derivative hits the symbol and a $|\xi|^{-1}$ appears, we notice that $\nabla_{\eta} \varphi$ hasn't been differenciated (because its derivatives have no $|\xi|^{-1}$ singularity) so we only have to apply the same estimates as above in $L^{6/5}$ instead of $L^2$, which gives a decay of order $s^{-4/3}$. The remaining terms are~: 
\begin{subequations}
\begin{align}
&\int e^{i s \varphi(\xi, \eta)} \mu_0^{0, 0, 1}(\xi, \eta) \widehat{f}(s, \xi - \eta) \widehat{f}(s, \eta) ~ d\eta \label{est_B++_close_term2-1} \\
&+ \int e^{i s \varphi(\xi, \eta)} \mu_0^{0, 0, 2}(\xi, \eta) \nabla_{\xi} \widehat{f}(s, \xi - \eta) \widehat{f}(s, \eta) ~ d\eta \label{est_B++_close_term2-2} \\
&+ \int e^{i s \varphi(\xi, \eta)} \mu_0^{0, 0, 3}(\xi, \eta) \nabla^2_{\xi} \widehat{f}(s, \xi - \eta) \widehat{f}(s, \eta) ~ d\eta \label{est_B++_close_term2-3} \\
&+ \int e^{i s \varphi(\xi, \eta)} \mu_0^{0, 2, 0}(\xi, \eta) \nabla_{\xi} \widehat{f}(s, \xi - \eta) \widehat{f}(s, \eta) ~ d\eta \label{est_B++_close_term2-4} \\
&+ \int e^{i s \varphi(\xi, \eta)} \mu_0^{0, 2, 1}(\xi, \eta) \nabla^2_{\xi} \widehat{f}(s, \xi - \eta) \widehat{f}(s, \eta) ~ d\eta \label{est_B++_close_term2-5} 
\end{align}
\end{subequations}
Denote again $q_0$ such that $\frac{1}{q_0} = \frac{1}{2} + \frac{\kappa}{3}$. 
\begin{align*}
\Vert \Lambda^{-\kappa} \eqref{est_B++_close_term2-2} \Vert_{H^2} &\lesssim \Vert T_{\mu_0}(\Lambda^2 e^{\pm i s \Lambda} x f, u) \Vert_{W^{2, q_0}} \lesssim \Vert x f \Vert_{H^4} \Vert u \Vert_{W^{2, 3/\kappa}} \lesssim s^{-1+\delta+2\kappa/3} \Vert u \Vert_X^2 \\
\Vert \Lambda^{-\kappa} \eqref{est_B++_close_term2-3} \Vert_{H^2} &\lesssim \Vert T_{\mu_0}(\Lambda^3 e^{\pm i s \Lambda} x^2 f, u) \Vert_{W^{2, q_0}} \lesssim \Vert \Lambda x^2 f \Vert_{H^4} \Vert u \Vert_{W^{2, 3/\kappa}} \lesssim s^{-1+\delta+\gamma+2\kappa/3} \Vert u \Vert_X^2
\end{align*}
\eqref{est_B++_close_term2-1} is simpler, \eqref{est_B++_close_term2-4} is similar to \eqref{est_B++_close_term2-2} and \eqref{est_B++_close_term2-5} to \eqref{est_B++_close_term2-3}. 

Finally, let us treat the term without any weight \eqref{est_B++_close_term3}. We write it as~: 
\[ 
|\xi|^{-\kappa} \int e^{i s \varphi} \mu_0^{0, 0, 1} \nabla_{\eta} \varphi \widehat{f}(s, \xi - \eta) \widehat{f}(s, \eta) ~ d\eta \]
It is then a slight variation of lemma \ref{lem_decnonres} to prove a decay $s^{-2+}$~: we apply lemma \ref{lem_intfrac} and have then to control everything in $L^{q_0}$-norm instead of $L^2$, which leads to a small loss. 

This concludes the estimates on the $\pm \pm$ interactions. 

\subsection{\texorpdfstring{$\pm 0$}{+/- 0} and \texorpdfstring{$0 \pm$}{0 +/-} interactions}

In these cases, we separate three situations~: 
\begin{enumerate}
\item either the $\Lambda$ in the nonlinearity is on $u^0$, in which case it is easy to control even with a weight $x$~; 
\item or the $\Lambda$ in the nonlinearity is on $u$, in which case $u^0$ is subtler to control when weights $x$ appear. In this second case~: 
\begin{enumerate}
\item either $\varphi(\xi, \eta) = \pm (|\xi|_0 - |\xi - \eta|_0)$, which means we are in the $++0$ or $--0$ case, and we can use the identity of lemma \ref{identite_fond_0} to express $|\xi|_0 \nabla_{\xi} \varphi$ as a $\mu_0 |\eta|$ (and thus obtaining $\Lambda u^0$, less singular), by noticing that
\[ \frac{|\xi|_0 - |\xi - \eta|_0}{|\eta|_0} \]
is a symbol of order $0$. 
\item or $\varphi(\xi, \eta) = \pm (|\xi|_0 + |\xi - \eta|_0)$ and this allows to integrate by parts in time. One sensitive term to estimate will be $\partial_s u^0$, but for this one we may separate again $u^0$ into a quadratic expression. One last angular repartition is needed to control this case. 
\end{enumerate}
\end{enumerate}
In all this subsection, $\nabla_{\eta} \varphi$ is of the form $\pm g_0 \frac{\xi - \eta}{|\xi - \eta|_0}$ or $\pm g_0 \frac{\eta}{|\eta|_0}$ so we can always integrate by parts for free. 

\subsubsection{Case 1 (\texorpdfstring{$0 \pm$}{0 +/-})} 
In this case, there is one derivative on $u^0$ in the nonlinearity which removes its singularity. Let us apply the dispersive inequality followed by the inequality from lemma \ref{lem_L1L2} to obtain~: 
\[ \begin{aligned}
\left\Vert t e^{\pm i t \Lambda} \int_1^t e^{\mp i s \Lambda} T_{\mu_0}(\Lambda u^0, u) ~ ds \right\Vert_{\dot{B}^1_{\infty, 1}} 
&\lesssim \int_1^t \left( \Vert \Lambda^{-\kappa} |x| e^{\pm i s \Lambda} T_{\mu_0 |\xi|^2}(\Lambda u^0, u) \Vert_{H^2} + \Vert \Lambda^{-\kappa} T_{\mu_0}(\Lambda u^0, u) \Vert_{H^3} \right)^{1/2} \\
&~~~~~\left( \Vert \Lambda^{-\kappa} |x|^2 e^{\pm i s \Lambda} T_{\mu_0 |\xi|^2}(\Lambda u^0, u) \Vert_{H^2} + \Vert \Lambda^{-\kappa} T_{\mu_0}(\Lambda u^0, u) \Vert_{H^3} \right)^{1/2} ~ ds 
\end{aligned} \]
Denote 
\begin{subequations}
\begin{align}
&|x| e^{\pm i s \Lambda} T_{\mu_0 |\xi|^2}(\Lambda u^0, u) \label{est_B0+_term1} \\
&|x| e^{\pm i s \Lambda} T_{\mu_0 |\xi|^2}(\Lambda u^0, u) \label{est_B0+_term2} \\
&T_{\mu_0}(\Lambda u^0, u) \label{est_B0+_term3} 
\end{align}
\end{subequations}
to be estimated with a $\Lambda^{-\kappa}$ in front. 

For the term with weight $x$ \eqref{est_B0+_term1}, let us write~: 
\[ \begin{aligned}
\eqref{est_B0+_term1}
&= \int e^{i s \varphi(\xi, \eta)} \mu_0^{1, 0, 1}(\xi, \eta) \widehat{u^0}(s, \xi - \eta) \widehat{f}(s, \eta) ~ d\eta + \int e^{i s \varphi(\xi, \eta)} \mu_0^{2, 0, 0}(\xi, \eta) \widehat{u^0}(s, \xi - \eta) \widehat{f}(s, \eta) ~ d\eta \\
&+ \int e^{i s \varphi(\xi, \eta)} \mu_0^{2, 0, 1}(\xi, \eta) \nabla_{\xi} \widehat{u^0}(s, \xi - \eta) \widehat{f}(s, \eta) ~ d\eta + \int e^{i s \varphi(\xi, \eta)} s \mu_0^{2, 0, 1}(\xi, \eta) \widehat{u^0}(s, \xi - \eta) \widehat{f}(s, \eta) ~ d\eta
\end{aligned} \]
On the last term, we can apply an integration by part to get~: 
\begin{subequations}
\begin{align}
\eqref{est_B0+_term1}
&= \int e^{i s \varphi(\xi, \eta)} \mu_0^{1, 0, 1}(\xi, \eta) \widehat{u^0}(s, \xi - \eta) \widehat{f}(s, \eta) ~ d\eta \label{est_B0+_term1-1} \\
&+ \int e^{i s \varphi(\xi, \eta)} \mu_0^{2, 0, 0}(\xi, \eta) \widehat{u^0}(s, \xi - \eta) \widehat{f}(s, \eta) ~ d\eta \label{est_B0+_term1-2} \\
&+ \int e^{i s \varphi(\xi, \eta)} \mu_0^{2, 0, 1}(\xi, \eta) \nabla_{\xi} \widehat{u^0}(s, \xi - \eta) \widehat{f}(s, \eta) ~ d\eta \label{est_B0+_term1-3} \\
&+ \int e^{i s \varphi(\xi, \eta)} \mu_0^{2, -1, 1}(\xi, \eta) \widehat{u^0}(s, \xi - \eta) \widehat{f}(s, \eta) ~ d\eta \label{est_B0+_term1-4} \\
&+ \int e^{i s \varphi(\xi, \eta)} \mu_0^{2, 0, 1}(\xi, \eta) \widehat{u^0}(s, \xi - \eta) \nabla_{\eta} \widehat{f}(s, \eta) ~ d\eta \label{est_B0+_term1-5} 
\end{align}
\end{subequations}
\eqref{est_B0+_term1-1} and \eqref{est_B0+_term1-2} are easy to estimate. Then, for the following, we estimate by~: 
\[ \begin{aligned}
\Vert \Lambda^{-\kappa} \eqref{est_B0+_term1-3} \Vert_{H^2} = \Vert T_{\mu_0 |\xi|^{2 - \kappa}}(\Lambda x u^0, u) \Vert_{H^2} &\lesssim \Vert \Lambda x u^0 \Vert_{H^4} \Vert u \Vert_{W^{5, \infty-}} \lesssim s^{-2+\delta+a} \Vert u \Vert_X^2 \\
\Vert \Lambda^{-\kappa} \eqref{est_B0+_term1-5} \Vert_{H^2} = \Vert T_{\mu_0 |\xi|^{2 - \kappa}}(\Lambda u^0, e^{\pm i s \Lambda} x f) \Vert_{H^2} &\lesssim \Vert u^0 \Vert_{W^{6, \infty-}} \Vert xf \Vert_{H^4} \lesssim s^{-2+a+\delta} \Vert u \Vert_X^2
\end{aligned} \]
\eqref{est_B0+_term1-4} is similar to \eqref{est_B0+_term1-5} by Hardy's inequality. 

For the term with a weight $|x|^2$ \eqref{est_B0+_term2}, we differentiate the terms above. If the derivative hits the exponential, a factor $s$ appears and we can apply the same estimates as before to reach a decay $s^{-1}$. It then only remains~: 
\begin{align*}
&\int e^{i s \varphi(\xi, \eta)} \mu_0^{0, 0, 1}(\xi, \eta) \widehat{u^0}(s, \xi - \eta) \widehat{f}(s, \eta) ~ d\eta + \int e^{i s \varphi(\xi, \eta)} \mu_0^{1, 0, 0}(\xi, \eta) \widehat{u^0}(s, \xi - \eta) \widehat{f}(s, \eta) ~ d\eta \\
&+ \int e^{i s \varphi(\xi, \eta)} \mu_0^{1, 0, 1}(\xi, \eta) \nabla_{\xi} \widehat{u^0}(s, \xi - \eta) \widehat{f}(s, \eta) ~ d\eta + \int e^{i s \varphi(\xi, \eta)} \mu_0^{2, 0, -1}(\xi, \eta) \widehat{u^0}(s, \xi - \eta) \widehat{f}(s, \eta) ~ d\eta \\
&+ \int e^{i s \varphi(\xi, \eta)} \mu_0^{2, 0, 0}(\xi, \eta) \nabla_{\xi} \widehat{u^0}(s, \xi - \eta) \widehat{f}(s, \eta) ~ d\eta + \int e^{i s \varphi(\xi, \eta)} \mu_0^{2, 0, 1}(\xi, \eta) \nabla^2_{\xi} \widehat{u^0}(s, \xi - \eta) \widehat{f}(s, \eta) ~ d\eta \\
&+ \int e^{i s \varphi(\xi, \eta)} \mu_0^{1, -1, 1}(\xi, \eta) \widehat{u^0}(s, \xi - \eta) \widehat{f}(s, \eta) ~ d\eta + \int e^{i s \varphi(\xi, \eta)} \mu_0^{2, -1, 1}(\xi, \eta) \nabla_{\xi} \widehat{u^0}(s, \xi - \eta) \widehat{f}(s, \eta) ~ d\eta \\
&+ \int e^{i s \varphi(\xi, \eta)} \mu_0^{1, 0, 1}(\xi, \eta) \widehat{u^0}(s, \xi - \eta) \nabla_{\eta} \widehat{f}(s, \eta) ~ d\eta + \int e^{i s \varphi(\xi, \eta)} \mu_0^{2, 0, 0}(\xi, \eta) \widehat{u^0}(s, \xi - \eta) \nabla_{\eta} \widehat{f}(s, \eta) ~ d\eta \\
&+ \int e^{i s \varphi(\xi, \eta)} \mu_0^{2, 0, 1}(\xi, \eta) \nabla_{\xi} \widehat{u^0}(s, \xi - \eta) \nabla_{\eta} \widehat{f}(s, \eta) ~ d\eta
\end{align*}
Most of the terms are simple to estimate, except~: 
\[ \begin{aligned}
&\int e^{i s \varphi(\xi, \eta)} \mu_0^{2, 0, -1}(\xi, \eta) \widehat{u^0}(s, \xi - \eta) \widehat{f}(s, \eta) ~ d\eta + \int e^{i s \varphi(\xi, \eta)} \mu_0^{2, 0, 0}(\xi, \eta) \nabla_{\xi} \widehat{u^0}(s, \xi - \eta) \widehat{f}(s, \eta) ~ d\eta \\
&+ \int e^{i s \varphi(\xi, \eta)} \mu_0^{2, 0, 1}(\xi, \eta) \nabla^2_{\xi} \widehat{u^0}(s, \xi - \eta) \widehat{f}(s, \eta) ~ d\eta \\
&+ \int e^{i s \varphi(\xi, \eta)} \mu_0^{2, 0, 1}(\xi, \eta) \nabla_{\xi} \widehat{u^0}(s, \xi - \eta) \nabla_{\eta} \widehat{f}(s, \eta) ~ d\eta + \int e^{i s \varphi(\xi, \eta)} \mu_0^{2, -1, 1}(\xi, \eta) \nabla_{\xi} \widehat{u^0}(s, \xi - \eta) \widehat{f}(s, \eta) ~ d\eta
\end{aligned} \]
To treat the term $\nabla_{\xi}^2 \widehat{u^0}(s, \xi - \eta)$, we apply an integration by parts in frequency to recover terms we already have. The remaining terms are estimated by~: 
\[ \begin{aligned}
\Vert T_{\mu_0 |\xi|^{2 - \kappa}}(\Lambda^{-1} u^0, u) \Vert_{H^2} &\lesssim \Vert \Lambda^{-1} u^0 \Vert_{W^{4, 6}} \Vert u \Vert_{W^{4, 3}} \lesssim \Vert u^0 \Vert_{H^4} \Vert u \Vert_{W^{4, 3}} \lesssim s^{-4/3+a+\delta} \Vert u \Vert_X^2 \\
\Vert T_{\mu_0 |\xi|^{2 - \kappa}}(\Lambda x u^0, e^{\pm i s \Lambda} x f) \Vert_{H^2} &\lesssim \Vert \Lambda x u^0 \Vert_{W^{4, 4}} \Vert e^{\pm i s \Lambda} x f \Vert_{W^{4, 4}} \lesssim \Vert \Lambda x u^0 \Vert_{H^5} \Vert x f \Vert_{H^5} \lesssim s^{-1+a} \Vert u \Vert_X^2
\end{aligned} \]
The missing terms are similar to these ones. The term $\Lambda^{-\kappa} T_{\mu_0}(\Lambda u^0, u)$ to be estimated in $H^3$ is also simpler. 

\subsubsection{Case 2a (\texorpdfstring{$++0$}{++0} or \texorpdfstring{$--0$}{--0})} In this case, we apply lemma \ref{identite_fond_0} under the form~: 
\[ |\xi|_0 \nabla_{\xi} \varphi(\xi, \eta) = \mu_0^{0, 1, 0} \]
We can then apply as before the dispersive inequality. For the term with weight $x$, we have for instance~: 
\[ \begin{aligned}
&\nabla_{\xi} |\xi|^2 \int e^{i s \varphi(\xi, \eta)} \mu_0^{0, 0, 1}(\xi, \eta) \widehat{f}(s, \xi - \eta) \widehat{u^0}(s, \eta) ~ d\eta \\
&= \int e^{i s \varphi(\xi, \eta)} \mu_0^{1, 0, 1}(\xi, \eta) \widehat{f}(s, \xi - \eta) \widehat{u^0}(s, \eta) ~ d\eta + \int e^{i s \varphi(\xi, \eta)} \mu_0^{2, 0, 0}(\xi, \eta) \widehat{f}(s, \xi - \eta) \widehat{u^0}(s, \eta) ~ d\eta \\
&+ \int e^{i s \varphi(\xi, \eta)} \mu_0^{2, 0, 1}(\xi, \eta) \nabla_{\xi} \widehat{f}(s, \xi - \eta) \widehat{u^0}(s, \eta) ~ d\eta + \int e^{i s \varphi(\xi, \eta)} s \mu_0^{1, 1, 1}(\xi, \eta) \widehat{f}(s, \xi - \eta) \widehat{u^0}(s, \eta) ~ d\eta
\end{aligned} \]
We can apply an integration by parts on the last term, and estimate just as in the case $0 \pm$ just above. 

\subsubsection{Case 2b (\texorpdfstring{$+-0$}{+-0} or \texorpdfstring{$-+0$}{-+0})} This case is more complicated, and we start by sketching the proof to give the general idea. 

\begin{Intuit} The idea is to notice that, up to a sign, $\varphi(\xi, \eta) = |\xi|_0 + |\xi - \eta|_0$, so that we can always apply an integration by parts in time. If the derivative hits the $\pm$-type term, lemma \ref{lem_decnonres} offers the desired decay. However, if the time derivative hits $u^0$, we need to control $\partial_s u^0$, which does not have a better decay. Indeed, when developing $\partial_t u^0 = \partial_t \Lambda^{-1} \mathcal{N'}(\Lambda U, U)$, we get terms of the form $\Lambda^{-1} T_{\mu_0}(e^{\pm i s \Lambda} \Lambda \partial_s f, u)$ or $\Lambda^{-1} T_{\mu_0}(\Lambda u, e^{\pm i s \Lambda} \partial_s f)$, but also $\Lambda^{-1} T_{\mu_0 \varphi'}(\Lambda u, u)$ where $\varphi'$ is one of the wave interactions in $\mathcal{N}'$. 

Again, we will exhaust the cases, knowing that $\varphi'(\xi, \eta) = \pm |\xi - \eta|_0 \pm |\eta|_0$. If the signs here are opposite, then $\varphi'(\xi, \eta) = |\xi| \mu_0$, and the last term above from $\partial_s u^0$ is of the form $T_{\mu_0}(\Lambda u, u)$ and thus we removed the singularity. It can be controlled by lemma \ref{lem_decnonres} and this treats the cases $+-+-$, $+--+$, $-++-$, $-+-+$. 

If the signs are the same, that is the cases $+-++$, $+---$, $-+++$ or $-+--$, we know that
\[ \{ \nabla_{\rho} \varphi'(\eta, \rho) = 0 \} = \{ \rho = \lambda \eta, 0 \leq \lambda \leq 1 \} \]
and so we will apply an angular repartition given by 
\[ \chi(\rho, \eta) = \widetilde{\chi}\left( \frac{\rho}{|\rho|_0} \cdot g_0 \frac{\rho - \eta}{|\rho - \eta|_0} \right) \]
On the support of $\chi$, we are away from the space-resonant set, so we can apply integrations by parts in $\rho$ to gain more decay on $u^0$, without risking giving it a weight $x$ (which would have happened if we integrated by parts in $\eta$). On the support of $1 - \chi$, we have that $\rho \cdot g_0 (\rho - \eta) \leq \frac{1}{4} |\rho|_0 |\rho - \eta|_0$, so
\[ |\eta|_0^2 = |\rho|_0^2 + |\rho - \eta|_0^2 - 2 \rho \cdot g_0 (\rho - \eta) \geq |\rho|_0^2 + |\rho - \eta|_0^2 - \frac{1}{2} |\rho|_0 |\rho - \eta|_0 \geq \frac{3}{4} \max\left( |\rho|_0^2, |\rho - \eta|_0^2 \right) \]
so that $|\eta|^{-1}$ in front of the singularity of $u^0$ can be replaced by $|\rho - \eta|^{-1}$ and this absorbs the singularity. 

The $+- \pm 0$, $+- 0 \pm$, $+-00$ cases and their symmetries are essentially simpler. 
\end{Intuit}

Let us write the nonlinearity as~: 
\begin{equation} e^{\pm i t \Lambda} \int_1^t \int \int e^{i s \varphi(\xi, \eta)} e^{i s \varphi'(\eta, \rho)} \mu_0^{0, -1, 1}(\xi, \eta) \mu_0^{0, 0, 1}(\eta, \rho) \widehat{f}(s, \xi - \eta) \widehat{f}(s, \eta - \rho) \widehat{f}(s, \rho) ~ d\rho d\eta ds \label{est_B+-0_init} \end{equation}

\paragraph{Case $+-+-$ and analogous} In this case, $\varphi'(\eta, \rho) = \pm \left( |\eta - \rho|_0 - |\rho|_0 \right)$ while $\varphi(\xi, \eta) = \pm \left( |\xi|_0 + |\xi - \eta|_0 \right)$. Let us apply an integration by parts in time by artificially creating a $\varphi(\xi, \eta)$~: 
\begin{subequations}
\begin{align}
\eqref{est_B+-0_init}
&= e^{\pm i t \Lambda} \int_1^t \int \int e^{i s \varphi(\xi, \eta)} e^{i s \varphi'(\eta, \rho)} \frac{1}{\varphi(\xi, \eta)} \mu_0^{0, -1, 1}(\xi, \eta) \mu_0^{0, 0, 1}(\eta, \rho) \partial_s \widehat{f}(s, \xi - \eta) \widehat{f}(s, \eta - \rho) \widehat{f}(s, \rho) ~ d\rho d\eta ds \\
&+ e^{\pm i t \Lambda} \int_1^t \int \int e^{i s \varphi(\xi, \eta)} e^{i s \varphi'(\eta, \rho)} \frac{1}{\varphi(\xi, \eta)} \mu_0^{0, -1, 1}(\xi, \eta) \mu_0^{0, 0, 1}(\eta, \rho) \widehat{f}(s, \xi - \eta) \partial_s \widehat{f}(s, \eta - \rho) \widehat{f}(s, \rho) ~ d\rho d\eta ds \\
&+ e^{\pm i t \Lambda} \int_1^t \int \int e^{i s \varphi(\xi, \eta)} e^{i s \varphi'(\eta, \rho)} \frac{1}{\varphi(\xi, \eta)} \mu_0^{0, -1, 1}(\xi, \eta) \mu_0^{0, 0, 1}(\eta, \rho) \widehat{f}(s, \xi - \eta) \widehat{f}(s, \eta - \rho) \partial_s \widehat{f}(s, \rho) ~ d\rho d\eta ds \\
&+ e^{\pm i t \Lambda} \int_1^t \int \int e^{i s \varphi(\xi, \eta)} e^{i s \varphi'(\eta, \rho)} \frac{\varphi'(\eta, \rho)}{\varphi(\xi, \eta)} \mu_0^{0, 0, 1}(\xi, \eta) \mu_0^{-1, 0, 1}(\eta, \rho) \widehat{f}(s, \xi - \eta) \widehat{f}(s, \eta - \rho) \widehat{f}(s, \rho) ~ d\rho d\eta ds  \\
&+ e^{\pm i t \Lambda} \int \int e^{i t \varphi(\xi, \eta)} e^{i t \varphi'(\eta, \rho)} \frac{1}{\varphi(\xi, \eta)} \mu_0^{0, -1, 1}(\xi, \eta) \mu_0^{0, 0, 1}(\eta, \rho) \widehat{f}(t, \xi - \eta) \widehat{f}(t, \eta - \rho) \widehat{f}(t, \rho) ~ d\rho d\eta \label{est_B+-+-_term_final} \\
&+ e^{\pm i t \Lambda} \int \int e^{i \varphi(\xi, \eta)} e^{i \varphi'(\eta, \rho)} \frac{1}{\varphi} \mu_0^{0, -1, 1}(\xi, \eta) \mu_0^{0, 0, 1}(\eta, \rho) \widehat{f}(1, \xi - \eta) \widehat{f}(1, \eta - \rho) \widehat{f}(1, \rho) ~ d\rho d\eta \label{est_B+-+-_term_init}
\end{align}
\end{subequations}
Notice that $\frac{\mu_0^{0, -1, 1}(\xi, \eta)}{\varphi(\xi, \eta)} = \mu_0^{0, -1, 0}(\xi, \eta)$, and $\varphi'(\eta, \rho) \mu_0^{-1, 0, 1}(\eta, \rho) = \mu_0^{0, 0, 1}(\eta, \rho)$. The estimate on the initial time term \eqref{est_B+-+-_term_init} is a consequence of the hypothesis. For the final time term \eqref{est_B+-+-_term_final}, it can be written as~: 
\[ \Vert T_{\mu_0}(u, \Lambda^{-1} T_{\mu_0'}(\Lambda u, u)) \Vert_{\dot{B}^m_{\infty, 1}} \lesssim \Vert T_{\mu_0}(u, \Lambda^{-1} T_{\mu_0'}(\Lambda u, u)) \Vert_{W^{2, \infty-}} \lesssim \Vert u \Vert_{W^{2, \infty-}} \Vert \Lambda u \Vert_{W^{2, 3}} \Vert u \Vert_{W^{2, \infty-}} \lesssim t^{-7/3+3\delta} \Vert u \Vert_X^3 \]
For the other terms, we apply the dispersive inequality and lemmas \ref{lem_L1L2}, \ref{interp_Besov} to get back to $H^2$ norms with a weight $x$ or $x^2$, or a $H^3$ norm without any weight. Denote : 
\begin{subequations} \label{est_B+-+-_term1}
\begin{align}
&x e^{\mp is \Lambda} \Lambda^2 T_{\mu_0}(e^{\pm is \Lambda} \partial_s f, \Lambda^{-1} T_{\mu_0}(\Lambda u, u)) \label{est_B+-+-_term1-a} \\
&x^2 e^{\mp is \Lambda} \Lambda^2 T_{\mu_0}(e^{\pm is \Lambda} \partial_s f, \Lambda^{-1} T_{\mu_0}(\Lambda u, u)) \label{est_B+-+-_term1-b} \\
&\Lambda^2 T_{\mu_0}(e^{\pm is \Lambda} \partial_s f, \Lambda^{-1} T_{\mu_0}(\Lambda u, u)) \label{est_B+-+-_term1-c} \\
\nextParentEquation[est_B+-+-_term2]
&x e^{\mp is \Lambda} \Lambda^2 T_{\mu_0}(u, \Lambda^{-1} T_{\mu_0}(\Lambda e^{\pm is \Lambda} \partial_s f, u)) \label{est_B+-+-_term2-a} \\
&x^2 e^{\mp is \Lambda} \Lambda^2 T_{\mu_0}(u, \Lambda^{-1} T_{\mu_0}(\Lambda e^{\pm is \Lambda} \partial_s f, u)) \label{est_B+-+-_term2-b} \\
&\Lambda^2 T_{\mu_0}(u, \Lambda^{-1} T_{\mu_0}(\Lambda e^{\pm is \Lambda} \partial_s f, u)) \label{est_B+-+-_term2-c} \\
\nextParentEquation[est_B+-+-_term3]
&x e^{\mp is \Lambda} \Lambda^2 T_{\mu_0}(u, \Lambda^{-1} T_{\mu_0}(\Lambda u, e^{\pm is \Lambda} \partial_s f)) \label{est_B+-+-_term3-a} \\
&x^2 e^{\mp is \Lambda} \Lambda^2 T_{\mu_0}(u, \Lambda^{-1} T_{\mu_0}(\Lambda u, e^{\pm is \Lambda} \partial_s f)) \label{est_B+-+-_term3-b} \\
&\Lambda^2 T_{\mu_0}(u, \Lambda^{-1} T_{\mu_0}(\Lambda u, e^{\pm is \Lambda} \partial_s f)) \label{est_B+-+-_term3-c} \\
\nextParentEquation[est_B+-+-_term4] 
&x e^{\mp is \Lambda} \Lambda^2 T_{\mu_0}(u, T_{\mu_0}(\Lambda u, u)) \label{est_B+-+-_term4-a} \\
&x^2 e^{\mp is \Lambda} \Lambda^2 T_{\mu_0}(u, T_{\mu_0}(\Lambda u, u)) \label{est_B+-+-_term4-b} \\
&\Lambda^2 T_{\mu_0}(u, T_{\mu_0}(\Lambda u, u)) \label{est_B+-+-_term4-c}
\end{align}
\end{subequations}
All of these terms have to be estimated with a $\Lambda^{-\kappa}$ in front. 


For \eqref{est_B+-+-_term1-a}~: 
\begin{subequations}
\begin{align}
\eqref{est_B+-+-_term1-a}
&= \int \int e^{i s \varphi(\xi, \eta)} e^{i s \varphi'(\eta, \rho)} \mu_0^{1, -1, 0}(\xi, \eta) \mu_0^{0, 0, 1}(\eta, \rho) \partial_s \widehat{f}(s, \xi - \eta) \widehat{f}(s, \eta - \rho) \widehat{f}(s, \rho) ~ d\rho d\eta \label{est_B+-+-_term1-a-1} \\
&+ \int \int e^{i s \varphi(\xi, \eta)} e^{i s \varphi'(\eta, \rho)} \mu_0^{2, -1, -1}(\xi, \eta) \mu_0^{0, 0, 1}(\eta, \rho) \partial_s \widehat{f}(s, \xi - \eta) \widehat{f}(s, \eta - \rho) \widehat{f}(s, \rho) ~ d\rho d\eta \label{est_B+-+-_term1-a-2} \\
&+ \int \int e^{i s \varphi(\xi, \eta)} e^{i s \varphi'(\eta, \rho)} \mu_0^{2, -1, 0}(\xi, \eta) \mu_0^{0, 0, 1}(\eta, \rho) \partial_s \nabla_{\xi} \widehat{f}(s, \xi - \eta) \widehat{f}(s, \eta - \rho) \widehat{f}(s, \rho) ~ d\rho d\eta \label{est_B+-+-_term1-a-3} \\
&+ \int \int e^{i s \varphi(\xi, \eta)} e^{i s \varphi'(\eta, \rho)} s \mu_0^{2, -1, 0}(\xi, \eta) \mu_0^{0, 0, 1}(\eta, \rho) \partial_s \widehat{f}(s, \xi - \eta) \widehat{f}(s, \eta - \rho) \widehat{f}(s, \rho) ~ d\rho d\eta \label{est_B+-+-_term1-a-4} 
\end{align}
\end{subequations}
Denote $q_0$ such that $\frac{1}{q_0} = \frac{1}{2} + \frac{\kappa}{3}$ and let us estimate~: 
\begin{align*}
\Vert \Lambda^{-\kappa} \eqref{est_B+-+-_term1-a-2} \Vert_{H^2} &\lesssim \Vert T_{\mu_0 |\xi|^2}(\Lambda^{-1} e^{\pm i s \Lambda} \partial_s f, \Lambda^{-1} T_{\mu_0}(\Lambda u, u)) \Vert_{W^{2, q_0}} \lesssim \Vert \Lambda^{-1} \partial_s f \Vert_{H^4} \Vert \Lambda^{-1} T_{\mu_0}(\Lambda u, u) \Vert_{W^{4, 3/\kappa}} \\
&\lesssim \Vert x \partial_s f \Vert_{H^4} \Vert u \Vert_{W^{5, 3}} \Vert u \Vert_{W^{4, 3/\kappa}} \lesssim s^{-7/3+\tau+2\delta+2\kappa/3} \Vert u \Vert_X^4 \\
\Vert \Lambda^{-\kappa} \eqref{est_B+-+-_term1-a-3} \Vert_{H^2} &\lesssim \Vert T_{\mu_0 |\xi|^2}(e^{\pm i s \Lambda} x \partial_s f, \Lambda^{-1} T_{\mu_0}(\Lambda u, u)) \Vert_{W^{2, q_0}} \lesssim \Vert x \partial_s f \Vert_{H^4} \Vert u \Vert_{W^{5, 3}} \Vert u \Vert_{W^{4, 3/\kappa}} \\
&\lesssim s^{-7/3+\tau+2\delta+2\kappa/3} \Vert u \Vert_X^4 \\
\Vert \Lambda^{-\kappa} \eqref{est_B+-+-_term1-a-4} \Vert_{H^2} &\lesssim s \Vert T_{\mu_0 |\xi|^2}(e^{\pm i s \Lambda} \partial_s f, \Lambda^{-1} T_{\mu_0}(\Lambda u, u)) \Vert_{W^{2, q_0}} \lesssim s \Vert \partial_s f \Vert_{H^4} \Vert u \Vert_{W^{5, 3}} \Vert u \Vert_{W^{4, 3/\kappa}} \\
&\lesssim s^{-7/3+\tau+2\delta+2\kappa/3} \Vert u \Vert_X^4
\end{align*}
by applying lemma \ref{lem_decnonres}. The term \eqref{est_B+-+-_term1-a-1} is simpler. 

If we consider the weight $x^2$ \eqref{est_B+-+-_term1-b}, we keep the same decomposition as above and add a $\nabla_{\xi}$. If this derivative hits the exponential, we can apply the same estimates as for \eqref{est_B+-+-_term1-a} by only losing a factor $s$. If the derivative adds a factor $|\xi|^{-1}$, the estimate is simpler than above. It only remains~: 
\begin{subequations}
\begin{align}
&\int \int e^{i s \varphi(\xi, \eta)} e^{i s \varphi'(\eta, \rho)} \mu_0^{2, -1, -2}(\xi, \eta) \mu_0^{0, 0, 1}(\eta, \rho) \partial_s \widehat{f}(s, \xi - \eta) \widehat{f}(s, \eta - \rho) \widehat{f}(s, \rho) ~ d\rho d\eta \label{est_B+-+-_term1-b-1} \\
&+ \int \int e^{i s \varphi(\xi, \eta)} e^{i s \varphi'(\eta, \rho)} \mu_0^{2, -1, -1}(\xi, \eta) \mu_0^{0, 0, 1}(\eta, \rho) \partial_s \nabla_{\xi} \widehat{f}(s, \xi - \eta) \widehat{f}(s, \eta - \rho) \widehat{f}(s, \rho) ~ d\rho d\eta \label{est_B+-+-_term1-b-2} \\
&+ \int \int e^{i s \varphi(\xi, \eta)} e^{i s \varphi'(\eta, \rho)} \mu_0^{2, -1, 0}(\xi, \eta) \mu_0^{0, 0, 1}(\eta, \rho) \partial_s \nabla^2_{\xi} \widehat{f}(s, \xi - \eta) \widehat{f}(s, \eta - \rho) \widehat{f}(s, \rho) ~ d\rho d\eta \label{est_B+-+-_term1-b-ipp} 
\end{align}
On the last term, we apply an integration by parts in $\eta$~: 
\begin{align}
\eqref{est_B+-+-_term1-b-ipp}
&= \int \int e^{i s \varphi(\xi, \eta)} e^{i s \varphi'(\eta, \rho)} s \mu_0^{2, -1, 0}(\xi, \eta) \mu_0^{0, 0, 1}(\eta, \rho) \partial_s \nabla_{\xi} \widehat{f}(s, \xi - \eta) \widehat{f}(s, \eta - \rho) \widehat{f}(s, \rho) ~ d\rho d\eta \label{est_B+-+-_term1-b-simple1} \\
&+ \int \int e^{i s \varphi(\xi, \eta)} e^{i s \varphi'(\eta, \rho)} \mu_0^{2, -1, -1}(\xi, \eta) \mu_0^{0, 0, 1}(\eta, \rho) \partial_s \nabla_{\xi} \widehat{f}(s, \xi - \eta) \widehat{f}(s, \eta - \rho) \widehat{f}(s, \rho) ~ d\rho d\eta \label{est_B+-+-_term1-b-simple2} \\
&+ \int \int e^{i s \varphi(\xi, \eta)} e^{i s \varphi'(\eta, \rho)} \mu_0^{2, -2, 0}(\xi, \eta) \mu_0^{0, 0, 1}(\eta, \rho) \partial_s \nabla_{\xi} \widehat{f}(s, \xi - \eta) \widehat{f}(s, \eta - \rho) \widehat{f}(s, \rho) ~ d\rho d\eta \label{est_B+-+-_term1-b-3} \\
&+ \int \int e^{i s \varphi(\xi, \eta)} e^{i s \varphi'(\eta, \rho)} \mu_0^{2, -1, 0}(\xi, \eta) \mu_0(\eta, \rho) \partial_s \nabla_{\xi} \widehat{f}(s, \xi - \eta) \widehat{f}(s, \eta - \rho) \widehat{f}(s, \rho) ~ d\rho d\eta \label{est_B+-+-_term1-b-simple3} \\
&+ \int \int e^{i s \varphi(\xi, \eta)} e^{i s \varphi'(\eta, \rho)} \mu_0^{2, -1, 0}(\xi, \eta) \mu_0^{0, 0, 1}(\eta, \rho) \partial_s \nabla_{\xi} \widehat{f}(s, \xi - \eta) \nabla_{\eta} \widehat{f}(s, \eta - \rho) \widehat{f}(s, \rho) ~ d\rho d\eta \label{est_B+-+-_term1-b-4} 
\end{align}
\end{subequations}
\eqref{est_B+-+-_term1-b-simple1} and \eqref{est_B+-+-_term1-b-simple2} are analogous to already present terms~; \eqref{est_B+-+-_term1-b-simple3} is simpler. 
We then estimate~:
\begin{align*}
\Vert \Lambda^{-\kappa} \eqref{est_B+-+-_term1-b-1} \Vert_{H^2} &\lesssim \Vert T_{\mu_0 |\xi|^{2 - \kappa}}(\Lambda^{-2} e^{\pm i s \Lambda} \partial_s f, \Lambda^{-1} T_{\mu_0}(\Lambda u, u)) \Vert_{H^2} \lesssim \Vert \Lambda^{-2} e^{\pm i s \Lambda} \partial_s f \Vert_{W^{4, 6}} \Vert \Lambda^{-1} T_{\mu_0}(\Lambda u, u) \Vert_{H^4} \\
&\lesssim \Vert \Lambda^{-1} \partial_s f \Vert_{H^4} \Vert x T_{\mu_0}(\Lambda u, u) \Vert_{H^4} \lesssim \Vert x \partial_s f \Vert_{H^4} \Vert u \Vert_X^2 s^{-1+\tau} \lesssim s^{-2+2\tau} \Vert u \Vert_X^4 \\
\Vert \Lambda^{-\kappa} \eqref{est_B+-+-_term1-b-3} \Vert_{H^2} &\lesssim \Vert T_{\mu_0 |\xi|^{2 - \kappa}}(e^{\pm i s \Lambda} x \partial_s f, \Lambda^{-2} T_{\mu_0}(\Lambda u, u)) \Vert_{H^2} \lesssim \Vert x \partial_s f \Vert_{H^4} \Vert \Lambda^{-2} T_{\mu_0}(\Lambda u, u) \Vert_{W^{4, 6}} \\
&\lesssim \Vert x \partial_s f \Vert_{H^4} \Vert x T_{\mu_0}(\Lambda u, u) \Vert_{H^4} \lesssim s^{-2+2\tau} \Vert u \Vert_X^4 \\
\Vert \Lambda^{-\kappa} \eqref{est_B+-+-_term1-b-4} \Vert_{H^2} &\lesssim \Vert T_{\mu_0 |\xi|^{2 - \kappa}}(e^{\pm i s \Lambda} x \partial_s f, \Lambda^{-1} T_{\mu_0}(\Lambda e^{\pm i s \Lambda} x f, u)) \Vert_{H^2} \lesssim \Vert x \partial_s f \Vert_{H^4} \Vert \Vert x f \Vert_{H^5} \Vert u \Vert_{W^{5, \infty-}} \\
&\lesssim s^{-2+\tau+\delta} \Vert u \Vert_X^4
\end{align*}
because $T_{\mu_0}(\Lambda u, u)$ satisfies the non-resonance condition, so that we can apply lemma \ref{lem_decnonres}. \eqref{est_B+-+-_term1-b-2} can be treated as \eqref{est_B+-+-_term1-b-1}, controlling the first factor in $L^6$ and then applying a fractional integration inequality to win a derivative, and changing the $\Lambda^{-1}$ of the other factor into a $x$ by Hardy's inequality, in order to recover a term controlled by lemma \ref{lem_decnonres}. 

The terms \eqref{est_B+-+-_term2}, \eqref{est_B+-+-_term3} are very similar to \eqref{est_B+-+-_term1}. Finally, \eqref{est_B+-+-_term4} can be seen as a term similar to 
\[ \int_1^t \int e^{i s \varphi(\xi, \eta)} \mu_0^{0, 1, 1}(\xi, \eta) \widehat{f}(s, \xi - \eta) \widehat{u^0}(s, \eta) ~ d\eta \]
and thus can be treated as in the case $0 \pm$. 

The terms \eqref{est_B+-+-_term1-c}, \eqref{est_B+-+-_term2-c}, \eqref{est_B+-+-_term3-c}, \eqref{est_B+-+-_term4-c} without weight are simpler to estimate. 

\paragraph{Case $+-++$ and analogous} This time, we consider the interactions in $u^0$ where signs are identical. Therefore, we have $\varphi'(\eta, \rho) = \pm \left( |\eta - \rho|_0 + |\rho|_0 \right)$ so that
\[ \{ \nabla_{\rho} \varphi'(\eta, \rho) = 0 \} = \left\{ \frac{\eta - \rho}{|\eta - \rho|_0} = \frac{\rho}{|\rho|_0} \right\} = \{ \rho = \lambda (\eta - \rho), \lambda \geq 0 \} \]
Let us choose the following angular repartition~: 
\[ \chi(\eta, \rho) = \widetilde{\chi}\left( \frac{\rho}{|\rho|_0} \cdot g_0 \frac{\rho - \eta}{|\eta - \rho|_0} \right) \]
where $\widetilde{\chi}$ is the same smooth function already used to treat the $\pm \pm$ interactions earlier. 

On the support of $\chi$, $\frac{\rho}{|\rho|_0} \cdot g_0 \frac{\eta - \rho}{|\eta - \rho|_0} \leq \frac{1}{4}$, so in particular $\nabla_{\rho} \varphi'$ does not vanish~: 
\[ |\nabla_{\rho} \varphi'(\eta, \rho)|_{0'}^2 \chi(\eta, \rho) = 2 \chi(\eta, \rho) \left( 1 + \frac{\rho}{|\rho|_0} \cdot g_0 \frac{\rho - \eta}{|\rho - \eta|_0} \right) \geq \frac{3}{2} \chi(\eta, \rho) \]
On the support of $1 - \chi$, as we already saw, $|\eta|_0 \geq c |\rho - \eta|_0$, so that $\frac{|\rho - \eta|}{|\eta|} (1 - \chi(\eta, \rho))$ is a symbol of order $0$. 

Let us write $1 = (1 - \chi) + \chi$ and thus separating the interactions into two contributions.  

For the part containing $\chi$, we apply the dispersive inequality and we have to estimate 
\begin{equation} e^{\pm i s \Lambda} \Lambda^2 T_{\mu_0}(\Lambda u, \Lambda^{-1} T_{\mu_0' \chi}(\Lambda u, u)) \label{est_B+-++_chi_init} \end{equation}
with a weight $x$ or $x^2$, both in $H^2$, or without any weight in $H^3$. 

For the weight $x$, let us compute in Fourier~:  
\begin{align*}
\nabla_{\xi} \eqref{est_B+-++_chi_init} 
&= \int \int e^{i s \varphi(\xi, \eta)} e^{i s \varphi'(\eta, \rho)} \mu_0^{1, -1, 1}(\xi, \eta) \mu_0^{0, 0, 1}(\eta, \rho) \chi(\eta, \rho) \widehat{f}(s, \xi - \eta) \widehat{f}(s, \eta - \rho) \widehat{f}(s, \rho) ~ d\rho d\eta \\
&+ \int \int e^{i s \varphi(\xi, \eta)} e^{i s \varphi'(\eta, \rho)} \mu_0^{2, -1, 0}(\xi, \eta) \mu_0^{0, 0, 1}(\eta, \rho) \chi(\eta, \rho) \widehat{f}(s, \xi - \eta) \widehat{f}(s, \eta - \rho) \widehat{f}(s, \rho) ~ d\rho d\eta \\
&+ \int \int e^{i s \varphi(\xi, \eta)} e^{i s \varphi'(\eta, \rho)} \mu_0^{2, -1, 1}(\xi, \eta) \mu_0^{0, 0, 1}(\eta, \rho) \chi(\eta, \rho) \nabla_{\xi} \widehat{f}(s, \xi - \eta) \widehat{f}(s, \eta - \rho) \widehat{f}(s, \rho) ~ d\rho d\eta \\
&+ \int \int e^{i s \varphi(\xi, \eta)} e^{i s \varphi'(\eta, \rho)} s \mu_0^{2, -1, 1}(\xi, \eta) \mu_0^{0, 0, 1}(\eta, \rho) \chi(\eta, \rho) \widehat{f}(s, \xi - \eta) \widehat{f}(s, \eta - \rho) \widehat{f}(s, \rho) ~ d\rho d\eta
\end{align*}
We apply an integration by parts in $\rho$ on the last term, by using the fact that $\chi$ always authorize it. We get~: 
\begin{align*}
\nabla_{\xi} \eqref{est_B+-++_chi_init}
&= \int \int e^{i s \varphi(\xi, \eta)} e^{i s \varphi'(\eta, \rho)} \mu_0^{1, -1, 1}(\xi, \eta) \mu_0^{0, 0, 1}(\eta, \rho) \chi(\eta, \rho) \widehat{f}(s, \xi - \eta) \widehat{f}(s, \eta - \rho) \widehat{f}(s, \rho) ~ d\rho d\eta \\
&+ \int \int e^{i s \varphi(\xi, \eta)} e^{i s \varphi'(\eta, \rho)} \mu_0^{2, -1, 0}(\xi, \eta) \mu_0^{0, 0, 1}(\eta, \rho) \chi(\eta, \rho) \widehat{f}(s, \xi - \eta) \widehat{f}(s, \eta - \rho) \widehat{f}(s, \rho) ~ d\rho d\eta \\
&+ \int \int e^{i s \varphi(\xi, \eta)} e^{i s \varphi'(\eta, \rho)} \mu_0^{2, -1, 1}(\xi, \eta) \mu_0^{0, 0, 1}(\eta, \rho) \chi(\eta, \rho) \nabla_{\xi} \widehat{f}(s, \xi - \eta) \widehat{f}(s, \eta - \rho) \widehat{f}(s, \rho) ~ d\rho d\eta \\
&+ \int \int e^{i s \varphi(\xi, \eta)} e^{i s \varphi'(\eta, \rho)} \mu_0^{2, -1, 1}(\xi, \eta) \mu_0(\eta, \rho) \chi(\eta, \rho) \widehat{f}(s, \xi - \eta) \widehat{f}(s, \eta - \rho) \widehat{f}(s, \rho) ~ d\rho d\eta \\
&+ \int \int e^{i s \varphi(\xi, \eta)} e^{i s \varphi'(\eta, \rho)} \mu_0^{2, -1, 1}(\xi, \eta) \mu_0^{0, -1, 1}(\eta, \rho) \chi(\eta, \rho) \widehat{f}(s, \xi - \eta) \widehat{f}(s, \eta - \rho) \widehat{f}(s, \rho) ~ d\rho d\eta \\
&+ \int \int e^{i s \varphi(\xi, \eta)} e^{i s \varphi'(\eta, \rho)} \mu_0^{2, -1, 1}(\xi, \eta) \mu_0^{0, 0, 1}(\eta, \rho) \chi(\eta, \rho) \widehat{f}(s, \xi - \eta) \nabla_{\rho} \widehat{f}(s, \eta - \rho) \widehat{f}(s, \rho) ~ d\rho d\eta \\
&+ \int \int e^{i s \varphi(\xi, \eta)} e^{i s \varphi'(\eta, \rho)} \mu_0^{2, -1, 1}(\xi, \eta) \mu_0^{0, 0, 1}(\eta, \rho) \chi(\eta, \rho) \widehat{f}(s, \xi - \eta) \widehat{f}(s, \eta - \rho) \nabla_{\rho} \widehat{f}(s, \rho) ~ d\rho d\eta
\end{align*}
Now we apply to every term an integration by parts in $\rho$, allowed by the presence of $\chi$ (otherwise, the structure of the nonlinearity would only allow one). We get~:
\begin{subequations}
\begin{align}
\nabla_{\xi} \eqref{est_B+-++_chi_init}
&= \int \int e^{i s \varphi(\xi, \eta)} e^{i s \varphi'(\eta, \rho)} s^{-1} \mu_0^{1, -1, 1}(\xi, \eta) \widehat{f}(s, \xi - \eta) \nabla_{\rho} \left( \mu_0^{0, 0, 1}(\eta, \rho) \chi(\eta, \rho) \widehat{f}(s, \eta - \rho) \widehat{f}(s, \rho) \right) ~ d\rho d\eta \label{est_B+-++_chi_term1} \\
&+ \int \int e^{i s \varphi(\xi, \eta)} e^{i s \varphi'(\eta, \rho)} s^{-1} \mu_0^{2, -1, 0}(\xi, \eta) \widehat{f}(s, \xi - \eta) \nabla_{\rho} \left( \mu_0^{0, 0, 1}(\eta, \rho) \chi(\eta, \rho) \widehat{f}(s, \eta - \rho) \widehat{f}(s, \rho) \right) ~ d\rho d\eta \label{est_B+-++_chi_term2} \\
&+ \int \int e^{i s \varphi(\xi, \eta)} e^{i s \varphi'(\eta, \rho)} s^{-1} \mu_0^{2, -1, 1}(\xi, \eta) \nabla_{\xi} \widehat{f}(s, \xi - \eta) \nabla_{\rho} \left( \mu_0^{0, 0, 1}(\eta, \rho) \chi(\eta, \rho) \widehat{f}(s, \eta - \rho) \widehat{f}(s, \rho) \right) ~ d\rho d\eta \label{est_B+-++_chi_term3} \\
&+ \int \int e^{i s \varphi(\xi, \eta)} e^{i s \varphi'(\eta, \rho)} s^{-1} \mu_0^{2, -1, 1}(\xi, \eta) \mu_0^{0, 0, -1}(\eta, \rho) \chi(\eta, \rho) \widehat{f}(s, \xi - \eta) \widehat{f}(s, \eta - \rho) \widehat{f}(s, \rho) ~ d\rho d\eta \label{est_B+-++_chi_term4} \\
&+ \int \int e^{i s \varphi(\xi, \eta)} e^{i s \varphi'(\eta, \rho)} s^{-1} \mu_0^{2, -1, 1}(\xi, \eta) \mu_0(\eta, \rho) \chi(\eta, \rho) \widehat{f}(s, \xi - \eta) \nabla_{\rho} \widehat{f}(s, \eta - \rho) \widehat{f}(s, \rho) ~ d\rho d\eta \label{est_B+-++_chi_term5} \\
&+ \int \int e^{i s \varphi(\xi, \eta)} e^{i s \varphi'(\eta, \rho)} s^{-1} \mu_0^{2, -1, 1}(\xi, \eta) \mu_0^{0, -2, 1}(\eta, \rho) \chi(\eta, \rho) \widehat{f}(s, \xi - \eta) \widehat{f}(s, \eta - \rho) \widehat{f}(s, \rho) ~ d\rho d\eta \label{est_B+-++_chi_term6} \\
&+ \int \int e^{i s \varphi(\xi, \eta)} e^{i s \varphi'(\eta, \rho)} s^{-1} \mu_0^{2, -1, 1}(\xi, \eta) \mu_0^{0, -1, 1}(\eta, \rho) \chi(\eta, \rho) \widehat{f}(s, \xi - \eta) \nabla_{\rho} \widehat{f}(s, \eta - \rho) \widehat{f}(s, \rho) ~ d\rho d\eta \label{est_B+-++_chi_term7} \\
&+ \int \int e^{i s \varphi(\xi, \eta)} e^{i s \varphi'(\eta, \rho)} s^{-1} \mu_0^{2, -1, 1}(\xi, \eta) \mu_0^{0, -1, 1}(\eta, \rho) \chi(\eta, \rho) \widehat{f}(s, \xi - \eta) \widehat{f}(s, \eta - \rho) \nabla_{\rho} \widehat{f}(s, \rho) ~ d\rho d\eta \label{est_B+-++_chi_term8} \\
&+ \int \int e^{i s \varphi(\xi, \eta)} e^{i s \varphi'(\eta, \rho)} s^{-1} \mu_0^{2, -1, 1}(\xi, \eta) \mu_0^{0, 0, 1}(\eta, \rho) \chi(\eta, \rho) \widehat{f}(s, \xi - \eta) \nabla^2_{\rho} \widehat{f}(s, \eta - \rho) \widehat{f}(s, \rho) ~ d\rho d\eta \label{est_B+-++_chi_term9} \\
&+ \int \int e^{i s \varphi(\xi, \eta)} e^{i s \varphi'(\eta, \rho)} s^{-1} \mu_0^{2, -1, 1}(\xi, \eta) \mu_0^{0, 0, 1}(\eta, \rho) \chi(\eta, \rho) \widehat{f}(s, \xi - \eta) \nabla_{\rho} \widehat{f}(s, \eta - \rho) \nabla_{\rho} \widehat{f}(s, \rho) ~ d\rho d\eta \label{est_B+-++_chi_term10} \\
&+ \int \int e^{i s \varphi(\xi, \eta)} e^{i s \varphi'(\eta, \rho)} s^{-1} \mu_0^{2, -1, 1}(\xi, \eta) \mu_0^{0, 0, 1}(\eta, \rho) \chi(\eta, \rho) \widehat{f}(s, \xi - \eta) \widehat{f}(s, \eta - \rho) \nabla^2_{\rho} \widehat{f}(s, \rho) ~ d\rho d\eta \label{est_B+-++_chi_term11}
\end{align}
\end{subequations}
\eqref{est_B+-++_chi_term1}, \eqref{est_B+-++_chi_term2}, \eqref{est_B+-++_chi_term3} are simple to estimate. In \eqref{est_B+-++_chi_term6}, \eqref{est_B+-++_chi_term8} and \eqref{est_B+-++_chi_term11}, we distribute the derivative using \eqref{distrib_deriv_symbols} to avoid having too many singularities. When we have a $|\rho|$ appearing that way, we recover terms already present above (by symmetry). Therefore, we finally have objects of the form $\Lambda^2 T_{\mu_0}(\Lambda u, \Lambda^{-1} A)$ or $\Lambda^2 T_{\mu_0}(\Lambda u, B)$, which we can estimate by
\[ \begin{aligned}
&\Vert s^{-1} \Lambda^{-\kappa} T_{\mu_0 |\xi|^2}(\Lambda u, \Lambda^{-1} A) \Vert_{H^2} \lesssim s^{-1} \Vert u \Vert_{W^{5, 3}} \Vert \Lambda^{-1} A \Vert_{W^{4, 6}} \lesssim s^{-4/3+\delta} \Vert u \Vert_X \Vert A \Vert_{H^4} \\
&\Vert s^{-1} \Lambda^{-\kappa} T_{\mu_0 |\xi|^2}(\Lambda u, B) \Vert_{H^2} \lesssim s^{-1} \Vert u \Vert_{W^{5, \infty-}} \Vert B \Vert_{H^4} \lesssim s^{-2+\delta} \Vert u \Vert_X \Vert B \Vert_{H^4} 
\end{aligned} \]
Then, for $A$, we have the following possible expressions~: 
\[ \begin{aligned}
&\Vert T_{\mu_0}(e^{\pm i s \Lambda} x f, u) \Vert_{H^4} \lesssim \Vert xf  \Vert_{H^4} \Vert u \Vert_{W^{5, \infty-}} \lesssim s^{-1+\delta} \Vert u \Vert_X^2 \\
&\Vert T_{\mu_0}(e^{\pm i s \Lambda} \Lambda x f, e^{\pm i s \Lambda} x f) \Vert_{H^4} \lesssim \Vert e^{\pm i s \Lambda} x f \Vert_{W^{5, 4}}^2 \lesssim s^{-1} \Vert \Lambda x f \Vert_{W^{5, 4/3}} \lesssim \Vert \langle x \rangle \Lambda x f \Vert_{H^5} \lesssim s^{-1+2\gamma} \Vert u \Vert_X^2 \\
&\Vert T_{\mu_0}(\Lambda e^{\pm i s \Lambda} x^2 f, u) \Vert_{H^4} \lesssim \Vert \Lambda x^2 f \Vert_{H^4} \Vert u \Vert_{W^{5, \infty-}} \lesssim s^{-1+\delta+\gamma} \Vert u \Vert_X^2
\end{aligned} \]
or simpler or analogous forms~; and for $B$~: 
\[ \begin{aligned}
&\Vert T_{\mu_0}(\Lambda u, \Lambda^{-2} u) \Vert_{H^4} \lesssim \Vert u \Vert_{W^{5, 3}} \Vert \Lambda^{-2} u \Vert_{W^{4, 6}} \lesssim s^{-1/3+\delta} \Vert u \Vert_X \Vert x f \Vert_{H^4} \lesssim s^{-1/3+\delta} \Vert u \Vert_X^2 \\
&\Vert T_{\mu_0}(\Lambda u, e^{\pm i s \Lambda} x^2 f) \Vert_{H^4} \lesssim \Vert u \Vert_{W^{5, 3}} \Vert e^{\pm i s \Lambda} x^2 f \Vert_{W^{4, 6}} \lesssim \Vert u \Vert_{W^{5, 3}} \Vert \Lambda x^2 f \Vert_{H^4} \lesssim s^{-1/3+\delta+\gamma} \Vert u \Vert_X^2
\end{aligned} \]
and the others possibilities are simpler or analogous. Summing up, we obtain a decay of order $s^{-7/3}$. 

For the weight $x^2$, we simply add a derivative in $\xi$ and estimate as before. More precisely, if the derivative hits the exponential, we only lose a factor $s$ and can procede exactly as above~; if the derivative adds a $|\xi|^{-1}$ or a $|\xi - \eta|^{-1}$, we can also keep the same estimates~; if the derivative adds a $\nabla_{\xi}$ on $\widehat{f}(s, \xi - \eta)$, either we estimate as above, or we get a $|\xi - \eta| \nabla_{\xi}^2 \widehat{f}(s, \xi - \eta)$ that we also estimate in $L^3$, then going back to $L^2$ by Sobolev's inequality, which means losing only $s^{1/3+\gamma}$ with respect to the estimate with weight $x$. In either case, we have a decay of order $s^{-4/3}$, up to small parameters. 

The term without any weight is simpler. 

For the part with $(1 - \chi)$, we obtain a term of the form
\[ \int_1^t \int e^{i s \varphi(\xi, \eta)} \mu_0^{0, 0, 1}(\xi, \eta) \widehat{f}(s, \xi - \eta) \widehat{T_{\mu_0}(u, u)}(\eta) ~ d\eta ds \]
It can be estimated the same way as \eqref{est_B+-+-_term4} from case $+-+-$. 

\paragraph{Case $+-\pm 0$ and analogous} In this case, we have an object of the form~:
\begin{equation} \int_1^t \int e^{i s \varphi(\xi, \eta)} \mu_0^{0, -1, 1}(\xi, \eta) \widehat{f}(s, \xi - \eta) \widehat{A}(\eta) ~ d\eta ds \label{est_B+-+0_term} \end{equation}
where $A = T_{\mu_0}(\Lambda u, u^0)$ or $T_{\mu_0}(\Lambda u^0, u)$. We can apply again an integration by parts in $\rho$~:
\[ \begin{aligned}
s A &= T_{\mu_0}(u, u^0), ~~~ T_{\mu_0}(\Lambda e^{\pm i s \Lambda} x f, u^0), ~~~ T_{\mu_0}(\Lambda u, \Lambda^{-1} u^0), ~~~ T_{\mu_0}(\Lambda u, x u^0), \\
&~~~~~~~~ T_{\mu_0}(\Lambda u^0, \Lambda^{-1} u), ~~~ T_{\mu_0}(\Lambda u^0, e^{\pm i s \Lambda} x f), ~~~ T_{\mu_0}(\Lambda x u^0, u) 
\end{aligned} \]
Let us write $|\eta - \rho| = \mu_0 |\eta| + \mu_0' |\rho|$ (as in \eqref{distrib_deriv_symbols_back}) whenever we have $\Lambda^{-1} u^0$ or $x u^0$ and bring us back to~:
\begin{equation} \begin{aligned}
s A &= T_{\mu_0}(u, u^0), ~~~ T_{\mu_0}(\Lambda e^{\pm i s \Lambda} x f, u^0), ~~~ \Lambda T_{\mu_0}(u, \Lambda^{-1} u^0), ~~~ \Lambda T_{\mu_0}(u, x u^0), \\
&~~~~~~~~ T_{\mu_0}(\Lambda u^0, \Lambda^{-1} u), ~~~ T_{\mu_0}(\Lambda u^0, e^{\pm i s \Lambda} x f), ~~~ T_{\mu_0}(\Lambda x u^0, u) 
\end{aligned} \label{est_B+-+0_decomposition} \end{equation}
We now apply the dispersive inequality of lemma \ref{lem_disp} together with lemmas \ref{lem_L1L2}, \ref{interp_Besov} to get~: 
\[ \begin{aligned}
\Vert \eqref{est_B+-+0_term} \Vert_{\dot{B}^1_{1, \infty}} \lesssim t^{-1} \int_1^t &\left( \Vert \Lambda^{-\kappa} |x| \Lambda^2 T_{\mu_0}(\Lambda u, \Lambda^{-1} A) \Vert_{H^1} + \Vert T_{\mu_0}(\Lambda u, \Lambda^{-1} A) \Vert_{H^4} \right)^{1/2} \\
&\quad \left( \Vert \Lambda^{-\kappa} |x|^2 \Lambda^2 T_{\mu_0}(\Lambda u, \Lambda^{-1} A) \Vert_{H^1} + \Vert T_{\mu_0}(\Lambda u, \Lambda^{-1} A) \Vert_{H^4} \right)^{1/2} ds 
\end{aligned} \]
Denote 
\begin{subequations}
\begin{align}
&|x| \Lambda^2 T_{\mu_0}(\Lambda u, \Lambda^{-1} A) \label{est_B+-+0_terma} \\
&|x|^2 \Lambda^2 T_{\mu_0}(\Lambda u, \Lambda^{-1} A) \label{est_B+-+0_termb} \\
&T_{\mu_0}(\Lambda u, \Lambda^{-1} A) \label{est_B+-+0_termc} 
\end{align}
\end{subequations}
We expand : 
\[ \eqref{est_B+-+0_terma} = \Lambda T_{\mu_0}(\Lambda u, \Lambda^{-1} A) + \Lambda^2 T_{\mu_0}(u, \Lambda^{-1} A) + \Lambda^2 T_{\mu_0}(e^{\pm i s \Lambda} \Lambda x f, \Lambda^{-1} A) + \Lambda^2 T_{\mu_0}(\Lambda u, \Lambda^{-1} s A) \]
For the first terms, we have the direct estimate
\[ \Vert A \Vert_{H^3} \lesssim \Vert T_{\mu_0}(\Lambda u, u^0) \Vert_{H^3} + \Vert T_{\mu_0}(\Lambda u^0, u) \Vert_{H^3} \lesssim \Vert u^0 \Vert_{H^4} \Vert u \Vert_{W^{5, \infty-}} \lesssim s^{-2+\delta} \Vert u \Vert_X^2 \]
so that we can easily bound for instance
\[ \Vert \Lambda^{2-kappa} T_{\mu_0}(e^{\pm i s \Lambda} \Lambda x f, \Lambda^{-1} A) \Vert_{H^1} \lesssim \Vert e^{\pm i s \Lambda} x f \Vert_{W^{4, 3}} \Vert \Lambda^{-1} A \Vert_{W^{3, 6}} \lesssim s^{-1/3} \Vert x f \Vert_{W^{5, 3/2}} \Vert A \Vert_{H^3} \lesssim s^{-2-1/3+\gamma+\delta} \Vert u \Vert_X^3 \]
and similarly for the other ones. Finally, for the last term containing $sA$, we apply the decomposition \eqref{est_B+-+0_decomposition} to get 
\[ \Vert \Lambda^{2-\kappa} T_{\mu_0}(\Lambda u, \Lambda^{-1} s A) \Vert_{H^1} \lesssim s^{-7/3+\delta+a+\gamma} \Vert u \Vert_X \]
by similar estimates. In all cases, up to choosing the parameters of the $X$-norm small enough, we get for \eqref{est_B+-+0_terma} a better decay than $s^{-2}$. 

On the other hand, if we expand 
\[ \eqref{est_B+-+0_termb} = x \eqref{est_B+-+0_terma} \sim \Lambda^{-1} \eqref{est_B+-+0_terma} + s \eqref{est_B+-+0_terma} + \eqref{est_B+-+0_terma}' + \eqref{est_B+-+0_terma}'' \]
where $\eqref{est_B+-+0_terma}' := |x| \Lambda^2 T_{\mu_0}(u, \Lambda^{-1} A)$ and $\eqref{est_B+-+0_terma}'' := |x| \Lambda^2 T_{\mu_0}(e^{\pm i s \Lambda} \Lambda xf, \Lambda^{-1} A)$. $s \eqref{est_B+-+0_terma}$ has already been estimated ; $\Lambda^{-1} \eqref{est_B+-+0_terma}$ is simpler, as well as $\eqref{est_B+-+0_terma}'$, while 
\[ \eqref{est_B+-+0_terma}'' = \Lambda T_{\mu_0}(e^{\pm i s \Lambda} \Lambda x f, \Lambda^{-1} A) + \Lambda^2 T_{\mu_0}(e^{\pm i s \Lambda} x f, \Lambda^{-1} A) + \Lambda^2 T_{\mu_0}(e^{\pm i s \Lambda} \Lambda xf, \Lambda^{-1} s A) + \Lambda^2 T_{\mu_0}(e^{\pm i s \Lambda} \Lambda x^2 f, \Lambda^{-1} A) \]
Here above, everything except the last term is already present in one of the previous terms. Finally, 
\[ \begin{aligned}
\Vert \Lambda^{2-\kappa} T_{\mu_0}(e^{\pm is \Lambda} \Lambda x^2 f, \Lambda^{-1} A) \Vert_{H^1} &\lesssim \Vert \Lambda x^2 f \Vert_{H^3} \Vert \Lambda^{-1} e^{\pm i s \Lambda} A \Vert_{W^{3, \infty-}} \\
&\lesssim s^{\gamma} \Vert u \Vert_X \left( \Vert T_{\mu_0}(\Lambda u, u^0) \Vert_{W^{3, 3-}} + \Vert T_{\mu_0}(\Lambda u^0, u) \Vert_{W^{3, 3-}} \right) \\
&\lesssim s^{\gamma} \Vert u \Vert_X \Vert u \Vert_{W^{4, 6-}} \Vert u^0 \Vert_{W^{4, 6-}} \lesssim s^{-2+a+\gamma+} \Vert u \Vert_X^3
\end{aligned} \]

Finally, we skip the estimate of \eqref{est_B+-+0_termc} which decays faster than the precedent terms. Putting everything together, we get a strong enough decay to conclude. 

%

\paragraph{Case $+-00$ and analogous} This time, we have~: 
\[ \int_1^t \int e^{i s \varphi(\xi, \eta)} \mu_0^{0, -1, 1}(\xi, \eta) \widehat{f}(s, \xi - \eta) \widehat{T_{\mu_0}(\Lambda u^0, u^0)}(\eta) ~ d\eta \]
We apply again the dispersive inequality and obtain terms with weight $x$ or $x^2$, on which we lose at most $s$ or $s^2$ respectively. We estimate $\widehat{f}(s, \xi - \eta)$ in $L^3$ to win $s^{-1/3}$, which leaves $T_{\mu_0}(\Lambda u^0, u^0)$ in $L^2$, in which we win $s^{-3+2a}$, and this is largely enough. 

\subsection{\texorpdfstring{$0 0$}{0 0} interactions}

Let us apply again the dispersive inequality. We need to estimate~: 
\[ \Vert \Lambda^{-\kappa} |x| T_{\mu_0 |\xi|^2}(\Lambda u^0, u^0) \Vert_{H^2}, ~~~ \Vert \Lambda^{-\kappa} |x|^2 T_{\mu_0 |\xi|^2}(\Lambda u^0, u^0) \Vert_{H^2}, ~~~ \Vert \Lambda^{-\kappa} T_{\mu_0}(\Lambda u^0, u^0) \Vert_{H^3} \]
so that their mean decay be integrable in time. Denote
\begin{subequations}
\begin{align}
&|x| T_{\mu_0 |\xi|^2}(\Lambda u^0, u^0) \label{est_B00_term1} \\
&|x|^2 T_{\mu_0 |\xi|^2}(\Lambda u^0, u^0) \label{est_B00_term2} \\
&T_{\mu_0}(\Lambda u^0, u^0) \label{est_B00_term3}
\end{align}
\end{subequations}

For the term with weight $x$ \eqref{est_B00_term1}~: 
\begin{subequations}
\begin{align}
\eqref{est_B00_term1}
&= \int e^{i s \varphi(\xi, \eta)} \mu_0^{1, 0, 1}(\xi, \eta) \widehat{u^0}(s, \xi - \eta) \widehat{u^0}(s, \eta) ~ d\eta \label{est_B00_term1-1} \\
&+ \int e^{i s \varphi(\xi, \eta)} \mu_0^{2, 0, 0}(\xi, \eta) \widehat{u^0}(s, \xi - \eta) \widehat{u^0}(s, \eta) ~ d\eta \label{est_B00_term1-2} \\
&+ \int e^{i s \varphi(\xi, \eta)} \mu_0^{2, 0, 1}(\xi, \eta) \nabla_{\xi} \widehat{u^0}(s, \xi - \eta) \widehat{u^0}(s, \eta) ~ d\eta \label{est_B00_term1-3} \\
&+ \int e^{i s \varphi(\xi, \eta)} s \mu_0^{2, 0, 1}(\xi, \eta) \widehat{u^0}(s, \xi - \eta) \widehat{u^0}(s, \eta) ~ d\eta \label{est_B00_term1-4} 
\end{align}
\end{subequations}
We can estimate everything by making use of the strong decay of $u^0$~: 
\[ \begin{aligned}
\Vert \Lambda^{-\kappa} \eqref{est_B00_term1-3} \Vert_{H^2} \lesssim \Vert \Lambda x u^0 \Vert_{H^4} \Vert u^0 \Vert_{W^{4, 3/\kappa}} \lesssim s^{3+2a+\delta+2\kappa/3} \Vert u \Vert_X^2 \\
\Vert \Lambda^{-\kappa} \eqref{est_B00_term1-4} \Vert_{H^2} \lesssim s \Vert u^0 \Vert_{H^5} \Vert u^0 \Vert_{W^{4, 3/\kappa}} \lesssim s^{-2+2a+\delta+2\kappa/3} \Vert u \Vert_X^2
\end{aligned} \]
The terms \eqref{est_B00_term1-1} and \eqref{est_B00_term1-2} are simpler. 

For the term with weight $x^2$ \eqref{est_B00_term2}~: 
\[ \begin{aligned}
\eqref{est_B00_term2}
&= \nabla_{\xi} \int e^{i s \varphi(\xi, \eta)} \mu_0^{1, 0, 1}(\xi, \eta) \widehat{u^0}(s, \xi - \eta) \widehat{u^0}(s, \eta) ~ d\eta + \nabla_{\xi} \int e^{i s \varphi(\xi, \eta)} \mu_0^{2, 0, 0}(\xi, \eta) \widehat{u^0}(s, \xi - \eta) \widehat{u^0}(s, \eta) ~ d\eta \\
&+ \nabla_{\xi} \int e^{i s \varphi(\xi, \eta)} \mu_0^{2, 0, 1}(\xi, \eta) \nabla_{\xi} \widehat{u^0}(s, \xi - \eta) \widehat{u^0}(s, \eta) ~ d\eta + \nabla_{\xi} \int e^{i s \varphi(\xi, \eta)} s \mu_0^{2, 0, 1}(\xi, \eta) \widehat{u^0}(s, \xi - \eta) \widehat{u^0}(s, \eta) ~ d\eta
\end{aligned} \]
If the second derivative hits the exponential, we estimate just like above with weight $x$, only losing a factor $s$. The remaining terms are~: 
\begin{subequations}
\begin{align}
&\int e^{i s \varphi(\xi, \eta)} \mu_0^{0, 0, 1}(\xi, \eta) \widehat{u^0}(s, \xi - \eta) \widehat{u^0}(s, \eta) ~ d\eta \\
&+ \int e^{i s \varphi(\xi, \eta)} \mu_0^{1, 0, 0}(\xi, \eta) \widehat{u^0}(s, \xi - \eta) \widehat{u^0}(s, \eta) ~ d\eta \\
&+ \int e^{i s \varphi(\xi, \eta)} \mu_0^{1, 0, 1}(\xi, \eta) \nabla_{\xi} \widehat{u^0}(s, \xi - \eta) \widehat{u^0}(s, \eta) ~ d\eta \\
&+ \int e^{i s \varphi(\xi, \eta)} \mu_0^{2, 0, -1}(\xi, \eta) \widehat{u^0}(s, \xi - \eta) \widehat{u^0}(s, \eta) ~ d\eta \\
&+ \int e^{i s \varphi(\xi, \eta)} \mu_0^{2, 0, 0}(\xi, \eta) \nabla_{\xi} \widehat{u^0}(s, \xi - \eta) \widehat{u^0}(s, \eta) ~ d\eta \label{est_B00_term2-5} \\
&+ \int e^{i s \varphi(\xi, \eta)} \mu_0^{2, 0, 1}(\xi, \eta) \nabla^2_{\xi} \widehat{u^0}(s, \xi - \eta) \widehat{u^0}(s, \eta) ~ d\eta \label{est_B00_term2-6}
\end{align}
\end{subequations}
For \eqref{est_B00_term2-6}, we apply an integration by parts to recover only already present terms. For \eqref{est_B00_term2-5}~: 
\[ \Vert \Lambda^{2-\kappa} T_{\mu_0}(x u^0, u^0) \Vert{H^2} \lesssim \Vert x u^0 \Vert_{W^{4, 6}} \Vert u^0 \Vert_{W^{4, 3}} \lesssim \Vert \Lambda x u^0 \Vert_{H^4} \Vert u^0 \Vert_{W^{4, 3}} \lesssim s^{-4/3+2a+\delta} \Vert u \Vert_X^2 \]
All the other terms are simpler or analogous. 

\section{Estimate of the \texorpdfstring{$L^2$}{L2} norm with weight \texorpdfstring{$x^2$}{x2}} \label{section_L2x2}

Here, we may see appear terms with a weight $x^2$ and one derivative too much with respect to what we control with the norm $\Vert u \Vert_X$ (because the nonlinearity is of order 1). Therefore, we will use an approach similar to the one already used in the energy estimate $H^N$ and show that~: 
\[ \int_1^t \partial_s \Vert \partial^{\alpha} |x|^2 f \Vert_{L^2}^2 ~ ds \lesssim t^{2 \gamma} \Vert u \Vert_X^3 \]
which will prove~: 
\[ \Vert \partial^{\alpha} |x|^2 f \Vert_{L^2} \leq \Vert \partial^{\alpha} |x|^2 f_0 \Vert_{L^2} + C t^{\gamma} \Vert u \Vert_X^{3/2} \]
By Parseval's identity, this is the same as estimating~: 
\[ \int_1^t \int \xi^{\alpha} \nabla_{\xi}^2 \widehat{f}(s, \xi) \xi^{\alpha} \partial_t \nabla_{\xi}^2 \widehat{f}(s, \xi) ~ d\xi ds \]

Since we chose $\gamma_k = \frac{k \gamma_5 + (5-k) \gamma_0}{5}$, the estimates for $k = 1, 2, 3, 4$ are a consequence for those for $k = 0$, $k = 5$ by interpolation. 

\subsection{Quasi-linear structure}

Let us fix $\alpha \in \mathbb{N}^3$, with $|\alpha| = 6$ or $|\alpha| = 1$. We write~: 
\begin{subequations} \label{est_L2x2_init}
\begin{align}
&\mathcal{F} D^{\alpha} |x|^2 \partial_t f \tag{\ref{est_L2x2_init}} \\
&= - \xi^{\alpha} \Delta_{\xi} \int e^{i s \varphi} \mu_0^{0, 0, 1}(\xi, \eta) \widehat{f}(s, \xi - \eta) \widehat{f}(s, \eta) ~ d\eta \nonumber \\
&= - \xi^{\alpha} \int e^{i s \varphi} s^2 \nabla_{\xi} \varphi \cdot \nabla_{\xi} \varphi \mu_0^{0, 0, 1}(\xi, \eta) \widehat{f}(s, \xi - \eta) \widehat{f}(s, \eta) ~ d\eta \\
&- \xi^{\alpha} \int e^{i s \varphi} is \Delta_{\xi} \varphi \mu_0^{0, 0, 1}(\xi, \eta) \widehat{f}(s, \xi - \eta) \widehat{f}(s, \eta) ~ d\eta \\
&- \xi^{\alpha} \int e^{i s \varphi} \mu_0^{-2, 0, 1}(\xi, \eta) \widehat{f}(s, \xi - \eta) \widehat{f}(s, \eta) ~ d\eta \\
&- \xi^{\alpha} \int e^{i s \varphi} \mu_0^{0, 0, -1}(\xi, \eta) \widehat{f}(s, \xi - \eta) \widehat{f}(s, \eta) ~ d\eta \\ 
&- \xi^{\alpha} \int e^{i s \varphi} \mu_0^{0, 0, 1}(\xi, \eta) \Delta_{\xi} \widehat{f}(s, \xi - \eta) \widehat{f}(s, \eta) ~ d\eta \label{est_L2x2_term4} \\
&- \xi^{\alpha} \int e^{i s \varphi} i s \nabla_{\xi} \varphi \mu_0(\xi, \eta) \widehat{f}(s, \xi - \eta) \widehat{f}(s, \eta) ~ d\eta \\
&- \xi^{\alpha} \int e^{i s \varphi} i s \nabla_{\xi} \varphi \mu_0^{0, 0, 1}(\xi, \eta) \nabla_{\xi} \widehat{f}(s, \xi - \eta) \widehat{f}(s, \eta) ~ d\eta \\
&- \xi^{\alpha} \int e^{i s \varphi} i s \nabla_{\xi} \varphi \mu_0^{-1, 0, 1}(\xi, \eta) \widehat{f}(s, \xi - \eta) \widehat{f}(s, \eta) ~ d\eta \\
&- \xi^{\alpha} \int e^{i s \varphi} \mu_0^{-1, 0, 1}(\xi, \eta) \nabla_{\xi} \widehat{f}(s, \xi - \eta) \widehat{f}(s, \eta) ~ d\eta \\
&- \xi^{\alpha} \int e^{i s \varphi} \mu_0(\xi, \eta) \nabla_{\xi} \widehat{f}(s, \xi - \eta) \widehat{f}(s, \eta) ~ d\eta
\end{align}
\end{subequations}
On the other hand, 
\[ \mathcal{F} D^{\alpha} |x|^2 f = - \xi^{\alpha} \Delta_{\xi} \widehat{f}(s, \xi) \]

Let us only treat \eqref{est_L2x2_term4} for now, by writing that $\xi^{\alpha} = ((\xi - \eta) + \eta)^{\alpha}$, developing and isolating the term $(\xi - \eta)^{\alpha}$. This is~: 
\[ \int_1^t \int e^{i s A(D)} D^{\alpha} |x|^2 f^{\epsilon_1} \cdot P^{\epsilon_1} \mathcal{N}(\Lambda P^{\epsilon_2} e^{i s A(D)} D^{\alpha} |x|^2 f^{\epsilon_2}, u^{\epsilon_3}) ~ dx ds \]
The projection operators are orthogonal and symmetric, so we can transfer $P^{\epsilon_1}$ to the other term outside $B$, sum on $\epsilon_1 = +, -$ and $\epsilon_2 = +, -$, use the symmetric structure given in proposition \ref{prop_quasi_lin} and recover a product of the form~: 
\[ \int_1^t \int (e^{i s A(D)} D^{\alpha} g) (e^{i s A(D)} D^{\alpha} g) \Lambda u^{\epsilon_3} ~ dx ds, \quad \mbox{where } g = P^{+} |x|^2 f^{+} + P^{-} |x|^2 f^{-} \]
Then, we conclude as for the $H^N$ estimate, obtaining~: 
\[ \lesssim \int_1^t \Vert \Lambda |x|^2 f \Vert_{H^k}^2 \Vert u \Vert_{W^{1, \infty}} ~ ds \lesssim \Vert u \Vert_X^3 t^{2\gamma_k} \]
where $k = 0$ or $k = 5$. Then, in \eqref{est_L2x2_term4} we restricted our attention to, the only sub-term to estimate has the form~: 
\[ \int e^{i s \varphi} \mu_0^{|\alpha|-1, 1, 1}(\xi, \eta) \Delta_{\xi} \widehat{f}(s, \xi - \eta) \widehat{f}(s, \eta) ~ d\eta \]
Note that the whole 4th term remains in case of a $0 \pm$ or $0 0$ interaction. In the following, we fix $k = 0$ or $k = 5$. 

\subsection{\texorpdfstring{$\pm \pm$}{+/- +/-} interactions} \label{section-L2x2++}

In this case, the structure of the nonlinearity allows to factor the symbol by $\nabla_{\eta} \varphi$ once, and to use lemma \ref{identite_fond} to replace $\nabla_{\xi} \varphi$ by $\nabla_{\eta} \varphi$ and $\varphi$. Starting from the decomposition above, we get~: 
\begin{align*}
&\int e^{i s \varphi} s^2 \nabla_{\xi} \varphi \nabla_{\xi} \varphi \nabla_{\eta} \varphi \mu_0^{k+1, 0, 1}(\xi, \eta) \widehat{f}(s, \xi - \eta) \widehat{f}(s, \eta) ~ d\eta + \int e^{i s \varphi} s \nabla_{\eta} \varphi \mu_0^{k, 0, 1}(\xi, \eta) \widehat{f}(s, \xi - \eta) \widehat{f}(s, \eta) ~ d\eta \\
&+ \int e^{i s \varphi} s \nabla_{\eta} \varphi \mu_0^{k+1, 0, 0}(\xi, \eta) \widehat{f}(s, \xi - \eta) \widehat{f}(s, \eta) ~ d\eta + \int e^{i s \varphi} \nabla_{\eta} \varphi \mu_0^{k-1, 0, 1}(\xi, \eta) \widehat{f}(s, \xi - \eta) \widehat{f}(s, \eta) ~ d\eta \\
&+ \int e^{i s \varphi} \mu_0^{k+1, 0, -1}(\xi, \eta) \widehat{f}(s, \xi - \eta) \widehat{f}(s, \eta) ~ d\eta + \int e^{i s \varphi} \mu_0^{k, 1, 1}(\xi, \eta) \Delta_{\xi} \widehat{f}(s, \xi - \eta) \widehat{f}(s, \eta) ~ d\eta \\
&+ \int e^{i s \varphi} s \nabla_{\xi} \varphi \nabla_{\eta} \varphi \mu_0^{k+1, 0, 1}(\xi, \eta) \nabla_{\xi} \widehat{f}(s, \xi - \eta) \widehat{f}(s, \eta) ~ d\eta + \int e^{i s \varphi} s \nabla_{\xi} \varphi \mu_0^{k+1, 0, 0}(\xi, \eta) \widehat{f}(s, \xi - \eta) \widehat{f}(s, \eta) ~ d\eta \\
&+ \int e^{i s \varphi} \mu_0^{k, 0, 1}(\xi, \eta) \nabla_{\xi} \widehat{f}(s, \xi - \eta) \widehat{f}(s, \eta) ~ d\eta + \int e^{i s \varphi} \mu_0^{k+1, 0, 0}(\xi, \eta) \nabla_{\xi} \widehat{f}(s, \xi - \eta) \widehat{f}(s, \eta) ~ d\eta
\end{align*}
On the terms with $s \nabla_{\eta} \varphi$, we apply an integration by parts in frequency~: 
\begin{subequations} \label{est_L2x2_calculs1_part1}
\begin{align}
&= \int e^{i s \varphi} s \nabla_{\xi} \varphi \mu_0^{k+1, 0, 0}(\xi, \eta) \widehat{f}(s, \xi - \eta) \widehat{f}(s, \eta) ~ d\eta \\
&+ \int e^{i s \varphi} s \nabla_{\xi} \varphi \nabla_{\xi} \varphi \mu_0^{k+1, -1, 1}(\xi, \eta) \widehat{f}(s, \xi - \eta) \widehat{f}(s, \eta) ~ d\eta \\
&+ \int e^{i s \varphi} s \nabla_{\xi} \varphi \nabla_{\xi} \varphi \mu_0^{k+1, 0, 1}(\xi, \eta) \nabla_{\eta} \widehat{f}(s, \xi - \eta) \widehat{f}(s, \eta) ~ d\eta \\
&+ \int e^{i s \varphi} s \nabla_{\xi} \varphi \nabla_{\xi} \varphi \mu_0^{k+1, 0, 1}(\xi, \eta) \widehat{f}(s, \xi - \eta) \nabla_{\eta} \widehat{f}(s, \eta) ~ d\eta \\
&+ \int e^{i s \varphi} s \nabla_{\xi} \varphi \nabla_{\eta} \varphi \mu_0^{k+1, 0, 1}(\xi, \eta) \nabla_{\xi} \widehat{f}(s, \xi - \eta) \widehat{f}(s, \eta) ~ d\eta \\
\nextParentEquation[est_L2x2_calculs1_part2]
&+ \int e^{i s \varphi} \mu_0^{k, 0, 0}(\xi, \eta) \widehat{f}(s, \xi - \eta) \widehat{f}(s, \eta) ~ d\eta \\
&+ \int e^{i s \varphi} \mu_0^{k, -1, 1}(\xi, \eta) \widehat{f}(s, \xi - \eta) \widehat{f}(s, \eta) ~ d\eta \\
&+ \int e^{i s \varphi} \mu_0^{k, 0, 1}(\xi, \eta) \nabla_{\eta} \widehat{f}(s, \xi - \eta) \widehat{f}(s, \eta) ~ d\eta \\
&+ \int e^{i s \varphi} \mu_0^{k, 0, 1}(\xi, \eta) \widehat{f}(s, \xi - \eta) \nabla_{\eta} \widehat{f}(s, \eta) ~ d\eta \\
&+ \int e^{i s \varphi} \mu_0^{k+1, 0, -1}(\xi, \eta) \widehat{f}(s, \xi - \eta) \widehat{f}(s, \eta) ~ d\eta \\
&+ \int e^{i s \varphi} \mu_0^{k+1, 0, 0}(\xi, \eta) \nabla_{\eta} \widehat{f}(s, \xi - \eta) \widehat{f}(s, \eta) ~ d\eta \\
&+ \int e^{i s \varphi} \mu_0^{k, 1, 1}(\xi, \eta) \Delta_{\xi} \widehat{f}(s, \xi - \eta) \widehat{f}(s, \eta) ~ d\eta \\
\nextParentEquation[est_L2x2_calculs1_part3]
&+ \int e^{i s \varphi} s^{-1} \mu_0^{k-1, 0, 0}(\xi, \eta) \widehat{f}(s, \xi - \eta) \widehat{f}(s, \eta) ~ d\eta \\
&+ \int e^{i s \varphi} s^{-1} \mu_0^{k-1, -1, 1}(\xi, \eta) \widehat{f}(s, \xi - \eta) \widehat{f}(s, \eta) ~ d\eta \\
&+ \int e^{i s \varphi} s^{-1} \mu_0^{k-1, 0, 1}(\xi, \eta) \nabla_{\eta} \widehat{f}(s, \xi - \eta) \widehat{f}(s, \eta) ~ d\eta \\
&+ \int e^{i s \varphi} s^{-1} \mu_0^{k-1, 0, 1}(\xi, \eta) \widehat{f}(s, \xi - \eta) \nabla_{\eta} \widehat{f}(s, \eta) ~ d\eta
\end{align}
\end{subequations}
where we grouped all terms containing a factor $s$ in \eqref{est_L2x2_calculs1_part1}, then all terms without any $s$ factor in \eqref{est_L2x2_calculs1_part2}, then all terms with a $s^{-1}$ factor in \eqref{est_L2x2_calculs1_part3}. In \eqref{est_L2x2_calculs1_part1}, we have a $\nabla_{\xi} \varphi$, so we can apply lemma \ref{identite_fond} and an integration by parts in frequency whenever it is possible to obtain~: 
\begin{subequations} \label{est_L2x2_calculs2_part1}
\begin{align}
\eqref{est_L2x2_init} &= \int e^{i s \varphi} s \varphi \mu_0^{k, 0, 0}(\xi, \eta) \widehat{f}(s, \xi - \eta) \widehat{f}(s, \eta) ~ d\eta \\
&+ \int e^{i s \varphi} s \varphi \nabla_{\xi} \varphi \mu_0^{k, -1, 1}(\xi, \eta) \widehat{f}(s, \xi - \eta) \widehat{f}(s, \eta) ~ d\eta \\
&+ \int e^{i s \varphi} s \varphi \nabla_{\xi} \varphi \mu_0^{k, 0, 1}(\xi, \eta) \nabla_{\eta} \widehat{f}(s, \xi - \eta) \widehat{f}(s, \eta) ~ d\eta \\
&+ \int e^{i s \varphi} s \varphi \nabla_{\xi} \varphi \mu_0^{k, 0, 1}(\xi, \eta) \widehat{f}(s, \xi - \eta) \nabla_{\eta} \widehat{f}(s, \eta) ~ d\eta \\
&+ \int e^{i s \varphi} s \varphi \nabla_{\eta} \varphi \mu_0^{k, 0, 1}(\xi, \eta) \nabla_{\xi} \widehat{f}(s, \xi - \eta) \widehat{f}(s, \eta) ~ d\eta \\
\nextParentEquation[est_L2x2_calculs2_part2]
&+ \int e^{i s \varphi} \nabla_{\xi} \varphi \mu_0^{k, 1, 1}(\xi, \eta) \nabla_{\eta} \widehat{f}(s, \xi - \eta) \nabla_{\eta} \widehat{f}(s, \eta) ~ d\eta \\
\nextParentEquation[est_L2x2_calculs2_part3]
&+ \int e^{i s \varphi} \mu_0^{k, 0, 0}(\xi, \eta) \widehat{f}(s, \xi - \eta) \widehat{f}(s, \eta) ~ d\eta \\
&+ \int e^{i s \varphi} \mu_0^{k, -1, 1}(\xi, \eta) \widehat{f}(s, \xi - \eta) \widehat{f}(s, \eta) ~ d\eta \\
&+ \int e^{i s \varphi} \mu_0^{k, 0, 1}(\xi, \eta) \nabla_{\eta} \widehat{f}(s, \xi - \eta) \widehat{f}(s, \eta) ~ d\eta \\
&+ \int e^{i s \varphi} \mu_0^{k, 0, 1}(\xi, \eta) \widehat{f}(s, \xi - \eta) \nabla_{\eta} \widehat{f}(s, \eta) ~ d\eta \\
&+ \int e^{i s \varphi} \mu_0^{k+1, 0, -1}(\xi, \eta) \widehat{f}(s, \xi - \eta) \widehat{f}(s, \eta) ~ d\eta \\
&+ \int e^{i s \varphi} \mu_0^{k+1, 0, 0}(\xi, \eta) \nabla_{\eta} \widehat{f}(s, \xi - \eta) \widehat{f}(s, \eta) ~ d\eta \\
&+ \int e^{i s \varphi} \mu_0^{k, 1, 1}(\xi, \eta) \nabla^2_{\xi} \widehat{f}(s, \xi - \eta) \widehat{f}(s, \eta) ~ d\eta \\
&+ \int e^{i s \varphi} s^{-1} \mu_0^{k-1, 0, 0}(\xi, \eta) \widehat{f}(s, \xi - \eta) \widehat{f}(s, \eta) ~ d\eta \\
&+ \int e^{i s \varphi} s^{-1} \mu_0^{k-1, -1, 1}(\xi, \eta) \widehat{f}(s, \xi - \eta) \widehat{f}(s, \eta) ~ d\eta \\
&+ \int e^{i s \varphi} s^{-1} \mu_0^{k-1, 0, 1}(\xi, \eta) \nabla_{\eta} \widehat{f}(s, \xi - \eta) \widehat{f}(s, \eta) ~ d\eta \\
&+ \int e^{i s \varphi} s^{-1} \mu_0^{k-1, 0, 1}(\xi, \eta) \widehat{f}(s, \xi - \eta) \nabla_{\eta} \widehat{f}(s, \eta) ~ d\eta
\end{align}
\end{subequations}
Again, on the only middle term remaining \eqref{est_L2x2_calculs2_part2}, we apply lemma \ref{identite_fond} and an integration by parts in frequency to obtain~: 
\begin{subequations} \label{est_L2x2_calculs3_part1}
\begin{align}
\eqref{est_L2x2_init} &= \int e^{i s \varphi} s \varphi \mu_0^{k, -1, 1}(\xi, \eta) \widehat{f}(s, \xi - \eta) \widehat{f}(s, \eta) ~ d\eta \\
&+ \int e^{i s \varphi} s \varphi \mu_0^{k, 0, 1}(\xi, \eta) \nabla_{\eta} \widehat{f}(s, \xi - \eta) \widehat{f}(s, \eta) ~ d\eta \\
&+ \int e^{i s \varphi} s \varphi \mu_0^{k, 0, 1}(\xi, \eta) \widehat{f}(s, \xi - \eta) \nabla_{\eta} \widehat{f}(s, \eta) ~ d\eta \\
\nextParentEquation[est_L2x2_calculs3_part2]
&+ \int e^{i s \varphi} \varphi \mu_0^{k-1, 1, 1}(\xi, \eta) \nabla_{\eta} \widehat{f}(s, \xi - \eta) \nabla_{\eta} \widehat{f}(s, \eta) ~ d\eta \label{est_L2x2_calculs3_part2-1} \\
&+ \int e^{i s \varphi} s^{-1} \mu_0^{k-1, 2, 0}(\xi, \eta) \nabla_{\eta} \widehat{f}(s, \xi - \eta) \nabla_{\eta} \widehat{f}(s, \eta) ~ d\eta \\
&+ \int e^{i s \varphi} s^{-1} \mu_0^{k-1, 2, 1}(\xi, \eta) \nabla^2_{\eta} \widehat{f}(s, \xi - \eta) \nabla_{\eta} \widehat{f}(s, \eta) ~ d\eta \\
&+ \int e^{i s \varphi} s^{-1} \mu_0^{k-1, 2, 1}(\xi, \eta) \nabla_{\eta} \widehat{f}(s, \xi - \eta) \nabla^2_{\eta} \widehat{f}(s, \eta) ~ d\eta \\
\nextParentEquation[est_L2x2_calculs3_part3]
&+ \int e^{i s \varphi} \mu_0^{k, -1, 1}(\xi, \eta) \widehat{f}(s, \xi - \eta) \widehat{f}(s, \eta) ~ d\eta \\
&+ \int e^{i s \varphi} \mu_0^{k, 0, 1}(\xi, \eta) \nabla_{\eta} \widehat{f}(s, \xi - \eta) \widehat{f}(s, \eta) ~ d\eta \\
&+ \int e^{i s \varphi} \mu_0^{k, 0, 1}(\xi, \eta) \widehat{f}(s, \xi - \eta) \nabla_{\eta} \widehat{f}(s, \eta) ~ d\eta \\
&+ \int e^{i s \varphi} \mu_0^{k, 1, 1}(\xi, \eta) \nabla^2_{\xi} \widehat{f}(s, \xi - \eta) \widehat{f}(s, \eta) ~ d\eta \\
&+ \int e^{i s \varphi} s^{-1} \mu_0^{k-1, -1, 1}(\xi, \eta) \widehat{f}(s, \xi - \eta) \widehat{f}(s, \eta) ~ d\eta \\
\nextParentEquation[est_L2x2_calculs3_part4]
&+ \int e^{i s \varphi} s^{-1} \mu_0^{k-1, 0, 1}(\xi, \eta) \nabla_{\eta} \widehat{f}(s, \xi - \eta) \widehat{f}(s, \eta) ~ d\eta \\
&+ \int e^{i s \varphi} s^{-1} \mu_0^{k-1, 0, 1}(\xi, \eta) \widehat{f}(s, \xi - \eta) \nabla_{\eta} \widehat{f}(s, \eta) ~ d\eta
\end{align}
\end{subequations}
Let us denote 
\begin{align}
\eqref{est_L2x2_calculs2_part1} + \eqref{est_L2x2_calculs3_part1} + \eqref{est_L2x2_calculs3_part2-1} \label{est_L2x2_term1} \\
\eqref{est_L2x2_calculs1_part2} + \eqref{est_L2x2_calculs1_part3} + \eqref{est_L2x2_calculs2_part3} + \left( \eqref{est_L2x2_calculs3_part2} - \eqref{est_L2x2_calculs3_part2-1} \right) + \eqref{est_L2x2_calculs3_part3} + \eqref{est_L2x2_calculs3_part4} \label{est_L2x2_term2}
\end{align}
that is \eqref{est_L2x2_term1} has all the terms with $\varphi$ and \eqref{est_L2x2_term2} all the others. 

Concerning $\pm \pm$ interactions, we showed~: 
\[ \partial_s \xi^{\alpha} \nabla_{\xi}^2 \widehat{f}(s, \xi) = \eqref{est_L2x2_term1} + \eqref{est_L2x2_term2} + \mathcal{F} \mathcal{N}(\Lambda D^{\alpha} e^{\pm i s \Lambda} |x|^2 f, u) \]
and the contribution of the last term has already been estimated. 

\begin{Intuit}
We will show that $\Vert \eqref{est_L2x2_term2} \Vert_{L^2} \lesssim s^{-1+\gamma_k} \Vert u \Vert_X^2$, so that
\[ \int_1^t \int \xi^{\alpha} \nabla_{\xi}^2 \widehat{f}(s, \xi) \eqref{est_L2x2_term2} ~ d\xi ds \lesssim \int_1^t \Vert \Lambda |x|^2 f \Vert_{H^k} \Vert \eqref{est_L2x2_term2} \Vert_{L^2} ~ ds \lesssim \int_1^t s^{-1+2\gamma_k} \Vert u \Vert_X^3 ~ ds \lesssim t^{2\gamma_k} \Vert u \Vert_X^3 \]

On the other hand, for \eqref{est_L2x2_term1}, the presence of $\varphi$ encourages to apply an integration by parts in time. If the derivative hits one of the terms from \eqref{est_L2x2_term1}, we can apply lemma \ref{lem_decnonres} and obtain a sufficiant decay to compensate the factor $s$. However, if the derivative hits the term across, it produces a $\partial_s \xi^{\alpha} \nabla_{\xi}^2 \widehat{f}(s, \xi)$. In this case, we develop again using the decomposition above and note that (a) the term $\mathcal{N}(\Lambda D^{\alpha} e^{\pm i s \Lambda} |x|^2 f, u)$ can be avoided (there is a $\Lambda^{-1}$ in front of \eqref{est_L2x2_term1} that we can freely transfer)~; (b) when we develop and get \eqref{est_L2x2_term2}, what remains of \eqref{est_L2x2_term1} can be estimated by $s^{\gamma_k}$ and we have the correct total decay~; (c) when we develop and get \eqref{est_L2x2_term1}, we can reapply an integration by parts using a symmetry argument. 
\end{Intuit}

\subsubsection{Estimate of \eqref{est_L2x2_term2}} The terms of \eqref{est_L2x2_term2} are~: 
\begin{subequations} \label{est_L2x2_term2-1}
\begin{align}
\eqref{est_L2x2_term2}
&= \int e^{i s \varphi} s^{-1} \mu_0^{k-1, 2, 0}(\xi, \eta) \nabla_{\eta} \widehat{f}(s, \xi - \eta) \nabla_{\eta} \widehat{f}(s, \eta) ~ d\eta \label{est_L2x2_term2-1-1} \\
&+ \int e^{i s \varphi} s^{-1} \mu_0^{k-1, 2, 1}(\xi, \eta) \nabla^2_{\eta} \widehat{f}(s, \xi - \eta) \nabla_{\eta} \widehat{f}(s, \eta) ~ d\eta \label{est_L2x2_term2-1-2} \\
&+ \int e^{i s \varphi} s^{-1} \mu_0^{k-1, 2, 1}(\xi, \eta) \nabla_{\eta} \widehat{f}(s, \xi - \eta) \nabla^2_{\eta} \widehat{f}(s, \eta) ~ d\eta \label{est_L2x2_term2-1-3} \\
&+ \int e^{i s \varphi} s^{-1} \mu_0^{k-1, -1, 1}(\xi, \eta) \widehat{f}(s, \xi - \eta) \widehat{f}(s, \eta) ~ d\eta \label{est_L2x2_term2-1-4} \\
&+ \int e^{i s \varphi} s^{-1} \mu_0^{k-1, 0, 1}(\xi, \eta) \nabla_{\eta} \widehat{f}(s, \xi - \eta) \widehat{f}(s, \eta) ~ d\eta \label{est_L2x2_term2-1-5} \\
&+ \int e^{i s \varphi} s^{-1} \mu_0^{k-1, 0, 1}(\xi, \eta) \widehat{f}(s, \xi - \eta) \nabla_{\eta} \widehat{f}(s, \eta) ~ d\eta \label{est_L2x2_term2-1-6} \\
\nextParentEquation[est_L2x2_term2-2]
&+ \int e^{i s \varphi} \mu_0^{k, -1, 1}(\xi, \eta) \widehat{f}(s, \xi - \eta) \widehat{f}(s, \eta) ~ d\eta \label{est_L2x2_term2-2-1} \\
&+ \int e^{i s \varphi} \mu_0^{k, 0, 1}(\xi, \eta) \nabla_{\eta} \widehat{f}(s, \xi - \eta) \widehat{f}(s, \eta) ~ d\eta \label{est_L2x2_term2-2-2} \\
&+ \int e^{i s \varphi} \mu_0^{k, 0, 1}(\xi, \eta) \widehat{f}(s, \xi - \eta) \nabla_{\eta} \widehat{f}(s, \eta) ~ d\eta \label{est_L2x2_term2-2-3} \\
&+ \int e^{i s \varphi} \mu_0^{k, 1, 1}(\xi, \eta) \nabla^2_{\xi} \widehat{f}(s, \xi - \eta) \widehat{f}(s, \eta) ~ d\eta \label{est_L2x2_term2-2-4}
\end{align}
\end{subequations}
In \eqref{est_L2x2_term2-1}, we have a factor $s^{-1}$ but possible a singularity $|\xi|^{-1}$ in the case $k = 0$. In \eqref{est_L2x2_term2-2}, we have no factor $s^{-1}$ but also no singularity. In the case $k = 0$~: 
\begin{align*}
\Vert \eqref{est_L2x2_term2-1-1} \Vert_{L^2} &\lesssim s^{-1} \Vert T_{\mu_0}(e^{\pm i s \Lambda} x f, \Lambda^2 e^{\pm i s \Lambda} x f) \Vert_{L^{6/5}} \lesssim s^{-1} \Vert e^{\pm i s \Lambda} x f \Vert_{L^3} \Vert x f \Vert_{H^2} \lesssim s^{-1} \Vert x f \Vert_{H^2}^2 \lesssim s^{-1+\gamma_0} \Vert u \Vert_X^2 \\
\Vert \eqref{est_L2x2_term2-1-2} \Vert_{L^2} &\lesssim s^{-1} \Vert T_{\mu_0}(e^{\pm i s \Lambda} \Lambda x^2 f, \Lambda^2 e^{\pm i s \Lambda} x f) \Vert_{L^{6/5}} \lesssim s^{-1} \Vert \Lambda x^2 f \Vert_{L^2} \Vert \Lambda^2 e^{\pm i s \Lambda} x f \Vert_{L^3} \lesssim s^{-4/3+ \gamma_0} \Vert u \Vert_X \Vert \Lambda^{8/3} x f \Vert_{L^{3/2}} \\
&\lesssim s^{-4/3+2\gamma} \Vert u \Vert_X^2 \lesssim s^{-1+\gamma_0} \Vert u \Vert_X^2 \\
\Vert \eqref{est_L2x2_term2-1-3} \Vert_{L^2} &\lesssim s^{-1} \Vert T_{\mu_0}(e^{\pm i s \Lambda} \Lambda x f, e^{\pm i s \Lambda} \Lambda^2 x^2 f) \Vert_{L^{6/5}} \lesssim s^{-1} \Vert \Lambda e^{\pm i s \Lambda} x f \Vert_{L^3} \Vert \Lambda x^2 f \Vert_{H^1} \lesssim s^{-4/3+\gamma_1} \Vert \Lambda^{5/3} x f \Vert_{L^{3/2}} \Vert u \Vert_X \\
&\lesssim s^{-4/3+2\gamma} \Vert u \Vert_X^2 \lesssim s^{-1+\gamma_0} \Vert u \Vert_X^2
\end{align*}
The other terms are simpler or similar. In the case $k = 5$~: 
\begin{align*}
\Vert \eqref{est_L2x2_term2-1-1} \Vert_{L^2} &\lesssim s^{-1} \Vert T_{\mu_0}(e^{\pm i s \Lambda} x f, \Lambda^2 e^{\pm i s \Lambda} x f) \Vert_{H^4} \lesssim s^{-1} \Vert e^{\pm i s \Lambda} x f \Vert_{W^{3, 6}} \Vert x f \Vert_{H^6} \lesssim s^{-1+b} \Vert x f \Vert_{H^4} \Vert u \Vert_X \lesssim s^{-1+\gamma} \Vert u \Vert_X^2 \\
\Vert \eqref{est_L2x2_term2-1-2} \Vert_{L^2} &\lesssim s^{-1} \Vert T_{\mu_0}(e^{\pm i s \Lambda} \Lambda x^2 f, \Lambda^2 e^{\pm i s \Lambda} x f) \Vert_{H^4} \lesssim s^{-1} \Vert \Lambda x^2 f \Vert_{H^4} \Vert \Lambda^2 e^{\pm i s \Lambda} x f \Vert_{W^{1, 6}} + s^{-1} \Vert \Lambda e^{\pm i s \Lambda} x^2 f \Vert_{W^{1, 6}} \Vert x f \Vert_{H^6} \\
&\lesssim s^{-1} \Vert \Lambda x^2 f \Vert_{H^4} \Vert x f \Vert_{H^4} + s^{-1} \Vert \Lambda x^2 f \Vert_{H^2} \Vert x f \Vert_{H^6} \lesssim \left( s^{-1+\gamma_4} + s^{-1+\gamma_2+b} \right) \lesssim s^{-1+\gamma} \Vert u \Vert_X^2 \\
\Vert \eqref{est_L2x2_term2-1-3} \Vert_{L^2} &\lesssim s^{-1} \Vert T_{\mu_0}(e^{\pm i s \Lambda} \Lambda x f, e^{\pm i s \Lambda} \Lambda^2 x^2 f) \Vert_{H^4} \lesssim s^{-1} \Vert \Lambda e^{\pm i s \Lambda} x f \Vert_{W^{1, 6}} \Vert \Lambda x^2 f \Vert_{H^5} + s^{-1} \Vert x f \Vert_{H^5} \Vert e^{\pm i s \Lambda} \Lambda x^2 f \Vert_{W^{2, 6}} \\
&\lesssim s^{-1} \Vert x f \Vert_{H^3} \Vert \Lambda x^2 f \Vert_{H^5} + s^{-1} \Vert x f \Vert_{H^5} \Vert \Lambda x^2 f \Vert_{H^3} \lesssim s^{-1+\gamma} \Vert u \Vert_X^2
\end{align*}
where we used that $\gamma_2 + b = \gamma_2 + \frac{\gamma_4 + \varepsilon}{2} = \frac{4}{5} \gamma + \frac{7}{10} \gamma_0 + \frac{\varepsilon}{2} \leq \gamma$ if $5\gamma_0 < \gamma$ and $\varepsilon$ is small enough with respect to $\gamma, \gamma_5-5\gamma_0$. Again, the other terms are similar or simpler. 

For \eqref{est_L2x2_term2-2}, let us estimate~: 
\begin{align*}
\Vert \eqref{est_L2x2_term2-2-1} \Vert_{L^2} &= \Vert T_{\mu_0}(\Lambda u, \Lambda^{-1} u) \Vert_{H^k} \lesssim \Vert u \Vert_{W^{k+2, \infty-}} \Vert x f \Vert_{H^k} \lesssim s^{-1+\delta} \Vert u \Vert_X^2 \lesssim s^{-1+\gamma_k} \Vert u \Vert_X^2 \\
\Vert \eqref{est_L2x2_term2-2-4}_{k = 0} \Vert_{L^2} &= \Vert T_{\mu_0}(e^{\pm i s \Lambda} \Lambda x^2 f, \Lambda u) \Vert_{L^2} \lesssim \Vert \Lambda x^2 f \Vert_{L^2} \Vert u \Vert_{\dot{B}^1_{\infty, 1}} \lesssim s^{-1+\gamma_0} \Vert u \Vert_X^2 \\
\Vert \eqref{est_L2x2_term2-2-4}_{k = 5} \Vert_{L^2} &\lesssim \Vert \Lambda x^2 f \Vert_{H^5} \Vert u \Vert_{\dot{B}^1_{\infty, 1}} + \Vert \Lambda x^2 f \Vert_{H^4} \Vert u \Vert_{W^{7, \infty-}} \lesssim \left( s^{-1+\gamma} + s^{-1+\delta+\gamma_4} \right) \Vert u \Vert_X^2 \lesssim s^{-1+\gamma} \Vert u \Vert_X^2
\end{align*}
where we used lemma \ref{lem_symb_Besov}, and then that when $\delta$ is small enough, $\gamma_4 + \delta \leq \gamma$. The other terms are simpler or similar. 

\subsubsection{Estimate of \eqref{est_L2x2_term1}} consists in~: 
\[ \begin{aligned}
\eqref{est_L2x2_term1} &= \int e^{i s \varphi} s \varphi \mu_0^{k, -1, 1}(\xi, \eta) \widehat{f}(s, \xi - \eta) \widehat{f}(s, \eta) ~ d\eta \\
&+ \int e^{i s \varphi} s \varphi \mu_0^{k, 0, 1}(\xi, \eta) \nabla_{\eta} \widehat{f}(s, \xi - \eta) \widehat{f}(s, \eta) ~ d\eta + \int e^{i s \varphi} s \varphi \mu_0^{k, 0, 1}(\xi, \eta) \widehat{f}(s, \xi - \eta) \nabla_{\eta} \widehat{f}(s, \eta) ~ d\eta \\
&+ \int e^{i s \varphi} |\eta| \varphi \mu_0^{k-1, 0, 1}(\xi, \eta) \nabla_{\eta} \widehat{f}(s, \xi - \eta) \nabla_{\eta} \widehat{f}(s, \eta) ~ d\eta 
\end{aligned} \]

Let us apply an integration by parts in time~: 
\begin{subequations} \label{est_L2x2_term1-1}
\begin{align}
&\int_1^t \int \xi^{\alpha} \Delta_{\xi}^2 \widehat{f}(s, \xi) \eqref{est_L2x2_term1} ~ d\xi ds \nonumber \\
&= \int_1^t \int \int \xi^{\alpha} \Delta_{\xi}^2 \partial_s \widehat{f}(s, \xi) e^{i s \varphi} s \mu_0^{k, -1, 1}(\xi, \eta) \widehat{f}(s, \xi - \eta) \widehat{f}(s, \eta) ~ d\eta d\xi ds \\
&+ \int_1^t \int \int \xi^{\alpha} \Delta_{\xi}^2 \partial_s \widehat{f}(s, \xi) e^{i s \varphi} s \mu_0^{k, 0, 1}(\xi, \eta) \nabla_{\eta} \widehat{f}(s, \xi - \eta) \widehat{f}(s, \eta) ~ d\eta d\xi ds \\
&+ \int_1^t \int \int \xi^{\alpha} \Delta_{\xi}^2 \partial_s \widehat{f}(s, \xi) e^{i s \varphi} s \mu_0^{k, 0, 1}(\xi, \eta) \widehat{f}(s, \xi - \eta) \nabla_{\eta} \widehat{f}(s, \eta) ~ d\eta d\xi ds \\
&+ \int_1^t \int \int \xi^{\alpha} \Delta_{\xi}^2 \partial_s \widehat{f}(s, \xi) e^{i s \varphi} \mu_0^{k-1, 1, 1}(\xi, \eta) \nabla_{\eta} \widehat{f}(s, \xi - \eta) \nabla_{\eta} \widehat{f}(s, \eta) ~ d\eta d\xi ds \\
\nextParentEquation[est_L2x2_term1-2-1]
&+ \int_1^t \int \int \xi^{\alpha} \Delta_{\xi}^2 \widehat{f}(s, \xi) e^{i s \varphi} s \mu_0^{k, -1, 1}(\xi, \eta) \partial_s \widehat{f}(s, \xi - \eta) \widehat{f}(s, \eta) ~ d\eta d\xi ds \label{est_L2x2_term1-2-1-1} \\
&+ \int_1^t \int \int \xi^{\alpha} \Delta_{\xi}^2 \widehat{f}(s, \xi) e^{i s \varphi} s \mu_0^{k, 0, 1}(\xi, \eta) \nabla_{\eta} \partial_s \widehat{f}(s, \xi - \eta) \widehat{f}(s, \eta) ~ d\eta d\xi ds \label{est_L2x2_term1-2-1-2} \\
&+ \int_1^t \int \int \xi^{\alpha} \Delta_{\xi}^2 \widehat{f}(s, \xi) e^{i s \varphi} s \mu_0^{k, 0, 1}(\xi, \eta) \partial_s \widehat{f}(s, \xi - \eta) \nabla_{\eta} \widehat{f}(s, \eta) ~ d\eta d\xi ds \label{est_L2x2_term1-2-1-3} \\
&+ \int_1^t \int \int \xi^{\alpha} \Delta_{\xi}^2 \widehat{f}(s, \xi) e^{i s \varphi} \mu_0^{k-1, 1, 1}(\xi, \eta) \nabla_{\eta} \partial_s \widehat{f}(s, \xi - \eta) \nabla_{\eta} \widehat{f}(s, \eta) ~ d\eta d\xi ds \label{est_L2x2_term1-2-1-4} \\
&+ \int_1^t \int \int \xi^{\alpha} \Delta_{\xi}^2 \widehat{f}(s, \xi) e^{i s \varphi} s \mu_0^{k, -1, 1}(\xi, \eta) \widehat{f}(s, \xi - \eta) \partial_s \widehat{f}(s, \eta) ~ d\eta d\xi ds \label{est_L2x2_term1-2-1-5} \\
&+ \int_1^t \int \int \xi^{\alpha} \Delta_{\xi}^2 \widehat{f}(s, \xi) e^{i s \varphi} s \mu_0^{k, 0, 1}(\xi, \eta) \nabla_{\eta} \widehat{f}(s, \xi - \eta) \partial_s \widehat{f}(s, \eta) ~ d\eta d\xi ds \label{est_L2x2_term1-2-1-6} \\
&+ \int_1^t \int \int \xi^{\alpha} \Delta_{\xi}^2 \widehat{f}(s, \xi) e^{i s \varphi} s \mu_0^{k, 0, 1}(\xi, \eta) \widehat{f}(s, \xi - \eta) \nabla_{\eta} \partial_s \widehat{f}(s, \eta) ~ d\eta d\xi ds \label{est_L2x2_term1-2-1-7} \\
\nextParentEquation[est_L2x2_term1-2-2]
&+ \int \int \xi^{\alpha} \Delta_{\xi}^2 \widehat{f}(t, \xi) e^{i t \varphi} t \mu_0^{k, -1, 1}(\xi, \eta) \widehat{f}(t, \xi - \eta) \widehat{f}(t, \eta) ~ d\eta d\xi \label{est_L2x2_term1-2-2-1} \\
&+ \int \int \xi^{\alpha} \Delta_{\xi}^2 \widehat{f}(t, \xi) e^{i t \varphi} t \mu_0^{k, 0, 1}(\xi, \eta) \nabla_{\eta} \widehat{f}(t, \xi - \eta) \widehat{f}(t, \eta) ~ d\eta d\xi \label{est_L2x2_term1-2-2-2} \\
&+ \int \int \xi^{\alpha} \Delta_{\xi}^2 \widehat{f}(t, \xi) e^{i t \varphi} t \mu_0^{k, 0, 1}(\xi, \eta) \widehat{f}(t, \xi - \eta) \nabla_{\eta} \widehat{f}(t, \eta) ~ d\eta d\xi \label{est_L2x2_term1-2-2-3} \\
&+ \int \int \xi^{\alpha} \Delta_{\xi}^2 \widehat{f}(t, \xi) e^{i t \varphi} \mu_0^{k-1, 1, 1}(\xi, \eta) \nabla_{\eta} \widehat{f}(t, \xi - \eta) \nabla_{\eta} \widehat{f}(t, \eta) ~ d\eta d\xi \label{est_L2x2_term1-2-2-4} \\
\nextParentEquation[est_L2x2_term1-2-3]
&+ \int \int \xi^{\alpha} \Delta_{\xi}^2 \widehat{f}(1, \xi) e^{i \varphi} \mu_0^{k, -1, 1}(\xi, \eta) \widehat{f}(1, \xi - \eta) \widehat{f}(1, \eta) ~ d\eta d\xi \\
&+ \int \int \xi^{\alpha} \Delta_{\xi}^2 \widehat{f}(1, \xi) e^{i \varphi} \mu_0^{k, 0, 1}(\xi, \eta) \nabla_{\eta} \widehat{f}(1, \xi - \eta) \widehat{f}(1, \eta) ~ d\eta d\xi \\
&+ \int \int \xi^{\alpha} \Delta_{\xi}^2 \widehat{f}(1, \xi) e^{i \varphi} \mu_0^{k, 0, 1}(\xi, \eta) \widehat{f}(1, \xi - \eta) \nabla_{\eta} \widehat{f}(1, \eta) ~ d\eta d\xi \\
&+ \int \int \xi^{\alpha} \Delta_{\xi}^2 \widehat{f}(1, \xi) e^{i \varphi} \mu_0^{k-1, 1, 1}(\xi, \eta) \nabla_{\eta} \widehat{f}(1, \xi - \eta) \nabla_{\eta} \widehat{f}(1, \eta) ~ d\eta d\xi
\end{align} 
\end{subequations}
\eqref{est_L2x2_term1-1} contains again $D^{\alpha} |x|^2 \partial_t f$. \eqref{est_L2x2_term1-2-1} contains the terms where the time derivative produces a term controlled by lemma \ref{lem_decnonres}. \eqref{est_L2x2_term1-2-2} contains final time terms, \eqref{est_L2x2_term1-2-3} initial time terms. 

We will first show that 
\begin{equation} \eqref{est_L2x2_term1-2-1} + \eqref{est_L2x2_term1-2-2} + \eqref{est_L2x2_term1-2-3} \label{est_L2x2_term1-2} \end{equation}
is controlled by $t^{2\gamma_k} \Vert u \Vert_X^4$. 

For \eqref{est_L2x2_term1-2-1}~:
\begin{align*}
\eqref{est_L2x2_term1-2-1-2} &\lesssim \int_1^t s \Vert \Lambda |x|^2 f \Vert_{H^k} \Vert T_{\mu_0}(e^{\pm i s \Lambda} \Lambda x \partial_s f, u) \Vert_{H^k} ~ ds \lesssim \int_1^t s^{1+\gamma_k} \Vert u \Vert_X \Vert x \partial_s f \Vert_{H^{k+1}} \Vert u \Vert_{W^{k+1, \infty-}} ~ ds \\
&\lesssim \int_1^t s^{-1+\frac{3\gamma_k}{2} + \delta + \tau} \Vert u \Vert_X^4 ~ ds \lesssim t^{2\gamma_k} \Vert u \Vert_X^4 \\
\eqref{est_L2x2_term1-2-1-3} &\lesssim \int_1^t s \Vert \Lambda |x|^2 f \Vert_{H^k} \Vert T_{\mu_0}(e^{\pm i s \Lambda} \Lambda \partial_s f, e^{\pm i s \Lambda} x f) \Vert_{H^k} ~ ds \\
&\lesssim \int_1^t s^{1+\gamma_k} \Vert u \Vert_X \left( \Vert \partial_s f \Vert_{H^{k+1}} \Vert e^{\pm i s \Lambda} x f \Vert_{W^{1, 6}} + \Vert e^{\pm i s \Lambda} \partial_s f \Vert_{W^{1, 6}} \Vert x f \Vert_{H^k} \right) ~ ds \\
&\lesssim \int_1^t s^{1+\gamma_k} \Vert u \Vert_X \left( \Vert \partial_s f \Vert_{H^{k+1}} \Vert x f \Vert_{H^2} + \Vert \partial_s f \Vert_{H^2} \Vert x f \Vert_{H^k} \right) ~ ds \lesssim s^{-1+\tau+\frac{3 \gamma_k}{2}} \Vert u \Vert_X^4 ~ ds \lesssim t^{2\gamma_k} \Vert u \Vert_X^4 \\
\eqref{est_L2x2_term1-2-1-4}_{k = 0} &\lesssim \int_1^t \Vert \Lambda |x|^2 f \Vert_{L^2} \Vert T_{\mu_0}(e^{\pm i s \Lambda} \Lambda \partial_s x f, e^{\pm i s \Lambda} \Lambda x f) \Vert_{L^{6/5}} ~ ds \lesssim \int_1^t s^{\gamma_0} \Vert u \Vert_X \Vert \partial_s x f \Vert_{H^1} \Vert x f \Vert_{H^2} ~ ds \\
&\lesssim \int_1^t s^{-1+\tau+\gamma_0} \Vert u \Vert_X^4 ~ ds \lesssim t^{2\gamma_0} \Vert u \Vert_X^4 \\
\eqref{est_L2x2_term1-2-1-4}_{k = 5} &\lesssim \int_1^t \Vert \Lambda |x|^2 f \Vert_{H^5} \Vert T_{\mu_0}(e^{\pm i s \Lambda} \Lambda \partial_s x f, e^{\pm i s \Lambda} \Lambda x f) \Vert_{H^4} ~ ds \lesssim \int_1^t s^{\gamma_5} \Vert u \Vert_X \Vert \partial_s x f \Vert_{H^5} \Vert x f \Vert_{H^7} ~ ds \\
&\lesssim \int_1^t s^{-1+\tau+\frac{3\gamma_5}{2} + \frac{\gamma_4}{2}} \Vert u \Vert_X^4 ~ ds \lesssim t^{2\gamma_5} \Vert u \Vert_X^4
\end{align*}
where we used that, by lemma \ref{lem_decnonres}, 
\[ \Vert x \partial_s f \Vert_{H^{k+1}} \lesssim s^{-1+\tau+\frac{\gamma_k}{2}} \Vert u \Vert_X^2, ~~~ \Vert \partial_s f \Vert_{H^{k+1}} \lesssim s^{-1+\tau+\frac{\gamma_k}{2}} \Vert u \Vert_X^2, ~~~ \Vert x f \Vert_{H^k} \lesssim s^{\frac{\gamma_k}{2}} \Vert u \Vert_X \]
and that all the other parameters ($\tau, \delta, \varepsilon, 1/N$) can be chosen small with respect to the $\gamma_k$ and to $\gamma_5 - \gamma_4$. The other terms are simpler or similar. 

For \eqref{est_L2x2_term1-2-2}, 
\begin{align*}
\eqref{est_L2x2_term1-2-2-2} &\lesssim t \Vert \Lambda |x|^2 f \Vert_{H^k} \Vert T_{\mu_0}(\Lambda e^{\pm i s \Lambda} x f, u) \Vert_{H^k} \lesssim t^{1+\gamma_k} \Vert u \Vert_X \Vert x f \Vert_{H^{k+1}} \Vert u \Vert_{W^{k+1, \infty-}} \lesssim t^{\frac{3\gamma_k}{2} + \delta} \Vert u \Vert_X^3 \lesssim t^{2\gamma_k} \Vert u \Vert_X^3 \\
\eqref{est_L2x2_term1-2-2-4}_{k = 0} &\lesssim \Vert \Lambda |x|^2 f \Vert_{L^2} \Vert T_{\mu_0}(\Lambda e^{\pm i s \Lambda} x f, \Lambda e^{\pm i s \Lambda} x f) \Vert_{L^{6/5}} \lesssim t^{\gamma_0} \Vert u \Vert_X \Vert x f \Vert_{H^2}^2 \lesssim t^{2\gamma_0} \Vert u \Vert_X^3 \\
\eqref{est_L2x2_term1-2-2-4}_{k = 5} &\lesssim \Vert \Lambda |x|^2 f \Vert_{H^5} \Vert T_{\mu_0}(\Lambda e^{\pm i s \Lambda} x f, \Lambda e^{\pm i s \Lambda} x f) \Vert_{H^4} \lesssim t^{\gamma_5} \Vert u \Vert_X \Vert x f \Vert_{H^5} \Vert x f \Vert_{H^3} \lesssim t^{2\gamma_5} \Vert u \Vert_X^3
\end{align*}
The other terms are simpler or similar. Finally, for \eqref{est_L2x2_term1-2-3}, the estimates are made as for \eqref{est_L2x2_term1-2-2} (and we even have no growth since we are at initial time). 

\subsubsection{A symmetry argument} \label{section_L2x2_symmetry}
Finally, let us treat \eqref{est_L2x2_term1-1}. It can be written as
\[ \int_1^t \int D^{\alpha} |x|^2 \partial_t f ~ (\ref{est_L2x2_term1}') ~ dx ds \]
where $(\ref{est_L2x2_term1}')$ corresponds in Fourier to \eqref{est_L2x2_term1} divided by $\varphi(\xi, \eta)$. We rewrite~: 
\[ D^{\alpha} |x|^2 \partial_t f \sim \eqref{est_L2x2_term1} + \eqref{est_L2x2_term2} + D^{\alpha} \mathcal{N}(\Lambda e^{\pm i s \Lambda} |x|^2 f, u) \]
plus some terms with $0 \pm$, $\pm 0$ or $0 0$ interactions, which will be treated in the corresponding subsections. Let us start by treating the last term. By integration by parts, we recover~: 
\[ \int_1^t \int \Lambda^{-1} D^{\alpha} \mathcal{N}(\Lambda e^{\pm i s \Lambda} |x|^2 f, u) \Lambda (\ref{est_L2x2_term1}') ~ dx ds \]
But $\Lambda^{-1} D^{\alpha} \mathcal{N}(\Lambda e^{\pm i s \Lambda} |x|^2 f)$ only contains terms already treated in \eqref{est_L2x2_term2}. Therefore, all we have to control is~: 
\[ \int_1^t \int \eqref{est_L2x2_term1} ~ (\ref{est_L2x2_term1}') ~ dx ds + \int_1^t \int \eqref{est_L2x2_term2} ~ (\ref{est_L2x2_term1}') ~ dx ds + \int_1^t \int \eqref{est_L2x2_term2} ~ \Lambda (\ref{est_L2x2_term1}') ~ dx ds \]
On the other hand, we already know that $\Vert \eqref{est_L2x2_term2} \Vert_{L^2} \lesssim s^{1+\gamma_k} \Vert u \Vert_X^2$. Therefore, it is enough to prove that $\Vert (\ref{est_L2x2_term1}') \Vert_{H^1} \lesssim s^{\gamma_k}$ to treat the two terms on the right. The first one will be treated by a symmetry argument. 

We write $(\ref{est_L2x2_term1}')$ as~: 
\begin{subequations}
\begin{align}
\widehat{(\ref{est_L2x2_term1}')} &= \int e^{i s \varphi} s \mu_0^{k, -1, 1}(\xi, \eta) \widehat{f}(s, \xi - \eta) \widehat{f}(s, \eta) ~ d\eta \label{est_L2x2_term1p-1} \\
&+ \int e^{i s \varphi} s \mu_0^{k, 0, 1}(\xi, \eta) \nabla_{\eta} \widehat{f}(s, \xi - \eta) \widehat{f}(s, \eta) ~ d\eta \label{est_L2x2_term1p-2} \\
&+ \int e^{i s \varphi} s \mu_0^{k, 0, 1}(\xi, \eta) \widehat{f}(s, \xi - \eta) \nabla_{\eta} \widehat{f}(s, \eta) ~ d\eta \label{est_L2x2_term1p-3} \\
&+ \int e^{i s \varphi} \mu_0^{k-1, 1, 1}(\xi, \eta) \nabla_{\eta} \widehat{f}(s, \xi - \eta) \nabla_{\eta} \widehat{f}(s, \eta) ~ d\eta \label{est_L2x2_term1p-4} 
\end{align} 
\end{subequations}
Then~:
\begin{align*}
\Vert \eqref{est_L2x2_term1p-2} \Vert_{H^1} &\lesssim s \Vert x f \Vert_{H^{k+2}} \Vert u \Vert_{W^{k+2, \infty-}} \lesssim s^{\delta+\varepsilon/2+\gamma_k/2} \Vert u \Vert_X^2 \lesssim s^{\gamma_k} \Vert u \Vert_X^2 \\
\Vert \eqref{est_L2x2_term1p-4}_{k=0} \Vert_{H^1} &\lesssim \Vert T_{\mu_0}(\Lambda e^{\pm i s \Lambda} x f, \Lambda e^{\pm i s \Lambda} x f) \Vert_{W^{1, 6/5}} \lesssim \Vert x f \Vert_{H^3}^2 \lesssim \Vert u \Vert_X^2 \\
\Vert \eqref{est_L2x2_term1p-4}_{k = 5} \Vert_{H^1} &\lesssim \Vert x f \Vert_{H^6} \Vert e^{\pm i s \Lambda} x f \Vert_{W^{2, 6}} \lesssim \Vert x f \Vert_{H^6} \Vert xf \Vert_{H^3} \lesssim s^b \Vert u \Vert_X^2 \lesssim s^{\gamma} \Vert u \Vert_X^2
\end{align*}
The other terms are simpler or analogous. We even have $\Vert (\ref{est_L2x2_term1}') \Vert_{H^1} \lesssim s^{3\gamma_k/4} \Vert u \Vert_X^2$ by chosing the parameters suitably. 

Finally, to treat 
\[ \int_1^t \int \eqref{est_L2x2_term1} ~ (\ref{est_L2x2_term1}') ~ dx ds \]
we notice that the structure of \eqref{est_L2x2_term1} encourages once again to apply an integration by parts in time. More precisely, we have terms of the form (in Fourier)~: 
\[ \sum_{i, j} \int_1^t s^2 \int \int \int e^{i s \varphi(\xi, \eta)} \varphi(\xi, \eta) \mu_0^{(i)}(\xi, \eta) \widehat{g_i}(\xi - \eta) \widehat{h_i}(\eta) e^{i s \varphi(\xi, \eta')} \mu_0^{(j)}(\xi, \eta') \widehat{g_j}(\xi - \eta') \widehat{h_j}(\eta') ~ d\eta' d\eta d\xi ds \]
for $\mu_0^{(i)}$ symbols of order $0$ depending only on $i$. But then, by regrouping $(i, j)$ and $(j, i)$, by symmetry between $\eta$ and $\eta'$ we obtain~: 
\[ \sum_{i, j} \int_1^t s^2 \int \int \int e^{i s \varphi(\xi, \eta)} e^{i s \varphi(\xi, \eta')} \frac{\varphi(\xi, \eta) + \varphi(\xi, \eta')}{2} \mu_0^{(i)}(\xi, \eta) \mu_0^{(i')}(\xi, \eta') \widehat{g_i}(\xi - \eta) \widehat{h_i}(\eta) \widehat{g_j}(\xi - \eta') \widehat{h_j}(\eta') ~ d\eta' d\eta d\xi ds \]
on which we can apply an integration by parts in time and recover objects controlled by lemma \ref{lem_decnonres}. To make the reasoning absolutely rigorous, we actually sum on all the interactions (because the symbols and the $g, h$ depends on the $\epsilon_2, \epsilon_3$, which can be distinct on each side). Once the derivative is present on $g$ or $h$, we estimate as we estimated $(1')$ above, by noticing that the presence of the time derivative allows a gain of a $s^{-1+\tau}$ decay in the case $k = 0$, $s^{-1+\tau+\gamma/2}$ in the case $k = 5$. Since we already estimated $(1')$ with a decay $s^{3\gamma_k/4}$, we obtain a total decay $s^{-1+3\gamma_k/2+\tau}$, so that once we integrated in time we have
\[ \int_1^t \int \eqref{est_L2x2_term1} ~ (\ref{est_L2x2_term1}') ~ dx ds \lesssim t^{2\gamma_k} \Vert u \Vert_X^4 \]

\subsection{\texorpdfstring{$\pm 0$}{+/- 0} interactions}

For these interactions, we rewrite~: 
\begin{subequations}
\begin{align}
\eqref{est_L2x2_init} &= \int e^{i s \varphi} s^2 \nabla_{\xi} \varphi \mu_0^{k+1, 0, 1}(\xi, \eta) \widehat{f}(s, \xi - \eta) \widehat{u^0}(s, \eta) ~ d\eta \label{est_L2x2+0_term1} \\
&+ \int e^{i s \varphi} is \mu_0^{k, 0, 1}(\xi, \eta) \widehat{f}(s, \xi - \eta) \widehat{u^0}(s, \eta) ~ d\eta \label{est_L2x2+0_term2} \\
&+ \int e^{i s \varphi} i s \mu_0^{k+1, 0, 0}(\xi, \eta) \widehat{f}(s, \xi - \eta) \widehat{u^0}(s, \eta) ~ d\eta \label{est_L2x2+0_term3} \\
&+ \int e^{i s \varphi} i s \mu_0^{k+1, 0, 1}(\xi, \eta) \nabla_{\xi} \widehat{f}(s, \xi - \eta) \widehat{u^0}(s, \eta) ~ d\eta \label{est_L2x2+0_term4} \\
&+ \int e^{i s \varphi} \mu_0^{k-1, 0, 1}(\xi, \eta) \widehat{f}(s, \xi - \eta) \widehat{u^0}(s, \eta) ~ d\eta \label{est_L2x2+0_term5} \\
&+ \int e^{i s \varphi} \mu_0^{k+1, 0, -1}(\xi, \eta) \widehat{f}(s, \xi - \eta) \widehat{u^0}(s, \eta) ~ d\eta \label{est_L2x2+0_term6} \\
&+ \int e^{i s \varphi} \mu_0^{k, 0, 1}(\xi, \eta) \nabla_{\xi} \widehat{f}(s, \xi - \eta) \widehat{u^0}(s, \eta) ~ d\eta \label{est_L2x2+0_term7} \\
&+ \int e^{i s \varphi} \mu_0^{k+1, 0, 0}(\xi, \eta) \nabla_{\xi} \widehat{f}(s, \xi - \eta) \widehat{u^0}(s, \eta) ~ d\eta \label{est_L2x2+0_term8} \\
&+ \xi^{\alpha} \int e^{i s \varphi} \mu_0^{0, 0, 1}(\xi, \eta) \Delta_{\xi} \widehat{f}(s, \xi - \eta) \widehat{u^0}(s, \eta) ~ d\eta \label{est_L2x2+0_term9}
\end{align}
\end{subequations}
Recall that, concerning the last term \eqref{est_L2x2+0_term9}, we already treated the case when all the derivatives hit the first factor (through the symmetric structure). Therefore, we can replace it by : 
\[ \eqref{est_L2x2+0_term9}' = \int e^{i s \varphi} \mu_0^{k, 1, 1}(\xi, \eta) \Delta_{\xi} \widehat{f}(s, \xi - \eta) \widehat{u^0}(s, \eta) ~ d\eta \]
For all these terms except \eqref{est_L2x2+0_term1}, we can estimate directly with a decay $s^{-1+a+\delta+\frac{\gamma_k}{2}}$ by estimating $u^0$ in $L^{\infty-}$, thus winning $s^{-2}$ (and possibly absorbing a factor $s$) and the other in $L^2$. More precisely, we have that~: 
\begin{align*}
\Vert \eqref{est_L2x2+0_term3} \Vert_{L^2} &\lesssim s \Vert u \Vert_{H^6} \Vert u^0 \Vert_{W^{7, \infty-}} \lesssim s^{-1+a+\delta} \Vert u \Vert_X^2 \\
\Vert \eqref{est_L2x2+0_term4} \Vert_{L^2} &\lesssim s \Vert x f \Vert_{H^{k+2}} \Vert u^0 \Vert_{W^{k+2, \infty-}} \lesssim s^{-1+a+\delta+\frac{\gamma_k}{2}} \Vert u \Vert_X^2 \\
\Vert \eqref{est_L2x2+0_term5}_{k = 0} \Vert_{L^2} &\lesssim \Vert T_{\mu_0}(\Lambda u, u^0) \Vert_{L^{6/5}} \lesssim \Vert u \Vert_{H^1} \Vert u^0 \Vert_{L^3} \lesssim s^{-4/3+a} \Vert u \Vert_X^2 \\
\Vert \eqref{est_L2x2+0_term9} \Vert_{L^2} &\lesssim \Vert \Lambda |x|^2 f \Vert_{H^k} \Vert u^0 \Vert_{W^{k+2, \infty-}} \lesssim s^{-2+a+\delta+\gamma_k} \Vert u \Vert_X^2 \lesssim s^{-1} \Vert u \Vert_X^2
\end{align*}
The terms \eqref{est_L2x2+0_term2}, \eqref{est_L2x2+0_term6}, \eqref{est_L2x2+0_term7}, \eqref{est_L2x2+0_term8} are simpler or similar. 

Only \eqref{est_L2x2+0_term1} remains, that is
\[ \int e^{i s \varphi} s^2 \nabla_{\xi} \varphi \mu_0^{k+1, 0, 1}(\xi, \eta) \widehat{f}(s, \xi - \eta) \widehat{u^0}(s, \eta) ~ d\eta \]
As for the $L^{\infty}$ estimate, let us distinguish several cases, first by sketching the proof. 

\begin{Intuit} If the interaction is $++0$ or $--0$, then by lemma \ref{identite_fond_0} $|\xi| \nabla_{\xi} \varphi = |\eta| \mu_0$ and we absorb the singularity in $u^0$. We can then apply an integration by parts in frequency (since $\nabla_{\eta} \varphi$ does not vanish) to absorb a factor $s$, and finally get a total decay better than $s^{-1+\gamma_k}$. The reasoning is very similar to the one made for the $L^{\infty}$ norm. 

In the case $+-0$ or $-+0$, $\varphi(\xi, \eta) = |\xi|_0 + |\xi - \eta|_0$ so we can freely integrate by parts in time. As for the $L^{\infty}$ norm, we separate $u^0$ as a quadratic term. In the $+-+-$ case or analogous ones, we can integrate by parts in time, knowing that if it hits $D^{\alpha} |x|^2 f$ we can procede by symmetry, that if it hits $u^0$ it creates either terms controlled by lemma \ref{lem_decnonres} or a $\varphi'(\eta, \rho) = |\eta| \mu_0$, so that it absorbs the singularity in $u^0$. 

In the $+-++$ case or analogous ones, we can apply an angular repartition as in for the $L^{\infty}$ estimate.
\end{Intuit}

\subsubsection{Case $++0$ or $--0$} In this case, we have~:  
\[ \eqref{est_L2x2+0_term1} = \int e^{i s \varphi} s^2 \mu_0^{k, 1, 1}(\xi, \eta) \widehat{f}(s, \xi - \eta) \widehat{u^0}(s, \eta) ~ d\eta \]
Let us integrate by parts in $\eta$~: 
\begin{subequations}
\begin{align}
\eqref{est_L2x2+0_term1} &= \int e^{i s \varphi} s \mu_0^{k, 1, 0}(\xi, \eta) \widehat{f}(s, \xi - \eta) \widehat{u^0}(s, \eta) ~ d\eta \label{est_L2x2++0_term1-1} \\
&+ \int e^{i s \varphi} s \mu_0^{k, 0, 1}(\xi, \eta) \widehat{f}(s, \xi - \eta) \widehat{u^0}(s, \eta) ~ d\eta \label{est_L2x2++0_term1-2} \\
&+ \int e^{i s \varphi} s \mu_0^{k, 1, 1}(\xi, \eta) \nabla_{\eta} \widehat{f}(s, \xi - \eta) \widehat{u^0}(s, \eta) ~ d\eta \label{est_L2x2++0_term1-3} \\
&+ \int e^{i s \varphi} s \mu_0^{k, 1, 1}(\xi, \eta) \widehat{f}(s, \xi - \eta) \nabla_{\eta} \widehat{u^0}(s, \eta) ~ d\eta \label{est_L2x2++0_term1-4} 
\end{align} 
\end{subequations}
The first three terms \eqref{est_L2x2++0_term1-1}, \eqref{est_L2x2++0_term1-2}, \eqref{est_L2x2++0_term1-3} have already been estimated above. For the last one \eqref{est_L2x2++0_term1-4}~: 
\[ s \Vert T_{\mu_0}(\Lambda u, \Lambda x u^0) \Vert_{H^k} \lesssim s \Vert u \Vert_{W^{k+2, \infty-}} \Vert \Lambda x u^0 \Vert_{H^k} \lesssim s^{-1+\delta+a+\frac{\gamma_k}{2}} \Vert u \Vert_X^2 \]

\subsubsection{Case $+-+-$ and analogous} In this case, $\varphi(\xi, \eta) = \pm \left( |\xi|_0 + |\xi - \eta|_0 \right)$ does not vanish so we can integrate by parts in time. We then have~: 
\begin{subequations}
\begin{align}
&\eqref{est_L2x2+0_term1} = \int_1^t \int \int \int \xi^{\alpha} \Delta_{\xi} \widehat{f}(s, \xi) e^{i s \varphi(\xi, \eta)} e^{i s \varphi'(\eta, \rho)} s^2 \mu_0^{k+1, -1, 1}(\xi, \eta) \mu_0^{0, 0, 1}(\eta, \rho) \widehat{f}(s, \xi - \eta) \widehat{f}(s, \eta - \rho) \widehat{f}(s, \rho) ~ d\rho d\eta d\xi ds \nonumber \\
&= \int_1^t \int \int \int \xi^{\alpha} \Delta_{\xi} \partial_s \widehat{f}(s, \xi) e^{i s \varphi(\xi, \eta)} e^{i s \varphi'(\eta, \rho)} s^2 \mu_0^{k+1, -1, 1}(\xi, \eta) \mu_0^{0, 0, 1}(\eta, \rho) \frac{1}{\varphi(\xi, \eta)} \widehat{f}(s, \xi - \eta) \widehat{f}(s, \eta - \rho) \widehat{f}(s, \rho) ~ d\rho d\eta d\xi ds \label{est_L2x2+-+-_term1-1} \\
&+ \int_1^t \int \int \int \xi^{\alpha} \Delta_{\xi} \widehat{f}(s, \xi) e^{i s \varphi(\xi, \eta)} e^{i s \varphi'(\eta, \rho)} s^2 \mu_0^{k+1, 0, 1}(\xi, \eta) \mu_0^{-1, 0, 1}(\eta, \rho) \frac{\varphi'(\eta, \rho)}{\varphi(\xi, \eta)}  \widehat{f}(s, \xi - \eta) \widehat{f}(s, \eta - \rho) \widehat{f}(s, \rho) ~ d\rho d\eta d\xi ds \label{est_L2x2+-+-_term1-2} \\
&+ \int_1^t \int \int \int \xi^{\alpha} \Delta_{\xi} \widehat{f}(s, \xi) e^{i s \varphi(\xi, \eta)} e^{i s \varphi'(\eta, \rho)} s^2 \mu_0^{k+1, -1, 1}(\xi, \eta) \mu_0^{0, 0, 1}(\eta, \rho) \frac{1}{\varphi(\xi, \eta)}  \partial_s \widehat{f}(s, \xi - \eta) \widehat{f}(s, \eta - \rho) \widehat{f}(s, \rho) ~ d\rho d\eta d\xi ds \\
&+ \int_1^t \int \int \int \xi^{\alpha} \Delta_{\xi} \widehat{f}(s, \xi) e^{i s \varphi(\xi, \eta)} e^{i s \varphi'(\eta, \rho)} s^2 \mu_0^{k+1, -1, 1}(\xi, \eta) \mu_0^{0, 0, 1}(\eta, \rho) \frac{1}{\varphi(\xi, \eta)}  \widehat{f}(s, \xi - \eta) \partial_s \widehat{f}(s, \eta - \rho) \widehat{f}(s, \rho) ~ d\rho d\eta d\xi ds  \\
&+ \int_1^t \int \int \int \xi^{\alpha} \Delta_{\xi} \widehat{f}(s, \xi) e^{i s \varphi(\xi, \eta)} e^{i s \varphi'(\eta, \rho)} s^2 \mu_0^{k+1, -1, 1}(\xi, \eta) \mu_0^{0, 0, 1}(\eta, \rho) \frac{1}{\varphi(\xi, \eta)}  \widehat{f}(s, \xi - \eta) \widehat{f}(s, \eta - \rho) \partial_s \widehat{f}(s, \rho) ~ d\rho d\eta d\xi ds  \\
&+ \int \int \int \xi^{\alpha} \Delta_{\xi} \widehat{f}(t, \xi) e^{i t \varphi(\xi, \eta)} e^{i t \varphi'(\eta, \rho)} t^2 \mu_0^{k+1, -1, 1}(\xi, \eta) \mu_0^{0, 0, 1}(\eta, \rho) \frac{1}{\varphi(\xi, \eta)}  \widehat{f}(t, \xi - \eta) \widehat{f}(t, \eta - \rho) \widehat{f}(t, \rho) ~ d\rho d\eta d\xi \\
&+ \int \int \int \xi^{\alpha} \Delta_{\xi} \widehat{f}(1, \xi) e^{i \varphi(\xi, \eta)} e^{i \varphi'(\eta, \rho)} \mu_0^{k+1, -1, 1}(\xi, \eta) \mu_0^{0, 0, 1}(\eta, \rho) \frac{1}{\varphi(\xi, \eta)}  \widehat{f}(1, \xi - \eta) \widehat{f}(1, \eta - \rho) \widehat{f}(1, \rho) ~ d\rho d\eta d\xi  
\end{align}
\end{subequations}
Yet $\frac{\mu_0^{k+1, m, l}(\xi, \eta)}{\varphi(\xi, \eta)} = \mu_0^{k, m, l}(\xi, \eta)$, and $\varphi'(\eta, \rho) \mu_0^{-1, 0, 1}(\eta, \rho) = \mu_0^{0, 0, 1}(\eta, \rho)$ for this choice of signs $\epsilon_i$. The second line \eqref{est_L2x2+-+-_term1-2} is therefore identical to the $++0$ case treated above. 

The first line \eqref{est_L2x2+-+-_term1-1} can be treated by the symmetry argument already used~: namely, we develop $\xi^{\alpha} \Delta_{\xi} \partial_s \widehat{f}(s, \xi)$ as a bilinear product and apply similar estimates as before, noting that the other side of \eqref{est_L2x2+-+-_term1-1}, that is $s^2 T_{\mu_0^{k, -1, 1}}(u, T_{\mu_0^{0, 0, 1}}(u, u))$ has only $s^{\varepsilon+\delta}$ growth in $L^2$, so behaves better than $\Lambda^{k+1} |x|^2 f$. The only term we need to use more structure is when we use the non-resonant structure to apply integrations by parts in time, in which case the argument is exactly similar as the one in section \ref{section_L2x2_symmetry} ; note also that the quasi-linear structure is not needed here since $s^2 T_{\mu_0^{k, -1, 1}}(u, T_{\mu_0^{0, 0, 1}}(u, u))$ can be estimated the same way with one more derivative. 

Finally, the other lines can be treated by lemma \ref{lem_decnonres} which ensures an additional decay on the terms that bear a time derivative. The initial time or final time terms are also well controlled thanks to the absence of a time integration. 

\subsubsection{Case $+-++$ and analogous} In this case, we reuse the angular partition already used in the treatment of the $L^{\infty}$ norm, between one area where $\nabla_{\rho} \varphi'(\eta, \rho)$ does not vanish and one area where $\frac{|\eta - \rho|}{|\rho|}$ acts like a symbol of order $0$. 

In the first area, we write~: 
\begin{subequations} \label{est_L2x2+-++_area1_init}
\begin{align}
&\xi^{\alpha} \int \int e^{i s \varphi(\xi, \eta)} e^{i s \varphi'(\eta, \rho)} s^2 \mu_0^{0, -1, 1}(\xi, \eta) \mu_0^{0, 0, 1}(\eta, \rho) \chi(\eta, \rho) \widehat{f}(s, \xi - \eta) \widehat{f}(s, \eta - \rho) \widehat{f}(s, \rho) ~ d\rho d\eta \tag{\ref{est_L2x2+-++_area1_init}} \\
&= \int \int e^{i s \varphi(\xi, \eta)} e^{i s \varphi'(\eta, \rho)} \mu_0^{k+1, -1, 1}(\xi, \eta) \mu_0^{0, 0, -1}(\eta, \rho) \chi(\eta, \rho) \widehat{f}(s, \xi - \eta) \widehat{f}(s, \eta - \rho) \widehat{f}(s, \rho) ~ d\rho d\eta \label{est_L2x2+-++_area1_term1} \\
&+ \int \int e^{i s \varphi(\xi, \eta)} e^{i s \varphi'(\eta, \rho)} \mu_0^{k+1, -1, 1}(\xi, \eta) \mu_0(\eta, \rho) \chi(\eta, \rho) \widehat{f}(s, \xi - \eta) \nabla_{\rho} \widehat{f}(s, \eta - \rho) \widehat{f}(s, \rho) ~ d\rho d\eta \label{est_L2x2+-++_area1_term2} \\
&+ \int \int e^{i s \varphi(\xi, \eta)} e^{i s \varphi'(\eta, \rho)} \mu_0^{k+1, 0, 1}(\xi, \eta) \mu_0^{-1, -2, 1}(\eta, \rho) \chi(\eta, \rho) \widehat{f}(s, \xi - \eta) \widehat{f}(s, \eta - \rho) \widehat{f}(s, \rho) ~ d\rho d\eta \label{est_L2x2+-++_area1_term3} \\
&+ \int \int e^{i s \varphi(\xi, \eta)} e^{i s \varphi'(\eta, \rho)} \mu_0^{k+1, -1, 1}(\xi, \eta) \mu_0^{0, -1, 1}(\eta, \rho) \chi(\eta, \rho) \widehat{f}(s, \xi - \eta) \nabla_{\rho} \widehat{f}(s, \eta - \rho) \widehat{f}(s, \rho) ~ d\rho d\eta \label{est_L2x2+-++_area1_term4} \\
&+ \int \int e^{i s \varphi(\xi, \eta)} e^{i s \varphi'(\eta, \rho)} \mu_0^{k+1, 0, 1}(\xi, \eta) \mu_0^{-1, -1, 1}(\eta, \rho) \chi(\eta, \rho) \widehat{f}(s, \xi - \eta) \widehat{f}(s, \eta - \rho) \nabla_{\rho} \widehat{f}(s, \rho) ~ d\rho d\eta \label{est_L2x2+-++_area1_term5} \\
&+ \int \int e^{i s \varphi(\xi, \eta)} e^{i s \varphi'(\eta, \rho)} \mu_0^{k+1, -1, 1}(\xi, \eta) \mu_0^{0, 0, 1}(\eta, \rho) \chi(\eta, \rho) \widehat{f}(s, \xi - \eta) \nabla_{\rho}^2 \widehat{f}(s, \eta - \rho) \widehat{f}(s, \rho) ~ d\rho d\eta \label{est_L2x2+-++_area1_term6} \\
&+ \int \int e^{i s \varphi(\xi, \eta)} e^{i s \varphi'(\eta, \rho)} \mu_0^{k+1, -1, 1}(\xi, \eta) \mu_0^{0, 0, 1}(\eta, \rho) \chi(\eta, \rho) \widehat{f}(s, \xi - \eta) \nabla_{\rho} \widehat{f}(s, \eta - \rho) \nabla_{\rho} \widehat{f}(s, \rho) ~ d\rho d\eta \label{est_L2x2+-++_area1_term7} \\
&+ \int \int e^{i s \varphi(\xi, \eta)} e^{i s \varphi'(\eta, \rho)} \mu_0^{k+1, 0, 1}(\xi, \eta) \mu_0^{-1, 0, 1}(\eta, \rho) \chi(\eta, \rho) \widehat{f}(s, \xi - \eta) \widehat{f}(s, \eta - \rho) \nabla^2_{\rho} \widehat{f}(s, \rho) ~ d\rho d\eta \label{est_L2x2+-++_area1_term8}
\end{align}
\end{subequations}
In \eqref{est_L2x2+-++_area1_term3}, \eqref{est_L2x2+-++_area1_term3}, \eqref{est_L2x2+-++_area1_term8}, we distribute $\mu_0^{-1, m, 1}(\eta, \rho) = \mu_0^{0, m, 0}(\eta, \rho) + \mu_0^{-1, m+1, 0}(\eta, \rho)$ to recover~: 
\begin{align*}
\eqref{est_L2x2+-++_area1_init}
&= \int \int e^{i s \varphi(\xi, \eta)} e^{i s \varphi'(\eta, \rho)} \mu_0^{k+1, -1, 1}(\xi, \eta) \mu_0^{0, 0, -1}(\eta, \rho) \chi(\eta, \rho) \widehat{f}(s, \xi - \eta) \widehat{f}(s, \eta - \rho) \widehat{f}(s, \rho) ~ d\rho d\eta \\
&+ \int \int e^{i s \varphi(\xi, \eta)} e^{i s \varphi'(\eta, \rho)} \mu_0^{k+1, -1, 1}(\xi, \eta) \mu_0(\eta, \rho) \chi(\eta, \rho) \widehat{f}(s, \xi - \eta) \nabla_{\rho} \widehat{f}(s, \eta - \rho) \widehat{f}(s, \rho) ~ d\rho d\eta \\
&+ \int \int e^{i s \varphi(\xi, \eta)} e^{i s \varphi'(\eta, \rho)} \mu_0^{k+1, 0, 1}(\xi, \eta) \mu_0^{0, -2, 0}(\eta, \rho) \chi(\eta, \rho) \widehat{f}(s, \xi - \eta) \widehat{f}(s, \eta - \rho) \widehat{f}(s, \rho) ~ d\rho d\eta \\
&+ \int \int e^{i s \varphi(\xi, \eta)} e^{i s \varphi'(\eta, \rho)} \mu_0^{k+1, -1, 1}(\xi, \eta) \mu_0^{0, -1, 1}(\eta, \rho) \chi(\eta, \rho) \widehat{f}(s, \xi - \eta) \nabla_{\rho} \widehat{f}(s, \eta - \rho) \widehat{f}(s, \rho) ~ d\rho d\eta \\
&+ \int \int e^{i s \varphi(\xi, \eta)} e^{i s \varphi'(\eta, \rho)} \mu_0^{k+1, 0, 1}(\xi, \eta) \mu_0^{0, -1, 0}(\eta, \rho) \chi(\eta, \rho) \widehat{f}(s, \xi - \eta) \widehat{f}(s, \eta - \rho) \nabla_{\rho} \widehat{f}(s, \rho) ~ d\rho d\eta \\
&+ \int \int e^{i s \varphi(\xi, \eta)} e^{i s \varphi'(\eta, \rho)} \mu_0^{k+1, -1, 1}(\xi, \eta) \mu_0^{0, 0, 1}(\eta, \rho) \chi(\eta, \rho) \widehat{f}(s, \xi - \eta) \nabla_{\rho}^2 \widehat{f}(s, \eta - \rho) \widehat{f}(s, \rho) ~ d\rho d\eta \\
&+ \int \int e^{i s \varphi(\xi, \eta)} e^{i s \varphi'(\eta, \rho)} \mu_0^{k+1, -1, 1}(\xi, \eta) \mu_0^{0, 0, 1}(\eta, \rho) \chi(\eta, \rho) \widehat{f}(s, \xi - \eta) \nabla_{\rho} \widehat{f}(s, \eta - \rho) \nabla_{\rho} \widehat{f}(s, \rho) ~ d\rho d\eta \\
&+ \int \int e^{i s \varphi(\xi, \eta)} e^{i s \varphi'(\eta, \rho)} \mu_0^{k+1, 0, 1}(\xi, \eta) \mu_0(\eta, \rho) \chi(\eta, \rho) \widehat{f}(s, \xi - \eta) \widehat{f}(s, \eta - \rho) \nabla^2_{\rho} \widehat{f}(s, \rho) ~ d\rho d\eta
\end{align*}
Now we can control everything. The terms above have the form $D^{\alpha} T_{\mu_0}(\Lambda u, \Lambda^{-1} A)$ or $D^{\alpha} T_{\mu_0}(\Lambda u, B)$ and we can estimate them by~: 
\begin{align*}
&\Vert T_{\mu_0}(\Lambda u, \Lambda^{-1} A) \Vert_{H^{k+1}} \lesssim \Vert u \Vert_{W^{k+2, 3}} \Vert \Lambda^{-1} A \Vert_{L^6} + \Vert u \Vert_{W^{k+2, \infty-}} \Vert A \Vert_{H^k} \lesssim s^{-1/3+\delta} \Vert u \Vert_X \Vert A \Vert_{L^2} + s^{-1+\delta} \Vert u \Vert_X \Vert A \Vert_{H^k} \\
&\Vert T_{\mu_0}(\Lambda u, B) \Vert_{H^{k+1}} \lesssim s^{-1+\delta} \Vert u \Vert_X \Vert B \Vert_{H^{k+1}} 
\end{align*}
Then, for $A$, we have the following possibilities~:
\[ T_{\mu_0}(\Lambda^{-1} u, u), ~~~ T_{\mu_0}(e^{\pm i s \Lambda} x f, u), ~~~ T_{\mu_0}(e^{\pm i s \Lambda} x f, \Lambda^{-1} u), ~~~ T_{\mu_0}(\Lambda e^{\pm i s \Lambda} x^2 f, u), ~~~ T_{\mu_0}(\Lambda e^{\pm i s \Lambda} x f, e^{\pm i s \Lambda} x f) \]
and for $B$~: 
\[ T_{\mu_0}(u, \Lambda^{-2} u), ~~~ T_{\mu_0}(u, \Lambda^{-1} e^{\pm i s \Lambda} x f), ~~~ T_{\mu_0}(u, e^{\pm i s \Lambda} x^2 f) \]
Let us estimate~: 
\begin{align*}
\Vert T_{\mu_0}(\Lambda^{-1} u, u) \Vert_{H^k} &\lesssim \Vert x f \Vert_{H^k} \Vert u \Vert_{W^{k+1, \infty-}} \lesssim s^{-1+\delta} \Vert u \Vert_X^2 \\
\Vert T_{\mu_0}(e^{\pm i s \Lambda} x f, u) \Vert_{H^k} &\lesssim \Vert x f \Vert_{H^k} \Vert u \Vert_{W^{k+1, \infty-}} \lesssim s^{-1+\delta} \Vert u \Vert_X^2 \\
\Vert T_{\mu_0}(e^{\pm i s \Lambda} x f, \Lambda^{-1} u) \Vert_{L^2} &\lesssim \Vert e^{\pm i s \Lambda} x f \Vert_{L^4} \Vert \Lambda^{-1} e^{\pm i s \Lambda} f \Vert_{L^4} \lesssim s^{-1} \Vert \Lambda x f \Vert_{L^{4/3}} \Vert f \Vert_{L^{4/3}} \\
&\lesssim s^{-1} \Vert \langle x \rangle \Lambda x f \Vert_{L^2} \Vert \langle x \rangle f \Vert_{L^2} \lesssim s^{-1+\gamma_0} \Vert u \Vert_X^2 \\
\Vert T_{\mu_0}(e^{\pm i s \Lambda} x f, \Lambda^{-1} u) \Vert_{H^k} &\lesssim \Vert x f \Vert_{H^k} \Vert \Lambda^{-1} u \Vert_{W^{k+1, 6}} \lesssim \Vert x f \Vert_{H^k} \Vert u \Vert_{H^{k+1}} \lesssim s^{\varepsilon} \Vert u \Vert_X^2 \\
\Vert T_{\mu_0}(\Lambda e^{\pm i s \Lambda} x^2 f, u) \Vert_{H^k} &\lesssim \Vert \Lambda x^2 f \Vert_{H^k} \Vert u \Vert_{W^{k+1, \infty-}} \lesssim s^{-1+\delta+\gamma_k} \Vert u \Vert_X^2 \\
\Vert T_{\mu_0}(\Lambda e^{\pm i s \Lambda} x f, e^{\pm i s \Lambda} x f) \Vert_{L^2} &\lesssim \Vert \Lambda e^{\pm i s \Lambda} x f \Vert_{L^4} \Vert e^{\pm i s \Lambda} x f \Vert_{L^4} \lesssim s^{-1} \Vert \Lambda x f \Vert_{W^{1, 4/3}}^2 \lesssim s^{-1} \Vert \langle x \rangle \Lambda x f \Vert_{H^1} \\
&\lesssim s^{-1+\gamma_1} \Vert u \Vert_X^2 \\
\Vert T_{\mu_0}(\Lambda e^{\pm i s \Lambda} x f, e^{\pm i s \Lambda} x f) \Vert_{H^k} &\lesssim \Vert x f \Vert_{H^{k+1}} \Vert e^{\pm i s \Lambda} x f \Vert_{W^{2, 6}} \lesssim \Vert x f \Vert_{H^{k+1}} \Vert x f \Vert_{H^3} \lesssim s^b \Vert u \Vert_X^2 \\
\Vert T_{\mu_0}(u, \Lambda^{-2} u) \Vert_{H^{k+1}} &\lesssim \Vert u \Vert_{W^{k+1, 3}} \Vert \Lambda^{-2} u \Vert_{W^{k+1, 6}} \lesssim s^{-1/3+\delta+b} \Vert u \Vert_X^2 \\
\Vert T_{\mu_0}(u, \Lambda^{-1} e^{\pm i s \Lambda} x f) \Vert_{H^{k+1}} &\lesssim \Vert u \Vert_{W^{k+1, 3}} \Vert \Lambda^{-1} e^{\pm i s \Lambda} x f \Vert_{W^{k+1, 6}} \lesssim s^{-1/3+\delta+b} \Vert u \Vert_X^2 \\
\Vert T_{\mu_0}(u, e^{\pm i s \Lambda} x^2 f) \Vert_{H^{k+1}} &\lesssim \Vert u \Vert_{W^{k+1, 3}} \Vert e^{\pm i s \Lambda} x^2 f \Vert_{L^6} + \Vert u \Vert_{W^{k+1, \infty-}} \Vert \Lambda x^2 f \Vert_{H^k} \\
&\lesssim s^{-1/3+\delta+\gamma_0} \Vert u \Vert_X^2 + s^{-1+\delta+\gamma_k} \Vert u \Vert_X^2
\end{align*}
Summing up, we get~: 
\[ \lesssim \left( s^{-4/3+2\delta+\gamma_k+b} + s^{-1+\varepsilon+\delta+b} \right) \Vert u \Vert_X^3 \lesssim s^{-1+\gamma_k} \Vert u \Vert_X^3 \]

In the second area, we get~:
\begin{equation} \int \int e^{i s \varphi(\xi, \eta)} e^{i s \varphi'(\eta, \rho)} s^2 \mu_0^{k+1, 0, 1}(\xi, \eta) \mu_0(\eta, \rho) \widehat{f}(s, \xi - \eta) \widehat{f}(s, \eta - \rho) \widehat{f}(s, \rho) ~ d\rho d\eta \label{est_L2x2+-++_area2_init} \end{equation}
Let us integrate by parts in $\rho$, using the structure of the nonlinearity. We obtain~: 
\[ \begin{aligned}
\eqref{est_L2x2+-++_area2_init} &= \int \int e^{i s \varphi(\xi, \eta)} e^{i s \varphi'(\eta, \rho)} s \mu_0^{k+1, 0, 1}(\xi, \eta) \mu_0^{0, 0, -1}(\eta, \rho) \chi(\eta, \rho) \widehat{f}(s, \xi - \eta) \widehat{f}(s, \eta - \rho) \widehat{f}(s, \rho) ~ d\rho d\eta \\
&+ \int \int e^{i s \varphi(\xi, \eta)} e^{i s \varphi'(\eta, \rho)} s \mu_0^{k+1, 0, 1}(\xi, \eta) \mu_0(\eta, \rho) \chi(\eta, \rho) \widehat{f}(s, \xi - \eta) \nabla_{\rho} \widehat{f}(s, \eta - \rho) \widehat{f}(s, \rho) ~ d\rho d\eta
\end{aligned} \]
Let us write these objects as~: 
\[ D^{\alpha} s T_{\mu_0}(\Lambda u, T_{\mu_0}(\Lambda^{-1} u, u)) + D^{\alpha} s T_{\mu_0}(\Lambda u, T_{\mu_0}(e^{\pm i s \Lambda} x f, u)) \]
Now we estimate~: 
\begin{align*}
&\Vert s T_{\mu_0}(\Lambda u, T_{\mu_0}(\Lambda^{-1} u, u)) \Vert_{H^{k+1}} \lesssim s \Vert u \Vert_{W^{k+2, \infty-}} \Vert x f \Vert_{H^{k+1}} \Vert u \Vert_{W^{k+2, \infty-}} \lesssim s^{-1+2\delta+\frac{\gamma_k+\varepsilon}{2}} \Vert u \Vert_X^3 \lesssim s^{-1+\gamma_k} \Vert u \Vert_X^3 \\
&\Vert s T_{\mu_0}(\Lambda u, T_{\mu_0}(e^{\pm i s \Lambda} x f, u)) \Vert_{H^{k+1}} \lesssim s \Vert u \Vert_{W^{k+2, \infty-}} \Vert x f \Vert_{H^{k+1}} \Vert u \Vert_{W^{k+2, \infty-}} \lesssim s^{-1+2\delta+\frac{\gamma_k+\varepsilon}{2}} \Vert u \Vert_X^3 \lesssim s^{-1+\gamma_k} \Vert u \Vert_X^3
\end{align*}

\subsubsection{Case $+-+0$ and analogous} We have~: 
\begin{equation} \int \int e^{i s \varphi(\xi, \eta)} e^{i s \varphi'(\eta, \rho)} s^2 \mu_0^{k+1, -1, 1}(\xi, \eta) \mu_0^{0, 0, 1}(\eta, \rho) \widehat{f}(s, \xi - \eta) \widehat{f}(s, \eta - \rho) \widehat{u^0}(s, \rho) ~ d\rho d\eta \label{est_L2x2+-+0_init} \end{equation}
Let us integrate by parts in $\rho$ (this time $\nabla_{\rho} \varphi'(\eta, \rho)$ does not vanish) and distribute some derivatives to get~: 
\begin{align*}
\eqref{est_L2x2+-+0_init} &= \int \int e^{i s \varphi(\xi, \eta)} e^{i s \varphi'(\eta, \rho)} s \mu_0^{k+1, -1, 1}(\xi, \eta) \mu_0(\eta, \rho) \widehat{f}(s, \xi - \eta) \widehat{f}(s, \eta - \rho) \widehat{u^0}(s, \rho) ~ d\rho d\eta \\
&+ \int \int e^{i s \varphi(\xi, \eta)} e^{i s \varphi'(\eta, \rho)} s \mu_0^{k+1, -1, 1}(\xi, \eta) \mu_0^{0, -1, 1}(\eta, \rho) \widehat{f}(s, \xi - \eta) \widehat{f}(s, \eta - \rho) \widehat{u^0}(s, \rho) ~ d\rho d\eta \\
&+ \int \int e^{i s \varphi(\xi, \eta)} e^{i s \varphi'(\eta, \rho)} s \mu_0^{k+1, -1, 1}(\xi, \eta) \mu_0^{0, 0, 1}(\eta, \rho) \widehat{f}(s, \xi - \eta) \nabla_{\rho} \widehat{f}(s, \eta - \rho) \widehat{u^0}(s, \rho) ~ d\rho d\eta \\
&+ \int \int e^{i s \varphi(\xi, \eta)} e^{i s \varphi'(\eta, \rho)} s \mu_0^{k+1, -1, 1}(\xi, \eta) \mu_0^{0, 0, 1}(\eta, \rho) \widehat{f}(s, \xi - \eta) \widehat{f}(s, \eta - \rho) \nabla_{\rho} \widehat{u^0}(s, \rho) ~ d\rho d\eta \\
&= \int \int e^{i s \varphi(\xi, \eta)} e^{i s \varphi'(\eta, \rho)} s \mu_0^{k+1, -1, 1}(\xi, \eta) \mu_0(\eta, \rho) \widehat{f}(s, \xi - \eta) \widehat{f}(s, \eta - \rho) \widehat{u^0}(s, \rho) ~ d\rho d\eta \\
&+ \int \int e^{i s \varphi(\xi, \eta)} e^{i s \varphi'(\eta, \rho)} s \mu_0^{k+1, 0, 1}(\xi, \eta) \mu_0^{0, -1, 0}(\eta, \rho) \widehat{f}(s, \xi - \eta) \widehat{f}(s, \eta - \rho) \widehat{u^0}(s, \rho) ~ d\rho d\eta \\
&+ \int \int e^{i s \varphi(\xi, \eta)} e^{i s \varphi'(\eta, \rho)} s \mu_0^{k+1, -1, 1}(\xi, \eta) \mu_0^{0, 0, 1}(\eta, \rho) \widehat{f}(s, \xi - \eta) \nabla_{\rho} \widehat{f}(s, \eta - \rho) \widehat{u^0}(s, \rho) ~ d\rho d\eta \\
&+ \int \int e^{i s \varphi(\xi, \eta)} e^{i s \varphi'(\eta, \rho)} s \mu_0^{k+1, 0, 1}(\xi, \eta) \mu_0(\eta, \rho) \widehat{f}(s, \xi - \eta) \widehat{f}(s, \eta - \rho) \nabla_{\rho} \widehat{u^0}(s, \rho) ~ d\rho d\eta \\
&+ \int \int e^{i s \varphi(\xi, \eta)} e^{i s \varphi'(\eta, \rho)} s \mu_0^{k+1, -1, 1}(\xi, \eta) \mu_0^{0, 1, 0}(\eta, \rho) \widehat{f}(s, \xi - \eta) \widehat{f}(s, \eta - \rho) \nabla_{\rho} \widehat{u^0}(s, \rho) ~ d\rho d\eta
\end{align*}
We estimate by~: 
\begin{align*}
&\Vert s T_{\mu_0}(\Lambda u, \Lambda^{-1} T_{\mu_0}(u, u^0)) \Vert_{H^{k+1}} \lesssim s \Vert u \Vert_{W^{k+2, 3}} \Vert u \Vert_{H^{k+1}} \Vert u^0 \Vert_{W^{k+2, \infty-}} \lesssim s^{-4/3+2\delta+a} \Vert u \Vert_X^3 \\
&\Vert s T_{\mu_0}(\Lambda u, T_{\mu_0}(u, \Lambda^{-1} u^0)) \Vert_{H^{k+1}} \lesssim s \Vert u \Vert_{W^{k+3, \infty-}} \Vert u \Vert_{W^{k+1, 3}} \Vert u^0 \Vert_{H^{k+1}} \lesssim s^{-4/3+2\delta+a} \Vert u \Vert_X^3 \\
&\Vert s T_{\mu_0}(\Lambda u, \Lambda^{-1} T_{\mu_0}(\Lambda e^{\pm i s \Lambda} x f, u^0)) \Vert_{H^{k+1}} \lesssim s \Vert u \Vert_{W^{k+2, 3}} \Vert x f \Vert_{H^{k+2}} \Vert u^0 \Vert_{W^{k+2, \infty-}} \lesssim s^{-4/3+\frac{\gamma_5+\varepsilon}{2} + 2\delta+a} \Vert u \Vert_X^3 \\
&\Vert s T_{\mu_0}(\Lambda u, T_{\mu_0}(u, x u^0)) \Vert_{H^{k+1}} \lesssim s \Vert u \Vert_{W^{k+3, \infty-}} \Vert u \Vert_{W^{k+1, 3}} \Vert \Lambda x u^0 \Vert_{H^{k+1}} \lesssim s^{-4/3+2\delta+b+a} \Vert u \Vert_X^3 \\
&\Vert s T_{\mu_0}(\Lambda u, \Lambda^{-1} T_{\mu_0}(u, \Lambda x u^0)) \Vert_{H^{k+1}} \lesssim s \Vert u \Vert_{W^{k+2, 3}} \Vert u \Vert_{W^{k+2, \infty-}} \Vert \Lambda x u^0 \Vert_{H^{k+1}} \lesssim s^{-4/3+a+b+2\delta} \Vert u \Vert_X^3
\end{align*}
Summing up, we control everytging with a decay of $s^{-1+\gamma_k}$. 

\subsubsection{Case $+-0+$ and analogous} We have~: 
\begin{equation} \int \int e^{i s \varphi(\xi, \eta)} e^{i s \varphi'(\eta, \rho)} s^2 \mu_0^{k+1, -1, 1}(\xi, \eta) \mu_0^{0, 0, 1}(\eta, \rho) \widehat{f}(s, \xi - \eta) \widehat{u^0}(s, \eta - \rho) \widehat{f}(s, \rho) ~ d\rho d\eta \label{est_L2x2+-0+_init} \end{equation}
Let us integrate by parts in $\rho$~: 
\begin{subequations}
\begin{align}
\eqref{est_L2x2+-0+_init} &= \int \int e^{i s \varphi(\xi, \eta)} e^{i s \varphi'(\eta, \rho)} s \mu_0^{k+1, -1, 1}(\xi, \eta) \mu_0(\eta, \rho) \widehat{f}(s, \xi - \eta) \widehat{u^0}(s, \eta - \rho) \widehat{f}(s, \rho) ~ d\rho d\eta \label{est_L2x2+-0+_term1} \\
&+ \int \int e^{i s \varphi(\xi, \eta)} e^{i s \varphi'(\eta, \rho)} s \mu_0^{k+1, -1, 1}(\xi, \eta) \mu_0^{0, -1, 1}(\eta, \rho) \widehat{f}(s, \xi - \eta) \widehat{u^0}(s, \eta - \rho) \widehat{f}(s, \rho) ~ d\rho d\eta \\
&+ \int \int e^{i s \varphi(\xi, \eta)} e^{i s \varphi'(\eta, \rho)} s \mu_0^{k+1, -1, 1}(\xi, \eta) \mu_0^{0, 0, 1}(\eta, \rho) \widehat{f}(s, \xi - \eta) \nabla_{\rho}\widehat{u^0}(s, \eta - \rho) \widehat{f}(s, \rho) ~ d\rho d\eta \label{est_L2x2+-0+_term3} \\
&+ \int \int e^{i s \varphi(\xi, \eta)} e^{i s \varphi'(\eta, \rho)} s \mu_0^{k+1, -1, 1}(\xi, \eta) \mu_0^{0, 0, 1}(\eta, \rho) \widehat{f}(s, \xi - \eta) \widehat{u^0}(s, \eta - \rho) \nabla_{\rho} \widehat{f}(s, \rho) ~ d\rho d\eta
\end{align}
\end{subequations}
\eqref{est_L2x2+-0+_term1} and \eqref{est_L2x2+-0+_term3} have already been estimated in $+-+0$. For the other ones~: 
\begin{align*}
&\Vert s T_{\mu_0}(\Lambda u, \Lambda^{-1} T_{\mu_0}(\Lambda u^0, \Lambda^{-1} u)) \Vert_{H^{k+1}} \lesssim s \Vert u \Vert_{W^{k+2, 3}} \Vert u^0 \Vert_{W^{k+3, \infty-}} \Vert x f \Vert_{H^{k+1}} \lesssim s^{-4/3+2\delta+a+b} \Vert u \Vert_X^3 \\
&\Vert s T_{\mu_0}(\Lambda u, \Lambda^{-1} T_{\mu_0}(\Lambda u^0, e^{\pm i s \Lambda} x f) \Vert_{H^{k+1}} \lesssim s \Vert u \Vert_{W^{k+2, 3}} \Vert u^0 \Vert_{W^{k+3, \infty-}} \Vert x f \Vert_{H^{k+1}} \lesssim s^{-4/3+2\delta+a+b} \Vert u \Vert_X^3
\end{align*}

\subsubsection{Case $+-00$ and analogous} We have~: 
\[ \int \int e^{i s \varphi(\xi, \eta)} e^{i s \varphi'(\eta, \rho)} s^2 \mu_0^{k+1, -1, 1}(\xi, \eta) \mu_0^{0, 0, 1}(\eta, \rho) \widehat{f}(s, \xi - \eta) \widehat{u^0}(s, \eta - \rho) \widehat{u^0}(s, \rho) ~ d\rho d\eta \]
We estimate it directly~: 
\[ \Vert s^2 T_{\mu_0}(\Lambda u, \Lambda^{-1} T_{\mu_0}(\Lambda u^0, u^0)) \Vert_{H^{k+1}} \lesssim s^2 \Vert u \Vert_{W^{k+2, 3}} \Vert u^0 \Vert_{H^{k+2}} \Vert u^0 \Vert_{W^{k+2, \infty-}} \lesssim s^{-4/3+2\delta+2a} \Vert u \Vert_X^3 \lesssim s^{-1+\gamma_k} \Vert u \Vert_X^3 \]

\subsection{\texorpdfstring{$0 \pm$}{0 +/-} interactions}

We rewrite the corresponding part of $\xi^{\alpha} \Delta_{\xi} \widehat{f}(s, \xi)$ as~: 
\begin{subequations}
\begin{align}
\eqref{est_L2x2_init} &= \int e^{i s \varphi} s^2 \mu_0^{k+1, 0, 1}(\xi, \eta) \widehat{u^0}(s, \xi - \eta) \widehat{f}(s, \eta) ~ d\eta \label{est_L2x20+_line1} \\
&+ \int e^{i s \varphi} s \mu_0^{k, 0, 1}(\xi, \eta) \widehat{u^0}(s, \xi - \eta) \widehat{f}(s, \eta) ~ d\eta \label{est_L2x20+_line2} \\
&+ \int e^{i s \varphi} s \mu_0^{k+1, 0, 0}(\xi, \eta) \widehat{u^0}(s, \xi - \eta) \widehat{f}(s, \eta) ~ d\eta \label{est_L2x20+_line3} \\
&+ \int e^{i s \varphi} s \mu_0^{k+1, 0, 1}(\xi, \eta) \nabla_{\xi} \widehat{u^0}(s, \xi - \eta) \widehat{f}(s, \eta) ~ d\eta \label{est_L2x20+_line4} \\
&+ \int e^{i s \varphi} \mu_0^{k-1, 0, 1}(\xi, \eta) \widehat{u^0}(s, \xi - \eta) \widehat{f}(s, \eta) ~ d\eta \label{est_L2x20+_line5} \\
&+ \int e^{i s \varphi} \mu_0^{k, 0, 1}(\xi, \eta) \nabla_{\xi} \widehat{u^0}(s, \xi - \eta) \widehat{f}(s, \eta) ~ d\eta \label{est_L2x20+_line6} \\
&+ \int e^{i s \varphi} \mu_0^{k+1, 0, -1}(\xi, \eta) \widehat{u^0}(s, \xi - \eta) \widehat{f}(s, \eta) ~ d\eta \label{est_L2x20+_line7} \\
&+ \int e^{i s \varphi} \mu_0^{k+1, 0, 0}(\xi, \eta) \nabla_{\xi} \widehat{u^0}(s, \xi - \eta) \widehat{f}(s, \eta) ~ d\eta \label{est_L2x20+_line8} \\
&+ \int e^{i s \varphi} \mu_0^{k+1, 0, 1}(\xi, \eta) \Delta_{\xi} \widehat{u^0}(s, \xi - \eta) \widehat{f}(s, \eta) ~ d\eta \label{est_L2x20+_line9}
\end{align}
\end{subequations}
On the 1st and last term \eqref{est_L2x20+_line1} and \eqref{est_L2x20+_line9}, we integrate by parts in frequency (note that $\nabla_{\eta} \varphi$ does not vanish) and get~: 
\begin{subequations}
\begin{align}
\eqref{est_L2x2_init} &= \int e^{i s \varphi} s \mu_0^{k+1, 0, 1}(\xi, \eta) \widehat{u^0}(s, \xi - \eta) \nabla_{\eta} \widehat{f}(s, \eta) ~ d\eta \label{est_L2x20+_term1} \\
&+ \int e^{i s \varphi} s \mu_0^{k+1, -1, 1}(\xi, \eta) \widehat{u^0}(s, \xi - \eta) \widehat{f}(s, \eta) ~ d\eta \label{est_L2x20+_term2} \\
&+ \int e^{i s \varphi} s \mu_0^{k, 0, 1}(\xi, \eta) \widehat{u^0}(s, \xi - \eta) \widehat{f}(s, \eta) ~ d\eta \label{est_L2x20+_term3} \\
&+ \int e^{i s \varphi} s \mu_0^{k+1, 0, 0}(\xi, \eta) \widehat{u^0}(s, \xi - \eta) \widehat{f}(s, \eta) ~ d\eta \label{est_L2x20+_term4} \\
&+ \int e^{i s \varphi} s \mu_0^{k+1, 0, 1}(\xi, \eta) \nabla_{\xi} \widehat{u^0}(s, \xi - \eta) \widehat{f}(s, \eta) ~ d\eta \label{est_L2x20+_term5} \\
&+ \int e^{i s \varphi} \mu_0^{k-1, 0, 1}(\xi, \eta) \widehat{u^0}(s, \xi - \eta) \widehat{f}(s, \eta) ~ d\eta \label{est_L2x20+_term6} \\
&+ \int e^{i s \varphi} \mu_0^{k, 0, 1}(\xi, \eta) \nabla_{\xi} \widehat{u^0}(s, \xi - \eta) \widehat{f}(s, \eta) ~ d\eta \label{est_L2x20+_term7} \\
&+ \int e^{i s \varphi} \mu_0^{k+1, 0, -1}(\xi, \eta) \widehat{u^0}(s, \xi - \eta) \widehat{f}(s, \eta) ~ d\eta \label{est_L2x20+_term8} \\
&+ \int e^{i s \varphi} \mu_0^{k+1, 0, 0}(\xi, \eta) \nabla_{\xi} \widehat{u^0}(s, \xi - \eta) \widehat{f}(s, \eta) ~ d\eta \label{est_L2x20+_term9} \\
&+ \int e^{i s \varphi} \mu_0^{k+1, -1, 1}(\xi, \eta) \nabla_{\xi} \widehat{u^0}(s, \xi - \eta) \widehat{f}(s, \eta) ~ d\eta \label{est_L2x20+_term10} \\
&+ \int e^{i s \varphi} \mu_0^{k+1, 0, 1}(\xi, \eta) \nabla_{\xi} \widehat{u^0}(s, \xi - \eta) \nabla_{\eta} \widehat{f}(s, \eta) ~ d\eta \label{est_L2x20+_term11}
\end{align}
\end{subequations}
Then, 
\begin{align*}
\Vert \eqref{est_L2x20+_term1} \Vert_{L^2} &= s \Vert T_{\mu_0}(\Lambda u^0, e^{\pm i s \Lambda} x f) \Vert_{H^{k+1}} \lesssim s \Vert u^0 \Vert_{W^{k+3, \infty-}} \Vert x f \Vert_{H^{k+1}} \lesssim s^{-1+a+\delta+\frac{\gamma_k}{2}} \Vert u \Vert_X^2 \lesssim s^{-1+\gamma_k} \Vert u \Vert_X^2 \\
\Vert \eqref{est_L2x20+_term5} \Vert_{L^2} &= s \Vert T_{\mu_0}(\Lambda x u^0, u) \Vert_{H^{k+1}} \lesssim s \Vert \Lambda x u^0 \Vert_{H^{k+1}} \Vert u \Vert_{W^{k+2, \infty-}} \lesssim s^{-1+a+\delta+\frac{\gamma_k}{2}} \Vert u \Vert_X^2 \lesssim s^{-1+\gamma_k} \Vert u \Vert_X^2 \\
\Vert \eqref{est_L2x20+_term6}_{k = 0} \Vert_{L^2} &\lesssim \Vert T_{\mu_0}(\Lambda u^0, u) \Vert_{L^{6/5}} \lesssim \Vert u^0 \Vert_{H^1} \Vert u \Vert_{L^3} \lesssim s^{-4/3+a} \Vert u \Vert_X^2 \\
\Vert \eqref{est_L2x20+_term9} \Vert_{L^2} &= \Vert T_{\mu_0}(x u^0, u) \Vert_{H^{k+1}} \lesssim \Vert x u^0 \Vert_{W^{k+1, 6}} \Vert u \Vert_{W^{k+1, 3}} \lesssim \Vert \Lambda x u^0 \Vert_{H^{k+1}} \Vert u \Vert_{W^{k+1, 3}} \lesssim s^{-4/3+a+\delta+\frac{\gamma_k}{2}} \Vert u \Vert_X^2 \\
\Vert \eqref{est_L2x20+_term11} \Vert_{L^2} &= \Vert T_{\mu_0}(\Lambda x u^0, e^{\pm i s \Lambda} x f) \Vert_{H^{k+1}} \lesssim \Vert \Lambda x u^0 \Vert_{H^{k+1}} \Vert x f \Vert_{H^2} + \Vert \Lambda x u^0 \Vert_{H^2} \Vert x f \Vert_{H^{k+1}} \lesssim s^{-1+a+\frac{\gamma_k}{2}} \Vert u \Vert_X^2
\end{align*}
The other are simpler or similar. 

\subsection{\texorpdfstring{$0 0$}{0 0} interactions}

We rewrite the corresponding part of $\xi^{\alpha} \Delta_{\xi} \widehat{f}(s, \xi)$ as~:  
\begin{align*}
\eqref{est_L2x2_init} &= \int e^{i s \varphi} s^2 \mu_0^{k+1, 0, 1}(\xi, \eta) \widehat{u^0}(s, \xi - \eta) \widehat{u^0}(s, \eta) ~ d\eta + \int e^{i s \varphi} s \mu_0^{k, 0, 1}(\xi, \eta) \widehat{u^0}(s, \xi - \eta) \widehat{u^0}(s, \eta) ~ d\eta \\
&+ \int e^{i s \varphi} s \mu_0^{k+1, 0, 0}(\xi, \eta) \widehat{u^0}(s, \xi - \eta) \widehat{u^0}(s, \eta) ~ d\eta + \int e^{i s \varphi} s \mu_0^{k+1, 0, 1}(\xi, \eta) \nabla_{\xi} \widehat{u^0}(s, \xi - \eta) \widehat{u^0}(s, \eta) ~ d\eta \\
&+ \int e^{i s \varphi} \mu_0^{k-1, 0, 1}(\xi, \eta) \widehat{u^0}(s, \xi - \eta) \widehat{u^0}(s, \eta) ~ d\eta + \int e^{i s \varphi} \mu_0^{k, 0, 1}(\xi, \eta) \nabla_{\xi} \widehat{u^0}(s, \xi - \eta) \widehat{u^0}(s, \eta) ~ d\eta \\
&+ \int e^{i s \varphi} \mu_0^{k+1, 0, -1}(\xi, \eta) \widehat{u^0}(s, \xi - \eta) \widehat{u^0}(s, \eta) ~ d\eta + \int e^{i s \varphi} \mu_0^{k+1, 0, 0}(\xi, \eta) \nabla_{\xi} \widehat{u^0}(s, \xi - \eta) \widehat{u^0}(s, \eta) ~ d\eta \\
&+ \int e^{i s \varphi} \mu_0^{k+1, 0, 1}(\xi, \eta) \Delta_{\xi} \widehat{u^0}(s, \xi - \eta) \widehat{u^0}(s, \eta) ~ d\eta
\end{align*}
We integrate by parts on the last term to get~:
\begin{subequations}
\begin{align}
\eqref{est_L2x2_init} &= \int e^{i s \varphi} s^2 \mu_0^{k+1, 0, 1}(\xi, \eta) \widehat{u^0}(s, \xi - \eta) \widehat{u^0}(s, \eta) ~ d\eta \label{est_L2x200_term1} \\
&+ \int e^{i s \varphi} s \mu_0^{k, 0, 1}(\xi, \eta) \widehat{u^0}(s, \xi - \eta) \widehat{u^0}(s, \eta) ~ d\eta \label{est_L2x200_term2} \\
&+ \int e^{i s \varphi} s \mu_0^{k+1, 0, 0}(\xi, \eta) \widehat{u^0}(s, \xi - \eta) \widehat{u^0}(s, \eta) ~ d\eta \label{est_L2x200_term3} \\
&+ \int e^{i s \varphi} s \mu_0^{k+1, 0, 1}(\xi, \eta) \nabla_{\xi} \widehat{u^0}(s, \xi - \eta) \widehat{u^0}(s, \eta) ~ d\eta \label{est_L2x200_term4} \\
&+ \int e^{i s \varphi} \mu_0^{k-1, 0, 1}(\xi, \eta) \widehat{u^0}(s, \xi - \eta) \widehat{u^0}(s, \eta) ~ d\eta \label{est_L2x200_term5} \\
&+ \int e^{i s \varphi} \mu_0^{k, 0, 1}(\xi, \eta) \nabla_{\xi} \widehat{u^0}(s, \xi - \eta) \widehat{u^0}(s, \eta) ~ d\eta \label{est_L2x200_term6} \\
&+ \int e^{i s \varphi} \mu_0^{k+1, 0, -1}(\xi, \eta) \widehat{u^0}(s, \xi - \eta) \widehat{u^0}(s, \eta) ~ d\eta \label{est_L2x200_term7} \\
&+ \int e^{i s \varphi} \mu_0^{k+1, 0, 0}(\xi, \eta) \nabla_{\xi} \widehat{u^0}(s, \xi - \eta) \widehat{u^0}(s, \eta) ~ d\eta \label{est_L2x200_term8} \\
&+ \int e^{i s \varphi} \mu_0^{k+1, -1, 1}(\xi, \eta) \nabla_{\xi} \widehat{u^0}(s, \xi - \eta) \widehat{u^0}(s, \eta) ~ d\eta \label{est_L2x200_term9} \\
&+ \int e^{i s \varphi} \mu_0^{k+1, 0, 1}(\xi, \eta) \nabla_{\xi} \widehat{u^0}(s, \xi - \eta) \nabla_{\eta} \widehat{u^0}(s, \eta) ~ d\eta \label{est_L2x200_term10}
\end{align}
\end{subequations}
Then~: 
\begin{align*}
\Vert \eqref{est_L2x200_term1} \Vert_{L^2} &= s^2 \Vert T_{\mu_0}(\Lambda u^0, u^0) \Vert_{H^{k+1}} \lesssim s^2 \Vert u^0 \Vert_{H^{k+2}} \Vert u^0 \Vert_{W^{k+2, \infty-}} \lesssim s^{-1+2a+\delta} \Vert u \Vert_X^2 \\
\Vert \eqref{est_L2x200_term4} \Vert_{L^2} &= s \Vert T_{\mu_0}(\Lambda x u^0, u^0) \Vert_{H^{k+1}} \lesssim s \Vert \Lambda x u^0 \Vert_{H^{k+1}} \Vert u^0 \Vert_{W^{k+2, \infty-}} \lesssim s^{-2+2a+\frac{\gamma_k}{2}} \Vert u \Vert_X^2 \\
\Vert \eqref{est_L2x200_term5}_{k = 0} \Vert_{L^2} &\lesssim \Vert T_{\mu_0}(\Lambda u^0, u^0) \Vert_{L^{6/5}} \lesssim \Vert u^0 \Vert_{H^1} \Vert u^0 \Vert_{L^3} \lesssim s^{7/3+2a} \Vert u \Vert_X^2 \\
\Vert \eqref{est_L2x200_term8} \Vert_{L^2} &= \Vert T_{\mu_0}(x u^0, u^0) \Vert_{H^{k+1}} \lesssim \Vert \Lambda x u^0 \Vert_{H^k} \Vert u^0 \Vert_{W^{k+2, \infty-}} + \Vert x u^0 \Vert_{L^6} \Vert u^0 \Vert_{W^{k+1, 3}} \lesssim s^{-7/3+2a+\delta+\frac{\gamma_k}{2}} \Vert u \Vert_X^2 \\
&\lesssim s^{-1+\gamma_k} \Vert u \Vert_X^2 \\
\Vert \eqref{est_L2x200_term10} \Vert_{L^2} &= \Vert T_{\mu_0}(\Lambda x u^0, x u^0) \Vert_{H^{k+1}} \lesssim \Vert \Lambda x u^0 \Vert_{H^{k+1}} \Vert x u^0 \Vert_{W^{1, 6}} + \Vert \Lambda x u^0 \Vert_{W^{1, 6}} \Vert \Lambda x u^0 \Vert_{H^k} \\
&\lesssim s^{-2+2a+\frac{\gamma_k}{2}} \Vert u \Vert_X^2 \lesssim s^{-1+\gamma_k} \Vert u \Vert_X^2
\end{align*}
The remaining terms are simpler or similar. 

\section*{Acknowledgments}

The author is very grateful to Yann Brenier and Frédéric Rousset for their constant support and numerous advices in the study of this problem. 

\printbibliography

\end{document}